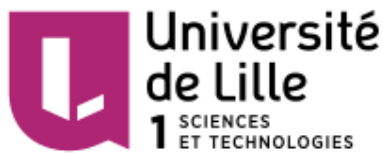
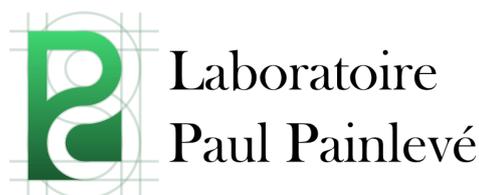

Université Lille 1

Laboratoire Paul Painlevé
CNRS UMR 8524

École Doctorale Sciences pour L'ingénieur
Université Lille Nord-de-France

# Thèse de doctorat

Discipline : Mathématiques

présentée par

# NGUYEN HUU KIEN

---

# La rationalité uniforme de la série Poincaré de relations d'équivalence p-adiques et la conjecture d'Igusa sur des sommes exponentielles

---

dirigée par Raf Cluckers

*Cette thèse est dédiée à mon père, ma mère
ma femme, Thị Thu et ma fille, Quỳnh Chi
sans lesquelles je n'aurais jamais vu le jour.*

# La rationalité uniforme de la série Poincaré de relations d'équivalence p-adiques et la conjecture d'Igusa sur des sommes exponentielles

## Résumé

Comprendre l'ensemble des solutions de système d'équation congruence est un problème important de la théorie des nombres. Un outil utile pour étudier le problème est la rationalité des séries de Poincaré associées avec les nombres des solutions d'une équation congruence modulo $p^m$, où $p$ est un nombre premier et $m$ varie dans l'ensemble des nombres naturels. Pour un premier fixé, la rationalité des series de Poincaré a été prouvée par Igusa, Meuser et Denef par quelques méthodes. Denef a aussi donné une extension de ce problème par replaçant le système ci-dessus par une famille définissable des relations équivalences indexés par éléments de group valué dans la structure de corps valués et considérant Poincaré séries de cette famille. Cette extension a une relation avec l'existence d'élimination des imaginaires des théories sur corps valués (voir [HMR]). La théorie d'intégration motivique nous aide pour montrer la dépendance uniforme dans $p$ de la rationalité des séries de Poincaré. Ce résultat a des applications dans fonctions zêta de groupe, etc. Dans le chapitre 1 de cette thèse, je donne une extension du résultat sur la rationalité uniforme dans $p$ des séries de Poincaré associées avec famille définissable des relations équivalences dans quelques théories sur corps valués dans lesquels élimination des imaginaires n'est pas prouvée comme théories sur structures analytiques. Ma méthode est que j'étends la théorie d'intégration motivique pour fonctions constructibles dans [CL08] et [CL15] aux fonctions constructibles rationnelles.

Autre problème important dans la théorie des nombres est estimation des sommes exponentielles. Sommes exponentielles modulo un nombre premier $p$ ont été recherchées par beaucoup de mathématicien comme Weil, Deligne, Katz, Laumon, etc. Sommes exponentielles modulo $p^m$ pour un nombre premier $p$ et un nombre naturel $m$ ont été étudiée par Igusa. Pour un nombre premier $p$ fixé, Igusa a découvert une relation profonde entre estimation des sommes exponentielles avec pôles de la fonction zêta local d'Igusa. Igusa a aussi montré que si on a une estimation uniforme dans $p$ et $m$ des sommes exponentielles, on peut obtenir une formule sommatoire de Poisson pour Adèle de type de Siegel-Weil. Il y a beaucoup de recherche dans cette direction pour obtenir la meilleur estimation uniforme dans $p$ et $m$ des sommes exponentielles. Dans les chapitres 2, 3, 4, on va prouver quelques versions uniformes dans $p$ et $m$ de borne supérieure des sommes exponentielles donnés par le seuil log canonique ou par polyèdre de Newton.

# Uniform rationality of Poincaré series of $p$-adic equivalence relations and Igusa's conjecture on exponential sums

## Abstract

Understanding the set of solutions of a system of congruence equation is an important problem of number theory. A useful tool to approach this problem is the rationality of Poincaré series associated with the numbers of solution of congruence equation modulo $p^m$, where $p$ is a prime and $m$ ranges over the set of natural numbers. The results in the rationality of Poincaré series for a fixed prime was proved by Igusa, Meuser and Denef by various methods. Denef gave an extension of this problem by considering the set of equivalence classes associated with a definable family of equivalence relations indexed by value group in the structure of valued field and defining Poincaré series of this family. This extension has relation with the existence of elimination of imaginaries theorem for theories of valued fields (see [HMR]). The development of motivic integration theory provides us with a new tool to show the uniform dependence of the rationality of Poincaré series on $p$. This can have many applications, for example, to the zeta function of groups. In the chapter 1 of this thesis, I extend the result on $p$-uniform rationality of Poincaré series associated with definable family of equivalence relations in some theories in which elimination of imaginaries has not been proved yet, for example theories on analytic structures. My method is that I extend the motivic integration theory for constructible motivic functions in [CL08] and [CL15] to rational constructible motivic functions.

Another classical problem of number theory is estimation of exponential sums. Exponential sums modulo a prime $p$ was studied by many mathematicians as Weil, Deligne, Katz, Laumon, etcetera. Exponential sums modulo $p^m$ was studied by Igusa, and for a fixed prime $p$, he gave a deep relation between estimation of exponential sums modulo $p^m$ and poles of Igusa local zeta function. Igusa also showed that an uniform estimation in $p$ and $m$ of exponential sums modulo $p^m$ could give a Poisson summation formula of Siegel-Weil type. By this motivation, many researchers tried to give the best uniform upper bound of exponential sums modulo $p^m$. In the chapters 2, 3, 4, we will try to obtain some uniform versions for upper bound of exponential sums modulo $p^m$ given by log-canonical threshold or Newton polyhedron due to Igusa's, Denef-Sperber's and Cluckers-Veys's conjectures.

## Keywords

# Remerciements

Je voudrais tout d'abord remercier mon directeur de thèse Monsieur Raf Cluckers pour m'avoir proposé un sujet de recherche passionnant avec de nombreuses perspectives d'évolution, pour m'avoir encadré et conseillé depuis quatre années, pour les nombreuses discussions mathématiques, pour toutes les heures passées à relire avec minutie mes articles et ma thèse et pour la liberté totale accordée durant cette thèse. C'est un vrai plaisir de pouvoir exprimer ici ma profonde gratitude ; sans sa patience, disponibilité et soutien permanent cette thèse n'aurait jamais vu le jour. Être son étudiant a été un vrai honneur.

Je remercie mes collègues Jorge Cely, Pablo Cubides Kovacsics pour beaucoup de discussion utile. En particulier, Je voudrais remercie Wouter Castryck, Saskia Chambilles qui ont coopération avec moi dans quelques parts de cette thèse. Je remercie Professeur Wim Veys, Professeur Immi Halupczok, Professeur Nero Budur, Professeur Johannes Nicaise pour quelques discussions intéressantes dans chapitres 3 et 4.

Je remercie chaleureusement Professeur Zoé Chatzidakis, Professeur Pierre Dèbes, Professeur Deirdre Haskell, Professeur François Loeser, Professeur Christopher Voll d'avoir accepté d'être les membres du jury de ma thèse.

Je voudrais remercier Fondation Sciences Mathématiques de Paris qui a m'aidé la finance pour étudier le Master 2. Je voudrais aussi remercier mes professeurs de l'institut de Math Jussieu qui m'ont fourni beaucoup de la connaissance des mathématiques quand j'ai étudié le Master 2. En particulier, j'exprime ma gratitude aux mon tuteur Professeur Pierre Charollois, Professeur François Loeser et Professeur Antoine Ducros.

Je tiens à exprimer ma gratitude à mes professeurs de l'École Normale Supérieure de Hanoi et l'Institut de Mathématiques du Vietnam, le programme de Master internationale en mathématiques à l'Institut de Mathématiques du Vietnam. En particulier, je ne peux pas me permettre d'oublier de remercier Professeur Ngô Việt Trung, Professeur Phùng Hồ Hải, Professeur Trần Văn Tấn, Professeur Nguyễn Việt Dũng, Professeur Đỗ Đức Thái m'ont aidé beaucoup quand j'ai étudié au Vietnam.

Je remercie tous les membres du séminaire Gémétrie et théorie des Modèles (2014-2018) qui ont m'aidé beaucoup dans mon étude.

Je tiens à remercier tous les membres de Laboratoire Paul Painleé, l'Université Lille 1 et École Doctorale Sciences pour L'ingénieur Université Lille Nord -de-France, CNRS Délégation Nord Pas-de-Calais et Picardie pour m'avoir accueilli. En particulier, je voudrais remercier Professeur Pierre Dèbes, Professeur Benoit Fresse, Madame Soledad Cuenca, Madame Fatima Soulaimani, Madame Fumery Ludivine, Madame Aurore Smets, Madame Sabine Hertsoen, Madame Frédérique Maréchal, Monsieur Taumi Lasri, Madame Thi Nguyen, Madame Sandrine Schwenck.

J'aimerais remercier aux mes amis vietnamiens en France : anh Thường, anh Hoá, anh Việt et Long. Un merci particulier aux anh Hà, anh Quang, anh Điệp, Kị, Nhật, Thi qui m'ont aidé beaucoup dans la vie.

Un grand merci chaleureusement à ma famille. Merci profondément à mes parents, mes beaux-parents, mes beaux-grands-parents, mon beau-grand-père, mes frères, mon beau-frère et mes soeurs d'être, malgré la distance, toujours là.

Finalement, je suis extrêmement reconnaissant à ma femme Nguyễn Thị Thu et ma fille Nguyễn Quỳnh Chi qui m'ont toujours encouragé et soutenu dans les moments difficiles à terminer cette thèse. Cette thèse est dédiée à elles.

# Acknowledgment

During the preparation of this thesis, the author supported by the European Research Council under the European Community's Seventh Framework Programme (FP7/2007-2013) with ERC Grant Agreement nr. 615722 MOTMELSUM, by the Labex CEMPI (ANR-11-LABX-0007-01).

# Table des matières

# Introduction

The thesis consists of a collection of published papers and preprints (each chapter representing one paper, see [Ngud], [CNa], [CNb], [CNc]) together with an extra general introduction. In this general introduction, I will give an overview of the topics of my thesis, it contains the sketch of history, the meaning of the problems and the ideas from the past to present.

This research is inspired by the work of my thesis advisor, R. Cluckers, it is a combination of model theory, number theory and algebraic geometry by the influence of J. Igusa, J. Denef and F. Loeser.

## 1. Uniform rationality of Poincaré series on local fields

Let $f_1, ..., f_k$ be polynomials in $\mathbb{Z}[x_1, ..., x_n]$ and $N$ be a positive integer. A classical problem of number theory is determining the solutions of congruence equation system

$$f_i(x_1, ..., x_n) \equiv 0 \mod N \tag{0.0.0.1}$$

where $i = \overline{1, k}$.

By Chinese remainder theorem, we can reduce the question to the case in which $N$ is a power $p^m$ of a prime $p$ to obtain system

$$f_i(x_1, ..., x_n) \equiv 0 \mod p^m \tag{0.0.0.2}$$

where $i = \overline{1, k}$.

One of main tools to research this problem is using Poincaré series associated with 0.0.0.2. More precisely, for each prime $p$ and a natural number $m$, let $N_{m,p}$ be the number of solutions modulo $p^m$ of 0.0.0.2 and consider the Poincaré series

$$P_{f_1,...,f_k,p}(T) = \sum_{m \geq 0} N_{m,p} T^m$$

We remark that Poincaré did not work on this series but we construct it by the spirit of Poincaré's idea. More precisely, from Poincaré's idea, to research a graded object, we try to construct a power series whose coefficient with respect to the $d$-exponent is a quantity (as dimension, rank, the number of points, etcetera) associated with the $d$-grade of this object. A well-known conjecture of Borevich and Shafarevich in [BS66] stated that $P_{f_1,...,f_k,p}(T)$ is a rational function, i.e. a quotient of two polynomials in $T$. This question was proved by Igusa (see [Igu74], [Igu78]) in case $k = 1$ and by Meuser (see [Meu81]) for arbitrary $k$, by using Hironaka's resolution of singularities and $p$-adic integration. From the viewpoint of model theory, in [Den84], Denef showed that this conjecture can be generalized by replacing system 0.0.0.2 by equivalence relations indexed by $m$ in the language of valued field and he proved the rationality of the associated Poincaré series in this generality by using cell decomposition theorem in place of Hironaka's resolution of singularities. Furthermore, Denef's proof suggests us that we can compare $P_{f_1,\ldots,f_k,p}(T)$ and $P_{f_1,\ldots,f_k,p'}(T)$ for primes $p$ and $p'$ because the equivalence relations are described independent with $p$ in the language of valued fields. Hence a motivic question is raised: How does $P_{f_1,\ldots,f_k,p}(T)$ change when $p$ runs in the set of primes?

In 1995, Maxim Kontsevich introduced the notation 'motivic integration' and it was developed by Denef and Loeser. This type of integral allows us to define measure of subsets of arc spaces with values in the Grothendieck ring of varieties (over field of characteristic zero). The rationality of Poincaré series has an important role in the motivic integration theory. More precisely, it helps us to control the convergence of integrals over value group. To see the rational result in Denef-Loeser's theory, we need some notations. For each algebraic variety $X$ over a field $k$ of characteristic zero, we denote by $\mathcal{L}(X)$ the scheme of germs of arcs on $X$, it is a scheme over $k$ such that for any extension $k \subset K$, there is a natural bijection

$$\mathcal{L}(X)(K) \simeq Mor_{k-schemes}(\operatorname{Spec} K[[t]], X)$$

between the set of $K$-rational points of $\mathcal{L}(X)$ and the set of $K[[t]]$-rational points of $X$. The scheme $\mathcal{L}(X)$ is defined as projective limit $\mathcal{L}(X) := \varprojlim_m \mathcal{L}_m(X)$, where $\mathcal{L}_m(X)$ is the space of $m$-jets on $X$ representing the functor

$$R \mapsto Mor_{k-schemes}(R[[t]]/(t^{m+1}R[[t]]), X)$$

defined on the category of $k$-algebras. We denote by $\pi_m$ the canonical map $\mathcal{L}(X) \to \mathcal{L}_m(X)$ corresponding to truncation of arcs.

Motivated by Denef's result in [Den84], Denef and Loeser showed in [DL99a] that the Poincaré series

$$P_X(T) = \sum_{m \geq 0} [\pi_m(\mathcal{L}(X)]T^m$$

considered as an element in $\mathcal{M}_{loc}[[T]]$ is rational and belongs to $\mathcal{M}[T]_{loc}$. Here $\mathcal{M}$ is the Grothendieck ring of varieties over $k$, $\mathcal{M}_{loc} = \mathcal{M}[\mathbb{L}^{-1}]$ with $\mathbb{L} := [\mathbb{A}_k^1]$ is the class of the line in $\mathcal{M}$ and $\mathcal{M}[T]_{loc}$ is the subring of $\mathcal{M}_{loc}[[T]]$ generated by $\mathcal{M}_{loc}[T]$ and elements of the form $(1 - \mathbb{L}^a T^b)^{-1}$ with $a \in \mathbb{Z}$ and $b \in \mathbb{N}\setminus\{0\}$. It means that we can write

$$P_X(T) = \frac{\sum_{i=0}^r [Y_i]T^i}{\mathbb{L}^c \prod_{j=1}^s (1 - \mathbb{L}^{a_j} T^{b_j})}$$

where $[Y_i] \in \mathcal{M}$ and $a_j \in \mathbb{Z}$, $b_j \in \mathbb{N}\setminus\{0\}$ and $c \in \mathbb{N}$.

When $k = \mathbb{Q}$, suppose that $f_i$ is reduced for all $i = \overline{1,k}$, and $X = \{f_1 = \ldots = f_k = 0\}$ is subvariety of $\mathbb{A}_\mathbb{Q}^n$ associated with $f_1, \ldots, f_k$. There is a positive integer $M$ such that we can view $Y_i$ as a $\mathbb{Z}[1/M]$-scheme so we can consider $\overline{Y}_{i,p}$ as the reduction modulo $p$ of $Y_i$ for $p > M$. By logical compactness we can show that

$$P_{f_1,\ldots,f_k,p}(T) = \frac{\sum_{i=0}^r \#\overline{Y}_{i,p}T^i}{p^c \prod_{j=1}^s (1 - p^{a_j} T^{b_j})}$$

for any big enough prime $p$ (see [DL02]). So the motivic question for $P_{f_1,\ldots,f_k,p}(T)$ reduces to the question: For a $\mathbb{Z}[1/M]$-scheme $Y$, how does $\#\overline{Y}_p$ change when $p$ runs in the set of primes. This question go to the Conjecture de Weil and the development of theory of motives with the contribution of so many mathematicians as Weil, Grothendieck, Deligne, Tate, etcetera.

To generalize the motivic question to equivalence relations, we need to fix an extension $\mathcal{L}$ of the language of valued fields and a theory $\mathcal{T}$ over this fixed language, and then we need to construct a motivic integration theory over $\mathcal{T}$. In [CL08], Cluckers and Loeser construct a theory of constructible motivic functions for language Denef-Pas $\mathcal{L}_{\text{DP}}$ and the $\mathcal{L}_{\text{DP}}$-theory of Henselian valued fields with residue field of characteristic zero by giving a motivic cell decomposition theorem. In [HK06], Hrushovski and Kazdan construct a motivic integration theory over theory of algebraically closed valued field with equicharacteristic zero $(\text{ACVF}_{(0,0)})$ over a compatible language. Cluckers and Loeser extend motivic integration theory to some theories of $(0, p, e)$-fields in [CL15] by giving some axioms to keep good properties as existence of cell decomposition theorem, orthogonality of residue field and valued group, etcetera. Even when we have a motivic integration theory, the question for rationality of Poincaré series of equivalence relations still remains hard. The problem is that the set of equivalence classes can not be described by a definable set in general case. In the original question, the set of equivalence classes can be identified with the set of closed points of a variety, so we can get around the difficulty. A solution of this problem in model theory is that we can add an imaginary sort for each equivalent relation to ensure that the set of equivalence classes is definable. But this way can cause a serious trouble with motivic integration theory. To see the trouble, we observe that each definable equivalence relation $E$ gives us a definable set $\mathcal{A}_E$ (moduli space of $E$) consisting equivalence classes of $E$ in the new language. In the general case, we may not ensure that $\mathcal{A}_E$ has good enough structure to define its volume and the integral of functions on it. This trouble suggests us to give another solution is that we try to prove that theory $\mathcal{T}$ eliminates imaginaries, it means that, we can definably associate a point to each equivalence class (imaginary element). In [HHM06], the authors gave a solution for the elimination of imaginaries by constructing a language called by geometric language $\mathcal{L}_{\mathcal{G}}$ and showed that the theory ACVF eliminates imaginaries in $\mathcal{L}_{\mathcal{G}}$. In fact, they added some imaginary sorts whose moduli space is simple to describe. Using this idea, in [HMR], the authors proved a uniform version in $p$ of the elimination of imaginaries in $\mathcal{L}_{\mathcal{G}}$ and it suffices to give a motivic version for rationality of Poincaré series of equivalence relations because of the special properties of moduli space associated with each imaginary sort. Unfortunately, in many cases, theory does not admit elimination of imaginaries, for instance, some theories over an analytique structure as in [HHM13].

In chapter 1, I will give an other way to solve this problem. I developed an idea in the appendix of [HMR] given by Cluckers. The idea is that, instead of trying to code an imaginary element by a point, we will code it by a definable set with two key properties: we can calculate the volume of this set by a very simple way, and there is a function of motivic type in this set whose integral is 1. However the motivic function that we constructed is not the usual motivic function because we need to divided the usual motivic function by some factors in the Grothendieck semi-ring of definable sets over the residue field. It suggests us to enlarge the set of motivic function. In order to follow the method for constructing motivic functions outlined in [CL15], I work with an extension $\mathcal{L}$ (with the same sorts) of Denef-Pas language $\mathcal{L}_{\mathrm{DP}}$ and an $\mathcal{L}$-theory $\mathcal{T}$ containing theory of henselian valued fields whose residue field of characteristic zero. I also require some futher axioms for including Jacobian property for functions from field sort to field sort, Split property for definable sets in a finite Cartesian product of auxiliary sorts, finite $b$-minimality for definable sets in field sort to ensure that we can construct a theory of motivic integration for theory $\mathcal{T}$. These axioms imply the existence of cell decomposition theorem and the orthogonality of residue sort and valued group sort) and hence we can enlarge the class of constructible motivic functions to the class of functions obtained by localizing the semi-ring of constructible functions by the multiplicatively closed set of 'non-zero' elements in the Grothendieck semi-ring of definable sets over residue field. We will call these functions by the rational constructible motivic functions. Note that we define an element $Z$ to be non-zero in the Grothendieck semi-ring of definable set over residue field is 'non-zero' if $Z$ is nonempty set. By this method we can answer the motivic question of Poincaré series for any theory containing such axioms. In particular, we can do it for many theories over analytic structures in which elimination of imaginaries have not been proved yet. Remark that this technique can be used to obtain rationality result of Poincaré series of equivalence relations in the field of Laurent series $K((t))$ for some field $K$ having minimality property as strong minimality, $o$-minimality, etcetera. But we need some extra small technique and I will not call it in this thesis.

This problem may have applications in the theory of Zeta function of groups and related objects (see e.g. [HMR] and [Vol10]).

## 2. Exponential sums modulo $p^m$

This section explains my work in chapter 2, 3, 4.

Exponential sums is a classical problem of number theory with many interesting consequences. Let $f(x_1, ..., x_n) \in \mathbb{Z}[x_1, ..., x_n]$ be a non-constant polynomial and $N > 1$ be a natural number, we consider the exponential sums:

$$E_{f,N} := \frac{1}{N^n} \sum_{(x_1,...,x_n) \in (\mathbb{Z}/N\mathbb{Z})^n} \exp(\frac{2\pi i f(x_1,...,x_n)}{N})$$

**Question:** Find the best upper bound of $|E_{f,N}|$

By Chinese remainder theorem, we can reduce the question to estimate

$$E_{f,p,m} := \int_{\mathbb{Z}_p^n} \exp(\frac{2\pi i f(x)}{p^m})|dx| = \frac{1}{p^{mn}} \sum_{(x_1,...,x_n) \in (\mathbb{Z}/p^m\mathbb{Z})^n} \exp(\frac{2\pi i f(x_1,...,x_n)}{p^m})$$

where $p$ is a prime and $m$ is a natural number.

In general case, let $L$ be a finite extension of $\mathbb{Q}_p$ or $\mathbb{F}_p((t))$, we will call $\mathcal{O}_L$ the valuation ring of $L$, $\mathcal{M}_L$ the maximal ideal of $\mathcal{O}_L$, $k_L$ the residue field of $L$ and $q_L = \#(k_L) = p_L^r$ where $p_L = char(k_L)$. Let $L$ be a finite extention of $\mathbb{Q}_p$ and $\Phi : \mathcal{O}_L^n \to \mathbb{C}$ be a Schwartz-Bruhat function we define

$$E_{f,L,z}^{\Phi} := \int_{x \in \mathcal{O}_L^n} \Phi(x) \Psi(zf(x)) |dx|$$

where $\Psi(z) := \exp(2\pi i \mathrm{Tr}_{L/\mathbb{Q}_p}(z))$ the standard additive character on $L$. Similarly, we can define exponential sums on an extension of $\mathbb{F}_p((t))$.

Research on exponential sums has a long history. It was studied by Gauss in the case $f = x^2$, $E_{x^2,N}$ is so-called the Gauss sum and is related to the law of quadratic reciprocity. $E_{f,p,1}$ was studied by Weil in [Wei48] if $f$ is a polynomial in one variable and by Deligne in [Del74] if $f$ is a polynomial in several variables. Deligne also explained the deep relation of exponential sums and conjecture de Weil. More precisely, the exponential sums in case of $m = 1$ depends on the eigenvalues of Frobenius operation acting on the groups of $\ell$-adic cohomology, Deligne's proof gave the best estimation in a class of polynomial (so-called them by Deligne polynomials). This idea have been developed by Katz, Laumon, etcetera to larger classes of polynomials.

In the general case, Igusa showed that if $char(L) = 0$ and the set of critical points of $f$ in $\mathrm{Supp}(\Phi)$ contained in $f^{-1}(0)$ then for $|z|$ big enough, we have

$$E_{f,L,z}^{\Phi} = \sum_{\chi,\lambda,\beta} a_{\chi,\lambda,\beta}(L) \chi(z) |z|^{\lambda} \mathrm{ord}_L(z)^{\beta} \qquad (0.0.0.3)$$

where $\chi$ runs in the set of multiplicative characters of $\mathcal{O}_L^*$, $\lambda$ runs in the set of poles of Igusa local zeta function $Z_\Phi(L, \chi, s, f)$ if $\chi \neq \chi_{trivial}$ and $(1 - q_L^{s+1})Z_\Phi(L, \chi, s, f)$ if $\chi = \chi_{trivial}$, $\beta \in \mathbb{N}$ and $\beta \leq$ (multiplicity of pole $\lambda$) $- 1$. We remark that there are only finite multiplicative characters $\chi$ such that $Z_\Phi(L, \chi, s, f)$ has pole, so the above sum is finite sum. Hence, for each $\sigma < -\lambda_{0,\Phi,L}$ with $\lambda_{0,\Phi,L}$ is the maximal real part of poles appearing in 0.0.0.3, there exists a constant $c(L, \Phi, \sigma)$ such that $|E_{f,L,z}^{\Phi}| \leq c(L, \Phi, \sigma) |z|^{-\sigma}$

Let $f$ be a non-constant homogeneous polynomial over number field $K$ and $\sigma \in \mathbb{R}$. Suppose that $1 < \sigma < \sigma^0$, where $\sigma^0$ is the minimum take over quotients $\frac{\nu_i}{N_i}$ with $(\nu_i, N_i) \neq (1, 1)$ is the numerical data associated with an exceptional divisor of a fixed embedded resolution of singularities of $\mathrm{Spec}(\mathbb{C}[x_1, ..., x_n]/(f(x)))$ in $\mathbb{A}_{\mathbb{C}}^n$. Because of property of Igusa zeta functions, if a local fields $L$ contains $K$ and $\Phi_L = \mathbf{1}_{\mathcal{O}_L^n}$ then $\sigma^0 \leq -\lambda_{0,\Phi_L,L}$. By formula 0.0.0.3, there exists a constant $c(L)$ such that $|E_{f,L,z}^{\Phi}| \leq c(L)|z|^{-\sigma}$. Igusa conjectured that one can take $c(L)$ independent of $L$. In particular, if we can choose $\sigma > 2$ then this upper bound implies an Adèlic Poisson summation formula of Siegel-Weil type.

It is obvious that $\sigma^0$ depends on an embedded resolution and it can tend to 1 by a suitable choice of resolution. Since 0.0.0.3, it is more natural to expect that Igusa's conjecture will be true for any $1 < \sigma < -\alpha$, where $\alpha = \limsup_{char(k_L) \longrightarrow +\infty} \lambda_{0,\Phi_L}$. Moreover we can also expect that Igusa's conjecture can be extended to all non-constant polynomials $f$ and $\sigma < -\alpha$.

Igusa proved his conjecture if $f$ defines a non-singular hypersurface in $\mathbb{P}^{n-1}$. In [DS01], Denef and Sperber proved Igusa's conjecture for a homogeneous polynomial which is non-degenerate with respect to its Newton polyhedron at origin (called *non-degenerate* polynomial) and has an extra condition on the Newton polyhedron. Remark that they only worked for $\sigma^0$ depended on a toric resolution of singularities of $\mathrm{Spec}(\mathbb{C}[x_1,...,x_n]/(f(x)))$ in $\mathbb{A}_\mathbb{C}^n$. Denef and Sperber also introduced local conjecture for exponential sums, i.e. we take $\Phi_L = \mathbf{1}_{\mathcal{M}_L^n}$, and they also proved it when a non-constant polynomial is *non-degenerate* and has the same condition on the Newton polyhedron as the homogeneous case. In [Clu08a], Cluckers extended the result of Denef and Sperber to *non-degenerate* polynomials by removing the condition on the Newton polygon of Denef and Sperber. In [Clu10], Cluckers continued his method to give a version for *non-degenerate* quasi-homogeneous polynomials. In [Wri], Wright removed the *non-degenerate* condition of Denef-Sperber and obtained a result in case of non-constant quasi-homogeneous polynomials in two variables. Lichtin reproved result of Wright by a different method and extended his result to case of non-constant homogeneous polynomials in three varibles (see [Lic13] and [Lic16]). In [CV16], Cluckers and Veys stated a conjecture for arbitrary non-constant polynomial by restricting to $\mathrm{ord}_L(z) \leq -2$ and working for $\sigma = \mathrm{lct}_0(f)$ (resp. $\sigma = \mathrm{lct}(f)$) in the local exponential sums conjecture of Denef and Sperber (resp. the global exponential sums conjecture of Igusa). Here $\mathrm{lct}_0(f)$ is the log-canonical threshold of $f$ at 0 and $\mathrm{lct}(f)$ is the minimum take over log-canonical thresholds $\mathrm{lct}_b(f - f(b))$, $b \in \mathbb{C}^n$. We remark that if $\Phi = \mathbf{1}_{\mathcal{M}_L^n}$ and $|\lambda_{0,\Phi,L}| < 1$ then $\lambda_{0,\Phi,L} = -\mathrm{lct}_0(f)$, so Cluckers-Veys gives us the best estimation of exponential sums in this case. Cluckers and Veys also conjectured a uniform version of local exponential sums of a non-constant polynomial f. More precisely, let $\mathcal{O}$ be the ring of integers of $K$, for each point $y \in \mathcal{O}^n$ and a non-archimedean local field $L$ containing $K$, we denote by $a_{y,L}(f) = \min_{y' \in y+\mathcal{M}_L^n} \mathrm{lct}_{y'}(f - f(y'))$. Cluckers-Veys stated that there exists a positive constant $c$ such that $|E_{f,L,z}^{\Phi_{y,L}}| \leq c|z|^{-a_{y,L}}|\mathrm{ord}(z)|^{n-1}$ for all $y \in \mathcal{O}^n$, all local field $L$ containing $K$ and all $\mathrm{ord}_L(z) \leq -2$, where $\Phi_{y,L} = \mathbf{1}_{y+\mathcal{M}_L^n}$.

By logical compactness, we can extend the above discussion to any local field $L$ of characteristic larger than a constant $M$ depended on $f$.

In chapter 2, my work (a joint work with Wouter Castryck, see [CNa]) contributes to remove the quasi-homogeneous condition in the proof of Cluckers in [Clu10] but we need to restrict to the case $\mathrm{ord}_L(z) \leq -2$.

In chapter 3, my work (a joint work with Saskia Chambille, see [CNb]) contributes to prove Cluckers-Veys's conjecture in the case $\mathrm{lct}_0(f)$ (resp. $\mathrm{lct}(f)$) less than one half.

In chapter 4, my work (a joint work with Raf Cluckers, see [CNc]) contributes to give a geometric variant of Cluckers-Veys conjecture. More precisely, we introduce conjectures on numerical data of embedded resolution of singularities and prove that the conjectures on numerical data and the conjectures on exponential sums for log-canonical threshold are equivalent. More results on the numerical data conjectures will be studied in future projects.

# Chapter 1

# Uniform rationality of the Poincaré series of definable, analytic equivalence relations on local fields

*This chapter will appear in Transactions of the American Mathematical Society, series B, see [Ngud].*

**Contents**

**Abstract**

Poincaré series of *p*-adic, definable equivalence relations have been studied in various cases since Igusa's and Denef's work related to counting solutions of polynomial equations modulo $p^n$ for prime *p*. General semi-algebraic equivalence relations on local fields have been studied uniformly in *p* recently by Hrushovski-Martin-Rideau. Here we generalize the rationality result of Hrushovski-Martin-Rideau to the analytic case, unifomly in *p*. In particular, the results hold for large positive characteristic local fields. We also introduce rational motivic constructible functions and their motivic integrals, as a tool to prove our main results.

## 1.1 Introduction

### 1.1.1

After observing that Igusa's and Denef's rationality results (see e.g. [Den84], [Den]) can be rephrased in terms of counting the number of equivalence classes of particular semi-algebraic equivalence relations, it becomes natural to consider more general definable equivalence relations in the *p*-adic context and study the number of equivalence classes and related Poincaré series. The study of uniform *p*-adic, semi-algebraic equivalence relations is one of the main themes of [HMR], with general rationality results of the associated Poincaré series as part of the main results, generalizing [Pas89] and [Mac90]. In the appendix of [HMR], a more direct way of obtaining such rationality results was developed in a different case, namely, in the subanalytic setting on $\mathbb{Q}_p$, generalizing the rationality results by Denef and van den Dries in [DD88]. A deep tool of [HMR] to study equivalence relations, called elimination of imaginaries, is very powerful but also problematic since it does not extend well to the analytic setting, see [HHM13]. In this paper we follow the more direct approach of the appendix of [HMR] to obtain rationality results in situations where elimination of imaginaries is absent; here, we make this approach uniform in non-archimedean local fields (non-archimedean locally compact field of any characteristic). The two main such situations where this applies come from analytic structures on the one hand, and from an axiomatic approach from [CL15] on the other hand. In the analytic case, our results generalize the uniform analytic results of [Dri92], [CLR06]. We heavily rely on cell decomposition, a tool which was not yet available at the time of [Dri92], and which is obtained more recently in an analytic context in [CL11], [CL] and [CLR06] uniformly, and in [Clu04] for any fixed *p*-adic field. In our approach we also need more general denominators than in previous studies, which we treat by introducing rational motivic constructible functions and their motivic integrals, a slight generalization of the motivic constructible functions from [CL15]. The adjective 'rational' reflects the extra localization of certain Grothendieck semi-rings as compared to [CL15]. For these integrals to be well-defined, a property called Jacobian property is used; also this property was not yet available at the time of [Dri92], and is shown in [CL11] for analytic structures.

**1.1.2**

Let us begin by rephrasing some of the classical results by e.g. Igusa and Denef in terms of definable equivalence relations. Let $\mathcal{L}_0$ be a first order, multi-sorted language such that $(\mathbb{Q}_p, \mathbb{Z}, \mathbb{F}_p)$ are $\mathcal{L}_0$-structures for all $p$. A basic example is the ring language on the first sort together with the valuation map $\mathrm{ord} : \mathbb{Q}_p^\times \to \mathbb{Z}$ and the ordering on $\mathbb{Z}$. We consider a formula $\varphi(x, y, n)$ in the language $\mathcal{L}_0$ with free variables $x$ and $y$ running over $\mathbb{Q}_p^m$ and $n \in \mathbb{Z}$.

Suppose that for each $n \geq 0$ and each prime $p$, the condition $\varphi(x, y, n)$ on $x, y$ yields an equivalence relation $\sim_{p,n}$ on $\mathbb{Q}_p^m$ (or on a uniformly definable subset $X_p$ of $\mathbb{Q}_p^m$) with finitely many, say $a_{\varphi, p, n}$, equivalence classes. Then we can consider, for each $p$, the Poincaré series associated to $\varphi$ and $p$:

$$P_{\varphi, p}(T) = \sum_{n \geq 0} a_{\varphi, p, n} T^n \tag{1.1.2.1}$$

When $\varphi(x, y, n)$ is the collection of equivalence relations $\sim_n$ based on the vanishing of a polynomial $f(x)$ modulo $p^n$, more precisely, when, for tuples $(x, y)$ with non-negative valuation, $\varphi$ is defined by

$$\mathrm{ord}(f(x), x - y) \geq n \tag{1.1.2.2}$$

with $f \in \mathbb{Z}[x_1, ..., x_m]$ and where the order of a tuple is the minimum of the orders, the question of the rationality of $P_{\varphi, p}(T)$ was conjectured by Borevich and Shafarevich in [BS66] and was proved by Igusa in [Igu74] and [Igu78]. The proof relied on Hironaka's resolution of singularities and used $p$-adic integration. In [Den84], still using $p$-adic integration but using cell decomposition instead of Hironaka's resolution of singularities, Denef proved the rationality of $P_{\varphi, p}$ for more general $\varphi$ than in (1.1.2.2) with its formula related to lifting solutions modulo $p^n$ to solutions in $\mathbb{Z}_p^m$, answering on the way a question given by Serre in [Ser81]. The idea of Denef (using [Mac76] and [Coh69]) is to represent a semi-algebraic set by a union of finitely many cells which have a simple description (and so does $\mathrm{ord}\, f$ on each cell, for $f$ a semi-algebraic function) so that we can integrate easily on them. This was made uniform in $p$ in [Pas89] and [Mac90]. The advantage of the approach via cell decomposition is that more general parameter situations can also be understood via parameter integrals, a feature heavily used to get Fubini theorems for motivic integrals [CLR06], [CL08], [CL15]; for us this approach leads to our main result Theorem 1.4.2 below, as a natural generalization of our rationality results Theorems 1.1.2 and 1.4.1.

When $f$ in (1.1.2.2) is given by a converging power series instead of a polynomial, then the rationality has been obtained in [DD88] for fixed $p$, in [Dri92] uniformly in $p$, and in [CLR06] uniformly in $p$ with the extra strength coming from the cell decomposition approach.

The rationality of $P_{\varphi, p}(T)$ as in (1.1.2.1) for more general $\varphi$, as well as the uniformity for large $p$ in $\mathbb{Q}_p$ and $\mathbb{F}_p((t))$, is the focus of this paper. A common feature of all the previously mentioned results is to bundle the information of $P_{\varphi, p}(T)$ into a $p$-adic integral of some kind, and then study these integrals qualitatively. Here, we bundle the information into slightly more general integrals than the ones occurring before.

**1.1.3**

Let us recall the precise result of [HMR] which states the rationality of $P_{\varphi,p}(T)$ in the semi-algebraic realm, with its uniformity in (non-archimedean) local fields. We recall that a non-archimedean local field is a locally compact topological field with respect to a non-discrete topology such that its topology defines a non-archimedean absolute value. Any such field is either a finite field extension of $\mathbb{Q}_p$ for some $p$ or isomorphic to $\mathbb{F}_q((t))$ for some prime power $q$; we will from now on say local field for non-archimedean local field.

Let $\mathcal{L}_{\mathrm{DP}}$ be the Denef-Pas language, namely with the ring language on the valued field sort and on the residue field sort, the language of ordered abelian groups on the value group, the valuation map, and an angular component map $\overline{\mathrm{ac}}$ from the valued field sort to the residue field sort (see Section 1.2.2). All local fields $K$ with a chosen uniformizer $\varpi$ are $\mathcal{L}_{\mathrm{DP}}$-structures, where the angular component map sends nonzero $x$ to the reduction of $x\varpi^{-\operatorname{ord} x}$ modulo the maximal ideal, and sends zero to zero. Let $\mathcal{O}_K$ denote the valuation ring of $K$ with maximal ideal $\mathcal{M}_K$ and residue field $k_K$ with $q_K$ elements and characteristic $p_K$.

Let $\varphi(x, y, n)$ be an $\mathcal{L}_{\mathrm{DP}}$-formula with free variables $x$ running over $K^m$, $y$ running over $K^m$ and $n$ running over $\mathbb{N}$, with $\mathbb{N}$ being the set of nonnegative integers. Suppose that for each local field $K$ and each $n$, $\varphi(x, y, n)$ defines an equivalence relation $\sim_{K,n}$ on $K^m$ with finitely many, say, $a_{\varphi,K,n}$, equivalence classes. (The situation that $\varphi(x, y, n)$ defines an equivalence relation on a uniformly definable subset $X_{K,n}$ of $K^m$ for each $n$ can be treated similarly, e.g. extending to a relation on $K^m$ by a single extra equivalence class.) For each local field $K$ consider the associated Poincaré series

$$P_{\varphi,K}(T) = \sum_{n \geq 0} a_{\varphi,K,n} T^n,$$

In [HMR], the authors proved the following (as well as a variant by adding constants of a ring of integers to $\mathcal{L}_{\mathrm{DP}}$, and by allowing $n$ and $T$ to be tuples of finite length instead of length one; these features are also captured in Theorem 1.4.2 below).

**Theorem 1.1.1.** *There exists $M > 0$ such that the power series $P_{\varphi,K}(T)$ is rational in $T$ for each local field $K$ whose residue field has characteristic at least $M$. Moreover, for such $K$, the series $P_{\varphi,K}(T)$ only depends on the residue field $k_K$ (namely, two local fields with isomorphic residue field give rise to the same Poincaré series).*

*More precisely, there exist nonnegative integers $a, N, M, k, b_j, e_i, q$, integers $a_j$, and formulas $X_i$ in the ring language for $i = 0, \ldots, N$ and $j = 0, \ldots, k$, such that for each $j$, $a_j$ and $b_j$ are not both $0$, $q$ is nonzero, and for all local fields $K$ with residue field of characteristic at least $M$ we have*

$$P_{\varphi,K}(T) = \frac{\sum_{i=0}^{N}(-1)^{e_i} \# X_i(k_K) T^i}{q \cdot q_K^a \prod_{j=1}^{k}(1 - q_K^{a_j} T^{b_j})},$$

*where $X_i(k_K)$ is the set of $k_K$-points satisfying $X_i$.*

This theorem is furthermore applied in [HMR] to the theory of zeta functions in group theory. Theorem 1.1.1 is shown in [HMR] by proving general elimination of imaginaries in a language $\mathcal{L}_\mathcal{G}$ called the geometrical language and which expands the language of valued fields. This elimination allows one to rewrite the data in terms of classical (Denef-Pas style) uniform *p*-adic integrals, from which rationality follows uniformly in the local field.

In the appendix of [HMR], a more direct but similar reduction to classical *p*-adic integrals is followed, and it is this reduction that is made uniform in the local field, here.

An interesting aspect of Theorem 1.1.1 is the appearance of the positive integer $q$ in the denominator. In more classical Poincaré series in this context (e.g. [Den84], [Den]), less general denominators suffice, namely without a factor $q > 0$. In this paper we use even more general denominators, namely, we may divide by the number of points on (nonempty and finite) definable subsets over the residue field. We develop a corresponding theory of *p*-adic and motivic integration, of what we call *rational* motivic constructible functions (altering the notion of motivic constructible functions from [CL08] and [CL15]). The benefits are that we need not restrict to the semi-algebraic case and that we don't rely on elimination of imaginaries. This allows us to obtain rationality in the uniform analytic contexts from [CL], [Dri92], and [CLR06], and also in the axiomatic context from [CL15].

Let us state our main result to indicate the more general nature of our denominators.

Let $\mathcal{T}$ be a theory of valued fields in a language $\mathcal{L}$ extending $\mathcal{L}_{\mathrm{DP}}$. Suppose that $\mathcal{T}$ has properties $(*)$ and $(**)$ as in the Definition 1.2.5 and 1.2.6 below (see Section 1.2.4 for concrete, analytic examples of such $\mathcal{T}$). Suppose for convenience here that every definable subset in the residue field sort is definable in the language of rings (this assumption is removed in the later form 1.4.1 of the main theorem of the introduction). Let $\varphi(x, y, n)$ be an $\mathcal{L}$-formula with free variables $x$ running over $K^m$, $y$ running over $K^m$ and $n$ running over $\mathbb{N}$. Suppose that for each local field $K$ and each $n$, $\varphi(x, y, n)$ gives an equivalence relation $\sim_{K,n}$ on $K^m$ with finitely many, say, $a_{\varphi,K,n}$, equivalence classes. For each local field $K$ consider the associated Poincaré series

$$P_{\varphi,K}(T) = \sum_{n \geq 0} a_{\varphi,K,n} T^n.$$

**Theorem 1.1.2.** *There exists $M > 0$ such that the power series $P_{\varphi,K}(T)$ is rational in $T$ for each local field $K$ whose residue field has characteristic at least $M$. Moreover, for such $K$, the series $P_{\varphi,K}(T)$ only depends on the residue field $k_K$ (namely, two local fields with isomorphic residue field give rise to the same Poincaré series).*

*More precisely, there exist nonnegative integers $N, M, k, b_j, e_i,$, integers $a_j$ and formulas $X_i$ and $Y$ in the ring language for $i = 0, \ldots, N$ and $j = 0, \ldots, k$, such that for each $j$, $a_j$ and $b_j$ are not both 0, and, for all local fields $K$ with residue field of characteristic at least $M$, $Y(k_K)$ is nonempty and*

$$P_{\varphi,K}(T) = \frac{\sum\limits_{i=0}^{N}(-1)^{e_i}\#X_i(k_K)T^i}{\#Y(k_K) \cdot \prod_{j=1}^{k}(1 - q_K^{a_j}T^{b_j})}.$$

As in [HMR], our theorem is related to zeta functions of groups, zeta functions of twist isoclasses of characters, the abscissa of convergence of Euler products, etc., but we do not give new applications in this direction as compared to [HMR].

In fact, we will give a more general theorem 1.4.2 which describes the dependence of the numbers $a_{\varphi,K,n}$ on $n$ (and on completely general parameters) by means of a rational motivic constructible function.

In Section 1.2 we recall the conditions on the language $\mathcal{L}$ from [CL15]. In Section 1.3 we introduce rational motivic constructible functions, their motivic integrals, and their specializations to local fields. In Section 1.4 we give some generalizations and the proofs of our main theorems.

## 1.2 Analytic languages

### 1.2.1

In Section 1.2.2 we recall the Denef-Pas language and quantifier elimination in its corresponding theory of henselian valued fields of characteristic zero. In Section 1.2.3 we develop axioms for expansions of the Denef-Pas language and its theory, following [CL15]. In Section 1.2.4 we recall that certain analytic structures satisfy the axioms from Section 1.2.3. Based on these axioms, we extend in Section 1.3 the motivic integration from [CL15] to a situation with more denominators.

### 1.2.2 The language of Denef-Pas

Let $K$ be a valued field, with valuation map $\mathrm{ord} : K^\times \to \Gamma_K$ for some additive ordered group $\Gamma_K$, $\mathcal{O}_K$ the valuation ring of $K$ with maximal ideal $\mathcal{M}_K$ and residue field $k_K$. We denote by $x \to \overline{x}$ the projection $\mathcal{O}_K \to k_K$ modulo $\mathcal{M}_K$. An angular component map (modulo $\mathcal{M}_K$) on $K$ is a multiplicative map $\overline{\mathrm{ac}} : K^\times \to k_K^\times$ extended by setting $\overline{\mathrm{ac}}(0) = 0$ and satisfying $\overline{\mathrm{ac}}(x) = \overline{x}$ for all $x$ with $\mathrm{ord}(x) = 0$.

The language $\mathcal{L}_{\mathrm{DP}}$ of Denef-Pas is the three-sorted language with as sorts:

(i) a sort VF for the valued field-sort,

(ii) a sort RF for the residue field-sort, and

(iii) a sort VG for the value group-sort,

and with the following languages, together with two maps between the sorts

$$(\mathcal{L}_{\mathrm{ring}}, \mathcal{L}_{\mathrm{ring}}, \mathcal{L}_{\mathrm{oag}}, \mathrm{ord}, \overline{\mathrm{ac}}).$$

The first copy of $\mathcal{L}_{\mathrm{ring}}$ is used for the sort VF, the second copy for RF, the language $\mathcal{L}_{\mathrm{oag}}$ is the language $(+, <)$ of ordered abelian groups for VG, ord denotes the valuation map on non-zero elements of VF, and $\overline{\mathrm{ac}}$ stands for an angular component map from VF to RF.

As usual for first order formulas, $\mathcal{L}_{\mathrm{DP}}$-formulas are built up from the $\mathcal{L}_{\mathrm{DP}}$-symbols together with variables, the logical connectives $\wedge$ (and), $\vee$ (or), $\neg$ (not), the quantifiers $\exists, \forall$, the equality symbol $=$, and possibly parameters (see [Mar02] for more details). If the parameters are in a set $A$, we will say '$A$-definable'. We say simply 'definable', if $A = \emptyset$.

Let us briefly recall the statement of the Denef-Pas theorem on elimination of valued field quantifiers in the language $\mathcal{L}_{\mathrm{DP}}$. Denote by $H_{\overline{\mathrm{ac}},0}$ the $\mathcal{L}_{\mathrm{DP}}$-theory of the above described structures whose valued field is henselian and whose residue field of characteristic zero. Then the theory $H_{\overline{\mathrm{ac}},0}$ admits elimination of quantifiers in the valued field sort, see [Pas89], Thm. 4.1 or [CL08], Thm. 2.1.1.

### 1.2.3 Expansions of the Denef-Pas language: an axiomatic approach

In this section we single out precise axioms needed to perform motivic integration, following [CL15]. Apart from cell decomposition, the axioms involve a Jacobian property for definable functions and a so-called property $(*)$ which requires at the same time orthogonality between the value group and residue field and that the value group has no other structure than that of an ordered abelian group. Although these theories are about equicharacteristic 0 valued fields, by logical compactness we will be able to use them for local fields of large residue field characteristic.

Let us fix a language $\mathcal{L}$ which contains $\mathcal{L}_{\mathrm{DP}}$ and which has the same sorts as $\mathcal{L}_{\mathrm{DP}}$. Let $\mathcal{T}$ be an $\mathcal{L}$-theory containing $H_{\overline{\mathrm{ac}},0}$. The requirements on $\mathcal{T}$ will be summarized in Definition 1.2.5 below.

**Definition 1.2.1.** (Jacobian property for a function). Let $K$ be a valued field. Let $F : B \to B'$ be a function with $B, B' \subset K$. We say that $F$ has the Jacobian property if the following conditions hold all together:
— $F$ is a bijection and $B, B'$ are open balls in $K$, namely of the form $\{x \mid \mathrm{ord}(x - a) > \gamma\}$ for some $a \in K$ and $\gamma \in \Gamma_K$,
— $F$ is $C^1$ on $B$ with derivative $F'$,
— $F'$ is nonvanishing and $\mathrm{ord}(F')$ and $\overline{\mathrm{ac}}(F')$ are constant on $B$,
— for all $x, y \in B$ we have

$$\mathrm{ord}(F') + \mathrm{ord}(x - y) = \mathrm{ord}(F(x) - F(y))$$

and

$$\overline{\mathrm{ac}}(F') \cdot \overline{\mathrm{ac}}(x - y) = \overline{\mathrm{ac}}(F(x) - F(y)).$$

**Definition 1.2.2.** (Jacobian property for $\mathcal{T}$). We say that the Jacobian property holds for the $\mathcal{L}$-theory $\mathcal{T}$ if for any model $\mathcal{K}$ the following holds.

Write $K$ for the VF-sort of $\mathcal{K}$. For any finite set $A$ in $\mathcal{K}$ and any $A$-definable function $F : K \to K$ there exists an $A$-definable function

$$f : K \to S$$

with $S$ a Cartesian product of the form $k_K^r \times \Gamma_K^t$ for some $r, t$ such that each infinite fiber $f^{-1}(s)$ is a ball on which $F$ is either constant or has the Jacobian property.

**Definition 1.2.3.** (Split). We say that $\mathcal{T}$ is split if the following conditions hold for any model $\mathcal{K}$. Write $K$ for the VF-sort of $\mathcal{K}$.
— any $\mathcal{K}$-definable subset of $\Gamma_K^r$ for any $r \geq 0$ is $\Gamma_K$-definable in the language of ordered abelian groups $(+, <)$,
— for any finite set $A$ in $\mathcal{K}$ and any $r, s \geq 0$, any $A$-definable subset $X \subset k_K^s \times \Gamma_K^r$ is equal to a finite disjoint union of $Y_i \times Z_i$ where the $Y_i$ are $A$-definable subsets of $k_K^s$ and the $Z_i$ are $A$-definable subsets of $\Gamma_K^r$.

**Definition 1.2.4.** (Finite $b$-minimality). The theory $\mathcal{T}$ is called finitely $b$-minimal if for any model $\mathcal{K}$ of $\mathcal{T}$ the following conditions hold. Write $K$ for the VF-sort of $\mathcal{K}$. Each locally constant $\mathcal{K}$-definable function $g : K^\times \to K$ has finite image and for any finite set $A$ in $K$ and any $A$-definable set $X \subset K$ there exist an $A$-definable function
$$f : X \to S$$
with $S$ a Cartesian product of the form $k_K^r \times \Gamma_K^t$ for some $r, t$ and an $A$-definable function
$$c : S \to K$$
such that each nonempty fiber $f^{-1}(s)$ of $s \in S$ is either the singleton $\{c(s)\}$ or the ball of the form
$$\{x \in K | \overline{\mathrm{ac}}(x - c(s)) = \xi(s), \mathrm{ord}(x - c(s)) = \eta(s)\}$$
for some $\xi(s)$ in $k_K$ and some $\eta(s) \in \Gamma_K$.

Recall that $\mathcal{T}$ is an $\mathcal{L}$-theory containing $H_{\overline{\mathrm{ac}},0}$, where $\mathcal{L}$ contains $\mathcal{L}_{\mathrm{DP}}$ and has the same sorts as $\mathcal{L}_{\mathrm{DP}}$.

**Definition 1.2.5.** We say that $\mathcal{T}$ has property $(*)$ if it is split, finitely $b$-minimal, and has the Jacobian property.

**Definition 1.2.6.** We say that $\mathcal{T}$ has property $(**)$ if all local fields can be equipped with $\mathcal{L}$-structure and such that, for any finite subtheory $\mathcal{T}'$ of $\mathcal{T}$, local fields with large enough residue field characteristic are models of $\mathcal{T}'$ when equipped with this $\mathcal{L}$-structure.

**Example 1.2.7.** The $\mathcal{L}_{\mathrm{DP}}$-theory $H_{\overline{\mathrm{ac}},0}$ of Henselian valued field with equicharacteristic $(0,0)$ has properties $(*)$ and $(**)$. It even has property $(*)$ in a resplendent way, namely, the theory $\mathcal{T}$ in an expansion $\mathcal{L}$ of $\mathcal{L}_{\mathrm{DP}}$ which is obtained from $\mathcal{L}_{\mathrm{DP}}$ by adding constant symbols from a substructure of a model of $H_{\overline{\mathrm{ac}},0}$ (and putting its diagram into $\mathcal{T}$) and by adding any collection of relation symbols on $\mathrm{RF}^n$ for $n \geq 0$ has property $(*)$, see [CL15], Thm. 3.10 or [CL08], Section 7 and Theorem 2.1.1.

Analytic examples of theories with properties $(*)$ and $(**)$ are given in the next section.

### 1.2.4 Analytic Expansions of the Denef-Pas language

Our main example is a uniform version (on henselian valued fields) of the $p$-adic subanalytic language of [DD88]. This uniform analytic structure is taken from [CL] and is a slight generalization of the uniform analytic structure introduced by van den Dries in [Dri92]; it also generalizes [CLR06]. While van den Dries obtained quantifier elimination results and Ax-Kochen principles, the full property $(*)$ is shown in the more recent work [CL11] and [CL]; see Remark 1.2.11 below for a more detailed comparison. Property $(**)$ will be naturally satisfied.

Fix a commutative noetherian ring $A$ (with unit $1 \neq 0$) and fix an ideal $I$ of $A$ with $I \neq A$. Suppose that $A$ is complete for the $I$-adic topology. By complete we mean that the inverse limit of $A/I^n$ for $n \in \mathbb{N}$ is naturally isomorphic to $A$. An already interesting example is $A = \mathbb{Z}[[t]]$ and $I = t\mathbb{Z}[[t]]$. For each $m$, write $A_m$ for

$$A[\xi_1, \ldots, \xi_m]\hat{\ },$$

namely the $I$-adic completion of the polynomial ring $A[\xi_1, \ldots, \xi_m]$, and put $\mathcal{A} = (A_m)_{m \in \mathbb{N}}$.

**Definition 1.2.8** (Analytic structure)**.** Let $K$ be a valued field. An analytic $\mathcal{A}$-structure on $K$ is a collection of ring homomorphisms

$$\sigma_m : A_m \to \text{ ring of } \mathcal{O}_K\text{-valued functions on } \mathcal{O}_K^m$$

for all $m \geq 0$ such that:
  (1) $I \subset \sigma_0^{-1}(\mathcal{M}_K)$,
  (2) $\sigma_m(\xi_i) =$ the $i$-th coordinate function on $\mathcal{O}_K^m, i = 1, ..., m$,
  (3) $\sigma_{m+1}$ extends $\sigma_m$ where we identify in the obvious way functions on $\mathcal{O}_K^m$ with functions on $\mathcal{O}_K^{m+1}$ that do not depend on the last coordinate.

Let us expand the example of $A = \mathbb{Z}[[t]]$, equipped with the $t$-adic topology. For any field $k$, the natural $\mathcal{L}_{\text{DP}}$-structure on $k((t'))$ with the $t'$-adic valuation has a unique $\mathcal{A}$-structure if one fixes $\sigma_0(t)$ (in the maximal ideal, as required by (1)). Likewise, for any finite field extension $K$ of $\mathbb{Q}_p$, for any prime $p$, say, with a chosen uniformizer $\varpi_K$ of $\mathcal{O}_K$ so that $\overline{\text{ac}}$ is also fixed, the natural $\mathcal{L}_{\text{DP}}$-structure has a unique $\mathcal{A}$-structure up to choosing $\sigma_0(t)$ (in the maximal ideal).

**Definition 1.2.9.** The $\mathcal{A}$-analytic language $\mathcal{L}_\mathcal{A}$ is defined as $\mathcal{L}_{\text{DP}} \cup (A_m)_{m \in \mathbb{N}}$. An $\mathcal{L}_\mathcal{A}$-structure is an $\mathcal{L}_{\text{DP}}$-structure which is equipped with an analytic $\mathcal{A}$-structure. Let $\mathcal{T}_\mathcal{A}$ be the theory $H_{\overline{\text{ac}},0}$ together with the axioms of such $\mathcal{L}_\mathcal{A}$-structures.

**Theorem 1.2.10** ([CL])**.** *The theory $\mathcal{T}_\mathcal{A}$ has property $(*)$. It does so in a resplendent way (namely, also expansions as in Example 1.2.7 have property $(*)$ ). If $A = \mathbb{Z}[[t]]$ with ideal $I = t\mathbb{Z}[[t]]$, then it also has property $(**)$ and every definable subset in the residue field sort is definable in the language of rings.*

*Proof.* By Theorem 3.2.5 of [CL], there is a separated analytic structure $\mathcal{A}'$ such that $\mathcal{L}_{\mathcal{A}'}$ is a natural definitial expansion of $\mathcal{L}_\mathcal{A}$, with natural corresponding theory $\mathcal{T}_{\mathcal{A}'}$, specified in [CL]. Now property $(*)$ follows from Theorem 6.3.7 of [CL11] for $\mathcal{T}_{\mathcal{A}'}$ (even resplendently). The statements when $A = \mathbb{Z}[[t]]$ and $I = t\mathbb{Z}[[t]]$ are clear (that every definable subset in the residue field sort is definable in the language of rings follows from quantifier elimination for $\mathcal{L}_{\mathcal{A}'}$ of Theorem 6.3.7 of [CL11]). □

Note that Theorem 1.2.10 includes Example 1.2.7 as a special case by taking $A = \mathbb{Z}$ with $I$ the zero ideal. Other examples of analytic theories that have property $(*)$ can be found in [CL], see also Section 4.4 of [CL11].

**Remark 1.2.11.** Let us highlight some of the differences with the uniform analytic structure from [Dri92]. In [Dri92], a variant of Definition 1.2.8 of analytic $\mathcal{A}$-structures is given which is slightly more stringent, see Definition (1.7) of [Dri92]. With this notion of (1.7), van den Dries proves quantifier elimination (resplendently) in Theorem (3.9) of [Dri92], which implies that the theory is split (see Definition 1.2.3 above). However, more recent work is needed in order to prove the Jacobian property and finite *b*-minimality (see Definitions 1.2.2 and 1.2.4), and that is done in [CL11], Theorem 6.3.7, for separated analytic structures. A reduction (with an expansion by definition) from an analytic $\mathcal{A}$-structure (as in Definition 1.2.8) to a separated analytic structure (as in [CL11]) is given in [CL].

## 1.3 Rational constructible motivic functions

### 1.3.1

We introduce rational constructible motivic functions and their motivic integrals, as a variant of the construction of motivic integration in [CL15]. We will use this variant to prove Theorem 1.1.2 and its generalizations 1.4.1, 1.4.2.

Let us fix a theory $\mathcal{T}$ (in a language $\mathcal{L}$) with property (∗). From Section 1.3.3 on, we will assume that $\mathcal{T}$ also has property (∗∗), to enable to specialize to local fields of large residue field characteristic, by logical compactness.

Up to Section 1.3.1, we recall terminology from [CL15]. From Section 1.3.1 on, we introduce our variant of rational constructible motivic functions.

**The category of definable subsets**

By a $\mathcal{T}$-field we mean a valued field $K$ with residue field $k_K$ and value group $\mathbb{Z}$, equipped with an $\mathcal{L}$-structure so that it becomes a model of $\mathcal{T}$. (For set-theoretical reasons, one may want to restrict this notion to valued fields $K$ living in a very large set, or, to consider the class of all $\mathcal{T}$-fields.)

For any integers $n, m, r \geq 0$, we denote by $h[n, m, r]$ the functor sending a $\mathcal{T}$-field $K$ to
$$h[n, m, r](K) := K^n \times k_K^m \times \mathbb{Z}^r$$

Here, the convention is that $h[0, 0, 0]$ is the definable subset of the singleton $\{0\}$, i.e. $h[0, 0, 0](K) = \{0\}$.

We call a collection of subsets $X(K)$ of $h[n, m, r](K)$ for all $\mathcal{T}$-fields $K$ a definable subset if there exists an $\mathcal{L}$-formula $\phi(x)$ with free variables $x$ corresponding to elements of $h[m, n, r]$ such that

$$X(K) = \{x \in h[m, n, r](K) | \phi(x) \text{ holds in the } \mathcal{L}\text{-structure } (K, k_K, \mathbb{Z})\}$$

for all $\mathcal{T}$-fields $K$.

A definable morphism $f : X \to Y$ between two definable subsets $X, Y$ is given by a definable subset $G$ such that $G(K)$ is the graph of a function $X(K) \to Y(K)$ for all $\mathcal{T}$-fields $K$.

Denote by $\mathrm{Def}(\mathcal{T})$ (or simply Def) the category of definable subsets with definable morphisms as morphisms. If $Z$ is a definable subset, we denote by $\mathrm{Def}_Z(\mathcal{T})$ (or simply $\mathrm{Def}_Z$) the category of definable subsets $X$ with a specified definable morphism $X \to Z$; a morphism between $X, Y \in \mathrm{Def}_Z$ is a definable morphism $X \to Y$ which makes a commutative diagram with the specified morphisms $X \to Z$ and $Y \to Z$. To indicate that we work over $Z$ for some $X$ in $\text{Def}_Z$, we will often write $X_{/Z}$.

For every morphism $f : Z \to Z'$ in Def, by composition with $f$, we can define a functor
$$f_! : \text{Def}_Z \to \text{Def}_{Z'}$$
sending $X_{/Z}$ to $X_{/Z'}$. Using the fiber product, we can define a functor
$$f^* : \text{Def}_{Z'} \to \text{Def}_Z$$
by sending $Y_{/Z'}$ to $(Y \otimes_{Z'} Z)_{/Z}$.

When $Y$ and $Y'$ are definable sets, we write $Y \times Y'$ for their Cartesian product. We also write $Y[m, n, r]$ for the product $Y \times h[m, n, r]$.

By a point on a definable subset $X$, we mean a tuple $x = (x_0, K)$ where $K$ is a $\mathcal{T}$-field and $x_0 \in X(K)$. We write $|X|$ for the collection of all points that lie on $X$.

**Constructible Presburger functions**

We follow [CL15, Section 5, 6]. Consider a formal symbol $\mathbb{L}$ and the ring
$$\mathbb{A} := \mathbb{Z}[\mathbb{L}, \mathbb{L}^{-1}, \{\frac{1}{1 - \mathbb{L}^{-i}}\}_{(i>0)}]$$

For every real number $q > 1$, there is a unique morphism of rings $\vartheta_q : \mathbb{A} \to \mathbb{R}$ mapping $\mathbb{L}$ to $q$, and it is obvious that $\vartheta_q$ is injective for $q$ transcendental. Define a partial ordering on $\mathbb{A}$ by setting $a \geq b$ if for every real number with $q > 1$ one has $\vartheta_q(a) \geq \vartheta_q(b)$. We denote by $\mathbb{A}_+$ the set $\{a \in \mathbb{A} | a \geq 0\}$.

**Definition 1.3.1.** Let $S$ be a definable subset in Def. The ring $\mathcal{P}(S)$ of constructible Presburger functions on $S$ is the subring of the ring of functions $|S| \to \mathbb{A}$ generated by all constant functions $|S| \to \mathbb{A}$, by all functions $\widehat{\alpha} : |S| \to \mathbb{A}$ corresponding to a definable morphism $\alpha : S \to h[0, 0, 1]$, and by all functions $\mathbb{L}^{\widehat{\beta}} : |S| \to \mathbb{A}$ corresponding to a definable morphism $\beta : S \to h[0, 0, 1]$. We denote by $\mathcal{P}_+(S)$ the semi-ring consisting of functions in $\mathcal{P}(S)$ wich take values in $\mathbb{A}_+$. Let $\mathcal{P}_+^0(S)$ be the sub-semi-ring of $\mathcal{P}_+(S)$ generated by the characteristic functions $\mathbf{1}_Y$ of definable subsets $Y \subset S$ and by the constant function $\mathbb{L} - 1$.

If $Z \to Y$ is a morphism in Def, composition with $f$ yields a natural pullback morphism $f^* : \mathcal{P}(Y) \to \mathcal{P}(Z)$ with restrictions $f^* : \mathcal{P}_+(Y) \to \mathcal{P}_+(Z)$ and $f^* : \mathcal{P}_+^0(Y) \to \mathcal{P}_+^0(Z)$.

**Rational constructible motivic functions**

Definition 1.3.2 is taken from [CL08] [CL15]. Right after this, we start our further localizations.

**Definition 1.3.2.** Let $Z$ be a definable subset in Def. Define the semi-group $\mathcal{Q}_+(Z)$ as the quotient of the free abelian semigroup over symbols $[Y]$ with $Y_{/Z}$ a definable subset of $Z[0, m, 0]$ with the projection to $Z$, for some $m \geq 0$, by relations
  (1) $[\emptyset \to Z] = 0$,

(2) $[Y] = [Y']$ if $Y \to Z$ is isomorphic to $Y' \to Z$,

(3) $[(Y \cup Y')] + [(Y \cap Y')] = [Y] + [Y']$ for $Y$ and $Y'$ definable subsets of a common $Z[0, m, 0] \to Z$ for some $m$.

The Cartesian fiber product over $Z$ induces a natural semi-ring structure on $\mathcal{Q}_+(Z)$ by setting

$$[Y] \times [Y'] = [Y \otimes_Z Y']$$

Now let $\mathcal{Q}_+^*(Z)$ be the sub-semi-ring of $\mathcal{Q}_+(Z)$ given by

$$\{[Y \xrightarrow{f} Z] \in \mathcal{Q}_+(Z) \mid \forall x \in Z, \ f^{-1}(x) \neq \emptyset\}.$$

Then, $\mathcal{Q}_+^*(Z)$ is a multiplicatively closed set of $\mathcal{Q}_+(Z)$. So, we can consider the localization $\tilde{\mathcal{Q}}_+(Z)$ of $\mathcal{Q}_+(Z)$ with respect to $\mathcal{Q}_+^*(Z)$.

Note that if $f : Z_1 \to Z_2$ is a morphism in Def then we have natural pullback morphisms:

$$f^* : \mathcal{Q}_+(Z_2) \to \mathcal{Q}_+(Z_1)$$

by sending $[Y] \in \mathcal{Q}_+(Z_2)$ to $[Y \otimes_{Z_2} Z_1]$ and

$$f^* : \tilde{\mathcal{Q}}_+(Z_2) \to \tilde{\mathcal{Q}}_+(Z_1)$$

by sending $\dfrac{[Y]}{[Y']} \in \tilde{\mathcal{Q}}_+(Z_2)$ to $\dfrac{[Y \otimes_{Z_2} Z_1]}{[Y' \otimes_{Z_2} Z_1]}$. One easily checks that these are well-defined. We write $\mathbb{L}$ for the class of $Z[0, 1, 0]$ in $\mathcal{Q}_+(Z)$, and, in $\tilde{\mathcal{Q}}_+(Z)$.

**Definition 1.3.3.** Let $Z$ be in Def. Using the semi-ring morphism $\mathcal{P}_+^0(Z) \to \mathcal{Q}_+(Z)$ which sends $\mathbf{1}_Y$ to $[Y]$ and $\mathbb{L} - 1$ to $\mathbb{L} - 1$, the semi-ring $\mathcal{C}_+(Z)$ is defined as follows in [CL15, Section 7.1]:

$$\mathcal{C}_+(Z) = \mathcal{P}_+(Z) \otimes_{\mathcal{P}_+^0(Z)} \mathcal{Q}_+(Z).$$

Elements of $\mathcal{C}_+(Z)$ are called (nonnegative) constructible motivic functions on $Z$. In the same way, we define the semi-ring of rational (nonnegative) constructible motivic functions as

$$\tilde{\mathcal{C}}_+(Z) = \mathcal{P}_+(Z) \otimes_{\mathcal{P}_+^0(Z)} \tilde{\mathcal{Q}}_+(Z),$$

by using the semi-ring morphism $\mathcal{P}_+^0(Z) \to \tilde{\mathcal{Q}}_+(Z)$ which sends $\mathbf{1}_Y$ to $[Y]$ and $\mathbb{L} - 1$ to $\mathbb{L} - 1$.

If $f : Z \to Y$ is a morphism in Def then there is a natural pullback morphism from [CL15, Section 7.1]:

$$f^* : \mathcal{C}_+(Y) \to \mathcal{C}_+(Z)$$

sending $a \otimes b$ to $f^*(a) \otimes f^*(b)$, where $a \in \mathcal{P}_+(Y)$ and $b \in \mathcal{Q}_+(Y)$. Likewise, we have the pullback morphism:

$$f^* : \tilde{\mathcal{C}}_+(Y) \to \tilde{\mathcal{C}}_+(Z)$$

sending $a \otimes b$ to $f^*(a) \otimes f^*(b)$, where $a \in \mathcal{P}_+(Y)$ and $b \in \tilde{\mathcal{Q}}_+(Y)$.

Since $\mathcal{T}$ is split, the canonical morphism

$$\mathcal{P}_+(Z[0,0,r]) \otimes_{\mathcal{P}_+^0(Z)} \tilde{\mathcal{Q}}_+(Z[0,m,0]) \to \tilde{\mathcal{C}}_+(Z[0,m,r]) \qquad (1.3.1.1)$$

is an isomorphism of semi-rings, where we view $\mathcal{P}_+(Z[0,0,r])$ and $\tilde{\mathcal{Q}}_+(Z[0,m,0])$ as $\mathcal{P}_+^0(Z)$-modules via the homomorphisms $p^* : \mathcal{P}_+^0(Z) \to \mathcal{P}_+(Z[0,0,r])$ and $q^* : \mathcal{P}_+^0(Z) \to \tilde{\mathcal{Q}}_+(Z[0,m,0])$ that come from the pullback homomorphism of the two projections $p : Z[0,0,r] \to Z$ and $q : Z[0,m,0] \to Z$.

**Proposition 1.3.4.** *For $F$ in $\tilde{\mathcal{C}}_+(\mathbb{Z}^r)$ there exist $[Y]$ in $\mathcal{Q}_+^*(h[0,0,0])$ and $G \in \mathcal{C}_+(\mathbb{Z}^r)$ such that one has the equality*

$$\overline{G} = \overline{[Y]} \cdot F$$

*in $\tilde{\mathcal{C}}_+(\mathbb{Z}^r)$, with $\overline{G}$ the image of $G$ under the natural map $\mathcal{C}_+(\mathbb{Z}^r) \to \tilde{\mathcal{C}}_+(\mathbb{Z}^r)$ and similarly for $\overline{[Y]}$.*

*Proof.* This follows directly from the isomorphism from (1.3.1.1) and the definition of $\tilde{\mathcal{Q}}_+(h[0,0,0])$. □

### 1.3.2 Integration of rational constructible motivic functions

In the next three sections, we give the definitions of rational constructible motivic functions and their integrals, which follows the same construction as for integration of constructible motivic functions in [CL15], similarly using property $(*)$.

**Integration over the residue field**

We adapt [CL15, Section 6.2] to our setting. Suppose that $Z \subset X[0,k,0]$ for some $k \geq 0$ and $a \in \tilde{\mathcal{Q}}_+(Z)$, we write $a = \dfrac{[Y]}{[Y']}$ for some $[Y \xrightarrow{f} Z] \in \mathcal{Q}_+(Z)$ and $[Y' \xrightarrow{f'} Z] \in \mathcal{Q}_+^*(Z)$. We write $\mu_{/X}$ for the corresponding formal integral in the fibres of the coordinate projection $Z \to X$

$$\mu_{/X} : \tilde{\mathcal{Q}}_+(Z) \to \tilde{\mathcal{Q}}_+(X), \frac{[Y]}{[Y']} \mapsto \frac{[Y]}{[Y'']}$$

where $[Y'' \to X] = [(Y' \sqcup (X \setminus \mathrm{Im} f')) \to X]$, where $\sqcup$ denotes the disjoint union (a disjoint union of definable sets can be realized as a definable set by using suitable piecewise definable bijections). Note that $Y''$ is built from $Y$ by using that the class of a definable singleton is the multiplicative unit to preserve the property that the fibers over $Z$ are never empty.

**Integration over the value group**

This section follows [CL15, Section 5] and will be combined with the integration from section 1.3.2 afterwards. Let $Z \in \mathrm{Def}$ and $f \in \mathcal{P}(Z[0,0,r])$. For any $\mathcal{T}$-field $K$ and any $q > 1$ we write $\vartheta_{q,K}(f) : Z(K) \to \mathbb{R}$ for the function sending $z \in Z(K)$ to $\vartheta_q(f(z,K))$. Then $f$ is called *$Z$-integrable* if for each $\mathcal{T}$-field $K$, each $q > 1$ and for each $z \in Z(K)$, the family $(\vartheta_{q,K}(f)(z,i))_{i \in \mathbb{Z}^r}$ is summable. The collection of $Z$-integrable functions in $\mathcal{P}(Z[0,0,r])$ is denoted by $I_Z\mathcal{P}(Z[0,0,r])$ and $I_Z\mathcal{P}_+(Z[0,0,r])$ is the collection of $Z$-integrable functions in $\mathcal{P}_+(Z[0,0,r])$.

We recall from [CL15, Theorem-Definition 5.1] that for each $\phi \in I_Z\mathcal{P}(Z[0,0,r])$, there exists a unique function $\varphi := \mu_{/Z}(\phi)$ in $\mathcal{P}_Z$ such that for all $q > 1$, all $\mathcal{T}$-fields $K$, all $z \in Z(K)$, one has

$$\vartheta_{q,K}(\varphi)(z) = \sum_{i \in \mathbb{Z}^r} \vartheta_{q,K}(\phi)(z,i)$$

and the mapping $\phi \mapsto \mu_{/Z}(\phi)$ yields a morphism of $\mathcal{P}(Z)$-modules

$$\mu_{/Z} : I_Z\mathcal{P}(Z \times \mathbb{Z}^r) \to \mathcal{P}(Z).$$

**Integration over one valued field variable**

We first follow [CL15, Section 8] and then use it for our setting of rational motivic constructible functions in Lemma-Definiton 1.3.8 below. For a ball $B = a + b\mathcal{O}_K$ and any real number $q > 1$, we call $\vartheta_q(B) := q^{-\operatorname{ord} b}$ the $q$-volume of $B$. A finite or countable collection of disjoint balls in $K$, each with different $q$-volume is called a step-domain; we will identify a step-domain $S$ with the union of the balls in $S$. Recall from [CL15] that a nonnegative real valued function $\varphi : K \to \mathbb{R}_{\geq 0}$ is a step-function if there exists a unique step-domain $S$ such that $\varphi$ is constant and nonzero on each ball of $S$ and zero outside $S \cup \{a\}$ for some $a \in K$.

Let $q > 1$ be a real number. A step-function $\varphi : K \to \mathbb{R}_{\geq 0}$ with step-domain $S$ is $q$-integrable over $K$ if and only if

$$\sum_{B \in S} \vartheta_q(B).\varphi(B) < \infty$$

Suppose that $Z = X[1,0,0]$ for some $X \in \operatorname{Def}$ and $\varphi \in \mathcal{P}_+(Z)$. We call $\varphi$ an $X$-integrable family of step-functions if for each $\mathcal{T}$-field $K$, for each $x \in X(K)$ and for each $q > 1$, the function

$$\vartheta_{q,K}(\varphi)(x,.) : K \to \mathbb{R}_{\geq 0}, t \mapsto \vartheta_{q,K}(\varphi)(x,t),$$

is a step-function which is $q$-integrable over $K$. For such $\varphi$ there exists a unique function $\phi$ in $\mathcal{P}_+(X)$ such that $\vartheta_{q,K}(\phi)(x)$ equals the $q$-integral over $K$ of $\vartheta_{q,K}(\varphi)(x,.)$ for each $\mathcal{T}$-field $K$, each $x \in X(K)$ and each $q > 1$ (see [CL15], Lemma-Definition 8.1); we write

$$\mu_{/X}(\varphi) := \phi,$$

the integral of $\varphi$ in the fibers of $Z \to X$.

**Lemma 1.3.5.** *Let $\phi \in \mathcal{Q}_+^*(Z)$ and suppose that $Z = X[1,0,0]$. Then there exist a partition $Z = \sqcup_{i=1}^l Z_i$, a natural number $n \in \mathbb{N}$ and for each $1 \leq i \leq l$ there exist a definable injection $\lambda_i : Z_i \to Z[0,n,0]$ and $\phi_i \in \mathcal{Q}_+^*(X[0,n,0])$ such that*

$$\mathbf{1}_{Z_i}\phi = \mathbf{1}_{Z_i}\lambda_i^* p^*(\phi_i)$$

*and $\tilde{\pi} \circ \lambda_i : Z_i \to Z_i$ is the identity on $Z_i$, where $p$ is the projection from $Z[0,n,0]$ to $X[0,n,0]$ and $\tilde{\pi}$ is the projection from $Z[0,n,0]$ to $Z$.*

*Proof.* The existence of $n, r \in \mathbb{N}$, $Z_i$, $\lambda_i : Z_i \to Z[0, n, r]$ and $\phi_i \in \mathcal{Q}_+(X[0, n, r])$ such that $\mathbf{1}_{Z_i}\phi = \mathbf{1}_{Z_i}\lambda_i^* p^*(\phi_i)$ is a consequence of cell decomposition theorem for b-minimal theory as in [CL07] and in section 7 of [CL08] (remark that a finitely b-minimal theory is in particular b-minimal, see corollary 3.8 in [CL15]). Because $\mathcal{T}$ is split then we can choose $r = 0$. If $Z_i$ and $\phi_i$ have been given, then we can take $\tilde{\phi}_i(x) = \phi_i(x)$ if $x \in p(\lambda_i(Z_i))$ and $\tilde{\phi}_i(x) = \{*\}$ if $x \in X \setminus p(\lambda_i(Z_i))$ for a definable point $*$ in residue field. Then by definition of $\mathcal{Q}_+^*(Z)$ and $\mathcal{Q}_+^*(X[0, n, 0])$ we see that $\phi_i \in \mathcal{Q}_+^*(X[0, n, 0])$ and $\{Z_i, \tilde{\phi}_i\}$ will satisfy the lemma. □

**Lemma-Definition 1.3.6.** Suppose that $Z = X[1, 0, 0]$. Let $\psi = \sum_{i=1}^k a_i \otimes \pi^*(b_i) \in \mathcal{P}_+(Z) \otimes \pi^* \tilde{\mathcal{Q}}_+(X) \subset \tilde{\mathcal{C}}_+(Z)$, where $\pi : Z \to X$ be the projection, $a_i \in \mathcal{P}_+(Z)$ and $b_i \in \tilde{\mathcal{Q}}_+(X)$. We say that $\psi$ is $X$-integrable if $a_i$ is $X$-integrable for $1 \leq i \leq k$ and then

$$\mu_{/X}(\psi) := \sum_{i=1}^k \mu_{/X}(a_i) \otimes b_i \in \tilde{\mathcal{C}}_+(X)$$

is independent of the choices and is called the integral of $\psi$ in the fibers of $Z \to X$.

*Proof.* Since $\mathcal{T}$ has property (*), the proof is similar with the Lemma- Definition 7.6 of [CL15]. □

**Lemma 1.3.7.** *Let $\psi \in \tilde{\mathcal{C}}_+(Z)$ and suppose that $Z = X[1, 0, 0]$. Then there exists $\phi$ in $\mathcal{P}_+(Z[0, m, 0]) \otimes \pi^* \tilde{\mathcal{Q}}_+(X[0, m, 0]) \subset \tilde{\mathcal{C}}_+(Z[0, m, 0])$ such that $\mu_{/Z}(\phi) = \psi$, where $\pi : Z[0, m, 0] \to X[0, m, 0]$ be the projection.*

*Proof.* We write $\psi = \dfrac{\theta}{\varphi}$ with $\theta \in \mathcal{C}_+(Z)$ and $\varphi \in \mathcal{Q}_+(Z)^*$. By lemma 7.8 of [CL15], there exist $m \in \mathbb{N}$ and $v \in \mathcal{P}_+(Z[0, m, 0])$ such that $\mu_{/Z}(v) = \theta$. Let $n, l, Z_i, \lambda_i, p, \tilde{\pi}, \phi_i$ be as in the lemma 1.3.5. If $m' > n'$ we can view $Z[0, n', 0]$ as the set $Z[0, n', 0] \times \{*\} \subset Z[0, m', 0]$ and $X[0, n', 0] \times \{*\} \subset X[0, m', 0]$ with $*$ a definable point in $h[0, m' - n', 0]$. So by fixing a definable point in $h[0, |m - n|, 0]$, we can assume that $m = n$. For each $1 \leq i \leq l$ we can extend $\lambda_i$ to a definable map $\tilde{\lambda}_i : Z \to Z[0, n, 0]$. We set

$$\phi = \sum_{i=1}^l \frac{\mathbf{1}_{\pi^{-1}(Z_i)} v}{\tilde{\pi}^* \tilde{\lambda}_i^* p^*(\phi_i)} = \sum_{i=1}^l \frac{\mathbf{1}_{\pi^{-1}(Z_i)} v}{\pi^*(\phi_i)} \in \mathcal{P}_+(Z[0, m, 0]) \otimes \pi^* \tilde{\mathcal{Q}}_+(X[0, m, 0]).$$

Now we have

$$\mu_{/Z}(\phi) = \mu_{/Z}\left(\sum_{i=1}^l \frac{\mathbf{1}_{\pi^{-1}(Z_i)} v}{\tilde{\pi}^* \tilde{\lambda}_i^* p^*(\phi_i)}\right)$$

$$= \sum_{i=1}^l \frac{\mu_{/Z}(\mathbf{1}_{\pi^{-1}(Z_i)} v)}{\tilde{\lambda}_i^* p^*(\phi_i)}$$

$$= \sum_{i=1}^l \frac{\mathbf{1}_{Z_i} \theta}{\tilde{\lambda}_i^* p^*(\phi_i)}$$

$$= \sum_{i=1}^l \frac{\mathbf{1}_{Z_i} \theta}{\phi}$$

$$= \psi.$$

Here, we remark that $\mu_{/Z}(v) = \theta$ implies $\mu_{/Z}(\mathbf{1}_{\pi^{-1}(Z_i)} v) = \mathbf{1}_{Z_i} \theta$. □

**Lemma-Definition 1.3.8.** Suppose that $Z = X[1,0,0]$.

Let $\psi \in \tilde{\mathcal{C}}_+(Z)$. We say that $\psi$ is $X$-integrable if there exists a $\phi$ in $\mathcal{P}_+(Z[0,m,0]) \otimes \pi^*\tilde{\mathcal{Q}}_+(X[0,m,0]) \subset \tilde{\mathcal{C}}_+(Z[0,m,0])$ with $\mu_{/Z}(\phi) = \psi$ such that $\phi$ is $X[0,m,0]$-integrable and then

$$\mu_{/X}(\psi) := \mu_{/X}(\mu_{/X[0,m,0]}(\phi)) \in \tilde{\mathcal{C}}_+(X)$$

is independent of the choices and is called the integral of $\varphi$ in the fibers of $Z \to X$.

*Proof.* The proof is similar with the discussion in section 9 of [CL08] and the Lemma-definition 8.2 of [CL15]. □

**Integration of rational constructible motivic functions in the general case**

Combining the three cases above, we define integrability and the integral $\mu_{/X}(\varphi)$ of an integrable rational constructible motivic function $\varphi \in \tilde{\mathcal{C}}_+(X[m,n,r])$ by Tonelli-Fubini iterated integration in a similar way as in Lemma-Definition 9.1 of [CL15]. More precisely, we will define the integrals in the fibers of a general coordinate projection $X[n,m,r] \to X$ by induction on $n \geq 0$.

First of all, based on (1.3.1.1) and the definition of integrability, for each $Z$-integrable function $\varphi \in \tilde{\mathcal{C}}_+(Z[0,m,r])$, one can write $\varphi = \sum_{i=1}^k a_i \otimes b_i$, where $a_i \in I_Z\mathcal{P}_+(Z[0,0,r])$ and $b_i \in \tilde{\mathcal{Q}}_+(Z[0,m,0])$ and

$$\mu_{/Z}(\varphi) = \sum_{i=1}^k \mu_{/Z}(a_i) \otimes \mu_{/Z}(b_i).$$

**Lemma-Definition 1.3.9.** Let $\varphi$ be in $\tilde{\mathcal{C}}_+(Z)$ and suppose that $Z = X[n,m,r]$ for some $X$ in Def.

If $n = 0$ we say that $\varphi$ is $X$-integrable if and only if $\varphi$ is $X[0,m,0]$-integrable. If this holds then

$$\mu_{/X}(\varphi) := \mu_{/X}(\mu_{/X[0,m,0]}(\varphi)) \in \tilde{\mathcal{C}}_+(X)$$

is called the integral of $\varphi$ in the fibers of $Z \to X$.

If $n \geq 1$, we say that $\varphi$ is $X$-integrable if there exists a definable subset $Z' \subset Z$ whose complement in $Z$ has relative dimension $< n$ over $X$ such that $\varphi' := \mathbf{1}_{Z'}\varphi$ is $X[n-1,m,r]$-integrable and $\mu_{/X[n-1,m,r]}(\varphi')$ is $X$-integrable. If this holds then

$$\mu_{/X}(\varphi) := \mu_{/X}(\mu_{/X[n-1,m,r]}(\varphi')) \in \tilde{\mathcal{C}}_+(X)$$

does not depend on the choices and is called the integral of $\varphi$ in the fibers of $Z \to X$.

Slightly more generally, let $\varphi \in \tilde{\mathcal{C}}_+(Z)$ and suppose that $Z \subset X[n,m,r]$. We say that $\varphi$ is $X$-integrable if the extension by zero of $\varphi$ to a function $\tilde{\varphi} \in \tilde{\mathcal{C}}_+(X[n,m,r])$ is $X$-integrable and we define $\mu_{/X}(\varphi)$ as $\mu_{/X}(\tilde{\varphi})$.

*Proof.* Since $\mathcal{T}$ has property $(*)$ the proof is similar to the proof for Lemma-Definition 9.1 of [CL15]. □

### 1.3.3 Interpretation of rational constructible motivic functions in non-archimedean local fields

In this section we show how rational constructible motivic functions can be specialized to real valued functions on local fields of large residue field characteristic, in the spirit of the specializations in [CL05] and Proposition 9.2 of [CL15]. Importantly, taking motivic integrals combines well with this specialization and integration over local fields.

Let $\mathcal{T}$ be a theory in a language $\mathcal{L}$ extending $\mathcal{L}_{\mathrm{DP}}$. Suppose that $\mathcal{T}$ has properties $(*)$ and $(**)$ from section 1.2.3. For a definable set $X$, a definable function $f$ and a rational motivic constructible function $\phi$, the objects $X_K = X(K)$, $f_K$ and $\phi_K$ make sense for every local field $K$ with large residue field characteristic. We make this explicit for rational motivic constructible functions, where we assume $K$ to be a local field with large residue field characteristic (depending on the data).

— For $a \in \mathcal{P}_+(X)$, we get $a_K : X_K \to \mathbb{Q}_{\geq 0}$ by replacing $\mathbb{L}$ by $q_K$.
— For $b = \dfrac{[Y]}{[Y']} \in \tilde{\mathcal{Q}}_+(X)$ with $Y$ a definable subset of $X[0, m, 0]$, $Y'$ a definable subset of $X[0, m', 0]$ such that $[Y'] \in \mathcal{Q}_+^*(X)$, if we write $p : Y \to X$, $p' : Y' \to X$ for the projections, one defines $b_K : X_K \to \mathbb{Q}_{\geq 0}$ by sending $x \in X_K$ to $\dfrac{\#(p_K^{-1}(x))}{\#((p'_K)^{-1}(x))}$
— For $\phi \in \mathcal{C}_+(X)$ or $\phi \in \tilde{\mathcal{C}}_+(X)$, writing $\phi$ as a finite sum $\sum_i a_i \otimes b_i$ with $a_i \in \mathcal{P}_+(X)$ and $b_i \in \mathcal{Q}_+(X)$ or $b_i \in \tilde{\mathcal{Q}}_+(X)$, we get the function

$$\phi_K : X_K \to \mathbb{Q}_{\geq 0}, x \mapsto \sum_i a_{iK}(x).b_{iK}(x),$$

which does not depend on the choices made for $a_i$ and $b_i$.

Taking motivic integrals commutes with taking specializations, as follows.

**Proposition 1.3.10.** *Let $\varphi$ be an $X$-integrable rational constructible motivic function in $\tilde{\mathcal{C}}_+(X[m, n, r])$ and let $\mu_{/X}(\varphi)$ be its motivic integral, in the fibers of the projection $X[m, n, r] \to X$. Then there exists $M > 0$ such that for all local fields $K$ whose residue field has characteristic at least $M$ one has for each $x \in X_K$*

$$\left(\mu_{/X}(\varphi)\right)_K(x) = \int_{y \in K^m \times k_K^n \times \mathbb{Z}^r} \varphi_K(x, y) |dy|$$

*where one puts the normalized Haar measure on $K$, the counting measure on $k_K$ and on $\mathbb{Z}$, and the product measure $|dy|$ on $K^m \times k_K^n \times \mathbb{Z}^r$.*

*Proof.* This follows naturally from the corresponding result for $\mathcal{C}_+$ instead of $\tilde{\mathcal{C}}_+$ (see [CL05] and Proposition 9.2 of [CL15]), and the concrete definitions of $\tilde{\mathcal{C}}_+$ and their integration . □

## 1.4 The uniform rationality for Poincaré series of definable equivalence relations

### 1.4.1

We will prove Theorem 1.4.1, which is a slight generalization of the Main Theorem 1.1.2.

Let $\mathcal{T}$ be a theory in a language $\mathcal{L}$ extending $\mathcal{L}_{\mathrm{DP}}$. Suppose that $\mathcal{T}$ has properties $(*)$ and $(**)$ from section 1.2.3. Let $\varphi(x, y, n)$ be an $\mathcal{L}$-formula with free variables $x$ running over $K^m$, $y$ running over $K^m$ and $n$ running over $\mathbb{N}$. Suppose that for each local field $K$ and each $n$, $\varphi(x, y, n)$ gives an equivalence relation $\sim_{K,n}$ on $K^m$ with finitely many, say, $a_{\varphi,K,n}$, equivalence classes. (The situation that $\varphi(x, y, n)$ gives an equivalence relation on a definable subset $X_{K,n}$ of $K^m$ for each $K$ and each $n$ is similar.) For each local field $K$ put

$$P_{\varphi,K}(T) = \sum_{n \geq 0} a_{\varphi,K,n} T^n.$$

**Theorem 1.4.1.** *There exists $M > 0$ such that the power series $P_{\varphi,K}(T)$ is rational in $T$ for each local field $K$ whose residue field has characteristic at least $M$. Moreover, for such $K$, the series $P_{\varphi,K}(T)$ only depends on the $\mathcal{L}$-structure induced on the residue field sort $k_K$.*

*More precisely, there exist nonnegative integers $N, k, b_j, e_i$, integers $a_j$ and $\mathcal{L}$-formulas $X_i$ and $Y$ for subsets of some power of the residue field for $i = 0, \ldots, N$ and $j = 0, \ldots, k$, such that for each $j$, $a_j$ and $b_j$ are not both $0$, and, for all local fields $K$ with residue field of characteristic at least $M$, $Y(k_K)$ is nonempty and*

$$P_{\varphi,K}(T) = \frac{\sum_{i=0}^{N}(-1)^{e_i} \# X_i(k_K) T^i}{\# Y(k_K) \cdot \prod_{j=1}^{k}(1 - q_K^{a_j} T^{b_j})}.$$

In fact, Theorem 1.4.1 is a consequence of the following more versatile theorem:

**Theorem 1.4.2.** *Suppose that $\mathcal{T}$ has properties $(*)$ and $(**)$. Let $\varphi(x, y, z)$ be an $\mathcal{L}$-formula with free variables $x$ running over $\mathrm{VF}^n$, $y$ running over $\mathrm{VF}^n$ and $z$ running over an arbitrary $\mathcal{L}$-definable set $Z$. Let $R$ be a definable subset of $\mathrm{VF}^n \times Z$. Suppose that for each local field $K$, and each $z \in Z_K$, $\varphi(x, y, z)$ gives an equivalence relation $\sim_{K,z}$ on $R_{K,z} := \{x \in K^n \mid (x, z) \in R_K\}$ with finitely many, say, $a_{\varphi,K,z}$, equivalence classes. Then there exist a rational motivic function $F$ in $\tilde{\mathcal{C}}_+(Z)$ and a constant $M > 0$ such that for each local field $K$ whose residue field has characteristic at least $M$ one has*

$$a_{\varphi,K,z} = F_K(z).$$

*Proof of Theorem 1.4.1.* Theorem 1.4.1 follows from Theorem 1.4.2 by Proposition 1.3.4, the formula 1.3.1.1 with $m = 0$, and the rationality result Theorem 7.1 of [Pas89] for the $p$-adic case or the rationality result as in theorem 14.4.1 and Theorem 5.7.1 of [CL08] for the motivic case. □

Before giving the proof of Theorem 1.4.2 we give a few more definitions and lemmas.

### 1.4.2 Multiballs, multiboxes, and their multivolumes

We give definitions which are inspired by concepts of the appendix of [HMR].

Fix a local field $K$. Recall that $q_K$ stands for the number of elements of the residue field $k_K$ of $K$. We implicitly use an ordering of the coordinates on $\mathcal{O}_K^n$ in the following definition.

**Definition 1.4.3.** Let $n \geq 1$, $0 \leq r_i \leq +\infty$ for $i = 1, ..., n$ and let a nonempty set $Y \subset \mathcal{O}_K^n$ be given.

If $n = 1$, we consider the closed ball $Y = \mathcal{B}(x, r_1) = \{y \in K \,|\, \operatorname{ord}_K(y - x) \geq r_1\}$ with convention that $Y$ is a singleton if $r_1 = +\infty$. The volume Vol is taken for the Haar measure on $K$ such that $\mathcal{O}_K$ has measure 1 and we consider $q_K^{-\infty}$ to be zero. With the notation on volume, we have $\operatorname{Vol}(Y) = q_K^{-r_1}$. We will say that $Y$ is a multiball of multivolume $q_K^{-r_1}$.

If $n \geq 2$, then $Y$ is called a multiball of multivolume $(q_K^{-r_1}, ..., q_K^{-r_n})$ if and only if $Y$ is of the form

$$\{(x_1, ..., x_n) \mid (x_1, ..., x_{n-1}) \in A,\ x_n \in B_{x_1, ..., x_{n-1}}\},$$

where $A \subset \mathcal{O}_K^{n-1}$ is a multiball of multivolume $(q_K^{-r_1}, ..., q_K^{-r_{n-1}})$, and, for each $(x_1, \ldots, x_{n-1}) \in A$, $B_{x_1, ..., x_{n-1}} \subset O_K$ is a multiball of multivolume $q_K^{-r_n}$. The multivolume of a multiball $Y$ is denoted by $multivol(Y)$.

**Definition 1.4.4.** We put the inverse lexicographical order (the colexicographical order) on $\mathbb{R}^n$, namely, $(a_1, ..., a_n) > (b_1, ..., b_n)$ if and only if there exists $1 \leq k \leq n$ such that $a_i = b_i$ for all $i > k$ and $a_k > b_k$. By this order, we can compare multivolumes. Let $X \subset \mathcal{O}_K^n$. The multibox of $X$, denoted by $MB(X)$, is the union of the multiballs $Y$ contained in $X$ and with maximal multivolume (for the colexicographical ordering), where maximality is among all multiballs contained in $X$. We write $multivol(X)$ for $multivol(Y)$ for any multiball $Y$ contained in $X$ with maximal multivolume.

Note that taking $MB$ and taking projections does not always commute and it may be that $p(MB(X))$ and $MB(p(X))$ are different, say, with $p : \mathcal{O}_K^n \to \mathcal{O}_K^{n-1}$ the coordinate projection to the first $n - 1$ coordinates.

**Definition 1.4.5.** Fix $1 \leq m \leq n$, a set $X \subset \mathcal{O}_K^n$, and $x = (x_1, ..., x_n)$ in $MB(X)$. Set $x_{\leq m} = (x_1, ..., x_m)$ and let $MB(X, m)$ be the image of $MB(X)$ under the projection from $\mathcal{O}_K^n$ to $\mathcal{O}_K^m$. Denote by $MB(X, x, m)$ the fiber of $MB(X, m)$ over $x_{\leq m-1}$. We write

$$multinumber_m(X, x)$$

for the number of balls with maximal volume contained in $MB(X, x, m)$. Write

$$Multinumber_m(X, x)$$

for the number of balls $B$ of minimal volume with $B \cap MB(X, x, m) \neq \emptyset$ and with $B \not\subseteq MB(X, x, m)$.

Note that the number of balls in $\mathcal{O}_K$ of any fixed volume is automatically finite.

### 1.4.3 Definable equivalence relations

From now on we suppose that $\mathcal{T}$ has properties $(*)$ and $(**)$. Let $X$ be a definable set. Of course, $multinumber_m(X_K, x)$ and $Multinumber_m(X_K, x)$ may vary with $K$ and $x$. In fact, it is difficult to give a uniform in $K$ estimate for $multinumber_m(X_K, x)$ but for $Multinumber_m(X_K, x)$ we can do it, even in definable families, see Lemma 1.4.6 and its corollary.

**Lemma 1.4.6.** *Let $Z$ and $X \subset \mathrm{VF} \times Z$ be $\mathcal{L}$-definable such that $X_K \subset \mathcal{O}_K \times Z_K$ for all local fields $K$ of large residue field characteristic. Then there exist positive integers $M$ and $Q$ such that, for all local fields $K$ with residue field characteristic at least $M$ and for all $z \in Z_K$, one has*

$$N_{K,z} \leq Q,$$

*where $N_{K,z}$ is the number of balls $B$ of minimal volume with*

$$B \cap MB(X_{K,z}) \neq \emptyset \text{ and } B \not\subseteq MB(X_{K,z}),$$

*and where $X_{K,z}$ is the set $\{x \in \mathcal{O}_K \mid (x,z) \in X_K\}$.*

*Proof.* Write $X'$ for the definable subset of $X$ such that $X'_{K,z} = MB(X_{K,z})$ for each $K$ with large residue field characteristic and each $z \in Z_K$. Since $\mathcal{T}$ is finitely b-minimal, has property $(**)$, and by logical compactness, there exist positive integers $M$, a Cartesian product $S$ of sorts not involving the valued field sort and $\mathcal{L}$-definable functions $f$, $c$, $\xi$, $t$, such that for all local fields $K$ with residue field characteristic at least $M$, we have
— $f_K : X'_K \to Z_K \times S_K$ is a function over $Z_K$ (meaning that $f_K$ makes a commutative diagram with the projections $Z_K \times S_K \to Z_K$ and $X'_K \to Z_K$);
— $c_K : Z_K \times S_K \to \mathcal{O}_K$, $\xi_K : Z_K \times S_K \to k_K$, and $t_K : Z_K \times S_K \to \mathbb{Z}$ are such that each nonempty fiber $f_K^{-1}(z,s)$, for $(z,s) \in Z_K \times S_K$, is either the singleton

$$\{c_K(z,s)\}$$

or the ball

$$\{y \in \mathcal{O}_K | \overline{\mathrm{ac}}(y - c_K(z,s)) = \xi_K(z,s), \ \mathrm{ord}(y - c_K(z,s)) = t_K(z,s)\}.$$

One derives from finite b-minimality and compactness (as in [CL15]) that there exists an integer $Q_0 > 0$ such that for for all local fields $K$ with large residue field characteristic, for each $z \in Z_K$, the range of $s \mapsto c_K(z,s)$ has no more than $Q_0$ elements. We will show that we can take $Q = Q_0$.

Suppose first that $X'_{K,z}$ is a disjoint union of balls of volume $q_K^{-\alpha(K,z)}$ where

$$q_K^{-\alpha(K,z)} = multivol(X_{K,z}).$$

Choose a ball $B$ with volume $q_K^{-\alpha(K,z)+1}$ and with $B \cap X'_{K,z} \neq \emptyset$, fix $y \in B \cap X'_{K,z}$ and write $f_K(y) = (z,s)$, so that $y$ belongs to the ball

$$B' = \{v \in \mathcal{O}_K | \overline{\mathrm{ac}}(v - c_K(z,s)) = \xi_K(z,s), \ \mathrm{ord}(v - c_K(z,s)) = t_K(z,s)\},$$

which has volume $q_K^{-t_K(z,s)-1}$. Since $X'_{K,z} = MB(X_{K,z})$, the ball

$$B(y, q_K^{-\alpha(K,z)})$$

around $y$ of volume $q_K^{-\alpha(K,z)}$ is a maximal ball contained in $X_{K,z}$. Hence, $y \in B' \subset X_{K,z}$ implies $B' \subset B(y, q_K^{-\alpha(K,z)})$. It follows that

$$t_K(z,s) + 1 \geq \alpha(K,z),$$

which proves that

$$\operatorname{ord}(y - c_K(z,s)) \geq \alpha(K,z) - 1. \tag{1.4.3.1}$$

Since $B$ has volume $q_K^{-\alpha(K,z)+1}$ and contains $y$, the inequality (1.4.3.1) implies that $c_K(z,s) \in B$ and it follows that $N_{K,z} \leq Q_0$.

If $X'_{K,z}$ is not a union of balls then it is contained in the range of $c_K(z,.) : s \mapsto c_K(z,s)$, and thus also in this case we find $0 = N_{K,z} \leq Q_0$. This shows we can take $Q = Q_0$. □

**Corollary 1.4.7.** *Let $Z$ and $X \subset \operatorname{VF}^n \times Z$ be $\mathcal{L}$-definable such that $X_K \subset \mathcal{O}_K^n \times Z_K$ for all local fields $K$ of large residue field characteristic. Fix $m$ with $1 \leq m \leq n$. Then there exist positive integers $M$ and $Q$ such that*

$$Multinumber_m(X_{K,z}, x) < Q$$

*for all local fields $K$ with residue field characteristic at least $M$, for all $z \in Z_K$ and all $x \in \mathcal{O}_K$ with $x \in MB(X_{K,z})$, where $X_{K,z}$ is the set $\{x \in \mathcal{O}_K^n \mid (x,z) \in X_K\}$.*

*Proof.* The corollary follows from Lemma 1.4.6 since the condition $x \in MB(X_{K,z})$ is an $\mathcal{L}$-definable condition on $(x,z)$. □

**Lemma 1.4.8.** *Let $Z$ and $X \subset \operatorname{VF}^n \times Z$ be $\mathcal{L}$-definable such that $X_K \subset \mathcal{O}_K^n \times Z_K$ for all local fields $K$ of large residue field characteristic. Then there exist $M > 0$ and a definable function $f : Z \to (\operatorname{VG} \cup \{+\infty\})^n$ such that*

$$(q_K^{-f_{i,K}(z)})_{i=1}^n = multivol(X_{K,z})$$

*for each $z \in Z_K$ and each local field $K$ with residue field characteristic at least $M$, where $X_{K,z}$ is the set $\{x \in \mathcal{O}_K^n \mid (x,z) \in X_K\}$.*

*Proof.* The lemma follows easily by the definability of multiboxes and of the valuative radius of the balls involved. □

**Remark 1.4.9.** We don't need $M$ in the lemma 1.4.8 because we can still define $f_K$ for residue field characteristic less than $M$.

*Proof of the main theorem 1.4.2.* Firstly, we can suppose that the sets $R_{K,z}$ for $z \in Z(K)$ are subsets of $\mathcal{O}_K^n$, up to increasing $n$ and mapping a coordinate $w \in K$ to $(w, 0) \in \mathcal{O}_K^2$ if $\operatorname{ord}(w) \geq 0$ and to $(0, w^{-1}) \in \mathcal{O}_K^2$ if $\operatorname{ord}(w) < 0$.

By Lemma 1.4.8, there is a definable function $f_{K,z} : \mathcal{O}_K^n \to (\mathbb{Z} \cup \{+\infty\})^n$ such that

$$(q_K^{-f_{K,z,i}(x)})_{i=1}^n = multivol(x/_{\sim_{K,z}}),$$

where $x/_{\sim_{K,z}}$ denotes the equivalence class of $x$ under the equivalence relation $\sim_{K,z}$.

For each subset $I$ of $\{1,...,n\}$, we set $R_{K,z,I} = \{x \in R_{K,z} | f_{K,z,i}(x) < +\infty \Leftrightarrow i \in I\}$. So, $R_{K,I} = (R_{K,z,I})_{z \in Z(K)}$ is defined by an $\mathcal{L}$-formula for each $I$, uniformly in $z$. It is easy to see that $R_{K,z} = \coprod_{I \subset \{1,...,n\}} R_{K,z,I}$ and $(x_{\sim_{K,z}} y) \wedge (x \in R_{K,z,I}) \Rightarrow y \in R_{K,z,I}$. So if we set $a_{K,z,I} = R_{K,z,I}/_{\sim_{K,z}}$ then $a_{K,z} = \sum_{I \subset \{1,...,n\}} a_{K,z,I}$. The proof will be followed by the claim for each $I \subset \{1,...,n\}$. We will consider two cases.

**Case 1: $I$ is the empty set.**

Because of the definition of $R_{K,I}$, we deduce that for each $z \in Z(K)$ the definable set $R_{K,z,I}$ will be a finite set. From the proof of lemma 1.4.6 we have that $a_{K,z,I}$ must be bounded uniformly on all local field $K$ with large residue field characteristic and $z \in Z(K)$. We can assume that $a_{K,z,I} \leq Q$ for all $z \in Z(K)$ and $\text{char}(k_K) > M$. For each $1 \leq d \leq Q$, we see that the condition $a_{K,z,I} = d$ will be an $\mathcal{L}$-definable condition in $K, z$. Write

$$Z_d(K) = \{z \in Z(K) | a_{K,z,I} = d\}$$

and set

$$F = \sum_{d=1}^{Q} d \mathbf{1}_{Z_d} \otimes [Z],$$

where $\mathbf{1}_{Z_d}$ stands for the characteristic function of $Z_d$ and $[Z] = [Z \xrightarrow{id} Z]$. Then it is obvious that $F_K(z) = a_{K,z,I}$ and that we can ensure that $F \in \tilde{\mathcal{C}}_+(Z)$.

**Case 2: $I$ is not the empty set.**

By a similar argument as for Case 1, we may and do suppose that $I = \{1,...,n\}$. To simplify the notation, $R_{K,I}$ will be rewritten as $R_K$.

The number $a_{K,z,I}$ does not change if we remove from $R_{K,z}$ any $x$ with $x \notin MB(x/_{\sim_{K,z}})$, and thus, we can assume that

$$R_{K,z} = \cup_{x \in R_{K,z}} MB(x/_{\sim_{K,z}}).$$

Consider a definable subset $D_K$ of $R_K \times k_K^n$, described as follows:

For each $z \in Z_K$ and each $x = (x_1,...,x_n) \in R_{K,z}$, the fiber $D_{K,x,z}$ of $D_K$ over $x$ and $z$ is

$$\{(\xi_1,...,\xi_n) \in k_K^n | \wedge_{m=1}^n \left((\xi_m = 0) \vee (\exists y \in x/_{\sim_{K,z}}(x,m)\right.$$
$$\left.[(\text{ord}(x_m - y) = f_{K,z,m} - 1) \wedge (\xi_m = \overline{ac}(x_m - y))])\right)\}.$$

From the definition of $D$ we observe that

$$\#(D_{K,x,z}) = \prod_{m=1}^{n} \#(D_{K,x,z,m})$$

where $D_{K,z,x,m}$ the image under the projection of $D_{K,z,x}$ to the $m$-th coordinate. Moreover, if $B$ is a ball of volume $q_K^{-f_{K,z,m}+1}$ such that $x_m \in B$ then $B \cap x/_{\sim_{K,z}}(x,m)$ is disjoint union of $\#(D_{K,z,x,m})$ balls of volume $q_K^{-f_{K,z,m}}$. It follows that

$$\#(D_{K,z,x,m}) = \#(D_{K,z,x',m})$$

for all $x, x'$ with $(x' \sim_{K,z} x) \wedge (x_{\leq m-1} = x'_{\leq m-1}) \wedge (x'_m \in B)$ and so we can denote this number by $N_{x/_{\sim_{K,z}}, x_{\leq m-1}, B}$.

By Lemma 1.4.8, the function $f_{K,z}(x) = \sum_{i=1}^{n} f_{K,z,i}(x)$ is an $\mathcal{L}$-definable function from $R_K$ to $\mathbb{Z}$, we denote by $f$ the definable function from $R$ to $\mathbb{Z}$ such that $f_K(x, z) = f_{K,z}(x)$ for all $(x, z) \in R_K$. By Corollary 1.4.7, there exist $M$ and $Q$ in $\mathbb{N}$ such that

$$\prod_{m=1}^{n} Multinumber_m(x/_{\sim_{K,z}}, x) \leq Q$$

for every $(x, z) \in R_K$ and all $K$ with $\text{char}(k_K) > M$. Hence, for each $d$ with $1 \leq d \leq Q$ the set

$$R_K(d) := \{(x, z) \in R_K | \prod_{m=1}^{n} Multinumber_m(x/_{\sim_{K,z}}, x) = d\}$$

is an $\mathcal{L}$-definable subset of $R_K$. Consider a rational motivic constructible function $g \in \tilde{\mathcal{C}}(R)$ such that

$$g_K(x, z) = \sum_{d=1}^{Q} \frac{\mathbf{1}_{R_K(d)}(x, z)}{d}.$$

Here, the denominator $d$ can be viewed as a set of $d$ $\emptyset$-definable points in the residue field sort. Next we define the rational motivic constructible function $\Phi \in \tilde{\mathcal{C}}_+(R)$ by

$$\Phi = \frac{g \cdot \mathbb{L}^f}{[D]},$$

where the map $D \to R$ comes from the coordinate projection. By the definition, $\Phi$ is $Z$-integrable if we have that $f_{L,z}$ is bounded above on $R_{L,z}$ for each $z \in Z(L)$ and each $\mathcal{T}$-field $L$. The fact that $f_{L,z}$ is bounded above is a definable condition, hence, up to replacing $f$ so that it is zero if $f_{L,z}$ is not bounded above, we may suppose that $f_{L,z}$ is bounded above on $R_{L,z}$ for each $z \in Z(L)$ and each $\mathcal{T}$-field $L$ and thus that $\Phi$ is $Z$-integrable. Set

$$F = \mu_{/Z}(\Phi) \text{ in } \tilde{\mathcal{C}}_+(Z).$$

Finally we prove that $a_{K,z} = F_K(z)$ for all $z \in Z(K)$ and all local fields $K$ with $\text{char}(k_K) > M$. If $\text{char}(k_K) > M$ for well-chosen $M$ we have

$$F_K(z) = \mu_{/\{z\}}(\Phi|_{R_{K,z}}) = \int_{x \in R_{K,z}} \frac{g_K(x, z) q_K^{f_K(x,z)}}{\#(D_{K,x,z})} |dx|$$

where $|dx|$ is the Haar measure on $K^n$. Let $J_z$ be a set of $a_{K,z}$ representatives of $\sim_{K,z}$. Thus,

$$\int_{x \in R_{K,z}} \frac{g_K(x, z) q_K^{f_K(x,z)}}{\#(D_{K,x,z})} |dx|$$
$$= \sum_{x \in J_z} \int_{y \in x/_{\sim_{K,z}}} \frac{g_K(y, z) q_K^{f_K(y,z)}}{\#(D_{K,y,z})} |dy|$$

Let us for a moment fix $x$ and write $X = x/{\sim_{K,z}}$ and $d_m(y) = Multinumber_m(y/{\sim_{K,z}}, y)$. Then, using Fubini's theorem, we have

$$\int_{y \in X} \frac{g_K(y,z) q_K^{f_K(y,z)}}{\#(D_{K,y,z})} |dy|$$

$$= \int_{\overline{y} \in X(n-1)} \int_{y_n \in X_{\overline{y}}} \frac{q_K^{f_{K,z,1}(x)+\ldots+f_{K,z,n}(x)}}{\prod_{m=1}^{n}[d_m(\overline{y}, y_n) \times \#(D_{K,z,(\overline{y},y_n),m})]} |dy_n| \cdot |d\overline{y}|$$

$$= \int_{\overline{y} \in X(n-1)} \frac{q_K^{f_{K,z,1}(x)+\ldots+f_{K,z,n-1}(x)}}{\prod_{m=1}^{n-1}[d_m(\overline{y}) \times \#(D_{K,z,(\overline{y},y_n),m})]} |d\overline{y}|,$$

where $f_{K,z,i}(y) = f_{K,z,i}(x)$ for all $y \in X$ and $1 \leq i \leq n$ and where $d_m(\overline{y}) := d_m(\overline{y}, y_n)$ does not depend on $y_n$ when $y = (\overline{y}, y_n)$ varies in $X$ for all $1 \leq m \leq n-1$. The last equality comes from:

$$\int_{y_n \in X_{\overline{y}}} \frac{q_K^{f_{K,z,n}(x)}}{d_n(\overline{y}) \times \#(D_{K,z,(\overline{y},y_n),n})} |dy_n|$$

$$= \frac{q_K^{f_{K,z,n}(x)}}{d_n(\overline{y})} \sum_{i=1}^{d_n(\overline{y})} \int_{y_n \in B_i \cap X_{\overline{y}}} \frac{1}{\#(D_{K,z,(\overline{y},y_n),n})} |dy_n|$$

$$= \frac{q_K^{f_{K,z,n}(x)}}{d_n(\overline{y})} \sum_{i=1}^{d_n(\overline{y})} \text{Vol}(B_i \cap X_{\overline{y}}) \frac{1}{N_{x/\sim_{K,z},\overline{y},B_i}}$$

$$= \frac{q_K^{f_{K,z,n}(x)}}{d_n(\overline{y})} \sum_{i=1}^{d_n(\overline{y})} q_K^{-f_{K,z,n}(x)}$$

$$= 1$$

where $B_1, \ldots, B_{d_n(\overline{y})}$ are the balls with volume $q_K^{-f_{K,z,n}(x)+1}$ which have nonempty intersection with $X_{\overline{y}}$; indeed, one sees that $B_i \cap X_{\overline{y}}$ is the union of $N_{x/\sim_{K,z},\overline{y},B_i}$ many disjoint balls of volume $q_K^{-f_{K,z,n}(x)}$. By applying Fubini's theorem with a similar calculation to each of the remaining $n-1$ variables we deduce that

$$\int_{y \in X} \frac{g_K(y,z) q_K^{f_K(y,z)}}{\#(D_{K,y,z})} |dy| = 1.$$

We conclude that

$$F_K(z) = \mu_{/\{z\}}(\Phi|_{R_{K,z}}) = \#J_z = a_{K,z}$$

for all local fields $K$ with $char(k_K) > M$ and $z \in Z(K)$. The main theorem is proved. □

# Chapter 2

# New bounds for exponential sums with a non-degenerate phase polynomial

*This chapter is a joint work with Wouter Castryck, see [CNa].*

**Contents**

# Abstract

We prove a recent conjecture due to Cluckers and Veys on exponential sums modulo $p^m$ for $m \geq 2$ in the special case where the phase polynomial $f$ is sufficiently non-degenerate with respect to its Newton polyhedron at the origin. Our main auxiliary result is an improved bound on certain related exponential sums over finite fields. This bound can also be used to settle a conjecture of Denef and Hoornaert on the candidate-leading Taylor coefficient of Igusa's local zeta function associated to a non-degenerate polynomial, at its largest non-trivial real candidate pole.

## 2.1 Introduction

### 2.1.1

Let $f \in \mathbb{Z}[x]$ be a non-zero polynomial in the variables $x = x_1, \ldots, x_n$ such that $f(0) = 0$. In this article we prove new bounds on the absolute value of exponential sums of the form

$$S_f(p, m) := \frac{1}{p^{mn}} \sum_{x \in \{0, \ldots, p^m-1\}^n} \exp\left(2\pi \mathbf{i} \frac{f(x)}{p^m}\right)$$

where $p$ is a prime number and $m \geq 1$ is an integer. We work under the assumption that $f$ is non-degenerate with respect to the faces of its Newton polyhedron $\Delta_0(f)$ at the origin, in the strong sense recalled in Section 2.2.1 below. Concretely, let $\sigma \in \mathbb{Q}_{>0}$ be such that $(1/\sigma, \ldots, 1/\sigma)$ is contained in a proper face of $\Delta_0(f)$, and let $\kappa$ denote the maximal codimension in $\mathbb{R}^n$ of such a face. Then we prove the existence of a constant $c \in \mathbb{R}_{>0}$ which only depends on $f$ such that

$$|S_f(p, m)| \leq c p^{-\sigma m} m^{\kappa - 1} \tag{2.1.1.1}$$

for all sufficiently large prime numbers $p$ and all integers $m \geq 2$. Moreover, if $f$ is supported on a hyperplane which does not contain the origin and which has a normal vector in $\mathbb{R}^n_{\geq 0}$, then we can include $m = 1$ in the foregoing statement.

It is known that $\sigma_0(f) = \min\{1, \sigma\}$ where $\sigma_0(f)$ denotes the log canonical threshold $\sigma_0(f)$ of $f$ at the origin [Laz04, §9.3.C]. So our work implies the existence of a constant $c \in \mathbb{R}_{>0}$ such that

$$|S_f(p, m)| \leq c p^{-\sigma_0(f) m} m^{n-1} \tag{2.1.1.2}$$

for all primes $p$ and integers $m \geq 2$. Likewise, if the support of $f$ is contained in a hyperplane not passing through the origin and having a normal vector in $\mathbb{R}^n_{\geq 0}$, then the same conclusion holds with $m \geq 1$. The fact that we can write 'all primes' rather than 'all sufficiently large primes' follows from an observation due to Igusa, which implies the bound (2.1.1.2) for all primes $p$ but for some constant $c$ that is a priori allowed to depend on $p$; see e.g. [Den, p. 364] and [Igu78, p. 78].

## 2.1.2

Let us give some context for these results. Igusa's conjecture on exponential sums [Igu78, p. 2] predicts that a bound of the form (2.1.1.2) should hold for all primes $p$ and all integers $m \geq 1$, regardless of any non-degeneracy assumption but under the condition that $f$ is a non-constant homogeneous polynomial. This is related to the integrability of certain functions over the adèles, which in turn is connected with the validity of a generalized Poisson summation formula of Siegel-Weil type [Igu78, Ch. 4]. Igusa's conjecture was proven in the non-degenerate case by Denef and Sperber [DS01, Thm. 1.2] subject to a certain combinatorial constraint on $\Delta_0(f)$. This constraint was later removed by Cluckers [Clu10, Thm. 3.1], who in fact proved the bound (2.1.1.1) for all $m \geq 1$ and all non-constant *quasi*-homogeneous non-degenerate polynomials $f$. Cluckers' result naturally leads to the following strengthening of the statement of Igusa's conjecture:

**Conjecture 2.1.1** (Igusa, generalization due to Cluckers)**.** *For all non-constant quasi-homogeneous polynomials $f \in \mathbb{Z}[x]$ there exists a constant $c \in \mathbb{R}_{>0}$ such that the bound (2.1.1.2) holds for all primes $p$ and all integers $m \geq 1$.*

In a recent paper Cluckers and Veys predict [CV16, Conj. 1.2] that assuming $m \geq 2$ should allow to drop the quasi-homogeneity condition from the statement of Igusa's conjecture. However, now it should be taken into account that there may exist points $\alpha = (\alpha_1, \ldots, \alpha_n) \in \mathbb{C}^n$ at which the log canonical threshold $\sigma_\alpha(f)$ of the hypersurface $f(x) - f(\alpha) = 0$ at $\alpha$ is strictly smaller than $\sigma_0(f)$. This phenomenon does not occur in the quasi-homogeneous case.

**Conjecture 2.1.2** (Cluckers, Veys)**.** *For all non-constant polynomials $f \in \mathbb{Z}[x]$ there exists a constant $c \in \mathbb{R}_{>0}$ such that the bound (2.1.1.2) holds for all primes $p$ and all integers $m \geq 2$, provided that we replace $\sigma_0(f)$ by $\min_{\alpha \in \mathbb{C}^n} \sigma_\alpha(f)$.*

Our contribution to this topic is twofold. Firstly, our result confirms Cluckers and Veys' conjecture in the special case where $f$ is non-degenerate; in this case the foregoing concern is void since the minimal log canonical threshold is always realized by $\sigma_0(f)$. Secondly, we raise the question whether in Conjecture 2.1.1 the assumption that $f$ is non-constant and quasi-homogeneous can be relaxed to the condition that $\mathrm{Supp}\, f$ is contained in a hyperplane not passing through the origin and having a normal vector in $\mathbb{R}^n_{\geq 0}$. Here as well, this is provided that we replace $\sigma_0(f)$ by $\min_{\alpha \in \mathbb{C}^n} \sigma_\alpha(f)$; we can in fact restrict to those $\alpha$ for which $f(\alpha) = 0$ by a version of Euler's identity. We give an affirmative answer in the non-degenerate case.

## 2.1.3

For convenience we have stated Igusa's conjecture and the Cluckers–Veys conjecture in terms of the log canonical threshold. However, several other versions have been put forward which may, in some cases, predict sharper bounds. In these versions the log canonical threshold is replaced by the motivic oscillation index [Clu08b, CV16], by the complex oscillation index [AGZV12, §13.1.5], or by minus the largest non-trivial real pole of Igusa's local zeta function associated with $f$ [Igu78]. In fact, in the Cluckers–Veys conjecture for the log canonical threshold, it is not entirely clear whether the condition $m \geq 2$ is absolutely necessary. The reason for including it comes from the other versions, where it is unavoidable in general (e.g., for $f = x^2y - x$ as explained in [Clu08b, Ex. 7.2]).

## 2.1.4

We work at the following level of generality. We fix a number field $K \supseteq \mathbb{Q}$, let $\mathbb{Z}_K$ denote its ring of integers, and consider a polynomial $f \in \mathbb{Z}_K[x]$, where as before $x$ abbreviates a list of $n \geq 1$ variables $x_1, \ldots, x_n$. For each non-zero prime ideal $\mathfrak{p} \subseteq \mathbb{Z}_K$ we consider the $\mathfrak{p}$-adic completion $K_\mathfrak{p}$ of $K$, along with its ring of integers $\mathbb{Z}_\mathfrak{p} = \{\, a \in K_\mathfrak{p} \mid \operatorname{ord}_\mathfrak{p}(a) \geq 0 \,\}$ and its residue field $\mathbb{F}_\mathfrak{p} = \mathbb{Z}_\mathfrak{p}/\mathfrak{p}$, whose cardinality we denote by $N\mathfrak{p}$. We denote by $|\cdot|_\mathfrak{p} = (N\mathfrak{p})^{-\operatorname{ord}_\mathfrak{p}(\cdot)}$ the corresponding non-archimedean norm on $K_\mathfrak{p}$.

Let $p$ be the prime number below $\mathfrak{p}$, and consider the additive character
$$\psi_\mathfrak{p} : K_\mathfrak{p} \to \mathbb{C}^* : a \mapsto \exp(2\pi\mathbf{i}\operatorname{Tr}_{K_\mathfrak{p}/\mathbb{Q}_p}(a))$$
where $\exp(2\pi\mathbf{i}\cdot)$ is evaluated on $p$-adic numbers as follows: given $a \in \mathbb{Q}_p$ we let $a'$ be a representant inside $\mathbb{Z}[1/p]$ of the residue class of $a$ modulo the ring of $p$-adic integers $\mathbb{Z}_p$ and we define $\exp(2\pi\mathbf{i}a)$ as $\exp(2\pi\mathbf{i}a')$. Then to each $y \in K_\mathfrak{p}$ we associate the integral
$$S_{f,\mathfrak{p}}(y) := \int_{\mathbb{Z}_\mathfrak{p}^n} \psi_\mathfrak{p}(yf(x))|dx|$$
where $|dx| = |dx_1 \wedge \ldots \wedge dx_n|$ denotes the Haar measure, normalized such that the volume of $\mathbb{Z}_\mathfrak{p}^n$ is 1. Notice that the sum $S_f(p, m)$ from Section 2.1.1 can be rewritten as $S_{f,(p)}(p^{-m})$.

Our main result is as follows:

**Theorem 2.1.3.** *Let $f \in \mathbb{Z}_K[x]$ be a non-zero polynomial such that $f(0) = 0$ and assume that it is non-degenerate with respect to the faces of its Newton polyhedron $\Delta_0(f)$ at the origin. Let $\sigma \in \mathbb{Q}_{>0}$ be such that $(1/\sigma, \ldots, 1/\sigma)$ is contained in a proper face of $\Delta_0(f)$, and let $\kappa$ denote the maximal codimension of such a face. Then there exists a constant $c \in \mathbb{R}_{>0}$ only depending on $\Delta_0(f)$ such that for all non-zero prime ideals $\mathfrak{p} \subseteq \mathbb{Z}_K$ for which $N\mathfrak{p}$ is sufficiently large and all $y \in K_\mathfrak{p}$ satisfying $\operatorname{ord}_\mathfrak{p}(y) \leq -2$, we have*
$$|S_{f,\mathfrak{p}}(y)| \leq c|y|_\mathfrak{p}^{-\sigma}|\operatorname{ord}_\mathfrak{p}(y)|^{\kappa-1}.$$

*If $\operatorname{Supp} f$ is contained in a hyperplane not passing through the origin and having a normal vector in $\mathbb{R}_{\geq 0}^n$, then moreover $c$ can be chosen such that the bound also applies to all $y \in K_\mathfrak{p}$ for which $\operatorname{ord}_\mathfrak{p}(y) = -1$.*

It is clear that Theorem 2.1.3 implies the claims made in Sections 2.1.1 and 2.1.2.

The condition that $N\mathfrak{p}$ is sufficiently large allows us to assume that $f$ is non-degenerate at $\mathfrak{p}$, by which we mean that

$$f_\mathfrak{p} := f \bmod \mathfrak{p} \in \mathbb{F}_\mathfrak{p}[x]$$

is non-degenerate with respect to the faces of its Newton polygon $\Delta_0(f_\mathfrak{p})$ at the origin and $\Delta_0(f_\mathfrak{p}) = \Delta_0(f)$. If $f \in \mathbb{Z}_K[x] \subseteq K[x]$ is non-degenerate with respect to the faces of $\Delta_0(f)$ to start from, then it is indeed non-degenerate at all but finitely many non-zero prime ideals $\mathfrak{p} \subseteq \mathbb{Z}_K$; this follows, for instance, from Hilbert's Nullstellensatz.

Although we avoid a detailed discussion, we note that our proof of Theorem 2.1.3 also applies to other global fields, such as the field of rational functions $K = \mathbb{F}_q(t)$ over a finite field $\mathbb{F}_q$. However, here we have the extra condition that $\operatorname{char} \mathbb{F}_q$ should not be contained in the set $P$ of bad primes associated with $\operatorname{Supp} f$, appearing in the statement of Theorem 2.1.4 below. (More generally, the discussion from [Clu08a, §9] in the context of Igusa's conjecture applies here, too.)

## 2.1.5

In Section 2.2 we will explain how Theorem 2.1.3 arises as a consequence of the following finite field exponential sum estimate, which is the central auxiliary result of this paper and which is proven in Section 2.4.

**Theorem 2.1.4.** *Let $\mathbb{F}_q$ be a finite field with $q$ elements and let $f \in \mathbb{F}_q[x]$ be non-degenerate with respect to the faces of its Newton polyhedron $\Delta_0(f)$ at the origin. Suppose that $\operatorname{Supp} f$ is contained in a hyperplane which does not contain the origin and which has a normal vector in $\mathbb{R}^n_{\geq 0}$. Let $\sigma \in \mathbb{Q}_{>0}$ be maximal such that $(1/\sigma, \ldots, 1/\sigma) \in \Delta_0(f)$ and let $\varphi : \mathbb{F}_q \to \mathbb{C}^*$ be a non-trivial additive character on $\mathbb{F}_q$. There exist a constant $c \in \mathbb{R}_{>0}$ which only depends on $\Delta_0(f)$ and a finite set of primes $P$ which only depends $\operatorname{Supp} f$ such that*

$$\left| \frac{1}{(q-1)^n} \sum_{x \in \mathbb{F}_q^{*n}} \varphi(f(x)) \right| < cq^{-\sigma}$$

*as soon as $p = \operatorname{char} \mathbb{F}_q \notin P$.*

If the hyperplane can be chosen such that it has a normal vector in $\mathbb{R}^n_{>0}$ then $f$ is quasi-homogeneous, in which case the foregoing result is due to Cluckers [Clu10, Prop. 6.2 & its proof].

## 2.1.6

Theorem 2.1.4 also implies a conjecture by Denef and Hoornaert, involving Igusa's local zeta function

$$Z_{f,\mathfrak{p}} : \{s \in \mathbb{C} \mid \Re s > 0\} \to \mathbb{C} : s \mapsto \int_{\mathbb{Z}_\mathfrak{p}^n} |f(x)|_\mathfrak{p}^s |dx|,$$

which one can associate to any non-zero prime ideal $\mathfrak{p} \subseteq \mathbb{Z}_K$ and any polynomial $f \in \mathbb{Z}_K[x]$. It is well-known that this is a rational function in $N\mathfrak{p}^{-s}$, so it admits a meromorphic continuation to the entire complex plane, where one may see poles showing up. These poles are believed to contain important arithmetic and geometric information about the hypersurface $f = 0$; see [Den, Seg06] for more background. In our setting where $f$ is non-constant, vanishing at 0, and non-degenerate with respect to the faces of its Newton polyhedron at the origin, Hoornaert proved [Hoo02, §4.3] that there is always at least one real pole and that the largest such pole is either $s = -\sigma$ or $s = -1$. If $s = -1$ is a pole, it is called trivial. The expected pole order of $-\sigma$ is $\kappa$, unless $-\sigma = -1$ in which case the expected pole order is $\kappa + 1$. Here 'expected' means that the actual pole order is bounded by the said quantity and typically equals it; a sufficient condition for equality is $\sigma < 1$, but see [DH01, Thms. 5.5, 5.17, 5.19] and [Hoo02, Thm. 4.10] for more precise statements.

In Section 2.5 we will demonstrate how Theorem 2.1.4 implies a certain uniformity in $\mathfrak{p}$ of (what is expected to be) the leading Taylor coefficient at $s = -\sigma$. More precisely, we will show that the real number

$$\lim_{s \to -\sigma} (N\mathfrak{p}^{s+\sigma} - 1)^{\kappa + \delta_{\sigma,1}} Z_{f,\mathfrak{p}}(s)$$

is in $O(N\mathfrak{p}^{1-\max\{1,\sigma\}})$ as $\mathfrak{p}$ varies; here $\delta_{\cdot,\cdot}$ denotes the Kronecker delta. This was proved by Denef and Hoornaert [DH01, §5] under a certain combinatorial assumption on $\Delta_0(f)$ which they conjectured to be superfluous. Our work confirms their conjecture. One important subtlety is that Denef and Hoornaert work under a considerably weaker notion of non-degeneracy than we do, so that their conjecture is a priori stronger. However as explained in Section 2.5.2 this is not a concern: we will see that our a priori weaker conclusion easily implies the Denef–Hoornaert conjecture in its full strength. This will rely on some conclusions made in Section 2.2.2, where we elaborate on the difference between both non-degeneracy notions.

## 2.2 Non-degenerate polynomials and the invariant $\sigma$

### 2.2.1

Let us recall what it means for a polynomial to be non-degenerate with respect to the faces of its Newton polyhedron at the origin, while fixing some notation. Let $k$ be a field, which from Section 2.3 on will be either a number field $K$ or a finite field $\mathbb{F}_q$. Let

$$f = \sum_{i \in \mathbb{Z}_{\geq 0}^n} c_i x^i \in k[x] \setminus \{0\}$$

where as before $x = x_1, \ldots, x_n$ and $x^i = x^{(i_1, \ldots, i_n)}$ abbreviates $x_1^{i_1} \cdots x_n^{i_n}$. The Newton polyhedron of $f$ at the origin is defined as

$$\Delta_0(f) = \operatorname{conv} \operatorname{Supp} f + \mathbb{R}_{\geq 0}^n,$$

with $\operatorname{Supp} f = \left\{ i \in \mathbb{Z}_{\geq 0}^n \,\middle|\, c_i \neq 0 \right\}$ the support of $f$. For all non-empty faces $\tau \subseteq \Delta_0(f)$ of any dimension, ranging from vertices to $\Delta_0(f)$ itself, we write

$$f_\tau = \sum_{i \in \tau \cap \mathbb{Z}_{\geq 0}^n} c_i x^i$$

and we say that $f$ is non-degenerate with respect to $\tau$ if the system of equations

$$\frac{\partial f_\tau}{\partial x_1} = \ldots = \frac{\partial f_\tau}{\partial x_n} = 0$$

has no solutions in $\overline{k}^{*n}$, or in other words if the map

$$\overline{k}^{*n} \to \overline{k} : \alpha \mapsto f_\tau(\alpha)$$

has no critical values. Here $\overline{k}$ denotes an algebraic closure of $k$. Self-evidently we call $f$ non-degenerate with respect to the faces (resp. the compact faces) of $\Delta_0(f)$ if it is non-degenerate with respect to all choices (resp. all compact choices) of $\tau$.

### 2.2.2

Our non-degeneracy notion arises naturally in the study of exponential sums, but is considerably stronger than some of its counterparts which can be found elsewhere in the existing literature. Most notably, often one merely imposes the generic condition that 0 is not a critical value, or in other words that the hypersurface $f_\tau = 0$ has no *singularities* in $\overline{k}^{*n}$, rather than critical points. We will refer to the latter notion of non-degeneracy as *weak* non-degeneracy. (This is not intended to become standard terminology: the only purpose it serves is to avoid confusion throughout the remainder of this paper.) Clearly, if $f$ is non-degenerate with respect to a face $\tau \subseteq \Delta_0(f)$, then it is also weakly non-degenerate with respect to that face. An important remark is that the converse holds as soon as $\tau$ is contained in a hyperplane of the form

$$H = \{ (i_1, \ldots, i_n) \,|\, c_1 i_1 + \ldots + c_n i_n = b \}$$

with $c_1, \ldots, c_n, b \in \mathbb{Z}$ satisfying $\operatorname{char} k \nmid b$. This can be seen using a weighted version of Euler's identity. As a consequence, we could have equally well formulated Theorem 2.1.4 assuming weak non-degeneracy, rather than non-degeneracy.

### 2.2.3

It is convenient to introduce the following notation: to each vector $a = (a_1, \ldots, a_n) \in \mathbb{R}_{\geq 0}^n$ we associate

$$\nu(a) = a_1 + \ldots + a_n, \qquad N(a) = \min_{x \in \Delta_0(f)} x \cdot a,$$

$$F(a) = \{ x \in \Delta_0(f) \,|\, x \cdot a = N(a) \}.$$

The latter set is a face of $\Delta_0(f)$ which is called the first meet locus of $a$. It is contained in the hyperplane $a_1 x_1 + \ldots + a_n x_n = N(a)$. Every face arises as the first meet locus of some vector, where we note that $\Delta_0(f)$ itself is given by $F(\vec{0}) = F((0,\ldots,0))$. Similarly writing $\vec{1} = (1,\ldots,1)$, we define

$$\kappa = \operatorname{codim} F(\vec{1}),$$
$$\sigma = \nu(\vec{1})/N(\vec{1}) = n/N(\vec{1})$$

as in the introduction; if $N(\vec{1}) = 0$ then we let $\sigma = +\infty$. Note that if $\sigma < +\infty$ then $\sigma N(a) \leq \nu(a)$ for all $a \in \mathbb{R}^n_{\geq 0}$. We will usually write $\sigma(f)$ rather than $\sigma$ to emphasize the dependence on $f$. It is natural to define $\sigma(0) = 0$.

### 2.2.4

For use in Section 2.4.4 we prove the following list of properties of $\sigma$, which we believe to be interesting in their own right:

**Lemma 2.2.1.** *Let $f \in k[x]$ and $g \in k[y]$ be polynomials vanishing at the origin, in the respective variables $x = x_1,\ldots,x_n$ and $y = y_1,\ldots,y_m$. Then:*
  *(i) $\sigma(f+g) = \sigma(f) + \sigma(g)$,*
  *(ii) $\sigma(fg) = \min\{\sigma(f),\sigma(g)\}$,*
  *(iii) if $n \geq 2$ then $\sigma(f) \geq \sigma(f(x_1,\ldots,x_{n-1},x_{n-1}))$.*
*Here $f + g$ and $fg$ are viewed as polynomials in $k[x,y]$, while $f(x_1,\ldots,x_{n-1},x_{n-1})$ is viewed as an element of $k[x_1,\ldots,x_{n-1}]$.*

Specifying the ambient ring is actually not needed: viewing $f(x_1,\ldots,x_{n-1},x_{n-1})$ as an element of $k[x_1,\ldots,x_n]$ does not change the corresponding value of $\sigma$.

*Proof of the lemma.* We assume that $f, g \neq 0$, since otherwise the listed properties are trivial.

(i) If $P \in \partial\Delta_0(f)$ and $Q \in \partial\Delta_0(g)$, where $\partial$ denotes the topological boundary, then every convex combination of $(P, 0)$ and $(0, Q)$ is contained in $\partial\Delta_0(f+g)$. Applying this to $P = (1/\sigma(f),\ldots,1/\sigma(f)) \in \partial\Delta_0(f)$ and $Q = (1/\sigma(g),\ldots,1/\sigma(g)) \in \partial\Delta_0(g)$ and considering the convex combination

$$\frac{\sigma(f)}{\sigma(f) + \sigma(g)}(P, 0) + \frac{\sigma(g)}{\sigma(f) + \sigma(g)}(0, Q),$$

we see that the desired conclusion follows.

(ii) Assume without loss of generality that the minimum is realized by $\sigma(f)$. Note that

$$\Delta_0(fg) = \Delta_0(f) \times \Delta_0(g),$$

from which one sees that $1/\sigma(fg) \geq 1/\sigma(f)$ or in other words that $\sigma(fg) \leq \sigma(f)$. For the converse inequality, we write $(1/\sigma(f),\ldots,1/\sigma(f)) \in \mathbb{R}^{n+m}_{\geq 0}$ as

$$\left(\frac{1}{\sigma(f)},\ldots,\frac{1}{\sigma(f)},\frac{1}{\sigma(g)},\ldots,\frac{1}{\sigma(g)}\right) + \left(0,\ldots,0,\frac{1}{\sigma(f)} - \frac{1}{\sigma(g)},\ldots,\frac{1}{\sigma(f)} - \frac{1}{\sigma(g)}\right),$$

which is seen to be an element of $\Delta_0(fg)$, proving that $\sigma(fg) \geq \sigma(f)$.

(iii) Write $f' = f(x_1, \ldots, x_{n-1}, x_{n-1})$. We claim that if $(i_1, \ldots, i_{n-2}, i) \in \Delta_0(f')$ then there exist $a, b \in \mathbb{R}_{\geq 0}$ such that $(i_1, \ldots, i_{n-2}, a, b) \in \Delta_0(f)$ and $a + b = i$. Indeed, this property is clear for any point in the support of $f'$, and the claim follows by convexity considerations. Now apply this claim to the point

$$\left(\frac{1}{\sigma(f')}, \ldots, \frac{1}{\sigma(f')}\right) \in \Delta_0(f')$$

to find $a, b \in \mathbb{R}_{\geq 0}$ for which

$$\left(\frac{1}{\sigma(f')}, \ldots, \frac{1}{\sigma(f')}, a, b\right) \in \Delta_0(f)$$

with $a + b = 1/\sigma(f')$. When adding $(0, \ldots, 0, 1/\sigma(f') - a, 1/\sigma(f') - b) \in \mathbb{R}_{\geq 0}^n$ to this, we stay inside $\Delta_0(f)$, from which we conclude the desired inequality. □

We remark that these properties are reminiscent of well-known facts on the log canonical threshold at the origin. Indeed, with the notation and assumptions from above (and the extra assumption that $f, g \neq 0$), one has that
  (a) $\sigma_0(f + g) = \min\{1, \sigma_0(f) + \sigma_0(g)\}$,
  (b) $\sigma_0(fg) = \min\{\sigma_0(f), \sigma_0(g)\}$,
  (c) $\sigma_0(f) \geq \sigma_0(f(x_1, \ldots, x_{n-1}, x_{n-1}))$.
See [DK01, Rmk. 2.10], [Laz04, Prop. 9.5.22] and [Kol97, Thm. 7.5 & proof of Thm. 8.20], respectively. In fact, statements (i – iii) could have also been settled as mere corollaries to (a – c), by using the following trick. Notice that properties (i – iii) are purely combinatorial, so we can replace $f$ and $g$ by non-degenerate polynomials over $\mathbb{C}$ having the same Newton polyhedron; for what follows it suffices to assume weak non-degeneracy, which is generically satisfied. Next let $f_d = f(x_1^d, \ldots, x_n^d)$ and $g_d = g(y_1^d, \ldots, y_m^d)$ for some large positive integer $d$. These polynomials are again non-degenerate, and moreover $\Delta_0(f_d) = d\Delta_0(f)$ and $\Delta_0(g_d) = d\Delta_0(g)$ so that $\sigma(f) = d\sigma(f_d)$ and $\sigma(g) = d\sigma(g_d)$. If $d$ is large enough such that $\sigma(f_d) + \sigma(g_d) \leq 1$, then one sees that properties (i – iii) follow from properties (a – c) and the fact that $\sigma_0(h) = \min\{1, \sigma(h)\}$ for any polynomial $h$ that is weakly non-degenerate with respect to the faces of its Newton polyhedron at the origin [Laz04, §9.3.C].

**2.2.5**

We conclude with a bound on $\sigma(f)$ in terms of the dimension of the affine critical locus $C_f$ of $f$, i.e. the set of geometric points of $\mathbb{A}_k^n$ at which all partial derivatives of $f$ vanish, where we adopt the convention that the dimension of the empty scheme is $-\infty$. If $k$ is a number field then the bound can be viewed as a corollary to a general observation due to Cluckers: see [Clu08b, Thm. 5.1] and how this is applied in [Clu08a, Lem. 6.3], for instance. Our more direct approach has the advantage of working over any field.

**Lemma 2.2.2.** *Let $f \in k[x] \setminus \{0\}$ be non-degenerate with respect to the faces of its Newton polyhedron $\Delta_0(f)$ at the origin, and denote by $\delta$ the dimension inside $\mathbb{A}_k^n$ of $C_f$. Then $\sigma(f) \leq (n - \delta)/2$.*

*Proof.* By non-degeneracy with respect to the entire Newton polyhedron $\Delta_0(f)$, every critical point of $f$ has at least one coordinate which is zero. Therefore it suffices to prove for all proper subsets $I$ of $\{1, \ldots, n\}$ that

$$\sigma(f) \leq (n - \delta_I)/2$$

where $\delta_I = \dim C_I$ with $C_I = \{(a_1, \ldots, a_n) \in C_f \mid a_i \neq 0 \text{ if and only if } i \in I\}$.

By reordering the variables if needed we can assume that $I = \{n - d + 1, \ldots, n\}$ for some integer $d$ which satisfies $0 \leq \delta_I \leq d < n$. If $C_I = \emptyset$ then there is nothing to prove. If $C_I \neq \emptyset$ then it contains at least one point $(a_1, \ldots, a_n)$, which by our assumption satisfies $a_1 = \ldots = a_{n-d} = 0$ and $a_{n-d+1}, \ldots, a_n \neq 0$. We first claim that

$$\Delta_0(f) \subseteq \Delta_0(x_1 + \ldots + x_{n-d}) \times \mathbb{R}^d_{\geq 0}. \qquad (2.2.5.1)$$

If we let $H$ denote the face of $\mathbb{R}^n_{\geq 0}$ defined by

$$i_1 = \ldots = i_{n-d} = 0, \quad i_{n-d+1}, \ldots, i_n \geq 0$$

and we define $\tau := \Delta_0(f) \cap H$, then (2.2.5.1) is equivalent to $\tau = \emptyset$. Now suppose that $\tau \neq \emptyset$: then it must concern a face of $\Delta_0(f)$. But $\partial f_\tau / \partial x_i$ vanishes identically for $i = 1, \ldots, n - d$, while it vanishes at $(1, \ldots, 1, a_{n-d+1}, \ldots, a_n)$ for $i = n - d + 1, \ldots, n$. This is in contradiction with the fact that $f$ is non-degenerate with respect to $\tau$, so our claim follows.

We now let $\ell \in \{0, \ldots, n - d\}$ be minimal such that

$$\Delta_0(f) \subseteq \Delta_0(x_1 + \ldots + x_\ell + x_{\ell+1}^2 + \ldots + x_{n-d}^2) \times \mathbb{R}^2_{\geq 0}$$

up to reordering the variables $x_1, \ldots, x_{n-d}$. In other words $f$ can be written as

$$f = x_1 \cdot g_1(x_{n-d+1}, \ldots, x_n) + \ldots + x_\ell \cdot g_\ell(x_{n-d+1}, \ldots, x_n) + \ldots$$

for non-zero polynomials $g_1, \ldots, g_\ell \in k[x_{n-d+1}, \ldots, x_n]$, where the last dots consist of terms that are at least quadratic in the variables $x_1, \ldots, x_{n-d}$. It is easy to check that $\sigma(f) \leq (n - d + \ell)/2$, so it suffices to show that $\delta_I \leq d - \ell$. This is trivial if $\ell = 0$, so assume that $\ell \geq 1$. By taking partial derivatives one observes that

$$C_I = \{(0, \ldots, 0, a_{n-d+1}, \ldots, a_n) \mid (a_{n-d+1}, \ldots, a_n) \in S_I\}$$

with $S_I \subseteq \overline{k}^{*d}$ the scheme defined by $(g_1, \ldots, g_\ell)$. We will establish the desired bound on $\delta_I$ by proving that $S_I$ is either empty or a smooth complete intersection.

By the Jacobian criterion this amounts to showing that for all $(a_{n-d+1}, \ldots, a_n) \in S_I$ the rows of the matrix

$$J = \begin{pmatrix} \frac{\partial g_1}{\partial x_{n-d+1}}(a_{n-d+1}, \ldots, a_n) & \ldots & \frac{\partial g_1}{\partial x_n}(a_{n-d+1}, \ldots, a_n) \\ \vdots & \ddots & \vdots \\ \frac{\partial g_\ell}{\partial x_{n-d+1}}(a_{n-d+1}, \ldots, a_n) & \ldots & \frac{\partial g_\ell}{\partial x_n}(a_{n-d+1}, \ldots, a_n) \end{pmatrix}$$

are linearly independent. Suppose this is not the case, then there exists a vector $(\alpha_1, \ldots, \alpha_\ell) \neq (0, \ldots, 0)$ such that $(\alpha_1, \ldots, \alpha_\ell) \cdot J = \vec{0}$. Assume without loss of generality that $\alpha_1, \ldots, \alpha_{\ell'} \neq 0$ and $\alpha_{\ell'+1} = \ldots = \alpha_\ell = 0$, where $0 < \ell' \leq \ell$. Now consider the face $\tau \subseteq \Delta_0(f)$ obtained by intersecting $\Delta_0(f)$ with

$$H : i_1 + \ldots + i_{\ell'} = 1, i_{\ell'+1} = \ldots = i_\ell = i_{\ell+1} = \ldots = i_{n-d} = 0, i_{n-d+1}, \ldots, i_n \geq 0.$$

Then
$$f_\tau = x_1 \cdot g_1(x_{n-d+1}, \ldots, x_n) + \ldots + x_{\ell'} \cdot g_{\ell'}(x_{n-d+1}, \ldots, x_n)$$

and one verifies that $(\alpha_1, \ldots, \alpha_{\ell'}, 1, \ldots, 1, 1, \ldots, 1, a_{n-d+1}, \ldots, a_n) \in \overline{k}^{*n}$ is a common root of its partial derivatives: a contradiction with the non-degeneracy assumption with respect to $\tau$. □

## 2.3 Reduction to estimating finite field exponential sums

### 2.3.1

In this section we explain how proving Theorem 2.1.3 reduces to proving the bound on finite field exponential sums stated in Theorem 2.1.4. This resorts to a well-known reasoning by Denef and Sperber [DS01, Prop. 2.1], a slight generalization of which was elaborated by Cluckers [Clu10, Prop. 4.1]. The idea is to partition the integration domain $\mathbb{Z}_\mathfrak{p}^n$ according to all possible valuations, ignoring the zero-measure set of points $x$ in which $0$ appears as a coordinate:

$$S_{f,\mathfrak{p}}(y) = \int_{\mathbb{Z}_\mathfrak{p}^n} \psi_\mathfrak{p}(yf(x))|dx| = \sum_{a \in \mathbb{Z}_{\geq 0}^n} \int_{\substack{x \in \mathbb{Z}_\mathfrak{p}^n \\ \mathrm{ord}_\mathfrak{p}(x) = a}} \psi_\mathfrak{p}(yf(x))|dx|.$$

Let $\pi_\mathfrak{p}$ be a uniformizing parameter of $\mathbb{Z}_\mathfrak{p}$; if $\mathfrak{p}$ is unramified then one can just take $\pi_\mathfrak{p} = p$. By introducing new variables $u = (u_1, \ldots, u_n)$ through the substitution $x_j \leftarrow \pi_\mathfrak{p}^{a_j} u_j$, with $a_j$ the $j$th coordinate of $a$, we can rewrite the above sum as

$$\sum_{a \in \mathbb{Z}_{\geq 0}^n} N\mathfrak{p}^{-\nu(a)} \int_{u \in \mathbb{Z}_\mathfrak{p}^{*n}} \psi_\mathfrak{p}(y \pi_\mathfrak{p}^{N(a)} (f_{F(a)}(u) + \pi_\mathfrak{p}(\ldots)))|du|$$

where the expression $(\ldots)$ takes values in $\mathbb{Z}_\mathfrak{p}$. As explained in Section 2.1.4 we can assume that $f$ is non-degenerate at $\mathfrak{p}$, in which case Hensel's lemma implies that the integral is zero whenever $\mathrm{ord}_\mathfrak{p}(y) + N(a) \leq -2$. On the other hand since $\psi_\mathfrak{p}$ is trivial on $\mathbb{Z}_\mathfrak{p}$, as soon as $\mathrm{ord}_\mathfrak{p}(y) + N(a) \geq 0$ the integral is just the measure of $\mathbb{Z}_\mathfrak{p}^{*n}$, which is $(1 - N\mathfrak{p}^{-1})^n$. The most interesting case is where $\mathrm{ord}_\mathfrak{p}(y) + N(a) = -1$, in which case the integral equals

$$N\mathfrak{p}^{-n} \sum_{x \in \mathbb{F}_\mathfrak{p}^{*n}} \varphi_y(\overline{f}_{F(a)}(x))$$

where
$$\varphi_y : \mathbb{F}_\mathfrak{p} \to \mathbb{C}^* : x \bmod \mathfrak{p} \mapsto \psi_\mathfrak{p}(y \pi_\mathfrak{p}^{-\mathrm{ord}_\mathfrak{p}(y)-1} x)$$

and $\overline{f}_{F(a)}$ denotes the reduction of $f_{F(a)}$ mod $\mathfrak{p}$. Note that $\varphi_y$ is a well-defined non-trivial additive character on $\mathbb{F}_\mathfrak{p}$.

We end up with

$$S_{f,\mathfrak{p}}(y) = (1 - N\mathfrak{p}^{-1})^n \sum_{\substack{\tau \text{ face} \\ \text{of } \Delta_0(f)}} \left( A_{\tau,\mathfrak{p}}(y) + \frac{B_{\tau,\mathfrak{p}}(y)}{(N\mathfrak{p} - 1)^n} \sum_{x \in \mathbb{F}_{\mathfrak{p}}^{*n}} \varphi_y(\overline{f}_\tau(x)) \right)$$

where

$$A_{\tau,\mathfrak{p}}(y) = \sum_{\substack{a \in \mathbb{Z}_{\geq 0}^n \text{ s.t. } F(a) = \tau \\ \text{and } N(a) \geq -\operatorname{ord}_\mathfrak{p}(y)}} N\mathfrak{p}^{-\nu(a)} \quad \text{and} \quad B_{\tau,\mathfrak{p}}(y) = \sum_{\substack{a \in \mathbb{Z}_{\geq 0}^n \text{ s.t. } F(a) = \tau \\ \text{and } N(a) = -\operatorname{ord}_\mathfrak{p}(y) - 1}} N\mathfrak{p}^{-\nu(a)}.$$

Here a trivial but important remark is that $B_{\tau,\mathfrak{p}}(y) = 0$ as soon as $\operatorname{ord}_\mathfrak{p}(y) \leq -2$ and the affine span of $\tau$ passes through the origin, because in this case $N(a) = 0$ while $-\operatorname{ord}_\mathfrak{p}(y) - 1 \geq 1$.

### 2.3.2

Using the estimates from [Clu10, Prop. 5.2] or [DS01, §3.4] one sees that there exists a constant $c \in \mathbb{R}_{>0}$ such that

$$A_{\tau,\mathfrak{p}}(y) \leq c|y|_\mathfrak{p}^{-\sigma}|\operatorname{ord}_\mathfrak{p}(y)|^{\kappa-1} \quad \text{and} \quad B_{\tau,\mathfrak{p}}(y) \leq c|y|_\mathfrak{p}^{-\sigma} N\mathfrak{p}^{\sigma(f_\tau)}|\operatorname{ord}_\mathfrak{p}(y)|^{\kappa-1}$$

for all choices of $\mathfrak{p}$ and $y$, where we recall that $\sigma(f_\tau) \leq \sigma$ is the minimal rational number such that $(1/\sigma(f_\tau), \ldots, 1/\sigma(f_\tau))$ is contained in $\Delta_0(f_\tau) = \tau + \mathbb{R}_{\geq 0}^n$. Using these bounds one verifies that in order to prove Theorem 2.1.3 it suffices to show that for all faces $\tau \subseteq \Delta_0(f)$ there exists a constant $c$ only depending on $\Delta_0(f)$ such that

$$\left| \frac{1}{(N\mathfrak{p} - 1)^n} \sum_{x \in \mathbb{F}_{\mathfrak{p}}^{*n}} \varphi_y(\overline{f}_\tau(x)) \right| \leq c \cdot N\mathfrak{p}^{-\sigma(f_\tau)}. \tag{2.3.2.1}$$

for all non-zero prime ideals $\mathfrak{p}$ having a sufficiently large norm. In fact, it suffices to establish this bound for the faces $\tau$ whose affine span does not contain the origin. Indeed, if $\operatorname{ord}_\mathfrak{p}(y) \leq -2$ then this claim is straightforward in view of the remark concluding Section 2.3.1. On the other hand, if $\operatorname{ord}_\mathfrak{p}(y) = -1$ and $\operatorname{Supp} f$ is contained in a hyperplane $H$ which does not contain the origin and which has a normal vector in $\mathbb{R}_{\geq 0}^n$, then in order to conclude (2.3.2.1) for a given face $\tau \subseteq \Delta_0(f)$, it suffices to prove it for the face $\tau \cap H$, whose affine span does not contain the origin.

## 2.4    New bounds for finite field exponential sums

### 2.4.1

Thus we are left with proving Theorem 2.1.4. Let $H$ be a hyperplane as in the statement; we can assume it to be of the form

$$H : c_1 i_1 + \ldots + c_n i_n = b$$

for $b \in \mathbb{Z}_{>0}$ and $c_1, \ldots, c_n \in \mathbb{Z}_{\geq 0}$. Without loss of generality we can order the variables such that $c_1, \ldots, c_{n-r} > 0$ and $c_{n-r+1}, \ldots, c_n = 0$. For simplicity we choose $H$ such that $r$ is maximal. We can assume that $r > 0$ because if $r = 0$ then $f$ is quasi-homogeneous, and as mentioned in this case Theorem 2.1.4 is due to Cluckers [Clu10, Prop. 6.2]. Define $P$ as the set of prime numbers dividing $b$ and assume throughout the rest of this section that $\operatorname{char} \mathbb{F}_q \notin P$. Note that $P$ clearly depends on $\operatorname{Supp} f$ only. In fact it even depends on $\Delta_0(f)$ only, but the reason for writing $\operatorname{Supp} f$ in the statement of Theorem 2.1.4 is that $P$ will be enlarged in Section 2.4.6, in a way which could a priori depend on the specific configuration of $\operatorname{Supp} f$.

## 2.4.2

Rename the variables $x_{n-r+1}, \ldots, x_n$ as $z = z_1, \ldots, z_r$ and view $f$ as a quasi-homogeneous polynomial over $\mathbb{F}_q[z]$ in the remaining variables. Note that none of the weights $c_j$ can exceed $b$, otherwise it would be possible to remove the term $c_j i_j$, which contradicts the maximality of $r$. We order the indices such that $x_1, \ldots, x_s$ are the variables for which the corresponding weights $c_1, \ldots, c_s$ are equal to $b$. These necessarily appear linearly with a non-zero coefficient in $\mathbb{F}_q[z]$, again by the maximality of $r$. The variables $x_{s+1}, \ldots, x_{n-r}$ which have a smaller non-zero weight $c_i$ are renamed $y = y_1, \ldots, y_t$. Note that $n = s + r + t$. Then we can write

$$f = h + g_1 x_1 + \ldots + g_s x_s,$$

where $g_1, \ldots, g_s \in \mathbb{F}_q[z] \setminus \{0\}$ and where $h \in \mathbb{F}_q[y, z]$ is quasi-homogeneous when considered over $\mathbb{F}_q[z]$. Moreover for each concrete value of $z \in \overline{\mathbb{F}}_q^r$ the polynomial $h(\cdot, z) \in \overline{\mathbb{F}}_q[y]$ admits $y = 0$ as a critical point. For the sake of clarity, we note that this includes the case where $h(\cdot, z)$ is identically zero, in which case every point is considered critical.

For each $I \subseteq \{1, \ldots, s\}$ and $0 \leq d \leq t$ we define

$$V_{I,d} = \{\, z \in \overline{\mathbb{F}}_q^{*r} \mid \dim C_{h(\cdot, z)} = d, g_i(z) = 0 \Leftrightarrow i \in I \,\}$$

and we write

$$\left| \frac{1}{(q-1)^n} \sum_{(x,y,z) \in \mathbb{F}_q^{*n}} \varphi(f(x,y,z)) \right|$$

$$= \left| \sum_{I,d} \sum_{z \in V_{I,d}(\mathbb{F}_q)} \frac{1}{(q-1)^n} [\sum_{x \in \mathbb{F}_q^{*s}} \varphi(\sum_{i \notin I} g_i(z) x_i)] \times [\sum_{y \in \mathbb{F}_q^{*t}} \varphi(h(y,z))] \right|$$

$$= \left| \sum_{I,d} \sum_{z \in V_{I,d}(\mathbb{F}_q)} \frac{1}{(q-1)^n} (-1)^{s-\#I} (q-1)^{\#I} \times [\sum_{y \in \mathbb{F}_q^{*t}} \varphi(h(y,z))] \right|$$

$$\leq \sum_{I,d} c_{I,d} q^{\dim V_{I,d}} (q-1)^{\#I - n} \left| \sum_{y \in \mathbb{F}_q^{*t}} \varphi(h(y,z)) \right|$$

where we recall our convention that the dimension of the empty scheme is $-\infty$. In the last step we used the Lang–Weil estimates, which introduce constants $c_{I,d} \in \mathbb{R}_{>0}$ that can be taken to depend on $I, d$ and $\Delta_0(f)$ only.

### 2.4.3

If $h$ is identically zero then $\dim C_{h(\cdot,z)} = t$ for all $z$ and the foregoing bound simplifies to

$$\sum_I c_I q^{\dim V_I}(q-1)^{t+\#I-n} \tag{2.4.3.1}$$

where $c_I = c_{I,t}$ and

$$V_I = V_{I,t} = \{\, z \in \overline{\mathbb{F}}_q^{*r} \mid g_i(z) = 0 \Leftrightarrow i \in I \,\}.$$

By the non-degeneracy of $f$, as in the proof of Lemma 2.2.2 we see that the scheme cut out by $(g_i)_{i \in I}$ is either empty or a smooth complete intersection. In particular $\dim V_I \leq r - \#I$. From this we see that (2.4.3.1) is bounded by $cq^{-s}$ where $c = \sum_I c_I$. Now from

$$\Delta_0(f) = \Delta_0(g_1 x_1 + \ldots + g_s x_s) \subseteq \Delta_0(x_1 + \ldots + x_s) \times \mathbb{R}_{\geq 0}^{t+r}$$

it is immediate that $\sigma(f) \leq s$, from which the desired bound follows.

### 2.4.4

Thus we can assume that $h$ is not identically zero. In this case we estimate the bound from Section 2.4.2 by

$$\sum_{I,d} c_{I,d} q^{\dim V_{I,d}}(q-1)^{\#I-n} a_d q^{\frac{t+d}{2}}$$
$$\leq \sum_{I,d} c_{I,d} a_d q^{\dim V_{I,d}+\frac{t+d}{2}}(q-1)^{\#I-n}.$$

Here we used an estimate due to Cluckers [Clu10, Thm. 7.4], which introduces constants $a_d \in \mathbb{R}_{>0}$ and excludes a finite set of field characteristics. From the proof of [Clu10, Thm. 7.4] and the references therein, mainly to [Clu08a, Cor. 6.1] which in turn invokes [Kat99, Thm. 4], it follows that $a_d$ can be taken to depend on $d$ and $\Delta_0(f)$ only. It also follows that the excluded field characteristics are contained in $P$, so this was already accounted for (by our assumption that $\operatorname{char} \mathbb{F}_q \notin P$).

Therefore it suffices to prove that for each $I, d$ we have

$$\dim V_{I,d} + \frac{t+d}{2} + \#I - n \leq -\sigma(f). \tag{2.4.4.1}$$

Now consider the following lemma:

**Lemma 2.4.1.** *Let $I \subseteq \{1, ..., s\}$ and define*

$$f_I = h + \sum_{i \in I} g_i x_i \in \mathbb{F}_q[(x_i)_{i \in I}, y, z].$$

*Then $\sigma(f) \leq \sigma(f_I) + s - \#I$. In particular we have that $\sigma(f) \leq \sigma(h) + s$.*

*Proof.* We proceed by induction on $s - \#I$. Note that if $s - \#I = 0$ then $f_I = f$ and there is nothing to prove. If $s - \#I \geq 1$ then consider an index $j$ which is not contained in $I$. By the induction hypothesis we know that $\sigma(f) \leq \sigma(f_{I \cup \{j\}}) + s - \#I - 1$. The lemma then follows by

$$\begin{aligned}
\sigma(f_{I \cup \{j\}}) &= \sigma(h + \sum_{i \in I} g_i x_i + g_j(z_1, \ldots, z_r) x_j) \\
&\leq \sigma(h + \sum_{i \in I} g_i x_i + g_j(u_1, \ldots, u_r) x_j) \\
&\leq \sigma(f_I) + \sigma(g_j(u_1, \ldots, u_r) x_j) \\
&\leq \sigma(f_I) + 1,
\end{aligned}$$

where $u_1, \ldots, u_r$ are new variables and the first two inequalities follow from properties (iii) resp. (i) stated in Lemma 2.2.1. □

From this we see that in order to prove (2.4.4.1), it is sufficient to demonstrate the following estimate:

**Theorem 2.4.2.** *If $h$ is a non-zero polynomial then we have*

$$\dim V_{I,d} + \frac{t+d}{2} + \#I - n \leq -\sigma(f_I) - s + \#I$$

*for all $I \subseteq \{1, \ldots, s\}$ and $0 \leq d \leq t$.*

The remainder of this section is devoted to proving this theorem. Note that the condition that $h$ is non-zero implies that $\sigma(f_I) < +\infty$ for each $I \subseteq \{1, \ldots, s\}$.

## 2.4.5

We first prove a weaker statement:

**Lemma 2.4.3** (weak version of Theorem 2.4.2)**.** *If $h$ is a non-zero polynomial then we have*

$$\frac{\dim V_{I,d}}{2} + \frac{t+d}{2} + \#I - n \leq -\sigma(f_I) - s + \#I$$

*for all $I \subseteq \{1, \ldots, s\}$ and $0 \leq d \leq t$.*

*Proof.* Consider the following algebraic subsets of $\mathbb{A}^{\#I + t + r}$:

$$\begin{aligned}
C_{f_I} &: \frac{\partial f_I}{\partial x_i} = \frac{\partial f_I}{\partial y_j} = \frac{\partial f_I}{\partial z_\ell} = 0 \text{ (i.e., the affine critical locus of } f_I\text{)}, \\
W_{I,d} &: \frac{\partial f_I}{\partial x_i} = \frac{\partial f_I}{\partial y_j} = 0, z \in \tilde{V}_{I,d}, \\
W_I &: \frac{\partial f_I}{\partial x_i} = \frac{\partial f_I}{\partial y_j} = 0.
\end{aligned}$$

Here $i \in I$, $1 \leq j \leq t$ and $1 \leq \ell \leq r$, and

$$\tilde{V}_{I,d} = \{\, z \in \mathbb{A}^r \mid \dim C_{h(\cdot, z)} = d, g_i(z) = 0 \Leftrightarrow i \in I \,\} \supseteq V_{I,d}.$$

It is easy to verify that $W_{I,d}$ has dimension $\dim \tilde{V}_{I,d} + \#I + d$ and is contained in $W_I$. On the other hand $\dim C_{f_I} \geq \dim W_I - r$, so we see that

$$\dim C_{f_I} \geq \dim \tilde{V}_{I,d} + \#I + d - r \geq \dim V_{I,d} + \#I + d - r.$$

Now because $f_I$ is non-degenerate with respect to the faces of its Newton polyhedron at the origin, by Lemma 2.2.2 we have

$$\frac{-r - t - \#I + \dim C_{f_I}}{2} \leq -\sigma(f_I).$$

Hence

$$\frac{-r - t - \#I + \dim V_{I,d} + \#I + d - r}{2} \leq -\sigma(f_I).$$

From this the lemma follows. □

This implies that Theorem 2.4.2 is true if $\dim V_{I,d} = 0$. In fact we will prove Theorem 2.4.2 by induction on $\dim V_{I,d}$, so this settles the base case. Note that the theorem is trivial if $\dim V_{I,d} = -\infty$.

## 2.4.6

The induction will rely on the following auxiliary lemma, which is the source of our enlargement of $P$, which we announced in Section 2.4.1.

**Lemma 2.4.4.** *Let $I \subseteq \{1, \ldots, s\}$. There exists a finite set of primes $P_{\text{comp}}$ which only depends on $\operatorname{Supp} f_I$, such that the following holds as soon as $\operatorname{char} \mathbb{F}_q \notin P_{\text{comp}}$. For each $a \in \overline{\mathbb{F}}_q$ let*

$$f_{I,a} = f_I((x_i)_{i \in I}, y_1, \ldots, y_t, z_1, \ldots, z_{r-1}, a)$$

*denote the polynomial obtained from $f_I$ by substituting $a$ for $z_r$. Then there exists a non-zero polynomial $\zeta \in \mathbb{F}_q[z_r]$ such that for all $a$ for which $\zeta(a) \neq 0$ we have that*
  *(i) $\sigma(f_I) \leq \sigma(f_{I,a})$,*
  *(ii) $g_i(z_1, \ldots, z_{r-1}, a)$ is not identically zero for each $i = 1, \ldots, s$,*
  *(iii) $f_{I,a}$ is non-degenerate with respect to the faces of its Newton polyhedron at the origin.*

*Proof.* Consider the following variation on the above assertion:

> There exists a non-zero polynomial $\zeta \in \mathbb{F}_q[z_r]$ such that for all $a \in \overline{\mathbb{F}}_q$ for which $\zeta(a) \neq 0$ we have that
> (a) each non-zero coefficient $c(z_r)$ of $f_I$, when viewed as a polynomial in $(x_i)_{i \in I}, y, z_1, \ldots, z_{r-1}$ over $\mathbb{F}_q[z_r]$, satisfies $c(a) \neq 0$,
> (b) $f_{I,a}$ is weakly non-degenerate with respect to the faces of $\Delta_0(f_{I,a})$.

These two properties imply (i–iii). Indeed, it is obvious that (a) implies (ii), while it also ensures that the Newton polyhedron of $f_{I,a}$ at the origin equals the image of $\Delta_0(f_I)$ under the projection $\pi_r : \mathbb{R}^{\#I+t+r} \to \mathbb{R}^{\#I+t+r-1}$ corresponding to dropping the last coordinate. In particular we have that

$$\Delta_0(f_I) \subseteq \Delta_0(f_{I,a}) \times \mathbb{R}_{\geq 0}$$

which implies (i). Finally, note that $f_I$ is supported on the hyperplane $H_I \subseteq \mathbb{R}^{\#I+t+r}$ defined by
$$\sum_{j \in I} c_j i_j + c_{s+1} i_{s+1} + \ldots + c_{s+t} i_{s+t} = b,$$
so that $f_{I,a}$ is supported on the hyperplane $\pi_r(H_I) \subseteq \mathbb{R}^{\#I+t+r-1}$, which is defined by the same equation. Therefore property (b) implies (iii), because it suffices to verify non-degeneracy with respect to the faces $\sigma \subseteq \Delta(f_{I,a})$ which are contained in $\pi_r(H_I)$, and as discussed in Section 2.2.2 the notions of non-degeneracy and weak non-degeneracy coincide with respect to such faces; here we used that $\operatorname{char} \mathbb{F}_q \notin P$.

Now it is easy to see that for each fixed choice of Supp $f_I$ the foregoing variation amounts to the validity of a first-order sentence. Therefore, if we can prove the direct analogue of this variation in characteristic zero, then by the compactness theorem (Robinson's principle) we know that it must be true in all characteristics outside a finite set $P_{\text{comp}}$, from which the desired conclusion follows.

So assume that we are working over a field $k$ of characteristic 0. That is, we consider a polynomial $f = h + g_1 x_1 + \ldots + g_s x_s \in k[x, y, z]$ satisfying the analogues of the properties listed in Section 2.4.1 and at the beginning of Section 2.4.2. As before we let $f_I = h + \sum_{i \in I} g_i x_i$ and for each $a \in \overline{k}$ we write $f_{I,a}$ for the polynomial obtained from $f_I$ by substituting $a$ for $z_r$. Our task is to show that properties (a) and (b) hold for all but finitely many $a \in \overline{k}$. Since for (a) this is immediate, from now on we assume that property (a) is satisfied, so that $\Delta_0(f_{I,a}) = \pi_r(\Delta_0(f_I))$. Our task is to prove that $f_{I,a}$ is weakly non-degenerate with respect to its Newton polyhedron at the origin, except possibly for another finite number of values of $a$. We note that a very similar problem was recently tackled by Esterov, Lemahieu and Takeuchi [ELT, Prop. 7.1], but their setting is more difficult since they assume weak non-degeneracy with respect to the compact faces, only.

Every face $\sigma \subseteq \Delta_0(f_{I,a})$ arises as the projection along $\pi_r$ of a face $\tau \subseteq \Delta_0(f_I)$ of the same codimension; see Figure 2.1 for some elementary examples. We will concentrate on the face $\sigma = \Delta_0(f_{I,a}) \cap \pi_r(H_I)$, which is the projection of $\tau = \Delta_0(f_I) \cap H_I$. One verifies that $f_{I,a}$ is weakly non-degenerate with respect to $\sigma$ if and only if the fiber of the projection
$$p_r : \{ ((x_i)_{i \in I}, y, z) \in \overline{k}^{*\#I+t+r} \mid f_I((x_i)_{i \in I}, y, z) = 0 \} \to \overline{k}$$
on the last coordinate is smooth. Note that the source is smooth, thanks to the non-degeneracy of $f_I$ with respect to $\tau$. By the generic smoothness theorem, which is only available in characteristic zero, this implies that all but finitely many fibers of $p_r$ are smooth. (More precisely, this conclusion follows by applying [Har77, Cor. III.10.7] to $p_r|_X$ for each irreducible component $X$ of the source.) We find that $f_{I,a}$ is weakly non-degenerate with respect to $\sigma$ for all but finitely many $a \in \overline{k}$. Repeating this argument for all other faces $\sigma \subseteq \Delta_0(f_{I,a})$ then leads to the desired conclusion. $\square$

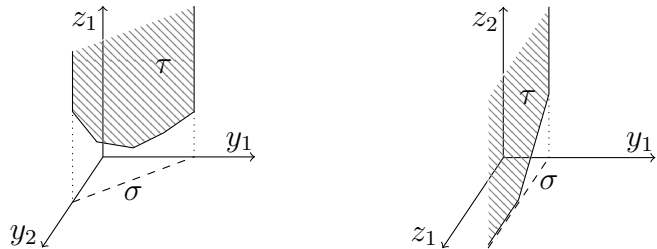

Figure 2.1 – Two examples of a face $\tau$ along with its projection $\sigma$

As announced in Section 2.4.1, we now enlarge the set $P$ of excluded primes by adjoining the sets $P_{\text{comp}}$ arising from multiple applications of Lemma 2.4.4. Indeed, the lemma will be applied for all possible choices of $I \subseteq \{1, \ldots, s\}$, but it will also be applied recursively to the polynomials $f_{I,a}$. Nevertheless, one easily sees that for a fixed choice of $\operatorname{Supp} f$ there are only a finite number of supports appearing: therefore $P$ remains finite.

### 2.4.7

We are now ready to prove Theorem 2.4.2. As announced, we will proceed by induction on $\dim V_{I,d}$, where the base case $\dim V_{I,d} = 0$ was taken care of in Section 2.4.5. So assume that $\dim V_{I,d} \geq 1$ and that Theorem 2.4.2 holds for all smaller dimensions.

If $r = 1$ then necessarily $\dim V_{I,d} = 1$, and because all $g_i(z) = g_i(z_1)$ are non-zero we have that $I = \emptyset$. In particular $f_I = h$. Now by Lemma 2.4.4 we can find an $a \in V_{I,d}$ such that $f_{I,a}$ is non-degenerate with respect to the faces of its Newton polyhedron at the origin, and such that $\sigma(f_I) \leq \sigma(f_{I,a})$. It is possible that we need to pick $a$ from $\overline{\mathbb{F}}_q \setminus \mathbb{F}_q$, but extending the coefficient field is of no harm to the statement we are trying to prove (or in other words, we can assume without loss of generality that $\#\mathbb{F}_q$ exceeds $\deg \zeta$). Our choice of $a$ implies that $C_{f_{I,a}}$ has dimension $d$, so by Lemma 2.2.2 we have that

$$\frac{-t+d}{2} \leq -\sigma(f_{I,a}) \leq \sigma(f_I),$$

which using $r = \dim V_{I,d} = 1$ and $\#I = 0$ implies that

$$\dim V_{I,d} + \frac{t+d}{2} + \#I - n \leq -\sigma(f_I) - s - \#I,$$

as wanted.

If $r \geq 2$ then we can proceed similarly. Indeed, since $\dim V_{I,d} \geq 1$ there exists at least one coordinate $z_j$ such that the image of the projection $\pi_j : V_{I,d} \to \mathbb{A}^1$ onto the $z_j$-coordinate is Zariski dense; we can suppose that $j = r$. Choose $a \in \pi_r(V_{I,d})$, again over an extension of $\mathbb{F}_q$ if needed, such that the fiber $V_{I,d,a}$ of $\pi_r$ over $a$ is of codimension 1 in $V_{I,d}$. By Lemma 2.4.4 we can moreover assume that $f_{I,a}$ is non-degenerate with respect to the faces of its Newton polyhedron at the origin and that $\sigma(f_I) \leq \sigma(f_{I,a})$. By applying our induction hypothesis we find that

$$\dim V_{I,d,a} + \frac{t+d}{2} + \#I - n + 1 \leq -\sigma(f_{I,a}) - s + \#I.$$

But we know that
$$\dim V_{I,d,a} = \dim V_{I,d} - 1,$$

therefore
$$\dim V_{I,d} + \frac{t+d}{2} + \#I - n \leq -\sigma(f_{I,a}) - s + \#I \leq -\sigma(f_I) - s + \#I,$$

as desired.

## 2.5 Denef and Hoornaert's conjecture

### 2.5.1

This final section is devoted to a proof of Denef and Hoornaert's conjecture:

**Theorem 2.5.1.** *Let $f \in \mathbb{Z}_K[x]$ be a non-zero polynomial such that $f(0) = 0$ and assume that it is weakly non-degenerate with respect to the faces of its Newton polyhedron at the origin. Let $\sigma$ and $\kappa$ be as in the statement of Theorem 2.1.3. Then*
$$\lim_{s \to -\sigma} (N\mathfrak{p}^{s+\sigma} - 1)^{\kappa + \delta_{\sigma,1}} Z_{f,\mathfrak{p}}(s) = O(N\mathfrak{p}^{1-\max\{1,\sigma\}})$$

*as $\mathfrak{p}$ varies over all non-zero prime ideals of $\mathbb{Z}_K$.*

By excluding finitely many prime ideals we can assume that $f$ is non-degenerate at $\mathfrak{p}$. Under this assumption Denef and Hoornaert proved the explicit formula
$$Z_{f,\mathfrak{p}}(s) = L_{\Delta_0(f)}(s) + \sum_{\substack{\tau \text{ proper face} \\ \text{of } \Delta_0(f)}} L_\tau(s) S_\tau(s)$$

where
$$L_\tau(s) = N\mathfrak{p}^{-n} \left( (N\mathfrak{p} - 1)^n - \#\{\, x \in \mathbb{F}_\mathfrak{p}^* \mid \overline{f}_\tau(x) = 0 \,\} \cdot \frac{N\mathfrak{p}^{s+1} - N\mathfrak{p}}{N\mathfrak{p}^{s+1} - 1} \right)$$

and
$$S_\tau(s) = \sum_{\substack{a \in \mathbb{Z}_{\geq 0}^n \\ \text{s.t. } F(a) = \tau}} N\mathfrak{p}^{-\nu(a) - N(a)s},$$

see [DH01, Thm. 4.2]. Moreover they showed that
$$\lim_{s \to -\sigma} (N\mathfrak{p}^{s+\sigma} - 1)^\kappa S_\tau(s) = 0$$

as soon as $\tau \not\subseteq \tau_0$, where $\tau_0 \subseteq \Delta_0(f)$ denotes the smallest face containing $(1/\sigma, \ldots, 1/\sigma)$, see [DH01, Lem. 5.4]. Thus it suffices to prove that
$$\lim_{s \to -\sigma} (N\mathfrak{p}^{s+\sigma} - 1)^{\kappa + \delta_{\sigma,1}} L_\tau(s) S_\tau(s) = O(N\mathfrak{p}^{1-\max\{1,\sigma\}}) \tag{2.5.1.1}$$

for all subfaces $\tau$ of $\tau_0$. We remark that Denef and Hoornaert restict their discussion to $K = \mathbb{Q}$, but everything readily generalizes to arbitrary number fields.

**2.5.2**

The face $\tau_0$ is contained in a hyperplane which does not contain the origin and which has a normal vector in $\mathbb{R}^n_{\geq 0}$. Necessarily the same is true for all subfaces $\tau \subseteq \tau_0$. This has an important consequence related to the subtlety mentioned in Section 2.1.6. Namely, Denef and Hoornaert work under the weak non-degeneracy assumption discussed in Section 2.2.2 (non-existence of singular points versus non-existence of critical points) and they also conjecture Theorem 2.5.1 in terms of this weaker hypothesis. But as mentioned at the end of Section 2.2.2, over fields of large enough characteristic, both non-degeneracy notions coincide with respect to faces that are contained in a hyperplane not passing through the origin. Therefore there is no ambiguity: proving (2.5.1.1) under the assumption that $f_{\tau_0}$ satisfies our stronger non-degeneracy assumption is sufficient to conclude the Denef–Hoornaert conjecture in its full strength.

**2.5.3**

The proof works by explicit computation along the lines of Denef and Hoornaert, making a distinction between the cases $\sigma < 1$, $\sigma > 1$ and $\sigma = 1$. We make use of two new plug-ins. One of these plug-ins is our finite field exponential sum estimate stated in Theorem 2.1.4, which implies:

**Lemma 2.5.2.** $N\mathfrak{p} \cdot \#\{\, x \in \mathbb{F}^*_\mathfrak{p} \mid \overline{f}_\tau(x) = 0 \,\} = (N\mathfrak{p} - 1)^n + O(N\mathfrak{p}^{n+1-\sigma(f_\tau)})$.

*Proof.* Let $\varphi : \mathbb{F}_\mathfrak{p} \to \mathbb{C}^*$ be a non-trivial additive character, then

$$N\mathfrak{p} \cdot \#\{\, x \in \mathbb{F}^*_\mathfrak{p} \mid \overline{f}_\tau(x) = 0 \,\} = (N\mathfrak{p} - 1)^n + \sum_{t \in \mathbb{F}^*_\mathfrak{p}} \sum_{x \in \mathbb{F}^{*n}_\mathfrak{p}} \varphi(t\overline{f}_\tau(x))$$

and because $\tau$ is contained in a hyperplane which does not contain the origin and which has a normal vector in $\mathbb{R}^n_{\geq 0}$, we can use Theorem 2.1.4 to find that

$$\left| \sum_{t \in \mathbb{F}^*_\mathfrak{p}} \sum_{x \in \mathbb{F}^{*n}_\mathfrak{p}} \varphi(t\overline{f}_\tau(x)) \right| \leq \sum_{t \in \mathfrak{p}^*} c(N\mathfrak{p} - 1)^n N\mathfrak{p}^{-\sigma(f_\tau)}$$

for some constant $c \in \mathbb{R}_{>0}$ that does not depend on $\mathfrak{p}$. The lemma follows. □

The second plug-in is a combinatorial inequality due to Cluckers [Clu08a, Thm. 4.1], which says that for any face $\tau \subseteq \Delta_0(f)$ and all $a \in \mathbb{R}^n_{\geq 0}$ such that $F(a) = \tau$ we have $\nu(a) \geq \sigma N(a) + \sigma - \sigma(f_\tau)$. This leads to the following statement:

**Lemma 2.5.3.** *If $\tau \subseteq \tau_0$ then*

$$\lim_{s \to -\sigma} (N\mathfrak{p}^{s+\sigma} - 1)^\kappa S_\tau(s) = O(N\mathfrak{p}^{\sigma(f_\tau)-\sigma}).$$

*Proof.* Denef and Hoornaert proved [DH01, Lem. 5.4] that

$$\lim_{s \to -\sigma} (N\mathfrak{p}^{s+\sigma} - 1)^\kappa S_{\tau_0}(s) = c_0$$

for some constant $c_0 \in \mathbb{R}_{>0}$ which does not depend on $\mathfrak{p}$. Since $\sigma(f_{\tau_0}) = \sigma$ this settles the case $\tau = \tau_0$. If $\tau \subsetneq \tau_0$ then we can redo the proof of [DH01, Lem. 5.16], but instead of invoking the Denef-Sperber inequality stated in [DH01, Lem. 5.15] we use Cluckers' result. □

## 2.5.4

We can now conclude. If $\sigma < 1$ then

$$\frac{N\mathfrak{p}^{-\sigma} - 1}{N\mathfrak{p}^{1-\sigma} - 1} \tag{2.5.4.1}$$

is a negative quantity which is in $O(N\mathfrak{p}^{\sigma-1})$. Together with Lemma 2.5.2 this implies that

$$L_\tau(-\sigma) = O(N\mathfrak{p}^{\sigma - \sigma(f_\tau)}),$$

which along with Lemma 2.5.3 implies that

$$\lim_{s \to -\sigma}(N\mathfrak{p}^{s+\sigma} - 1)^\kappa Z_{f,\mathfrak{p}}(s) = O(1)$$

as wanted. If $\sigma > 1$ then (2.5.4.1) becomes a positive quantity which equals $1 + o(1)$, leading to the estimate

$$L_\tau(-\sigma) = O(N\mathfrak{p}^{1-\sigma(f_\tau)}),$$

so here we find that

$$\lim_{s \to -\sigma}(N\mathfrak{p}^{s+\sigma} - 1)^\kappa Z_{f,\mathfrak{p}}(s) = O(N\mathfrak{p}^{1-\sigma}),$$

again as wanted. Finally if $\sigma = 1$ then using Lemma 2.5.2 one finds that

$$\lim_{s \to -1}(N\mathfrak{p}^{s+1} - 1)L_\tau(s) = O(N\mathfrak{p}^{1-\sigma(f_\tau)})$$

which together with Lemma 2.5.3 shows that

$$\lim_{s \to -1}(N\mathfrak{p}^{s+1} - 1)^{\kappa+1} Z_{f,\mathfrak{p}}(s) = O(1),$$

thereby concluding the proof of Theorem 2.5.1.

# Chapter 3

# Proof of Cluckers-Veys's conjecture on exponential sums for polynomials with log-canonical threshold at most a half

*This chapter is a joint work with Saskia Chambilles, see [CNb].*

## Contents

## Abstract

In this paper, we will give two proofs of the Cluckers-Veys conjecture on exponential sums for the case of polynomials in $\mathbb{Z}[x_1, \ldots, x_n]$ having log-canonical thresholds at most one half. In particular, these results imply Igusa's conjecture and Denef-Sperber's conjecture under the same restriction on the log-canonical threshold.

## 3.1 Introduction

Let $n \geq 1$ be a natural number and let $f \in \mathbb{Z}[x_1, \ldots, x_n]$ be a non-constant polynomial in $n$ variables, for which we assume that $f(0, \ldots, 0) = 0$. For homogeneous polynomials $f$, Igusa has formulated, on page 2 of [Igu78], a conjecture on the exponential sum

$$E_{m,p}(f) := \frac{1}{p^{mn}} \sum_{\overline{x} \in (\mathbb{Z}/p^m\mathbb{Z})^n} \exp\left(\frac{2\pi i f(x)}{p^m}\right),$$

where $p$ is a prime number and $m \in \mathbb{N}$. More precisely, he predicted that there exist a constant $\sigma$, which depends on the geometric properties of $f$, and a positive constant $C$, independent of $p$ and $m$, such that for all primes $p$ and for all $m \geq 1$,

$$|E_{m,p}(f)| \leq Cm^{n-1}p^{-m\sigma}.$$

In particular, his conjecture implies an adèlic Poisson summation formula.

A local version of this sum,

$$E_{m,p}^0(f) := \frac{1}{p^{mn}} \sum_{\overline{x} \in (p\mathbb{Z}/p^m\mathbb{Z})^n} \exp\left(\frac{2\pi i f(x)}{p^m}\right),$$

was considered by Denef and Sperber in [DS01]. Under certain conditions on the Newton polyhedron $\Delta$ of $f$, they proved that there exist constants $\sigma, \kappa$, depending only on $\Delta$, and a positive constant $C$, independent of $p$ and $m$, such that for all $m \geq 1$ and almost all $p$, we have

$$|E_{m,p}^0(f)| \leq Cm^{\kappa-1}p^{-m\sigma}.$$

In [Clu08a], Cluckers proved both conjectures in the case that $f$ is non-degenerate.

To generalise these facts, Cluckers and Veys formulated, in [CV16], a conjecture related to the log-canonical threshold of an arbitrary polynomial $f$. We will recall the definition of the log-canonical threshold in the next section. They also introduced the following local exponential sum, for each $y \in \mathbb{Z}^n$:

$$E_{m,p}^y(f) := \frac{1}{p^{mn}} \sum_{\overline{x} \in \overline{y} + (p\mathbb{Z}/p^m\mathbb{Z})^n} \exp\left(\frac{2\pi i f(x)}{p^m}\right).$$

We restate their conjecture here.

**Conjecture 3.1.1** (Cluckers-Veys)**.** *There exists a positive constant $C$ (that may depend on the polynomial $f$), such that for all primes $p$, for all $m \geq 2$ and for all $y \in \mathbb{Z}^n$, we have*

$$|E_{m,p}(f)| \leq Cm^{n-1}p^{-ma(f)}$$

*and*

$$|E_{m,p}^y(f)| \leq Cm^{n-1}p^{-ma_{y,p}(f)}.$$

*Here $a(f)$ is the minimum, over all $b \in \mathbb{C}$, of the log-canonical thresholds of the polynomials $f(x) - b$. And, for $y \in \mathbb{Z}^n$, $a_{y,p}(f)$ is the minimum of the log-canonical thresholds at $y'$ of the polynomials $f(x) - f(y')$, where $y'$ runs over $y + p\mathbb{Z}_p^n$.*

In this article, we will prove a special case of the Cluckers-Veys conjecture. More concretely, we will prove the case in which the log-canonical threshold of $f$ is at most a half. We will consider in detail the local sum where $y = 0$ and we will afterwards discuss how one can adapt the proofs to obtain uniform upper bounds for $|E_{m,p}^y(f)|$, for $y \in \mathbb{Z}^n$, and an upper bound for $|E_{m,p}(f)|$. Our main theorems will be the following.

**Theorem 3.1.2.** *Let $n \geq 1$ and let $f \in \mathbb{Z}[x_1, \ldots, x_n]$ be a non-constant polynomial with $f(0) = 0$. Put $\sigma = \min\left\{c_0(f), \frac{1}{2}\right\}$, where $c_0(f)$ is de log-canonical threshold of $f$ at $0$. Then there exists a positive constant $C$, not depending on $p$ and $m$, and a natural number $N$, such that for all $m \geq 1$ and all primes $p > N$, we have*

$$|E_{m,p}^0(f)| \leq C m^{n-1} p^{-m\sigma}.$$

**Theorem 3.1.3.** *Let $n \geq 1$ and let $f \in \mathbb{Z}[x_1, \ldots, x_n]$ be a non-constant polynomial. Put $\sigma = \min\left\{a(f), \frac{1}{2}\right\}$. Then there exists a positive constant $C$, not depending on $p$ and $m$, and a natural number $N$, such that for all $m \geq 2$ and all primes $p > N$, we have*

$$|E_{m,p}(f)| \leq C m^{n-1} p^{-m\sigma}.$$

Remark that by [DS01], [Igu74] and [Igu78], there exists, for each prime $p$, a positive constant $C_p$, such that

$$|E_{m,p}^0(f)| \leq C_p m^{n-1} p^{-mc_0(f)}$$

and

$$|E_{m,p}(f)| \leq C_p m^{n-1} p^{-ma(f)},$$

for all $m \geq 1$. Therefore we know that once the Main Theorems 3.1.2 and 3.1.3 are proven, they will hold for $N = 1$, possibly after enlarging the constant $C$.

Notice that the homogeneous polynomials $f$ in two variables that are not yet covered by Igusa in [Igu78], all satisfy that $a(f) \leq \frac{1}{2}$. Hence our results can be seen as a generalisation of a result of Lichtin from [Lic13] or of Wright from [Wri], in which they proved Igusa's conjecture for any homogeneous polynomial of two variables.

**Remark 3.1.4.** We observe that if $m = 1$, then $|E_{1,p}^0(f)| = \frac{1}{p^n}$. Hence the Main Theorem 3.1.2 is trivial for $m = 1$ and we only need to prove it for $m \geq 2$.

We will give two approaches to our main theorems and for the Main Theorem 3.1.2 we will give the details of these approaches. The first approach, in Section 3.3, will make use of model theory, an estimate of the dimension of arc spaces as in [Mus02], the Cluckers-Loeser motivic integration theory and an estimate of Weil on finite field exponential sums in one variable (see [Wei48]). We will also use an idea which is close to the construction of the local Artin map by Lubin-Tate theory. More concretely, we will prove that certain functions do not depend on the choice of a uniformiser in $\mathbb{Q}_p$, but only on the angular component of the chosen uniformiser. Hence, when varying uniformisers, we obtain orbits of points that have the same image under these functions. In fact, these orbits depend on actions of the group $\mu_{p-1}(\mathbb{Q}_p)$, the group of $(p-1)^{\text{th}}$ roots of unity of $\mathbb{Q}_p$, on the set of uniformisers of $\mathbb{Q}_p$ and on $\mathbb{Q}_p$. The second approach, in Section 3.4, will use a concrete expression of cohomology, as in [Den91]. Both of these approaches will use not only Lang-Weil estimates ([LW54]) for the number of points on varieties over finite fields, but also the theory of Igusa's local zeta functions. In Section 3.2 we will give some background on log-canonical thresholds, exponential sums and Igusa's local zeta functions. In Section 3.5, we will explain how the results from Section 3.4 can be used to prove the Main Theorem 3.1.3. We will end this paper by explaining, in Section 3.6, how to obtain uniform upper bounds for all local sums $E_{m,p}^y$. We will do this both from the geometric, as well as from the model theoretic point of view.

We remark that our results can be extended to the ring of integers $\mathcal{O}_K$ of any number field $K$, but we only work with $\mathbb{Z}$ and $\mathbb{Q}$ to simplify notation.

## 3.2 Log-canonical Thresholds and exponential sums

### 3.2.1 Log-canonical Threshold

In this section we will recall two possible definitions of the log-canonical threshold of a polynomial $f$.

**Definition 3.2.1.** Let $f$ be a non-constant polynomial in $n$ variables over an algebraically closed field $K$ of characteristic zero. Let $\pi : Y \to K^n$ be a proper birational morphism on a smooth variety $Y$. For any prime divisor $E$ on $Y$, we denote by $N$ and $\nu - 1$ the multiplicities along $E$ of the divisors of $\pi^* f$ and $\pi^*(dx_1 \wedge \ldots \wedge dx_n)$, respectively. For each $x \in Z(f) \subset K^n$, the *log-canonical threshold of $f$ at $x$*, denoted by $c_x(f)$, is the real number $\inf_{\pi, E} \left\{ \frac{\nu}{N} \right\}$, where $\pi$ runs over all proper birational morphisms to $K^n$ and $E$ runs over all prime divisors on $Y$ such that $x \in \pi(E)$. If we fix any embedded resolution $\pi$ of the germ of $f = 0$ at $x$, then

$$c_x(f) = \min_{E : x \in \pi(E)} \left\{ \frac{\nu}{N} \right\}.$$

Furthermore we always have $c_x(f) \leq 1$. We denote by $c(f) = \inf_{x \in Z(f)} c_x(f)$ the *log-canonical threshold of $f$*.

By the following theorem from [Mus02], which is true for any algebraically closed field $K$ of characteristic zero, there exists a description of the log-canonical threshold in terms of arc spaces and jet spaces.

**Theorem 3.2.2** ([Mus02], Corollaries 0.2 and 3.6). *Let $f$ be a non-constant polynomial over $K$ in $n$ variables and let $m$ be a natural number. We set*

$$\text{Cont}^{\geq m}(f) := \{ x \in K[[t]]^n \mid f(x) \equiv 0 \bmod t^m \}$$

*and*

$$\text{Cont}_0^{\geq m}(f) := \{ x \in (tK[[t]])^n \mid f(x) \equiv 0 \bmod t^m \}.$$

*We denote by $\pi_m$ the projection from $K[[t]]^n$ to $(K[t]/(t^m))^n$ and we consider the codimensions of $\pi_m(\text{Cont}^{\geq m}(f))$ and $\pi_m(\text{Cont}_0^{\geq m}(f))$ in $(K[t]/(t^m))^n \cong K^{nm}$. We* *denote these two values by* codim Cont$^{\geq m}(f)$ *and* codim Cont$_0^{\geq m}(f)$, *respectively. Then the log-canonical threshold of $f$ equals the real number*

$$c(f) = \inf_{m \geq 1} \frac{\text{codim Cont}^{\geq m}(f)}{m},$$

*and if $f(0) = 0$, then the log-canonical threshold of $f$ at 0 equals the real number*

$$c_0(f) = \inf_{m \geq 1} \frac{\text{codim Cont}_0^{\geq m}(f)}{m}.$$

### 3.2.2 Exponential sum and Igusa local zeta function

In this section we will discuss formulas for the exponential sums $E_{m,p}(f)$ and $E_{m,p}^0(f)$. These formulas can be found in the works of Igusa and Denef on Igusa local zeta functions. Most of the theory in this section comes from [Den]. We will just introduce the necessary notation here.

Let $K$ be a number field, $\mathcal{O}$ the ring of algebraic integers of $K$ and $\mathfrak{p}$ any maximal ideal of $\mathcal{O}$. We denote the completions of $K$ and $\mathcal{O}$ with respect to $\mathfrak{p}$ by $K_\mathfrak{p}$ and $\mathcal{O}_\mathfrak{p}$. Let $q = p^m$ be the cardinality of the residue field $k_\mathfrak{p}$ of the local ring $\mathcal{O}_\mathfrak{p}$, then $k_\mathfrak{p} = \mathbb{F}_q$. For $x \in K_\mathfrak{p}$, we denote by $\text{ord}(x) \in \mathbb{Z} \cup \{+\infty\}$ the $\mathfrak{p}$-valuation of $x$, $|x| = q^{-\text{ord}(x)}$ and $\text{ac}(x) = x\pi^{-\text{ord}(x)}$, where $\pi \in \mathcal{O}_\mathfrak{p}$ is a fixed uniformising parameter for $\mathcal{O}_\mathfrak{p}$.

Let $\chi : \mathcal{O}_\mathfrak{p}^\times \to \mathbb{C}^\times$ be a character on the group of units $\mathcal{O}_\mathfrak{p}^\times$ of $\mathcal{O}_\mathfrak{p}$, with finite image. By the *order* of such a character we mean the number of elements in its image. The *conductor* $c(\chi)$ of the character is the smallest $c \geq 1$ for which $\chi$ is trivial on $1 + \mathfrak{p}^c$. We formally put $\chi(0) = 0$. Let $f(x) \in K[x]$ be a polynomial in $n$ variables, $x = (x_1, \ldots, x_n)$, with $f \neq 0$, and let $\Phi : K_\mathfrak{p}^n \to \mathbb{C}$ be a *Schwartz-Bruhat function*, i.e., a locally constant function with compact support. We say that $\Phi$ is *residual* if $\text{Supp}(\Phi) \subset \mathcal{O}_\mathfrak{p}^n$ and $\Phi(x)$ only depends on $x$ mod $\mathfrak{p}$. Thus if $\Phi$ is residual, it induces a function $\overline{\Phi} : k_\mathfrak{p}^n \to \mathbb{C}$. Now we associate to these data *Igusa's local zeta function*

$$Z_\Phi(K_\mathfrak{p}, \chi, s, f) := \int_{K_\mathfrak{p}^n} \Phi(x)\chi\big(\text{ac}(f(x))\big)|f(x)|^s|dx|.$$

In [Igu78], Igusa showed that $Z_\Phi(K_\mathfrak{p}, \chi, s, f)$ is a rational function in $t = q^{-s}$. From now on we will write $Z_\Phi(\mathfrak{p}, \chi, s, f)$, whenever we have fixed $K$.

Let $\Psi$ be the standard additive character on $K_\mathfrak{p}$, i.e. for $z \in K_\mathfrak{p}$,

$$\Psi(z) := \exp(2\pi i \text{Tr}_{K_\mathfrak{p}/\mathbb{Q}_p}(z)),$$

where $\text{Tr}_{K_\mathfrak{p}/\mathbb{Q}_p}$ denotes the trace map. We set

$$E_\Phi(z, \mathfrak{p}, f) := \int_{K_\mathfrak{p}^n} \Phi(x)\Psi(zf(x))|dx|.$$

Whenever $\Phi = \mathbf{1}_{\mathcal{O}_\mathfrak{p}^n}$ or $\Phi = \mathbf{1}_{(\mathfrak{p}\mathcal{O}_\mathfrak{p})^n}$ and $K$ is fixed, we will simply denote this function by $E_\mathfrak{p}(z, f)$ or $E_\mathfrak{p}^0(z, f)$, respectively. When $K = \mathbb{Q}$, $\mathfrak{p} = p\mathbb{Z}$, $z = p^{-m}$ and $\Phi = \mathbf{1}_{\mathbb{Z}_p^n}$ or $\Phi = \mathbf{1}_{(p\mathbb{Z}_p)^n}$ we will simplify notation even more by writing $E_{m,p}(f)$ or $E_{m,p}^0(f)$, respectively, and this notation coincide with the notation in Section 3.1 by an easy calculation.

We recall the following proposition from [Den], that relates the exponential sums to Igusa's local zeta functions.

**Proposition 3.2.3** ([Den], Proposition 1.4.4). *Let $u \in \mathcal{O}_\mathfrak{p}^\times$ and $m \in \mathbb{Z}$. Then $E_\Phi(u\pi^{-m}, \mathfrak{p}, f)$ is equal to*

$$Z_\Phi(\mathfrak{p}, \chi_{\text{triv}}, 0, f) + \text{Coeff}_{t^{m-1}}\left(\frac{(t-q)Z_\Phi(\mathfrak{p}, \chi_{\text{triv}}, s, f)}{(q-1)(1-t)}\right)$$
$$+ \sum_{\chi \neq \chi_{\text{triv}}} g_{\chi^{-1}} \chi(u) \text{Coeff}_{t^{m-c(\chi)}}\left(Z_\Phi(\mathfrak{p}, \chi, s, f)\right),$$

*where $g_\chi$ is the Gaussian sum*

$$g_\chi = \frac{q^{1-c(\chi)}}{q-1} \sum_{\overline{v} \in (\mathcal{O}_\mathfrak{p}/\mathfrak{p}^{c(\chi)})^\times} \chi(v)\Psi(v/\pi^{c(\chi)}).$$

Now we will describe a formula for Igusa's local zeta function using resolution of singularities

Let $K$ and $f$ be as above. Put $X = \text{Spec}\, K[x]$ and $D = \text{Spec}\, K[x]/(f)$. We take an embedded resolution $(Y, h)$ for $f^{-1}(0)$ over $K$. This means that $Y$ is an integral smooth closed subscheme of projective space over $X$, $h : Y \to X$ is the natural map, the restriction $h : Y \backslash h^{-1}(D) \to X \backslash D$ is an isomorphism, and $(h^{-1}(D))_{\text{red}}$ has only normal crossings as subscheme of $Y$. Let $E_i, i \in T$, be the irreducible components of $(h^{-1}(D))_{\text{red}}$. For each $i \in T$, let $N_i$ be the multiplicity of $E_i$ in the divisor of $f \circ h$ on $Y$ and let $\nu_i - 1$ be the multiplicity of $E_i$ in the divisor of $h^*(dx_1 \wedge \ldots \wedge dx_n)$. The $(N_i, \nu_i)_{i \in T}$ are called the *numerical data* of the resolution. For each subset $I \subset T$, we consider the schemes

$$E_I := \cap_{i \in I} E_i \quad \text{and} \quad \overset{\circ}{E}_I := E_I \backslash \cup_{j \in T \backslash I} E_j.$$

In particular, when $I = \emptyset$ we have $E_\emptyset = Y$. We denote the critical locus of $f$ by $C_f$.

If $Z$ is a closed subscheme of $Y$, we denote the reduction mod $\mathfrak{p}$ of $Z$ by $\overline{Z}$ (see [Den87, Definition 2.2] and [Shi55] when $\overline{Z}$ is reduce). We say that the resolution $(Y, h)$ of $f$ has *good reduction modulo* $\mathfrak{p}$ if $\overline{Y}$ and all $\overline{E}_i$ are smooth, $\cup_{i \in T} \overline{E}_i$ has only normal crossings, and the schemes $\overline{E}_i$ and $\overline{E}_j$ have no common components whenever $i \neq j$. There exists a finite subset $S$ of $\text{Spec}\,\mathcal{O}$, such that for all $\mathfrak{p} \notin S$, we have $f \in \mathcal{O}_\mathfrak{p}[x]$, $f \not\equiv 0 \mod \mathfrak{p}$ and the resolution $(Y, h)$ for $f$ has good reduction mod $\mathfrak{p}$ (see [Den87], Theorem 2.4).

Let $\mathfrak{p} \notin S$ and $I \subset T$, then it is easy to prove that $\overline{E}_I = \cap_{i \in I} \overline{E}_i$. We put $\overset{\circ}{\overline{E}}_I := \overline{E}_I \backslash \cup_{j \notin I} \overline{E}_j$. Let $a$ be a closed point of $\overline{Y}$ and $T_a = \{i \in T | a \in \overline{E}_i\}$. In the local ring of $\overline{Y}$ at $a$ we can write

$$\overline{f} \circ \overline{h} = \overline{u} \prod_{i \in T_a} \overline{g}_i^{N_i},$$

where $\bar{u}$ is a unit, $(\bar{g}_i)_{i \in T_a}$ is a part of a regular system of parameters and $N_i$ is as above.

In two cases, depending on the conductor $c(\chi)$ of the character $\chi$, we will give a more explicit description of Igusa's zeta function $Z_\Phi(\mathfrak{p}, \chi, s, f)$. In the first case we consider a character $\chi$ on $\mathcal{O}_\mathfrak{p}^\times$ of order $d$, which is trivial on $1 + \mathfrak{p}\mathcal{O}_\mathfrak{p}$, i.e., $c(\chi) = 1$. Then $\chi$ induces a character (denoted also by $\chi$) on $k_\mathfrak{p}^\times$. We define a map

$$\Omega_\chi : \overline{Y}(k_\mathfrak{p}) \to \mathbb{C}$$

as follows. Let $a \in \overline{Y}(k_\mathfrak{p})$. If $d | N_i$ for all $i \in T_a$, then we put $\Omega_\chi(a) = \chi(\bar{u}(a))$, otherwise we put $\Omega_\chi(a) = 0$. This definition is independent of the choice of $\bar{g}_i$. In the following theorem we recall the formula of Igusa's local zeta function.

**Theorem 3.2.4** ([Den85], Theorem 2.2 or [Den], Theorem 3.4). *Let $\chi$ be a character on $\mathcal{O}_\mathfrak{p}^\times$ of order $d$, which is trivial on $1 + \mathfrak{p}\mathcal{O}_\mathfrak{p}$. Supose that $\mathfrak{p} \notin S$ and that $\Phi$ is residual, then we have*

$$Z_\Phi(\mathfrak{p}, \chi, s, f) = q^{-n} \sum_{\substack{I \subset T, \\ \forall i \in I : d | N_i}} c_{I,\Phi,\chi} \prod_{i \in I} \frac{(q-1)q^{-N_i s - \nu_i}}{1 - q^{-N_i s - \nu_i}},$$

*where*

$$c_{I,\Phi,\chi} = \sum_{a \in \overset{\circ}{\overline{E}}_I(k_\mathfrak{p})} \overline{\Phi}(\bar{h}(a)) \Omega_\chi(a).$$

If $\Phi = \mathbf{1}_{\mathcal{O}_\mathfrak{p}^n}$ or $\Phi = \mathbf{1}_{(\mathfrak{p}\mathcal{O}_\mathfrak{p})^n}$ we will denote $c_{I,\Phi,\chi}$ by $c_{I,\chi}$ or $c_{I,\chi}^0$, respectively.

We note that $c_{I,\Phi,\chi} = 0$, if there exists $i \in I$, such that $d \nmid N_i$. Therefore the number of characters $\chi$, for which $c(\chi) = 1$ and $c_{I,\Phi,\chi} \neq 0$ for some $I \subset T$, will have an upper bound $M$, which will only depend on the numerical data of $(Y, h)$, hence does not depend on $\text{char}(k_\mathfrak{p})$.

Now in the second case we consider a character $\chi$ on $\mathcal{O}_\mathfrak{p}^\times$, which is non-trivial on $1 + \mathfrak{p}\mathcal{O}_\mathfrak{p}$, i.e. $c(\chi) > 1$. Then we have the following theorem by Denef.

**Theorem 3.2.5** ([Den85], Theorem 2.1 or [Den], Theorem 3.3). *Let $\chi$ be a character on $\mathcal{O}_\mathfrak{p}^\times$, which is non-trivial on $1 + \mathfrak{p}\mathcal{O}_\mathfrak{p}$. Suppose that $\Phi$ is residual, $\mathfrak{p} \notin S$, $N_i \notin \mathfrak{p}$ for all $i \in T$, and $C_{\overline{f}} \cap \text{Supp}(\overline{\Phi}) \subset \overline{f}^{-1}(0)$. Then $Z_\Phi(\mathfrak{p}, \chi, s, f) = 0$*

As a consequence of these results, one can obtain the following description of the exponential sums $E_\Phi(z, \mathfrak{p}, f)$. This result and its proof are very similar to that of Corollary 1.4.5 from [Den].

**Corollary 3.2.6.** *Suppose that $\Phi$ is residual, $\mathfrak{p} \notin S$, $N_i \notin \mathfrak{p}$ for all $i \in T$, and $C_{\overline{f}} \cap \text{Supp}(\overline{\Phi}) \subset \overline{f}^{-1}(0)$. Then $E_\Phi(z, \mathfrak{p}, f)$ is a finite $\mathbb{C}$-linear combination of functions of the form $\chi(\text{ac}(z))|z|^\lambda (\log_q |z|)^\beta$ with coefficients independent of $z$, where $\lambda \in \mathbb{C}$ is a pole of $(q^{s+1} - 1)Z_\Phi(\mathfrak{p}, \chi_{\text{triv}}, s, f)$ or of $Z_\Phi(\mathfrak{p}, \chi, s, f)$, for $\chi \neq \chi_{\text{triv}}$, and $\beta \in \mathbb{N}$, such that $\beta \leq (\text{multiplicity of pole } \lambda) - 1$, provided that $|z|$ is big enough.*

*Proof.* It is easy to prove by combining the Theorems 3.2.3, 3.2.4 and 3.2.5. □

## 3.3 The first approach by Model theory

The first part of this section will contain some background on the theory of motivic integration. For the details we refer to [CL08] or [CL05]. In the second part we will use this theory to give our first proof of the Main Theorem 3.1.2.

### 3.3.1 Constructible Motivic Functions

**The language of Denef-Pas**

Let $K$ be a valued field, with valuation map $\operatorname{ord} : K^\times \to \Gamma_K$ for some additive ordered group $\Gamma_K$, and let $\mathcal{O}_K$ be the valuation ring of $K$ with maximal ideal $\mathcal{M}_K$ and residue field $k_K$. We denote by $x \to \overline{x}$ the projection $\mathcal{O}_K \to k_K$ modulo $\mathcal{M}_K$. An angular component map (modulo $\mathcal{M}_K$) on $K$ is any multiplicative map $\overline{\mathrm{ac}} : K^\times \to k_K^\times$ satisfying $\overline{\mathrm{ac}}(x) = \overline{x}$ for all $x$ with $\operatorname{ord}(x) = 0$. It can be extended to $K$ by putting $\overline{\mathrm{ac}}(0) = 0$.

The *language $\mathcal{L}_{\mathrm{DP}}$ of Denef-Pas* is the three-sorted language

$$(\mathcal{L}_{\mathrm{ring}}, \mathcal{L}_{\mathrm{ring}}, \mathcal{L}_{\mathrm{oag}}, \operatorname{ord}, \overline{\mathrm{ac}})$$

with as sorts:
   (i) a sort VF for the valued field-sort,
   (ii) a sort RF for the residue field-sort, and
   (iii) a sort VG for the value group-sort.

The first copy of $\mathcal{L}_{\mathrm{ring}}$ is used for the sort VF, the second copy for RF and the language $\mathcal{L}_{\mathrm{oag}}$, the language $(+, <)$ of ordered abelian groups, is used for VG. Furthermore ord denotes the valuation map from non-zero elements of VF to VG, and $\overline{\mathrm{ac}}$ stands for an angular component map from VF to RF.

As usual for first order formulas, $\mathcal{L}_{\mathrm{DP}}$-formulas are built up from the $\mathcal{L}_{\mathrm{DP}}$-symbols together with variables, the logical connectives $\wedge$ (and), $\vee$ (or), $\neg$ (not), the quantifiers $\exists, \forall$, the equality symbol $=$, and possibly parameters (see [Pas89] for more details).

Let us briefly recall the statement of the Denef-Pas theorem on elimination of valued field quantifiers in the language $\mathcal{L}_{\mathrm{DP}}$. Denote by $H_{\overline{\mathrm{ac}},0}$ the $\mathcal{L}_{\mathrm{DP}}$-theory of the above described structures whose valued field is Henselian and whose residue field is of characteristic zero. Then the theory $H_{\overline{\mathrm{ac}},0}$ admits elimination of quantifiers in the valued field sort, as stated in the following theorem.

**Theorem 3.3.1** (Pas, [Pas89])**.** *The theory $H_{\overline{\mathrm{ac}},0}$ admits elimination of quantifiers in the valued field sort. More precisely, every $\mathcal{L}_{\mathrm{DP}}$-formula $\phi(x, \xi, \alpha)$ (without parameters), with $x$ denoting variables in the VF-sort, $\xi$ variables in the RF-sort and $\alpha$ variables in the VG-sort, is $H_{\overline{\mathrm{ac}},0}$-equivalent to a finite disjunction of formulas of the form*

$$\psi\big(\overline{\mathrm{ac}} f_1(x), \ldots, \overline{\mathrm{ac}} f_k(x), \xi\big) \wedge \vartheta\big(\operatorname{ord} f_1(x), \ldots, \operatorname{ord} f_k(x), \alpha\big),$$

*where $\psi$ is an $\mathcal{L}_{\mathrm{ring}}$-formula, $\vartheta$ an $\mathcal{L}_{\mathrm{oag}}$-formula and $f_1, \ldots, f_k$ polynomials in $\mathbb{Z}[X]$.*

This theorem implies the following, useful corollary.

**Corollary 3.3.2** ([CL08], Corollary 2.1.2). *Let $(K, k, \Gamma)$ be a model of the theory $H_{\overline{ac},0}$ and $S$ a subring of $K$. Let $T_S$ be the set of atomic $\mathcal{L}_{\mathrm{DP}} \cup S$-sentences and negations of atomic sentences $\varphi$ such that $S \models \varphi$. We take $H_S$ to be the union of $H_{\overline{ac},0}$ and $T_S$. Then Theorem 3.3.1 holds with $H_{\overline{ac},0}$ replaced by $H_S$, $\mathcal{L}_{\mathrm{DP}}$ replaced by $\mathcal{L}_{\mathrm{DP}} \cup S$, and $\mathbb{Z}[X]$ replaced by $S[X]$.*

It is important to remark that by compactness, this theorem and its corollary are still true in the case of $\mathbb{Q}_p$ for $p$ sufficiently large.

We will need the following notion. Let $k$ be a fixed field of characteristic zero. We denote by $\mathcal{L}_{\mathrm{DP},k}$ the language obtained by adding constant symbols to the language $\mathcal{L}_{\mathrm{DP}}$ in the VF, resp. RF sort, for every element of $k((t))$, resp. $k$. Then for any field $K$ containing $k$, $(K((t)), K, \mathbb{Z})$ is an $\mathcal{L}_{\mathrm{DP},k}$-structure.

**Constructible motivic functions**

In this section we will recall very quickly the definition of constructible motivic functions. For the details we refer to [CL08].

We fix a field $k$ of characteristic zero. Denote by $\mathrm{Field}_k$ the category of all fields containing $k$. For any $\mathcal{L}_{\mathrm{DP},k}$-formula $\phi$, we denote by $h_\phi(K)$ the set of points in

$$h[m,n,r](K) := K((t))^m \times K^n \times \mathbb{Z}^r,$$

which satisfy $\phi$. We call the assignment $K \mapsto h_\phi(K)$ a *k-definable subassignment* and we define $\mathrm{Def}_k$ to be the category of $k$-definable subassignments. A point $x$ of $X \in \mathrm{Def}_k$ is a tuple $x = (x_0, K)$ where $x_0 \in X(K)$ and $K \in \mathrm{Field}_k$. In general, for $S \in \mathrm{Def}_k$ we define the category $\mathrm{Def}_S$ of definable subassignments $X$ with a definable map $X \to S$. We denote $\mathrm{RDef}_S$ for the category of definable subassignments of $S \times h[0,n,0]$ where $n \in \mathbb{N}$. We recall that the Grothendieck semigroup $SK_0(\mathrm{RDef}_S)$ is the quotient of the free abelian semigroup over symbols $[Y \to S]$, with $Y \to S$ in $\mathrm{RDef}_S$, by the relations
(1) $[\emptyset \to S] = 0$;
(2) $[Y \to S] = [Y' \to S]$, if $Y \to S$ is isomorphic to $Y' \to S$;
(3) $[(Y \cup Y') \to S] + [(Y \cap Y') \to S] = [Y \to S] + [Y' \to S]$, for $Y$ and $Y'$ definable subassignments of some $S[0,n,0] = S \times h[0,n,0] \to S$.

Similarly, we recall that the Grothendieck group $K_0(\mathrm{RDef}_S)$ is the quotient of the free abelian group over the symbols $[Y \to S]$, with $Y \to S$ in $\mathrm{RDef}_S$, by the relations (2) and (3). The Cartesian fiber product over $S$ induces a natural semi-ring (resp. ring) structure on $SK_0(\mathrm{RDef}_S)$ (resp. $K_0(\mathrm{RDef}_S)$) by setting

$$[Y \to S] \times [Y' \to S] = [Y \times_S Y' \to S].$$

We consider a formal symbol $\mathbb{L}$ and the ring

$$\mathbb{A} := \mathbb{Z}\left[\mathbb{L}, \mathbb{L}^{-1}, \left(\frac{1}{1-\mathbb{L}^{-i}}\right)_{i>0}\right].$$

For every real number $q > 1$, there is a unique morphism of rings $\vartheta_q : \mathbb{A} \to \mathbb{R}$ mapping $\mathbb{L}$ to $q$, and it is obvious that $\vartheta_q$ is injective for $q$ transcendental. We define a partial ordering on $\mathbb{A}$ by setting $a \geq b$ if, for every real number $q > 1$, $\vartheta_q(a) \geq \vartheta_q(b)$. Furthermore we denote by $\mathbb{A}_+$ the set $\{a \in \mathbb{A} | a \geq 0\}$.

**Definition 3.3.3.** Let $S$ be a definable subassignment in $\mathrm{Def}_k$ and denote by $|S|$ its set of points. The ring $\mathcal{P}(S)$ of *constructible Presburger functions on $S$* is defined as the subring of the ring of functions $|S| \to \mathbb{A}$, generated by
— the constant functions $|S| \to \mathbb{A}$;
— the functions $\widehat{\alpha} : |S| \to \mathbb{Z}$ that correspond to a definable morphism $\alpha : S \to h[0, 0, 1]$;
— the functions $\mathbb{L}^{\widehat{\beta}} : |S| \to \mathbb{A}$ that correspond to a definable morphism $\beta : S \to h[0, 0, 1]$.

We denote by $\mathcal{P}_+(S)$ the semiring of funtions in $\mathcal{P}(S)$ with values in $\mathbb{A}_+$.

**Definition 3.3.4.** Let $Z$ be in $\mathrm{Def}_k$. For $Y$ a definable subassignment of $Z$, we denote by $\mathbf{1}_Y$ the function in $\mathcal{P}(Z)$ with value 1 on $|Y|$ and 0 on $|Z \setminus Y|$. We denote by $\mathcal{P}_Z^0$ (resp. $\mathcal{P}_+^0(Z)$) the subring (resp. subsemiring) of $\mathcal{P}(Z)$ (resp. $\mathcal{P}_+(Z)$) generated by the functions $\mathbf{1}_Y$, for all definable subassignments $Y$ of $Z$, and by the constant function $\mathbb{L} - 1$. Notice that we have a canonical ring morphism $\mathcal{P}^0(Z) \to K_0(\mathrm{RDef}_Z)$ (resp. semiring morphism $\mathcal{P}_+^0(Z) \to SK_0(\mathrm{RDef}_Z)$) sending $\mathbf{1}_Y$ to the class of the inclusion morphism $[i : Y \to Z]$ and $\mathbb{L} - 1$ to $\mathbb{L}_Z - 1$. By $\mathbb{L}_Z$ we mean the class of the element $[Z \times h[0, 1, 0] \to Z]$ in $K_0(\mathrm{RDef}_Z)$ (resp. $SK_0(\mathrm{RDef}_Z)$).

**Definition 3.3.5.** We say that a function $\varphi \in \mathcal{P}(S \times \mathbb{Z}^r)$ is *S-integrable*, if for every $s \in S$, the family $(\varphi(s, i))_{i \in \mathbb{Z}^r}$ is summable. We denote by $I_S \mathcal{P}(S \times \mathbb{Z}^r)$ the $\mathcal{P}(S)$-module of $S$-integrable functions.

Now we define the semiring $\mathcal{C}_+(Z)$ of *positive constructible motivic functions on $Z$* as

$$\mathcal{C}_+(Z) = SK_0(\mathrm{RDef}_Z) \otimes_{\mathcal{P}_+^0(Z)} \mathcal{P}_+(Z)$$

and the ring $\mathcal{C}(Z)$ of *constructible motivic functions on $Z$* as

$$\mathcal{C}(Z) = K_0(\mathrm{RDef}_Z) \otimes_{\mathcal{P}^0(Z)} \mathcal{P}(Z).$$

Let $Z$ be a subassignment of $h[m, n, r]$. We denote by $\dim Z$ the dimension of Zariski closure of $p(Z)$ for $p$ the projection $h[m, n, r] \to h[m, 0, 0]$. For a natural number $d$, we denote by $\mathcal{C}^{\leq d}$ the ideal of $\mathcal{C}(Z)$ generated by all elements of the form $\mathbf{1}_Y$ with $Y$ a subassignment of $Z$ such that $\dim Y \leq d$. We set $\mathcal{C}^d = \mathcal{C}^{\leq d}/\mathcal{C}^{\leq d-1}$ and $C(Z) = \oplus_{d \geq 0} \mathcal{C}^d$.

For each $Y$ in $\mathrm{Def}_S$ we can define a graded subgroup $I_S(Y)$ of $C(Y)$, as in [CL10], together with a map $f_! : I_S(Y) \to I_S(Z)$, for any map $f : Y \to Z$ in $\mathrm{Def}_S$. When $S = h[0, 0, 0]$ and $f : Y \to h[0, 0, 0]$, the map $f_!$ is exactly the same as taking the integral over $Y$.

**The language $\mathcal{L}_\mathcal{O}$**

Now we suppose that $K$ is a number field with $\mathcal{O}$ its ring of integers. We denote by $\mathcal{F}_\mathcal{O}$ the set of all non-archimedean local fields over $\mathcal{O}$, which is endowed the structure of an $\mathcal{O}$-algebra. For $N \in \mathbb{N}$ we denote by $\mathcal{F}_{\mathcal{O},N}$ the set of all local fields in $\mathcal{F}_\mathcal{O}$ with residue field of characteristic at least $N$. The language $\mathcal{L}_\mathcal{O}$ is obtained from the language $\mathcal{L}_{\mathrm{DP},K}$ by restricting the constant symbols to $\mathcal{O}[[t]]$ for the valued field sort and to $\mathcal{O}$ for the residue field sort.

Let $F \in \mathcal{F}_\mathcal{O}$, we write $k_F$ for its residue field, $q_F$ for the number of elements in $k_F$, $\mathcal{O}_F$ for its valuation ring and $\mathcal{M}_F$ for its maximal ideal. For each choice of a uniformising element $\varpi_F$ of $\mathcal{O}_F$, there is a unique map $\overline{\mathrm{ac}}_{\varpi_F} : F^\times \to k_F^\times$, which extends the map $\mathcal{O}_F^\times \to k_F^\times$ and sends $\varpi_F$ to 1. Then $(F, k_F, \mathbb{Z})$ has an $\mathcal{L}_{\mathrm{DP}}$-structure with respect to $\varpi_F$. Moreover $F$ can be equipped with the structure of an $\mathcal{O}[[t]]$-algebra via the morphism:

$$\lambda_{\varpi_F} : \mathcal{O}[[t]] \to F;$$
$$\sum_{i \geq 0} a_i t^i \mapsto \sum_{i \geq 0} a_i \varpi_F^i.$$

By intepreting $a \in \mathcal{O}[[t]]$ as $\lambda_{\varpi_F}(a)$, an $\mathcal{L}_\mathcal{O}$-formula $\phi$ defines, for each $F \in \mathcal{F}_\mathcal{O}$, a definable subset $\phi(F)$ of $F^m \times k_F^n \times \mathbb{Z}^r$ for some $m, n, r \in \mathbb{N}$. If we have two $\mathcal{L}_\mathcal{O}$-formulas $\phi_1, \phi_2$ which define the same subassignment of $h[m, n, r]$, then, by compactness, $\phi_1(F) = \phi_2(F)$, for all $F \in \mathcal{F}_{\mathcal{O},N}$, for some large enough $N \in \mathbb{N}$, which does not depend on the choice of a uniformising element.

If a definable subassignment is defined in the language $\mathcal{L}_\mathcal{O}$, then we say that it belongs to $\mathrm{Def}_{\mathcal{L}_\mathcal{O}}$. In the same way we also say that a constructible function $\theta$ belongs to $\mathcal{C}(X, \mathcal{L}_\mathcal{O})$.

If $X \in \mathrm{Def}_{\mathcal{L}_\mathcal{O}}$, then $X$ is defined by a formula $\phi$ in $\mathcal{L}_\mathcal{O}$. By the above discussion we can define $X_F = \phi(F)$, for any $F \in \mathcal{F}_\mathcal{O}$. Also if $f : Y \to Z$ in $\mathrm{Def}_{\mathcal{L}_\mathcal{O}}$, then we can define a map $f_F : Y_F \to Z_F$, for any $F \in \mathcal{F}_\mathcal{O}$.

Now we will explain how to interprete a constructible function $\theta \in \mathcal{C}(X, \mathcal{L}_\mathcal{O})$ in a field $F \in \mathcal{F}_\mathcal{O}$. If $\theta \in \mathcal{P}(X)$ we will replace $\mathbb{L}$ by $q_F$ and a definable function $\alpha : X \to h[0, 0, 1]$ by a function $\alpha_F : X_F \to \mathbb{Z}$. If $\theta \in K_0(\mathrm{RDef}_{X, \mathcal{L}_\mathcal{O}})$ is of the form $[Y \xrightarrow{\pi} X]$ with $\pi : Y \to X$ defined by an $\mathcal{L}_\mathcal{O}$-formula, then we interpret $\theta$ by setting, for all $x \in X_F$,

$$\theta_F(x) := \#(\pi^{-1}(x)).$$

Notice that these interpretations can depend on the choice of formulas.

**Cell decomposition**

The structure of the sets appearing in a definable subassignment, can be better understood by decomposing the subassignment into 'cells'.

**Definition 3.3.6. Cells**. Let $S$ be in $\mathrm{Def}_K$ and $C$ a definable subassignment of $S$. Let $\alpha, \xi, c$ be definable morphisms $\alpha : C \to h[0, 0, 1], \xi : C \to h[0, 1, 0]$ and $c : C \to h[1, 0, 0]$. The 1-*cell* $Z_{C,\alpha,\xi,c}$ with *basis* $C$, *order* $\alpha$, *center* $c$, and *angular component* $\xi$, is the definable subassignment of $S[1, 0, 0]$, defined by the formula

$$y \in C \wedge \operatorname{ord}(z - c(y)) = \alpha(y) \wedge \overline{\operatorname{ac}}(z - c(y)) = \xi(y),$$

where $y$ belongs to $S$ and $z$ to $h[1,0,0]$. Similarly the 0-*cell* $Z_{C,c}$ with *basis* $C$ and *center* $c$, is the definable subassignment of $S[1,0,0]$, defined by the formula

$$y \in C \wedge z = c(y).$$

A definable subassignment $Z$ of $S[1,0,0]$ will be called a 1-*cell*, resp. a 0-*cell*, if there exists a definable isomorphism

$$\lambda : Z \to Z_C = Z_{C,\alpha,\xi,c} \subset S[1,s,r],$$

resp. a definable isomorphism

$$\lambda : Z \to Z_C = Z_{C,c} \subset S[1,s,0],$$

for some $r, s \geq 0$, some basis $C \subset S[0,s,r]$, resp. $S[0,s,0]$, and some 1-cell $Z_{C,\alpha,\xi,c}$, resp. 0-cell $Z_{C,c}$, such that the morphism $\pi \circ \lambda$, with $\pi$ the projection on the $S[1,0,0]$-factor, is the identity on $Z$. The data $(\lambda, Z_{C,\alpha,\xi,c})$, resp. $(\lambda, Z_{C,c})$, will be called a *presentation of the cell* $Z$ and denoted for short by $(\lambda, Z_C)$.

**Theorem 3.3.7** ([CL08], Thm 7.2.1). *Suppose that $K$ is a field of characteristic $0$. Let $X$ be a definable subassignment of $S[1,0,0]$ with $S$ in $\operatorname{Def}_K$.*
*(1) The subassignment $X$ is a finite disjoint union of cells.*
*(2) For every $\varphi \in \mathcal{C}(X)$, there exists a finite partition of $X$ into cells $Z_i$ with presentation $(\lambda_i, Z_{C_i})$, such that $\varphi|_{Z_i} = \lambda_i^* p_i^*(\psi_i)$, with $\psi_i \in \mathcal{C}(C_i)$ and $p_i : Z_{C_i} \to C_i$ the projection. Similar statements hold for $\varphi$ in $\mathcal{C}_+(C)$, in $\mathcal{P}(X)$, in $\mathcal{P}_+(X)$, in $K_0(\operatorname{RDef}_Z)$ and in $SK_0(\operatorname{RDef}_Z)$.*

**Corollary 3.3.8.** *Theorem 3.3.7 still holds, if we replace $\operatorname{Def}_K$ by $\operatorname{Def}_{\mathcal{L}_{\mathcal{O}}}$.*

*Proof.* The proof is the same as the proof of Theorem 3.3.7, but we replace $\mathcal{L}_{\operatorname{DP,K}}$ by $\mathcal{L}_{\mathcal{O}} \subset \mathcal{L}_{\operatorname{DP,K}}$. $\square$

### 3.3.2 Proof of the main theorem

We will give a proof of the Main Theorem 3.1.2 by splitting the exponential sum $E_{m,p}^0(f)$ into three subsums.

$$E_{p,m}^0(f) = \frac{1}{p^{nm}} \sum_{\substack{\overline{x} \in (p\mathbb{Z}/p^m\mathbb{Z})^n, \\ \operatorname{ord}_p(f(x)) \leq m-2}} \exp\left(\frac{2\pi i}{p^m} f(x)\right) +$$

$$\frac{1}{p^{nm}} \sum_{\substack{\overline{x} \in (p\mathbb{Z}/p^m\mathbb{Z})^n, \\ \operatorname{ord}_p(f(x)) = m-1}} \exp\left(\frac{2\pi i}{p^m} f(x)\right) + \frac{1}{p^{nm}} \sum_{\substack{\overline{x} \in (p\mathbb{Z}/p^m\mathbb{Z})^n, \\ \operatorname{ord}_p(f(x)) \geq m}} \exp\left(\frac{2\pi i}{p^m} f(x)\right).$$

In three different lemmas we will analyse each of these sums.

For the first subsum we will introduce a constructible function $G$, that expresses, for a certain input $z \in \mathbb{Z}_p$ with $\operatorname{ord}_p(z) \leq m-2$, how many $x \in (p\mathbb{Z}_p)^n$ are mapped close to $z$ by $f$. We will apply the Cell Decomposition Theorem to $G$ and with some further techniques like eliminiation of quantifiers, we will show that certain values $z$ of $f$ occur equally often. In the exponential sum these values will cancel out.

**Lemma 3.3.9.** *Let $f \in \mathbb{Z}[x_1, \ldots, x_n]$ be a non-constant polynomial such that $f(0) = 0$. There exists $N \in \mathbb{N}$ such that, for all $m \geq 1$ and all prime numbers $p > N$, we have*

$$\sum_{\substack{\overline{x} \in (p\mathbb{Z}/p^m\mathbb{Z})^n, \\ \mathrm{ord}_p(f(x)) \leq m-2}} \exp\left(\frac{2\pi i f(x)}{p^m}\right) = 0.$$

*Proof.* The statement is obvious when $m = 1$ or $m = 2$, so we can assume that $m > 2$. Let $\phi$ be the $\mathcal{L}_\mathbb{Z}$-formula given by

$$\phi(x_1, \ldots, x_n, z, m) = \bigwedge_{i=1}^{n}(\mathrm{ord}(x_i) \geq 1) \wedge (\mathrm{ord}(z) \leq m-2) \wedge (\mathrm{ord}(z - f(x_1, \ldots, x_n)) \geq m),$$

where $x_i, z$ are in the valued field-sort and $m$ is in the value group-sort. To shorten notation we set $x = (x_1, \ldots, x_n)$. For each prime $p$, we fix a uniformiser $\varpi_p$ of $\mathbb{Q}_p$, then $\phi$ defines, for each $p$, a definable set $X_p \subset (p\mathbb{Z}_p)^n \times \mathbb{Z}_p \times \mathbb{Z}$. More precisely, we have

$$X_p = \{(x, z, m) \in (p\mathbb{Z}_p)^n \times \mathbb{Z}_p \times \mathbb{Z} \mid \mathrm{ord}_p(f(x) - z) \geq m, \mathrm{ord}_p(z) \leq m - 2\}.$$

It is obvious that $X_p$ does not depend on $\varpi_p$.

We denote by $X \subset h[n+1, 0, 1]$ the definable subassignment defined by $\phi$. Let $F := \mathbf{1}_X \in I_{h[0,0,1]}(h[n+1, 0, 1])$ and $\pi$ the projection from $h[n+1, 0, 1]$ to $h[1, 0, 1]$. Then we have $G := \pi_!(F) \in I_{h[0,0,1]}(h[1, 0, 1])$. For each prime $p$ and each uniformiser $\varpi_p$ of $\mathbb{Q}_p$, there exist the following interpretations of $F$ and $G$ in $\mathbb{Q}_p$:

$$F_{\varpi_p} = \mathbf{1}_{X_p}$$

and

$$G_{\varpi_p}(z, m) = \int_{X_{p,z,m}} |dx| = p^{-mn} \#\{\overline{x} \in (p\mathbb{Z}/p^m\mathbb{Z})^n \mid f(x) \equiv z \bmod p^m\},$$

if $\mathrm{ord}_p(z) \leq m-2$, where $X_{p,z,m}$ is the fiber of $X_p$ over $(z, m)$, and

$$G_{\varpi_p}(z, m) = 0,$$

if $\mathrm{ord}_p(z) \geq m-1$. We can see that both $F_{\varpi_p}(x, z, m)$ and $G_{\varpi_p}(z, m)$ do not depend on $\varpi_p$.

Now we use Corollary 3.3.8 for $G \in I_{h[0,0,1]}(h[1, 0, 1])$. This means that there exists a finite partition of $h[1, 0, 1]$ into cells $Z_i$ (for $i$ in some finite set $I$) with presentation $(\lambda_i, Z_{C_i, \alpha_i, \xi_i, c_i})$, such that $G|_{Z_i} = \lambda_i^* p_i^*(G_i)$ with $G_i \in \mathcal{C}(C_i)$ and $p_i : Z_{C_i, \alpha_i, \xi_i, c_i} \to C_i$ the projection. Note that $C_i \subset h[0, r_i, s_i + 1]$ for some $r_i, s_i \in \mathbb{N}$. We denote by $\theta_i(z, \eta, \gamma, m)$ the $\mathcal{L}_\mathbb{Z}$-formula defining $c_i$, where $z \in h[1, 0, 0], \eta \in h[0, r_i, 0], \gamma \in h[0, 0, s_i]$ and $m \in h[0, 0, 1]$. By elimination of quantifiers (Corollary 3.3.2), there exist polynomials $f_1, \ldots, f_r$ in one variable $z$ with coefficients in $\mathbb{Z}[[t]]$, such that $\theta_i(z, \eta, \gamma, m)$ is equivalent to the formula

$$\bigvee_j \Big(\varsigma_{ij}\big(\overline{\mathrm{ac}} f_1(z), \ldots, \overline{\mathrm{ac}} f_r(z), \eta\big) \wedge \nu_{ij}\big(\mathrm{ord}\, f_1(z), \ldots, \mathrm{ord}\, f_r(z)\big), \gamma, m\Big),$$

where $\zeta_{ij}$ is an $\mathcal{L}_{\text{ring}}$-formula and $\nu_{ij}$ an $\mathcal{L}_{\text{oag}}$-formula. Since $c_i$ is a function, we know that, for each $(\eta, \gamma, m) \in C_i$, there exists a unique $z = c_i(\eta, \gamma, m)$ such that $\theta_i(z, \eta, \gamma, m)$ is true. We claim now that there exists $1 \leq l \leq r$ such that $f_l(z) = 0$. Indeed, if $f_l(z) \neq 0$, for all $l$, then there exists a small open neighborhood $V$ of $z$ and there exists an index $j$, such that, for all $y \in V$, $(y, \eta, \gamma, m)$ will satisfy the formulas $\zeta_{ij}, \eta_{ij}$. Since this would contradict the uniqueness of $z$, we must have that $f_l(z) = 0$ for some $l$. We deduce that $A := \{c_i(\eta, \gamma, m) \in h[1, 0, 0] \mid i \in I, (\eta, \gamma, m) \in C_i\} \subset \cup_{j=1}^r Z(f_j)$.

From the definition of $G_i$ we see that, if we fix $(\eta, \gamma, m) \in C_i$, then $G(\cdot, m)$ will be constant on the ball

$$\{y \in h[1, 0, 0] \mid \overline{\text{ac}}(y - c_i(\eta, \gamma, m)) = \xi_i(\eta, \gamma, m), \text{ord}(y - c_i(\eta, \gamma, m)) = \alpha_i(\eta, \gamma, m)\}.$$

Now, for each $m > 2$, we set

$$B_m := A \cap \{z \in h[1, 0, 0] \mid m - 2 \geq \text{ord}(z) \geq 1\}$$

and

$$U_m := \{y \in h[1, 0, 0] \mid \text{ord}(z - y) < m - 1, \forall z \in B_m\}.$$

So $U_m$ will be a union of balls of radius $m - 1$. Because $f(0) = 0$, we can see that $G(\cdot, m)$ will be zero on the set $\{z \in h[1, 0, 0] \mid \text{ord}(z) \leq 0\}$, if $m > 2$.

**Claim 3.3.10.** *If $m > 2$, $\text{ord}(z) \geq 1$ and $z \in U_m$, then $G(\cdot, m)$ will be constant on the ball $B(z, m - 1)$ (the ball with center $z$ and radius $m - 1$).*

From the cell decomposition of $h[1, 0, 1]$, we know that there exist $i \in I$ and $(\eta, \gamma) \in h[0, r_i, s_i]$, such that $(z, \eta, \gamma, m) \in Z_{C_i, \alpha_i, \xi_i, c_i}$. Hence $(\eta, \gamma, m) \in C_i$ and $z$ belongs to the ball

$$B = \{y \in h[1, 0, 0] \mid \overline{\text{ac}}(y - c_i(\eta, \gamma, m)) = \xi_i(\eta, \gamma, m), \text{ord}(y - c_i(\eta, \gamma, m)) = \alpha_i(\eta, \gamma, m)\}.$$

We will distinguish three cases, depending on the value of $c_i(\eta, \gamma, m)$. First of all, if $c_i(\eta, \gamma, m) \in B_m$, then we see that $\alpha_i(\eta, \gamma, m) = \text{ord}(z - c_i(\eta, \gamma, m)) < m - 1$. Therefore the ball $B$ will contain the ball $B(z, m-1)$, thus $G(\cdot, m)$ will be constant on $B(z, m - 1)$. Second of all, if $\text{ord}(c_i(\eta, \gamma, m)) \leq 0$, and since $\text{ord}(z) \geq 1$, we have $\alpha_i(\eta, \gamma, m) \leq 0 < m - 1$ so we have the same situation as above. Thirdly, if $\text{ord}(c_i(\eta, \gamma, m)) \geq m - 1$, then the case $\alpha_i(\eta, \gamma, m) < m - 1$ has already been treated above. Hence we can assume that $\alpha_i(\eta, \gamma, m) \geq m - 1$, in which case we have $B(z, m - 1) = B(0, m - 1)$. By definition of $G$ we have $G(\cdot, m)|_{B(0, m-1)} = 0$. This proves the claim.

Now there exists $N_0 \in \mathbb{N}$, independent of $m > 2$, for which we can interpret all of the above discussion in $\mathbb{Q}_p$, with any choice of uniformiser $\varpi_p \in \mathbb{Z}_p$ and for any $p > N_0$, by applying the map $\lambda_{\varpi_p}$ to the coefficients of the polynomials $f_1, \ldots, f_r$. Because $U_{m, \varpi_p}$ is an $\{m, \varpi_p\}$-definable set in the language $\mathcal{L}_{\text{DP}}$, it can vary when changing $\varpi_p$. This suggests us to set $\mathcal{U}_{m,p} := \cup_{\varpi_p} U_{m, \varpi_p}$ with $\varpi_p$ running over the set of all uniformisers of $\mathbb{Q}_p$. Then $\mathcal{U}_{m,p}$ is given by an $\mathcal{L}_{\text{DP}}$-formula.

**Claim 3.3.11.** *There exists $N \in \mathbb{N}$, such that $\mathcal{U}_{m,p} = \mathbb{Q}_p$, for all $m > 2$ and for all $p > N$.*

From the definition of $U_{m,\varpi_p}$ we see that $\mathcal{V}_{m,p} := \mathbb{Q}_p \setminus \mathcal{U}_{m,p}$ is a union of $d_{m,p}$ balls of radius $m - 1$, contained in $p\mathbb{Z}_p$, where $d_{m,p} \leq \sum_{j=1}^{r} \deg f_j$. Moreover, $\mathcal{V}_{m,p}$ will given by a $\mathcal{L}_{\text{DP}}$-formula. We use elimination of quantifiers (Theorem 3.3.1) for the formula defining $\mathcal{V}_{m,p}$. Hence there exist polynomials $q_1, \ldots, q_{\tilde{r}}$ of one variable $z$ with coefficients in $\mathbb{Z}$ and formulas $\varphi_j$ in $\mathcal{L}_{\text{ring}}$ and $\nu_j$ in $\mathcal{L}_{\text{oag}}$, for $1 \leq j \leq s$, such that

$$z \in \mathcal{V}_{m,p} \Leftrightarrow \bigvee_{j=1}^{s} \varphi_j\big(\overline{\text{ac}}_{\varpi_p}(q_1(z)), \ldots, \overline{\text{ac}}_{\varpi_p}(q_{\tilde{r}}(z))\big) \wedge \nu_j\big(\text{ord}_p(q_1(z)), \ldots, \text{ord}_p(q_{\tilde{r}}(z)), m\big),$$

for any $p > N_0$ (after enlarging $N_0$ if necessary) and any uniformiser $\varpi_p$.

We note that if $z \in \mathcal{V}_{m,p}$, then $\text{ord}_p(z) \geq 1$. Since $q_i$ has coefficients in $\mathbb{Z}$, we can assume, by possibly enlarging $N_0$, that $\overline{\text{ac}}_{\varpi_p}(q_i(z))$ only depends on $\overline{\text{ac}}_{\varpi_p}(z)$ and $\text{ord}_p(q_i(z))$ only depends on $\text{ord}_p(z)$, for any $p > N_0$. This follows from the $t$-adic version of this statement by a compactness argument. So if $z_1$ and $z_2$ satisfy that
— $\text{ord}_p(z_1) = \text{ord}_p(z_2) \geq 1$,
— there exist two uniformisers $\varpi_{1,p}$ and $\varpi_{2,p}$, such that $\overline{\text{ac}}_{\varpi_{1,p}}(z_1) = \overline{\text{ac}}_{\varpi_{2,p}}(z_2)$,
then we see that $z_1 \in \mathcal{V}_{m,p}$ if and only if $z_2 \in \mathcal{V}_{m,p}$. It implies that $\overline{\mathcal{V}}_{m,p} := \overline{\text{ac}}_{\varpi_p}(\mathcal{V}_{m,p})$ does not depend on $\varpi_p$, for any $p > N_0$. In particular, since $B(0, m-1) \not\subseteq \mathcal{V}_{m,p}$, we see that the number of elements in $\overline{\mathcal{V}}_{m,p}$ is at most $\sum_{j=1}^{r} \deg f_j$, for all $p > N_0$.

In what follows we will show that if $\mathcal{V}_{m,p}$ were not empty, then the set $\overline{\mathcal{V}}_{m,p}$ would grow with $p$. This will give the desired contradiction. We set

$$B_\infty = A \cap \{z \in h[1,0,0] \mid \infty > \text{ord}(z) \geq 1\} \subset \cup_{j=1}^{r} Z(f_j),$$

thus $B_\infty$ is a finite set with $0 \notin B_\infty$ and $B_m \subset B_\infty$ for all $m > 2$. Looking at the order of the coefficients of $f_j$ we see that there exists $M \in \mathbb{N}$ such that $\text{ord}_p(z) \leq M$ for all $z \in Z(f_{j,\varpi_p}) \setminus \{0\}$, for all $1 \leq j \leq r$, for all $p > N_0$ and for all uniformiser $\varpi_p$. So $\text{ord}_p(z) \leq M$ for all $z \in B_{\infty,\varpi_p}$, for all $\varpi_p$. It follows that $\text{ord}_p(z) \leq M$ for all $z \in \mathcal{V}_{m,p}$, for all $m > 2$ and $p > N_0$. Indeed, since $B(0, m-1) \not\subseteq \mathcal{V}_{m,p}$ we have $\text{ord}_p(z) < m-1$ for all $z \in \mathcal{V}_{m,p}$, so it is true if $m - 1 \leq M$. On the other hand, if $m - 1 > M$, then for each $z \in \mathcal{V}_{m,p}$ and each uniformiser $\varpi_p$, there exists $z_0 \in B_{\infty,\varpi_p}$ such that $\text{ord}_p(z - z_0) \geq m - 1 > M \geq \text{ord}_p(z_0)$, thus $\text{ord}_p(z) = \text{ord}_p(z_0) \leq M$. Now put $N := \max\{N_0, 1 + M \sum_{j=1}^{r} \deg f_j\}$. Suppose for a contradiction, that for some $p > N$, there exists $z \in \mathcal{V}_{p,m}$. Then $\overline{\text{ac}}_{\varpi_p}(z) \in \overline{\mathcal{V}}_{m,p}$, for every uniformiser $\varpi_p$, and so $\{\overline{\text{ac}}_{\varpi_p}(z) \mid \text{ord}_p(\varpi_p) = 1\} \subset \overline{\mathcal{V}}_{m,p}$. Suppose that $\overline{\text{ac}}_p(\varpi_p) = u$, then $u^{\text{ord}_p(z)}\overline{\text{ac}}_{\varpi_p}(z) = \overline{\text{ac}}_p(z)$, so we have $\{\overline{\text{ac}}_{\varpi_p}(z) \mid \text{ord}_p(\varpi_p) = 1\} = \{u^{-\text{ord}_p(z)}\overline{\text{ac}}_p(z) \mid u \in \mathbb{F}_p^\times\}$. Therefore $\#\{u^{-\text{ord}_p(z)}\overline{\text{ac}}_p(z) \mid u \in \mathbb{F}_p^\times\} \leq \sum_{j=1}^{r} \deg f_j$. But

$$\#\{u^{-\text{ord}_p(z)}\overline{\text{ac}}_p(z) \mid u \in \mathbb{F}_p^\times\} = \frac{p-1}{\gcd(\text{ord}_p(z), p-1)} \geq \frac{p-1}{\text{ord}_p(z)} \geq \frac{p-1}{M},$$

where $\gcd(a, b)$ is the greatest common divisor of $a$ and $b$. Then we have $p - 1 \leq M \sum_{j=1}^{r} \deg f_j \leq N - 1$. This is a contradiction, since $p > N$. So this proves the claim.

We know from Claim 3.3.10 that if $m > 2$ and $z \in U_{m,\varpi_p}$ such that $1 \leq \text{ord}_p(z) \leq m - 2$, then $G_{\varpi_p}(.,m)$ will be constant on the ball $B(z, m-1)$. Thus we have

$$\#\{\overline{x} \in (p\mathbb{Z}/p^m\mathbb{Z})^n \mid f(x) \equiv y \bmod p^m\} = \#\{\overline{x} \in (p\mathbb{Z}/p^m\mathbb{Z})^n \mid f(x) \equiv z \bmod p^m\}$$

for all $\overline{y} \in p\mathbb{Z}/p^m\mathbb{Z}$ with $y \equiv z \bmod p^{m-1}$. Hence

$$\sum_{\substack{\overline{x} \in (p\mathbb{Z}/p^m\mathbb{Z})^n, \\ f(x) \equiv z \bmod p^{m-1}}} p^{-mn} \exp\left(\frac{2\pi i f(x)}{p^m}\right) = G(z,m) \cdot \sum_{\substack{\overline{y} \in p\mathbb{Z}/p^m\mathbb{Z}, \\ y \equiv z \bmod p^{m-1}}} \exp\left(\frac{2\pi i y}{p^m}\right) = 0.$$

This implies that

$$\sum_{\substack{\overline{x} \in (p\mathbb{Z}/p^m\mathbb{Z})^n, \\ \overline{f(x)} \in \overline{\mathcal{U}}_{m,p}}} p^{-mn} \exp\left(\frac{2\pi i f(x)}{p^m}\right) = 0,$$

where $\overline{\mathcal{U}}_{m,p} := \{\overline{z} \in p\mathbb{Z}/p^{m-1}\mathbb{Z} \mid z \in \mathcal{U}_{m,p}, m-2 \geq \text{ord}_p(z) \geq 1\}$. For all $m > 2$ and $p > N$, we have $\mathcal{U}_{m,p} = \mathbb{Q}_p$, so $\overline{\mathcal{U}}_{m,p} = \{\overline{z} \in p\mathbb{Z}/p^{m-1}\mathbb{Z} \mid \text{ord}_p(z) \leq m-2\}$. Therefore we have

$$\sum_{\substack{\overline{x} \in (p\mathbb{Z}/p^m\mathbb{Z})^n, \\ \text{ord}_p(f(x)) \leq m-2}} p^{-mn} \exp\left(\frac{2\pi i f(x)}{p^m}\right) = 0. \qquad \square$$

In the proof of the following lemma we will introduce again a constructible function $G$, similar to the one from the previous proof. For this exponential sum the different values $z$ of $f$ do not cancel out completely. By using the Lang-Weil estimation (see [LW54]) and Theorem 3.2.2 we obtain the following upper bound for the second subsum.

**Lemma 3.3.12.** *Let $f \in \mathbb{Z}[x_1, \ldots, x_n]$ be a non-constant polynomial, such that $f(0) = 0$. Put $\sigma = \min\{c_0(f), \frac{1}{2}\}$, where $c_0(f)$ is the log-canonical threshold of $f$ at $0$. Then there exist, for each integer $m > 1$, a natural number $N_m$ and a positive constant $D_m$, such that, for all $p > N_m$, we have*

$$\left| \sum_{\substack{\overline{x} \in (p\mathbb{Z}/p^m\mathbb{Z})^n, \\ \text{ord}_p(f(x))=m-1}} p^{-mn} \exp\left(\frac{2\pi i f(x)}{p^m}\right) \right| \leq D_m p^{-m\sigma}.$$

*Proof.* Let $\phi, \overline{\phi}$ be two $\mathcal{L}_\mathbb{Z}$-formulas given by

$$\phi(x_1, \ldots, x_n, z, m) =$$
$$\bigwedge_{i=1}^{n} (\text{ord}(x_i) \geq 1) \wedge (\text{ord}(z) = m-1) \wedge (\text{ord}(z - f(x_1, \ldots, x_n)) \geq m),$$
$$\overline{\phi}(x_1, \ldots, x_n, \xi, m) =$$
$$\bigwedge_{i=1}^{n} (\text{ord}(x_i) \geq 1) \wedge (\text{ord}(f(x_1, \ldots, x_n) = m-1) \wedge (\overline{\text{ac}}(f(x_1, \ldots, x_n)) = \xi),$$

where $x_i, z$ are in the valued field-sort, $m$ is in the valued group-sort and $\xi$ is in the residue field-sort. To shorten notation we set $x = (x_1, \ldots, x_n)$. For each prime $p$, we fix a uniformiser $\varpi_p$ of $\mathbb{Q}_p$, then $\phi, \overline{\phi}$ define, for each $p$, two definable sets $X_p \subset (p\mathbb{Z}_p)^n \times \mathbb{Z}_p \times \mathbb{Z}$ and $\overline{X}_p \subset (p\mathbb{Z}_p)^n \times \mathbb{F}_p \times \mathbb{Z}$ given by

$$X_p = \{(x, z, m) \in (p\mathbb{Z}_p)^n \times \mathbb{Z}_p \times \mathbb{Z} \mid \mathrm{ord}_p(f(x) - z) \geq m, \mathrm{ord}_p(z) = m - 1\}$$

and

$$\overline{X}_p = \{(x, \xi, m) \in (p\mathbb{Z}_p)^n \times \mathbb{F}_p \times \mathbb{Z} \mid \mathrm{ord}_p(f(x)) = m - 1, \overline{\mathrm{ac}}_{\varpi_p}(f(x)) = \xi\}.$$

It is obvious that $X_p$ does not depend on $\varpi_p$.

We denote by $X \subset h[n+1, 0, 1]$, resp. $\overline{X} \subset h[n, 1, 1]$, the definable subassignments defined by $\phi$, resp. $\overline{\phi}$. Let $F := \mathbf{1}_X \in I_{h[0,0,1]}(h[n+1, 0, 1])$ and $\pi$ the projection from $h[n+1, 0, 1]$ to $h[1, 0, 1]$. Then we have $G := \pi_!(F) \in I_{h[0,0,1]}(h[1, 0, 1])$. For each prime $p$ and each uniformiser $\varpi_p$ of $\mathbb{Q}_p$, there exist the following interpretations of $F$ and $G$ in $\mathbb{Q}_p$.

$$F_{\varpi_p} = \mathbf{1}_{X_p}$$

and

$$G_{\varpi_p}(z, m) = \int_{X_{p,z,m}} |dx| = p^{-mn} \#\{\overline{x} \in (p\mathbb{Z}/p^m\mathbb{Z})^n \mid f(x) \equiv z \bmod p^m\},$$

if $\mathrm{ord}_p(z) = m - 1$, where $X_{p,z,m}$ is the fiber of $X_p$ over $(z, m)$, and

$$G_{\varpi_p}(z, m) = 0,$$

if $\mathrm{ord}_p(z) \neq m - 1$. We can see that both $F_{\varpi_p}(x, z, m)$ and $G_{\varpi_p}(z, m)$ do not depend on $\varpi_p$. So we can set $G(z, m, p) := G_{\varpi_p}(z, m)$. The idea is to partition $p^{m-1}\mathbb{Z}_p \backslash p^m \mathbb{Z}_p$ into sets on which $G(\cdot, m, p)$ is constant. First of all, we can see that $G(\cdot, m, p)$ is constant on balls of the form

$$\{z \in \mathbb{Z}_p \mid \mathrm{ord}_p(z) = m - 1, \overline{\mathrm{ac}}_{\varpi_p}(z) = \xi_0\},$$

with $\xi_0 \in \mathbb{F}_p^\times$. Now we will look more closely on which of these balls $G(\cdot, m, p)$ takes the same value. In what follows we will show is that for $p$ big enough, if $\varpi_p, \varpi_p'$ are two uniformiser, then $G(\cdot, m, p)$ will be the same on the sets $\{z \in \mathbb{Z}_p \mid \mathrm{ord}_p(z) = m - 1, \overline{\mathrm{ac}}_{\varpi_p}(z) = \xi_0\}$ and $\{z \in \mathbb{Z}_p \mid \mathrm{ord}_p(z) = m - 1, \overline{\mathrm{ac}}_{\varpi_p'}(z) = \xi_0\}$. When this holds, we can see that $G$ will be constant on the orbits of an action of the group $\mu_{p-1}(\mathbb{Q}_p)$ on $\mathbb{Q}_p$.

We take $\overline{F} := \mathbf{1}_{\overline{X}} \in I_{h[0,0,1]}(h[n, 1, 1])$ and $\overline{\pi}$ the projection from $h[n, 1, 1]$ to $h[0, 1, 1]$. Then we have $\overline{G} := \overline{\pi}_!(\overline{F}) \in I_{h[0,0,1]}(h[0, 1, 1])$. For each prime $p$ and each uniformiser $\varpi_p$ of $\mathbb{Q}_p$, there exist the following interpretations of $\overline{F}$ and $\overline{G}$ in $\mathbb{Q}_p$.

$$\overline{F}_{\varpi_p} = \mathbf{1}_{\overline{X}_p}$$

and

$$\overline{G}_{\varpi_p}(\xi, m) = \int_{\overline{X}_{p,\xi,m}} |dx|$$
$$= p^{-mn} \#\{\overline{x} \in (p\mathbb{Z}/p^m\mathbb{Z})^n \mid \mathrm{ord}_p(f(x)) = m - 1, \overline{\mathrm{ac}}_{\varpi_p}(f(x)) = \xi\},$$

where $\overline{X}_{p,\xi,m}$ the fiber of $\overline{X}_p$ over $(\xi, m)$.

Since $\overline{G} \in I_{h[0,0,1]}(h[0,1,1])$, we can write $\overline{G}$ in the form

$$\overline{G}(\xi, m) = \sum_{i \in I} n_i \alpha_i(\xi, m) \mathbb{L}^{\beta_i(\xi,m)}[V_i],$$

where $n_i \in \mathbb{Z}$, $\alpha_i, \beta_i$ are $\mathcal{L}_\mathbb{Z}$-definable functions from $h[0,1,1]$ to $h[0,0,1]$ and $[V_i] \in K_0(\mathrm{RDef}_{h[0,1,1],\mathcal{L}_\mathbb{Z}})$. We use elimition of quantifiers (Corollary 3.3.2) for the formulas defining $\alpha_i, \beta_i, V_i$, hence there exist $N \in \mathbb{N}$, and $(\mathcal{L}_{\mathrm{ring}} \cup \mathbb{Z})$-formulas $\phi_{ij}, \theta_{ij}, \varsigma_{ij}$ and $(\mathcal{L}_{\mathrm{oag}} \cup \mathbb{Z})$-formulas $\eta_{ij}, \nu_{ij}, \tau_{ij}$, where $j \in J$, such that for all $p > N$ and all uniformiser $\varpi_p$, we have

$$\alpha_{i,\varpi_p}(\xi, m) = \eta \Leftrightarrow \vee_{j \in J}(\phi_{ij}(\xi) \wedge \eta_{ij}(\eta, m));$$
$$\beta_{i,\varpi_p}(\xi, m) = \nu \Leftrightarrow \vee_{j \in J}(\theta_{ij}(\xi) \wedge \nu_{ij}(\nu, m));$$
$$(\xi, m, \varsigma) \in V_{i,\varpi_p} \Leftrightarrow \vee_{j \in J}(\varsigma_{ij}(\varsigma, \xi) \wedge \tau_{ij}(m)).$$

From these formulas we can see that $\overline{G}_{\varpi_p}(\xi, m)$ does not depend on the uniformiser $\varpi_p$, so we will write $\overline{G}(\xi, m, p)$ instead of $\overline{G}_{\varpi_p}(\xi, m)$. But by definition of $G$ and $\overline{G}$ we can see that $G(z, m, p) = G_{\varpi_p}(z, m) = \overline{G}_{\varpi_p}(\xi, m) = \overline{G}(\xi, m, p)$, if $\overline{\mathrm{ac}}_{\varpi_p}(z) = \xi$ and $\mathrm{ord}_p(z) = m - 1$. Therefore, for $m > 1$, $p > N$ and $\mathrm{ord}_p(z_1) = \mathrm{ord}_p(z_2) = m - 1$, we have $G(z_1, m, p) = G(z_2, m, p)$, if there exist two uniformisers $\varpi_{1,p}, \varpi_{2,p}$ such that $\overline{\mathrm{ac}}_{\varpi_{1,p}}(z_1) = \overline{\mathrm{ac}}_{\varpi_{2,p}}(z_2) \in \mathbb{F}_p^\times$. Let $d = \gcd(m-1, p-1)$, then by the same reasoning as in Lemma 3.3.9 we have that $G(\cdot, m, p)$ will be constant on the sets

$$\left\{z \mid \mathrm{ord}_p(z) = \mathrm{ord}_p(z_0) = m - 1 \wedge \overline{\mathrm{ac}}_p\left(\frac{z}{z_0}\right)^{\frac{p-1}{d}} = 1\right\},$$

for any $z_0 \in \mathbb{Z}_p$ with $\mathrm{ord}_p(z_0) = m - 1$. So we can decompose $p^{m-1}\mathbb{Z}_p \backslash p^m \mathbb{Z}_p$ into $d$ of these sets, each of them will consist of $\frac{p-1}{d}$ disjoint balls of volume $p^{-m}$ and $G(\cdot, m, p)$ will be constant on these sets. We denote these sets by $Y_1, \ldots, Y_d$ and the values of $G(\cdot, m, p)$ on these sets by $G_1, \ldots, G_d$ respectively. We remark that if $\mathrm{ord}_p(z) = m - 1$, then

$$\exp\left(\frac{2\pi i z}{p^m}\right) = \exp\left(\frac{2\pi i\, \overline{\mathrm{ac}}_p(z)}{p}\right),$$

so

$$\left|\sum_{\overline{y} \in Y_i/p^m \mathbb{Z}_p} \exp\left(\frac{2\pi i y}{p^m}\right)\right| = \left|\sum_{\xi \in \overline{\mathrm{ac}}_p(Y_i)} \exp\left(\frac{2\pi i \xi}{p}\right)\right|$$
$$= \left|\sum_{u \in \mathbb{F}_p^\times} \exp\left(\frac{2\pi i u^d \xi_0}{p}\right)\right|$$

for any $\xi_0 \in \overline{\mathrm{ac}}(Y_i)$. By the last result from [Wei48] we have

$$\left|\sum_{u \in \mathbb{F}_p^\times} \exp\left(\frac{2\pi i u^d \xi_0}{p}\right)\right| = \left|\sum_{u \in \mathbb{F}_p} \exp\left(\frac{2\pi i u^d \xi_0}{p}\right) - 1\right| \leq (d-1)p^{\frac{1}{2}} + 1 \leq dp^{\frac{1}{2}},$$

hence

$$\left| \sum_{\substack{\overline{x} \in (p\mathbb{Z}/p^m\mathbb{Z})^n, \\ \operatorname{ord}_p(f(\overline{x}))=m-1}} p^{-mn} \exp\left(\frac{2\pi i f(x)}{p^m}\right) \right| =$$

$$\left| \sum_{0 \neq \overline{z} \in p^{m-1}\mathbb{Z}_p/p^m\mathbb{Z}_p} G(z,m,p) \exp\left(\frac{2\pi i z}{p^m}\right) \right| =$$

$$\left| \sum_{i=1}^d G_i \sum_{\overline{y} \in Y_i/p^m\mathbb{Z}_p} \exp\left(\frac{2\pi i y}{p^m}\right) \right| \leq \left| \sum_{i=1}^d G_i d p^{\frac{1}{2}} \right|.$$

We also have

$$\sum_{i=1}^d \frac{p-1}{d} G_i = \sum_{0 \neq \overline{z} \in p^{m-1}\mathbb{Z}_p/p^m\mathbb{Z}_p} G(z,m,p)$$
$$= p^{-mn} \#\{\overline{x} \in (p\mathbb{Z}/p^m\mathbb{Z})^n \mid \operatorname{ord}_p(f(x)) = m-1\}$$
$$= p^{-mn} \# A_{p,m}$$

where $A_{p,m} := \{\overline{x} \in (p\mathbb{Z}_p/p^m\mathbb{Z}_p)^n \mid \operatorname{ord}_p(f(x)) = m-1\}$. When we view $A_{p,m}$ as a subvariety of $\mathbb{F}_p^{mn}$, then, by the Lang-Weil estimation (see [LW54]), there exists a constant $D'_m$, not depending on $p$, such that

$$\# A_{p,m} = D'_m p^{\dim_{\mathbb{F}_p}(A_{p,m})} + O(p^{\dim_{\mathbb{F}_p}(A_{p,m}) - \frac{1}{2}}).$$

By Theorem 3.2.2 we have

$$c_0(f) \leq \frac{(m-1)n - \dim_{\mathbb{F}_p}(\tilde{A}_{p,m})}{m-1},$$

where $\tilde{A}_{p,m}$ is the image of $A_{p,m}$ under the projection $\pi_m : (\mathbb{Z}_p/p^m\mathbb{Z}_p)^n \to (\mathbb{Z}_p/p^{m-1}\mathbb{Z}_p)^n$, viewed as a subvariety of $\mathbb{F}_p^{mn-n}$. Then we have

$$\dim_{\mathbb{F}_p}(A_{m,p}) \leq n + \dim_{\mathbb{F}_p}(\tilde{A}_{m,p}) \leq mn - (m-1)c_0(f).$$

And now we finish the proof by showing that for all $p$ big enough,

$$\left| \sum_{\substack{\overline{x} \in (p\mathbb{Z}/p^m\mathbb{Z})^n, \\ \operatorname{ord}_p(f(x))=m-1}} p^{-mn} \exp\left(\frac{2\pi i f(x)}{p^m}\right) \right| \leq \left| \sum_{i=1}^d G_i d p^{\frac{1}{2}} \right|$$

$$= d^2 \frac{p^{-mn+\frac{1}{2}}}{p-1} \# A_{p,m}$$
$$\leq 2d^2 p^{-mn-\frac{1}{2}} D'_m p^{mn-(m-1)c_0(f)}$$
$$\leq D_m p^{-m\sigma},$$

because $\sigma = \min\left\{\frac{1}{2}, c_0(f)\right\}$. Here $D_m = 2(m-1)^2 D'_m$. □

The last subsum can be easily estimated by use of the Lang-Weil estimation (see [LW54]) and Theorem 3.2.2.

**Lemma 3.3.13.** *Let $f \in \mathbb{Z}[x_1, \ldots, x_n]$ be a non-constant polynomial, such that $f(0) = 0$. Put $\sigma = \min\{c_0(f), \frac{1}{2}\}$, where $c_0(f)$ is the log-canonical threshold of $f$ at $0$. Then there exist, for each integer $m > 1$, a natural number $N_m$ and a positive constant $D_m$, such that, for all $p > N_m$, we have*

$$\left| \sum_{\substack{\overline{x} \in (p\mathbb{Z}/p^m\mathbb{Z})^n, \\ \operatorname{ord}_p(f(x)) \geq m}} p^{-mn} \exp\left(\frac{2\pi i f(x)}{p^m}\right) \right| \leq D_m p^{-m\sigma}.$$

*Proof.* If $\operatorname{ord}_p(f(x)) \geq m$, then $\exp\left(\frac{2\pi i f(x)}{p^m}\right) = 1$ so we have

$$\left| \sum_{\substack{\overline{x} \in (p\mathbb{Z}/p^m\mathbb{Z})^n, \\ \operatorname{ord}_p(f(x)) \geq m}} p^{-mn} \exp\left(\frac{2\pi i f(x)}{p^m}\right) \right| = p^{-mn} \# B_{p,m},$$

where $B_{p,m} := \{\overline{x} \in (p\mathbb{Z}_p/p^m\mathbb{Z}_p)^n \mid \operatorname{ord}_p(f(x)) \geq m\}$. We can view $B_{p,m}$ as a subvariety of $\mathbb{F}_p^{mn}$. Then by the Lang-Weil estimation (see [LW54]), there exists a number $D_m$, which does not depend on $p$, such that

$$\# B_{p,m} = D_m p^{\dim_{\mathbb{F}_p}(B_{m,p})} + O(p^{\dim_{\mathbb{F}_p}(B_{m,p}) - \frac{1}{2}}).$$

By Theorem 3.2.2 we have $c_0(f) \leq \dfrac{mn - \dim_{\mathbb{F}_p}(B_{m,p})}{m}$, so $\dim_{\mathbb{F}_p}(B_{m,p}) \leq mn - mc_0(f)$. Hence, for all $p$ big enough,

$$\left| \sum_{\substack{\overline{x} \in (p\mathbb{Z}/p^m\mathbb{Z})^n, \\ \operatorname{ord}_p(f(x)) \geq m}} p^{-mn} \exp\left(\frac{2\pi i f(x)}{p^m}\right) \right| \leq p^{-mn} D_m p^{mn - mc_0(f)}$$

$$\leq D_m p^{-m\sigma}. \qquad \square$$

We will now put the three lemmas together to prove one of our main theorems. The essential ingredient in this proof is the expression that was obtained in Corollary 3.2.6.

*Proof of the Main Theorem 3.1.2.* From the Lemmas 3.3.9, 3.3.12 and 3.3.13 it follows that, for each $m > 1$, there exists a natural number $N_m$ and a positive constant $C_m$, such that for all $p > N_m$, we have

$$|E^0_{m,p}(f)| \leq C_m p^{-\sigma m}. \tag{3.3.2.1}$$

By Corollary 3.2.6 (with $\operatorname{Supp}(\overline{\Phi}) = \{0\}$), there exist constants $s, M', N' \in \mathbb{N}$, and for each $1 \leq i \leq s$, there exist constants $\beta_i \in \mathbb{N}$, $\lambda_i \in \mathbb{Q}$ and a definable set $A_i \subset \mathbb{N}$ in the Presburger language $\mathcal{L}_{\operatorname{Pres}}$, such that for all $p > N'$ and for all $1 \leq i \leq s$, there exists $a_{i,p} \in \mathbb{C}$ for which the formula

$$E^0_{m,p}(f) = \sum_{i=1}^{s} a_{i,p} m^{\beta_i} p^{-\lambda_i m} \mathbf{1}_{A_i}(m)$$

holds, for all $m > M'$. Moreover from the results in Section 3.2 we can deduce that $0 \leq \beta_i \leq n - 1$ and $c_0(f) \leq \lambda_i$ for all $1 \leq i \leq s$. After enlarging $M'$ and removing some small elements from $A_i$, we can assume that, for each subset $I \subset \{1, \ldots, s\}$, the set $\cap_{i \in I} A_i \setminus \cup_{i \notin I} A_i$ is either empty or infinite. Notice that for each $m > M'$, there is a unique subset $I \subset \{1, \ldots, s\}$, such that $m \in \cap_{i \in I} A_i \setminus \cup_{i \notin I} A_i$.

**Claim 3.3.14.** *There exist $M_0 > M', N_0 > N'$ and a positive constant $C_0$, such that for all $m > M_0, p > N_0$ and $1 \leq i \leq s$, we have*

$$|a_{i,p} p^{-\lambda_i m}| \leq C_0 p^{-\sigma m}.$$

Since there are only finitely many subsets $I \subset \{1, \ldots, s\}$, it is sufficient to fix a subset $I$ and prove the claim for $m$ restricted to the set $\cap_{i \in I} A_i \setminus \cup_{i \notin I} A_i$. Without loss of generality, we can assume that $I = \{1, \ldots, r\}$. If $p > N'$, $m \in \cap_{i \in I} A_i \setminus \cup_{i \notin I} A_i$ and $m > M'$, then we have

$$E_{m,p}^0(f) = \sum_{i=1}^{r} a_{i,p} m^{\beta_i} p^{-\lambda_i m}.$$

From Equation 3.3.2.1 we can see that, for such $m$ and for all $p > \max\{N', N_m\}$, we have

$$|E_{m,p}^0| = \left| \sum_{i=1}^{r} a_{i,p} m^{\beta_i} p^{-\lambda_i m} \right| \leq C_m p^{-\sigma m}.$$

This implies that

$$\left| \sum_{i=1}^{r} a_{i,p} m^{\beta_i} p^{(\sigma - \lambda_i) m} \right| \leq C_m.$$

It is easy to see that there exist $m_1, \ldots, m_r \in \cap_{i \in I} A_i \setminus \cup_{i \notin I} A_i$, all bigger than $M'$, and $N_I > \max\{N', N_{m_1}, \ldots, N_{m_r}\}$, such that all of the determinants of the size $r$ and $r-1$ submatrices of the matrix $B_p = (m_j^{\beta_i} p^{(\sigma - \lambda_i) m_j})_{1 \leq j, i \leq r}$ are different from zero for every $p > N_I$. We set

$$C_I := \max\{C_{m_i} \mid 1 \leq i \leq r\};$$
$$c_{j,p} := \sum_{i=1}^{r} a_{i,p} m_j^{\beta_i} p^{(\sigma - \lambda_i) m_j}, \text{ for } 1 \leq j \leq r;$$
$$D_p := \det(B_p);$$
$$D_{k,l,p} := (-1)^{k+l} \det\left( (m_j^{\beta_i} p^{(\sigma - \lambda_i) m_j})_{j \neq k, i \neq l} \right), \text{ for } 1 \leq k, l \leq r.$$

If we write $x_p = (a_{1,p}, \ldots, a_{r,p})^T$ and $c_p = (c_{1,p}, \ldots, c_{r,p})^T$, then $x_p$ is a solution of the equation $B_p x = c_p$. By our assumption on $m_1, \ldots, m_r$ we see that $D_p \neq 0$ and $D_{k,l,p} \neq 0$ for every $1 \leq k, l \leq r$ and for $p > N_I$. Using Cramer's rule we have

$$a_{i,p} = \frac{\sum_{j=1}^{r} c_{j,p} D_{j,i,p}}{D_p},$$

for all $1 \leq i \leq r$ and $p > N_I$. We remark that $|c_{j,p}| \leq C_I$, for all $p > N_I$, and that $\lambda_i \geq \sigma$, for all $1 \leq i \leq r$. This gives us

$$|a_{i,p}| \leq \frac{\sum_{j=1}^{r} |c_{j,p} D_{j,i,p}|}{|D_p|} \leq C_I \frac{\sum_{j=1}^{r} |D_{j,i,p}|}{|D_p|},$$

for all $1 \leq i \leq r$ and $p > N_I$. Then, by the definition of determinant, there exists $\alpha$, such that, for $1 \leq i \leq r$ and $p > N_I$ we have $|a_{i,p}| \leq p^\alpha$. Let $1 \leq i_0 \leq r$, we will now distinguish two cases.

If $\lambda_{i_0} > \sigma$, then there exists $M_{i_0} > M'$ such that, for every $m > M_{i_0}$ and $p > N_I$ we have
$$|a_{i_0,p}p^{-m\lambda_{i_0}}| \leq p^{\alpha-m\lambda_{i_0}} \leq p^{-m\sigma}.$$

If $\lambda_{i_0} = \sigma$, we observe that
$$D_p = \sum_{j=1}^{r} m_j^{\beta_{i_0}} p^{(\sigma-\lambda_{i_0})m_j} D_{j,i_0,p} = \sum_{j=1}^{r} m_j^{\beta_{i_0}} D_{j,i_0,p}.$$

By the definition of determinant, there exist $\gamma_j, d_j$, for each $1 \leq j \leq r$, such that $D_{j,i_0,p} = d_j p^{\gamma_j}$, when $p \to \infty$. By changing $m_1, \ldots, m_r$ if necessary, we can assume that there exists $d > 0$, such that $|D_p| = dp^\gamma$, when $p \to \infty$, where $\gamma = \max\{\gamma_j \mid 1 \leq j \leq r\}$. Thus there exist $C_0 > 0$ and $N_{i_0} > N_I$, such that
$$|a_{i_0,p}| \leq C_I \frac{\sum_{j=1}^{r} |D_{j,i_0,p}|}{|D_p|} \leq C_0,$$

for all $p > N_{i_0}$. And so
$$|a_{i_0,p}p^{-m\lambda_{i_0}}| \leq C_0 p^{-m\sigma},$$

for all $p > N_{i_0}$ and all $m > 1$. This proves the claim.

Hence we have
$$|E^0_{p,m}(f)| = \left|\sum_{i=1}^{s} a_{i,p} m^{\beta_i} p^{-\lambda_i m} \mathbf{1}_{A_i}(m)\right| \leq sC_0 m^{n-1} p^{-m\sigma},$$

for all $m > M_0, p > N_0$. By Equation 3.3.2.1 we also have, for each $1 < m \leq M_0$, an upper bound for $|E^0_{m,p}(f)|$ in terms of some constant $C_m$. Now let $N := \max\{N_i \mid i \in \{0, 2, \ldots, M_0\}\}$ and $C := \max\{sC_0, C_2, \ldots, C_{M_0}\}$, then we have
$$|E^0_{p,m}| = \left|\sum_{i=1}^{s} a_{i,p} m^{\beta_i} p^{-\lambda_i m} \mathbf{1}_{A_i}(m)\right| \leq Cm^{n-1} p^{-m\sigma},$$

for all $m > 1, p > N$. □

## 3.4 The second approach by geometry

We take $f \in \mathbb{Z}[x_1, \ldots, x_n]$ a nonconstant polynomial with $f(0) = 0$ and we put $\sigma = \min\{c_0(f), \frac{1}{2}\}$, where $c_0(f)$ is the log-canonical threshold of $f$ at 0. We use the notation of Section 3.2.2 with $(Y, h)$ an embedded resolution of $f^{-1}(0)$, $K = \mathbb{Q}$ and $\mathcal{O}_K = \mathbb{Z}$. Then by Theorem 3.2.5 and the discussion preceding that theorem, there exist $M_0, N_0 \in \mathbb{N}$, such that for all $p > N_0$, there exist at most $M_0$ non-trivial characters $\chi$ of $\mathbb{Z}_p^\times$ with $Z_{\Phi_p}(p, \chi, s, f) \neq 0$, where $\Phi_p = \mathbf{1}_{(p\mathbb{Z}_p)^n}$, i.e., $\mathrm{Supp}(\overline{\Phi}_p) = \{0\}$. Moreover any such character has conductor $c(\chi) = 1$. To simplify we will omit $\Phi_p$ and $f$ in the notation of Igusa's local zeta functions.

We can suppose that $f$ has good reduction mod $p$ for all $p > N_0$ (after enlarging $N_0$ if necessary). Let $p > N_0$ and let $E$ be an irreducible component of $h^{-1}(Z(f))$, such that $0 \in \overline{h}(E)$, then $h(E) \cap p\mathbb{Z}_p^n \neq \emptyset$. Remark that $h$ is proper, so $h(E)$ is a closed subvariety of $\mathbb{A}^n$. Therefore, after possibly enlarging $N_0$ again, we can assume that if $0 \notin h(E)$, then $h(E) \cap p\mathbb{Z}_p^n = \emptyset$, for all $p > N_0$. Hence, for $p > N_0$, $0 \in \overline{h}(\overline{E})$ implies $0 \in h(E)$. So the map $E \mapsto \overline{E}$ is a bijection between

$$\{E_i \mid i \in T, \ 0 \in h(E_i)\} \quad \text{and} \quad \{\overline{E}_i \mid i \in T, \ 0 \in \overline{h}(\overline{E}_i)\},$$

where $T$ is as in Section 3.2, hence

$$c_0(f) = \min_{i \in T : 0 \in \overline{h}(\overline{E}_i(\mathbb{F}_p))} \left\{ \frac{\nu_i}{N_i} \right\}. \tag{3.4.0.1}$$

Now to prove the Main Theorem 3.1.2, we use Proposition 3.2.3 for $p > N_0$, $u = 1$, $\pi = p$ and $m > 1$. This tells us that $E_{p,m}^0(f)$ is equal to

$$Z(p, \chi_{\text{triv}}, 0) + \operatorname{Coeff}_{t^{m-1}} \frac{(t-p)Z(p, \chi_{\text{triv}}, s)}{(p-1)(1-t)} + \sum_{\chi \neq \chi_{\text{triv}}} g_{\chi^{-1}} \operatorname{Coeff}_{t^{m-1}} Z(p, \chi, s).$$

**Lemma 3.4.1.** *There exist a positive constant $C$ and a natural number $N$, such that for all $m > 1$, $p > N$, we have*

$$Z(p, \chi_{\text{triv}}, 0) + \operatorname{Coeff}_{t^{m-1}} \frac{(t-p)Z(p, \chi_{\text{triv}}, s)}{(p-1)(1-t)} \leq Cm^{n-1} p^{-mc_0(f)}.$$

*Proof.* We use Theorem 3.2.4 which tells us that there exists a natural number $N'$, such that for all $p > N'$,

$$Z(p, \chi_{\text{triv}}, 0) = p^{-n} \sum_{I \subset T} c_{I, \chi_{\text{triv}}}^0 \prod_{i \in I} \frac{(p-1)p^{-\nu_i}}{1 - p^{-\nu_i}}; \tag{3.4.0.2}$$

$$Z(p, \chi_{\text{triv}}, s) = p^{-n} \sum_{I \subset T} c_{I, \chi_{\text{triv}}}^0 \prod_{i \in I} \frac{(p-1)t^{N_i} p^{-\nu_i}}{1 - t^{N_i} p^{-\nu_i}}.$$

From the formula $\frac{(t-p)}{(p-1)(1-t)} = -\frac{1}{p-1} - \frac{1}{1-t}$ we get

$$\operatorname{Coeff}_{t^{m-1}} \frac{(t-p)Z(p, \chi_{\text{triv}}, s)}{(p-1)(1-t)} = -\operatorname{Coeff}_{t^{m-1}} \frac{Z(p, \chi_{\text{triv}}, s)}{p-1} - \operatorname{Coeff}_{t^{m-1}} \frac{Z(p, \chi_{\text{triv}}, s)}{1-t},$$

where

$$\operatorname{Coeff}_{t^{m-1}} \frac{Z(p, \chi_{\text{triv}}, s)}{p-1} = \sum_{I \subset T} p^{-n} c_{I, \chi_{\text{triv}}}^0 (p-1)^{\#I} \frac{1}{p-1} \operatorname{Coeff}_{t^{m-1}} \prod_{i \in I} \frac{t^{N_i} p^{-\nu_i}}{1 - t^{N_i} p^{-\nu_i}}; \tag{3.4.0.3}$$

$$\operatorname{Coeff}_{t^{m-1}} \frac{Z(p, \chi_{\text{triv}}, s)}{1-t} = \sum_{I \subset T} p^{-n} c_{I, \chi_{\text{triv}}}^0 (p-1)^{\#I} \operatorname{Coeff}_{t^{m-1}} \frac{1}{1-t} \prod_{i \in I} \frac{t^{N_i} p^{-\nu_i}}{1 - t^{N_i} p^{-\nu_i}}. \tag{3.4.0.4}$$

Notice that if $I \subset T$, such that $\overset{\circ}{\overline{E}}_I \cap \overline{h}^{-1}(0) = \emptyset$, then $c_{I, \chi_{\text{triv}}}^0 = 0$. Hence we can assume that $\overset{\circ}{\overline{E}}_I \cap \overline{h}^{-1}(0) \neq \emptyset$. For such $I \subset T$ we have

$$\operatorname{Coeff}_{t^{m-1}} \prod_{i \in I} \frac{t^{N_i} p^{-\nu_i}}{1 - t^{N_i} p^{-\nu_i}} = \sum_{(a_i)_{i \in I} \in J_{I,m}} p^{-\sum_{i \in I} \nu_i (a_i + 1)}; \tag{3.4.0.5}$$

$$\operatorname{Coeff}_{t^{m-1}} \frac{1}{1-t} \prod_{i \in I} \frac{t^{N_i} p^{-\nu_i}}{1 - t^{N_i} p^{-\nu_i}} = \sum_{(a_i)_{i \in I} \in J'_{I,m}} p^{-\sum_{i \in I} \nu_i (a_i + 1)}, \tag{3.4.0.6}$$

where $J_{I,m} := \{(a_i)_{i \in I} \in \mathbb{N}^{\#I} \mid \sum_{i \in I} N_i(a_i+1) = m-1\}$ and $J'_{I,m} := \{(a_i)_{i \in I} \in \mathbb{N}^{\#I} \mid \sum_{i \in I} N_i(a_i+1) \leq m-1\}$. When $(a_i)_{i \in I} \in J_{I,m}$ and $p > N_0$, we can use Equation 3.4.0.1 for the following estimate:

$$-\sum_{i \in I} \nu_i(a_i+1) = -\sum_{i \in I} N_i \sigma_i(a_i+1)$$
$$= -\sum_{i \in I} N_i(a_i+1)(\sigma_i - c_0(f)) - (m-1)c_0(f)$$
$$\leq -(m-1)c_0(f),$$

where $\sigma_i = \frac{\nu_i}{N_i} \geq c_0(f)$, since we assumed that $\overset{\circ}{E}_I \cap \overline{h}^{-1}(0) \neq \emptyset$. We also deduce from this assumption that $\#I \leq n$, thus by Equation 3.4.0.5,

$$\text{Coeff}_{t^{m-1}} \prod_{i \in I} \frac{t^{N_i} p^{-\nu_i}}{1 - t^{N_i} p^{-\nu_i}} \leq \#(J_{I,m}) p^{-(m-1)c_0(f)} \leq m^{n-1} p^{-(m-1)c_0(f)}, \quad (3.4.0.7)$$

for all $p > N_0$. Using Equation 3.4.0.6 we can see that, in order to find an upper bound for the difference of 3.4.0.2 and 3.4.0.4, we need to analyse the expression

$$\prod_{i \in I} \frac{p^{-\nu_i}}{1 - p^{-\nu_i}} - \sum_{(a_i)_{i \in I} \in J'_{I,m}} p^{-\sum_{i \in I} \nu_i(a_i+1)} =$$
$$\sum_{(a_i)_{i \in I} \in \mathbb{N}^{\#I}} p^{-\sum_{i \in I} \nu_i(a_i+1)} - \sum_{(a_i)_{i \in I} \in J'_{I,m}} p^{-\sum_{i \in I} \nu_i(a_i+1)} =$$
$$\sum_{(a_i)_{i \in I} \in J''_{I,m}} p^{-\sum_{i \in I} \nu_i(a_i+1)},$$

where $J''_{I,m} := \{(a_i)_{i \in I} \in \mathbb{N}^{\#I} \mid \sum_{i \in I} N_i(a_i+1) \geq m\}$. Let $m_I := m + \max\{N_i \mid i \in I\}$ and $\overline{J}_{I,m} := \{(a_i)_{i \in I} \in \mathbb{N}^{\#I} \mid m \leq \sum_{i \in I} N_i(a_i+1) \leq m_I\}$. Afters some calculations we find that

$$\sum_{(a_i)_{i \in I} \in J''_{I,m}} p^{-\sum_{i \in I} \nu_i(a_i+1)} \leq \left(1 + \prod_{i \in I} \frac{1}{1 - p^{-\nu_i}}\right) \sum_{(a_i)_{i \in I} \in \overline{J}_{I,m}} p^{-\sum_{i \in I} \nu_i(a_i+1)}.$$

But if $(a_i)_{i \in I} \in \overline{J}_{I,m+1}$, then, for all $p > N_0$,

$$-\sum_{i \in I} \nu_i(a_i+1) = -\sum_{i \in I} N_i \sigma_i(a_i+1)$$
$$\leq -\sum_{i \in I} N_i(a_i+1)(\sigma_i - c_0(f)) - mc_0(f)$$
$$\leq -mc_0(f).$$

Therefore, for all $p > N_0$, we have,

$$\sum_{(a_i)_{i \in I} \in J''_{I,m}} p^{-\sum_{i \in I} \nu_i(a_i+1)} \leq (1 + 2^{\#(I)}) \#(\overline{J}_{I,m}) p^{-mc_0(f)} \leq C_I m^{n-1} p^{-mc_0(f)}, \quad (3.4.0.8)$$

where $C_I$ is a constant which does not depend on $m$ and $p$, for example $C_I = (1 + 2^{\#(I)})(\max\{N_i \mid i \in I\} + 1)$.

Now if $I \subset T$, then, by the Lang-Weil estimate (see [LW54]), there exists a constant $D_I$ and a natural number $N_I$, depending only on $I$, such that for all $p > N_I$, we have

$$c^0_{I,\chi_{\text{triv}}} = \sum_{a \in \overset{\circ}{\overline{E}}_I \cap \overline{h}^{-1}(0)} \Omega_{\chi_{\text{triv}}}(a) = \#\left(\overset{\circ}{\overline{E}}_I \cap \overline{h}^{-1}(0)\right) \leq \#\left(\overset{\circ}{\overline{E}}_I\right) \leq D_I p^{n-\#I}. \quad (3.4.0.9)$$

Putting together the inequalities 3.4.0.7, 3.4.0.8 and 3.4.0.9 with the formulas 3.4.0.2, 3.4.0.3 and 3.4.0.4, we find that there exists a natural number $N > \max\{N_0, N', (N_I)_{I \subset T}\}$, such that for all $p > N$, we have

$$Z(p, \chi_{\text{triv}}, 0) + \text{Coeff}_{t^{m-1}} \frac{(t-p)Z(p, \chi_{\text{triv}}, s)}{(p-1)(1-t)}$$

$$\leq \sum_{I \subset T: \overset{\circ}{\overline{E}}_I \cap \overline{h}^{-1}(0) \neq \emptyset} p^{-n} c^0_{I,\chi_{\text{triv}}} (p-1)^{\#I} m^{n-1} \left(\frac{p^{-mc_0(f)+c_0(f)}}{p-1} + C_I p^{-mc_0(f)}\right)$$

$$\leq \sum_{I \subset T: \overset{\circ}{\overline{E}}_I \cap \overline{h}^{-1}(0) \neq \emptyset} p^{-n} D_I p^{n-\#I} (p-1)^{\#I} m^{n-1} \left(\frac{p^{-mc_0(f)+c_0(f)}}{p-1} + C_I p^{-mc_0(f)}\right)$$

$$\leq \sum_{I \subset T: \overset{\circ}{\overline{E}}_I \cap \overline{h}^{-1}(0) \neq \emptyset} D_I(C_I + 2) m^{n-1} p^{-mc_0(f)} \leq C m^{n-1} p^{-mc_0(f)},$$

where $C = \sum_{I \subset T} D_I(C_I + 2)$ is a constant that is independent of $p$ and $m$ and where we have used the fact that $c_0(f) \leq 1$. $\square$

**Lemma 3.4.2.** *There exist a positive constant $C$ and a natural number $N$, such that for all $m > 1$, $p > N$, we have*

$$\left| \sum_{\chi \neq \chi_{\text{triv}}} g_{\chi^{-1}} \text{Coeff}_{t^{m-1}} Z(p, \chi, s) \right| \leq C m^{n-1} p^{-m\sigma}.$$

*Proof.* We continue to use Theorem 3.2.4, hence there exists a natural number $N'$, such that for all $p > N'$,

$$Z(p, \chi, s) = p^{-n} \sum_{\substack{I \subset T, \\ \forall i \in I: d | N_i}} c^0_{I,\chi} \prod_{i \in I} \frac{(p-1) t^{N_i} p^{-\nu_i}}{1 - t^{N_i} p^{-\nu_i}},$$

with $\chi$ a character of order $d$ on $\mathbb{Z}_p^\times$ with conductor $c(\chi) = 1$.

For a subset $I \subset T$, such that $d | N_i, \forall i \in I$, and $\overset{\circ}{\overline{E}}_I \cap \overline{h}^{-1}(0) \neq \emptyset$, we have

$$\text{Coeff}_{t^{m-1}} \prod_{i \in I} \frac{t^{N_i} p^{-\nu_i}}{1 - t^{N_i} p^{-\nu_i}} = \sum_{(a_i)_{i \in I} \in J_{I,m}} p^{-\sum_{i \in I} \nu_i(a_i+1)},$$

where $J_{I,m}$ is as in the proof of Lemma 3.4.1. By the equations 3.4.0.5 and 3.4.0.7 we have

$$\sum_{(a_i)_{i \in I} \in J_{I,m}} p^{-\sum_{i \in I} \nu_i(a_i+1)} \leq m^{n-1} p^{-mc_0(f)+c_0(f)}.$$

We use the Lang-Weil estimate ([LW54]) again, as we did in Lemma 3.4.1. So there exist a constant $D_I$ and a natural number $N_I$, depending only on $I$, such that for all $p > N_I$, we have

$$|c_{I,\chi}^0| = \left| \sum_{a \in \overset{\circ}{\overline{E}}_I \cap \overline{h}^{-1}(0)} \Omega_\chi(a) \right| \leq \sum_{a \in \overset{\circ}{\overline{E}}_I \cap \overline{h}^{-1}(0)} 1 = \#\left(\overset{\circ}{\overline{E}}_I \cap \overline{h}^{-1}(0)\right)$$

$$\leq \#\left(\overset{\circ}{\overline{E}}_I\right) \leq D_I p^{n-\#(I)}.$$

If we take $N'' > \max_{I \subset T} N_I$, then we find that for all $p > N''$,

$$|\text{Coeff}_{t^{m-1}} Z(p,\chi,s)| \leq \sum_{\substack{I \subset T, \\ \forall i \in I : d|N_i}} p^{-n} D_I p^{n-\#(I)} (p-1)^{\#(I)} m^{n-1} p^{-mc_0(f)+c_0(f)}$$

$$\leq \sum_{\substack{I \subset T, \\ \forall i \in I : d|N_i}} D_I m^{n-1} p^{-mc_0(f)+c_0(f)}$$

$$\leq C' m^{n-1} p^{-mc_0(f)+c_0(f)},$$

where $C' := \sum_{I \subset T} D_I$. Furthermore, by a standard result on Gauss sums, we can see that, if $\chi \neq \chi_{\text{triv}}$, then $|g_{\chi^{-1}}| \leq Dp^{-\frac{1}{2}}$, for some constant $D$, that does not depend on $\chi$ and $p$. By Theorem 3.2.5 and the discussiong preceding this theorem, we know that for $p > N_0$, the set $\Upsilon_p$ of non-trivial characters $\chi$ such that $Z(p,\chi,s)$ is not zero, has atmost $M_0$ elements, for some positive integer $M_0$. Moreover, all these characters have conductor 1. So there exists a natural number $N > \max\{N_0, N''\}$, such that for all $p > N$ and $m > 1$, we have

$$\left| \sum_{\chi \neq \chi_{\text{triv}}} g_{\chi^{-1}} \text{Coeff}_{t^{m-1}} Z(p,\chi,s) \right| = \left| \sum_{\chi \in \Upsilon_p} g_{\chi^{-1}} \text{Coeff}_{t^{m-1}} Z(p,\chi,s) \right|$$

$$\leq \sum_{\chi \in \Upsilon_p} |g_{\chi^{-1}}| C' m^{n-1} p^{-mc_0(f)+c_0(f)}$$

$$\leq \sum_{\chi \in \Upsilon_p} C' D m^{n-1} p^{-m\sigma + \sigma - \frac{1}{2}}$$

$$\leq C m^{n-1} p^{-m\sigma}$$

where $C = M_0 C' D$ is a constant that is independent of $p$ and $m$ and where we have used the fact that $\sigma = \min\{c_0(f), \frac{1}{2}\}$. □

**Remark 3.4.3.** These two proofs still work if we take $\Phi_p = \mathbf{1}_{U_p}$ instead of $\mathbf{1}_{(p\mathbb{Z}_p)^n}$, where $U_p$ is a union of some multiballs $y + (p\mathbb{Z}_p)^n$ in $\mathbb{Z}_p^n$, such that $C_{\overline{f}} \cap \overline{U}_p \subset \overline{f}^{-1}(0)$ (this is needed in the proof of Lemma 3.4.2 to apply Theorem 3.2.5). We have to replace $c_0(f)$ for example by $c(f)$, $\overset{\circ}{\overline{E}}_I \cap \overline{h}^{-1}(0)$ by $\overset{\circ}{\overline{E}}_I \cap \overline{h}^{-1}(\overline{U}_p)$ and $c_{I,\chi_{\text{triv}}}^0$ by $c_{I,\mathbf{1}_{U_p},\chi_{\text{triv}}}$. The constant $C$ and the natural number $N$ that are found in these proofs, do not depend on $U_p$. They do depend however on $f$ and on the embedded resolution $(Y, h)$ of $f$.

*Proof of the Main Theorem 3.1.2.* The proof follows by combining the two Lemmas 3.4.1 and 3.4.2 and using the fact that $\sigma \leq c_0(f)$. □

## 3.5 Proof of the Main Theorem 3.1.3.

In this section we will prove the Main Theorem 3.1.3 by adapting the proofs from Section 3.4. First, we need the following lemma.

**Lemma 3.5.1.** *Let $f \in \mathbb{Z}[x_1, \ldots, x_n]$ and $V_{f,p}$ be the set of critical values $z$ of $f$ in $\mathbb{Q}_p$. Then $\#(V_{f,p})$ has an upper bound $d$, that does not depend on $p$. Furthermore, there exists $N$, such that for all $p > N$, the following holds:*

1. *for all $z \in V_{f,p}$, we have $\mathrm{ord}_p(z) = 0$;*

2. *for any two distinct points $z_1, z_2$ in $V_{f,p}$, we have $\mathrm{ord}_p(z_1 - z_2) = 0$;*

3. *if $x \in \mathbb{Z}_p^n$ such that $\mathrm{ord}_p(f(x) - z) = 0$ for all $z \in V_{f,p}$, then $x$, resp. $\overline{x}$, is a regular point of $f$, resp. $\overline{f} := (f \mod p)$.*

*Proof.* Remark that we can uniquely extend the valuation $\mathrm{ord}_p$ to $\overline{\mathbb{Q}}_p$ (the algebraic closure of $\mathbb{Q}_p$). We denote by $\mathcal{O}_p = \{z \in \overline{\mathbb{Q}}_p | \mathrm{ord}_p(z) \geq 0\}$ the ring of integers of $\overline{\mathbb{Q}}_p$ and by $\mathcal{M}_p = \{z \in \overline{\mathbb{Q}}_p | \mathrm{ord}_p(z) > 0\}$ its maximal ideal.

The set of critical values $V_f$ of $f$ is a definable set in $\mathcal{L}_{\mathrm{ring}}$ given by

$$z \in V_f \Leftrightarrow \exists y \left[ z = f(y) \wedge \frac{\partial f}{\partial x_1}(y) = 0 \wedge \ldots \wedge \frac{\partial f}{\partial x_n}(y) = 0 \right].$$

By elimination of quantifiers in the $\mathrm{ACF}_0$-theory, i.e., the theory of algebraically closed fields of characteristic 0, and because of the fact that $V_f$ is a finite set, there exist non-zero polynomials $T(z) \in \mathbb{Z}[z]$ and $R(z) \in \overline{\mathbb{Q}}[z]$, such that $V_f = Z(R) \subset Z(T)$. Moreover, we can assume that $T(z)$ and $R(z)$ only have simple roots in $\overline{\mathbb{Q}}$. By logical compactness, there exists $N_0$, such that for all $p > N_0$, $T_p(z) \in \mathbb{F}_p(z)$ and $R_p(z) \in \overline{\mathbb{F}}_p(z)$ also only have simple roots in $\overline{\mathbb{F}}_p$ and $V_{\overline{f}} = Z(R_p) \subset Z(T_p) \subset \overline{\mathbb{F}}_p$, where $T_p := (T \mod p)$ and $R_p := (R \mod \mathcal{M}_p)$. Since $V_{f,p} \subset V_f$, we have $\#(V_{f,p}) \leq \#(V_f) = \deg(R) =: d$. Because $Z(T) \subset \overline{\mathbb{Q}}$ is a finite set of algebraic numbers, there exists $N \geq N_0$, such that for all $p > N$, the conditions (1) and (2) are satisfied, not only for $V_{f,p}$, but for $Z(T)$ and $Z(R)$ as well.

To prove condition (3), we take $p > N$ and $x \in \mathbb{Z}_p^n$ such that $\mathrm{ord}_p(f(x) - z) = 0$ for all $z \in V_{f,p}$. Then $f(x) \notin V_{f,p}$, so $x$ is a regular point of $f$. Suppose, for a contradiction, that $\overline{x}$ is a critical point of $\overline{f}$, then $z' := \overline{f}(\overline{x}) \in V_{\overline{f}} = Z(R_p) \subset Z(T_p)$. From the facts that $T_p$ has only simple roots in $\overline{\mathbb{F}}_p$, $z' \in \mathbb{F}_p$ and $T_p(z') = 0$, it follows by Hensel's lemma that there exists $z_1 \in \mathbb{Z}_p$ such that $T(z_1) = 0$ and $\overline{z}_1 = z'$. Hence $\mathrm{ord}_p(f(x) - z_1) > 0$, and therefore $z_1 \notin V_{f,p}$. On the other hand, $R_p$ has also only simple roots in $\overline{\mathbb{F}}_p$ and $z' \in Z(R_p)$, so, by Hensel's lemma, there exists $z_2 \in \mathcal{O}_p$ such that $R(z_2) = 0$ and $\overline{z}_2 = z'$. From the facts that $z_1$ and $z_2$ are both roots of $T$, $\overline{z}_1 = z' = \overline{z}_2$ and the conditions (1) and (2) are true for $Z(T)$, it follows that $z_1 = z_2$. Hence $z_1 \in Z(R) = V_f$, and we knew already that $z_1 \in \mathbb{Z}_p$ so $z_1 \in V_{f,p}$. This contradiction proves that condition (3) also holds. □

*Proof of the Main Theorem 3.1.3.* Let $N, d$ be as in Lemma 3.5.1 and write $V_f = \{z_1, \ldots, z_d\}$. We fix $p > N$, then we can assume that $V_{f,p} = \{z_1, \ldots, z_r\}$ with $r \leq d$. For each $1 \leq i \leq r$, we put $\Phi_{i,p} := \mathbf{1}_{\{x \in \mathbb{Z}_p^n | \mathrm{ord}_p(f(x) - z_i) > 0\}} : \mathbb{Q}_p^n \to \mathbb{C}$. Because $f \in \mathbb{Z}[x_1, \ldots, x_n]$ and by Lemma 3.5.1 we see that $\Phi_{i,p}$ is residual, for all $1 \leq i \leq r$, and that $\mathrm{Supp}(\Phi_{i,p}) \cap \mathrm{Supp}(\Phi_{j,p}) = \emptyset$, if $i \neq j$. We denote $\Phi_{0,p} := \mathbf{1}_{\mathbb{Z}_p^n} - \sum_{i=1}^{r} \Phi_{i,p}$, then $\Phi_{0,p}$ will also be residual. Now we have

$$\begin{aligned}
E_p(z,f) &= \int_{\mathbb{Z}_p^n} \Psi(zf(x))|dx| \\
&= \sum_{i=0}^{r} \int_{\mathbb{Z}_p^n} \Phi_{i,p}(x)\Psi(zf(x))|dx| \\
&= \sum_{i=1}^{r} \int_{\mathbb{Z}_p^n} \Phi_{i,p}(x)\Psi(z(f(x)-z_i))\exp(2\pi i z z_i)|dx| + E_{\Phi_{0,p}}(z,p,f) \\
&= \sum_{i=1}^{r} \exp(2\pi i z z_i) E_{\Phi_{i,p}}(z,p,f_i) + E_{\Phi_{0,p}}(z,p,f),
\end{aligned}$$

where $f_i(x) = f(x) - z_i$ for $1 \leq i \leq r$.

If $1 \leq i \leq r$, we note that $C_{\overline{f}_i} \cap \mathrm{Supp}(\overline{\Phi}_{i,p}) \subset \overline{f}_i^{-1}(0)$. So we can use Theorems 3.2.4 and 3.2.5 for $f_i$. According to Remark 3.4.3, the Main Theorem 3.1.2 is still true for the exponential sum $E_{\Phi_{i,p}}(z,p,f_i)$, where we take $\sigma_i = \min\left\{c(f_i), \frac{1}{2}\right\}$. In the proofs from Section 3.4 we need to replace $c_{I,\chi}^0$ by $c_{I,\Phi_{i,p},\chi}$ and $\overset{\circ}{\overline{E}}_I \cap \overline{h}^{-1}(0)$ by $\overset{\circ}{\overline{E}}_I \cap \overline{h}_i^{-1}(Z(\overline{f}_i))$, with $h_i : Y_i \to \mathbb{Q}(z_i)^n$ an embedded resolution for $Z(f_i)$. For each $1 \leq i \leq r$, there exist a constant $C_i$ and a natural number $N_i > N$, only depending on the critical value $z_i \in V_f$ and the chosen resolution $h_i$ of $f_i$, such that, if $p > N_i$ and $z_i \in V_{f,p}$, then we have

$$|E_{\Phi_{i,p}}(z,p,f_i)| \leq C_i m^{n-1} p^{-m\sigma_i}.$$

We remark that by definition of $\Phi_{0,p}$ and by condition (3) from Lemma 3.5.1, we have $C_{\overline{f}} \cap \mathrm{Supp}(\overline{\Phi}_{0,p}) = \emptyset$, for all $p > N$, and thus it is well known that $E_{\Phi_{0,p}}(z,p,f) = 0$, for $|z| > p$ (see [Den], Remark 4.5.3).

We recall that $a(f)$ is the minimum, over all $b \in \mathbb{C}$, of the log-canonical thresholds of the polynomials $f(x) - b$. Therefore, if we set $\sigma = \min\left\{a(f), \frac{1}{2}\right\}$, then $\sigma \leq \min_{1 \leq i \leq d} \sigma_i$, hence there exist a constant $C > \max_{1 \leq i \leq d} C_i$ and a natural number $N' > \max_{1 \leq i \leq d} N_i$, such that for all $p > N'$ and $m \geq 2$, we have

$$|E_{m,p}| \leq C m^{n-1} p^{-m\sigma}. \qquad \square$$

## 3.6 The uniform version of the Main Theorem 3.1.2.

In this section $f \in \mathbb{Z}[x_1, \ldots, x_n]$ is a nonconstant polynomial. We will describe how to adapt the Sections 3.3 and 3.4 to obtain a constant $C$ and a natural number $N$, such that for all $y \in \mathbb{Z}^n$ and for all $m \geq 1$, $p > N$, we have

$$|E_{m,p}^y(f)| := \left| \frac{1}{p^{mn}} \sum_{\overline{x} \in \overline{y} + (p\mathbb{Z}/p^m\mathbb{Z})^n} \exp\left(\frac{2\pi i f(x)}{p^m}\right) \right| \leq C m^{n-1} p^{-m\sigma_{y,p}}. \qquad (3.6.0.1)$$

Here we take $\sigma_{y,p} = \min\{a_{y,p}(f), \frac{1}{2}\}$. We recall that $a_{y,p}(f)$ is the minimum of the log canonical thresholds at $y'$ of the polynomials $f(x) - f(y')$, where $y'$ runs over $y + (p\mathbb{Z}_p)^n$. Notice that the case $m = 1$ is covered by Remark 3.1.4. Hence we can assume that $m \geq 2$.

Let $V_f = \{z_1, \ldots, z_d\} \subset \overline{\mathbb{Q}}$ be the set of critical values of $f$ in $\overline{\mathbb{Q}}$, where $d$ is as in Lemma 3.5.1. For each $1 \leq j \leq d$ we put $f_j(x) := f(x) - z_j$ and we fix an embedded resolution $(Y_j, h_j)$ of $f_j^{-1}(0)$ over $\mathbb{Q}(z_i)$. Let $N_j'$ be a natural number, such that for all $p > N_j'$, $(Y_j, h_j)$ has good reduction modulo $p$. Furthermore, let $N_0'$ be a natural number, such that for all $p > N_0'$, we have $V_{f,p} = V_f \cap \mathbb{Q}_p \subset \mathbb{Z}_p$, any two distinct points $z, z'$ in $V_{f,p}$ satisfy $\mathrm{ord}_p(z - z') = 0$ and if $x \in \mathbb{Z}_p^n$ such that $\mathrm{ord}_p(f(x) - z) = 0$ for all $z \in V_{f,p}$, then $x$, resp. $\overline{x}$, is a regular point of $f$, resp. $\overline{f}$ (see Lemma 3.5.1). We put $N' := \max_{0 \leq i \leq d} N_i'$ and for each $p > N'$ we consider a partition of $\mathbb{Z}^n = \bigcup_{j=0}^d A_{j,p} \cup \bigcup_{j=1}^d B_{j,p}$, where

$$A_{j,p} := \{y \in \mathbb{Z}^n \mid \mathrm{ord}_p(f_j(y)) > 0 \text{ and } f \text{ has a critical point in } y + (p\mathbb{Z}_p)^n\},$$
$$B_{j,p} := \{y \in \mathbb{Z}^n \mid \mathrm{ord}_p(f_j(y)) > 0 \text{ and } f \text{ has no critical points in } y + (p\mathbb{Z}_p)^n\},$$

for $1 \leq j \leq d$, and

$$A_{0,p} := \mathbb{Z}^n \setminus \bigcup_{j=1}^d (A_{j,p} \cup B_{j,p}).$$

First of all, for $p > N'$, we observe that if $y \in A_{0,p}$, then $\mathrm{ord}_p(f(y) - z_j) \leq 0$ for all $1 \leq j \leq d$. In particular, $\mathrm{ord}_p(f(y) - z_j) = 0$ for all $z_j \in V_f \cap \mathbb{Z}_p = V_{f,p}$. So $\overline{y}$ is a regular point of $\overline{f}$, by Lemma 3.5.1, hence the condition $C_{\overline{f}} \cap \mathrm{Supp}(\overline{\Phi}_{y,p}) = \emptyset$, with $\Phi_{y,p} := \mathbf{1}_{y+(p\mathbb{Z}_p)^n}$, is satisfied. Thus, by Remark 4.5.3 from [Den], we get that $E_{m,p}^y(f) = 0$, for all $m \geq 2$, $p > N'$ and $y \in A_{0,p}$.

Secondly, if $1 \leq j \leq d$, $p > N'$, and $y \in B_{j,p}$, then $f_j$ has no critical points in $y + (p\mathbb{Z}_p)^n$. So by 1.4.1 from [Den], we have $E_{m,p}^y(f_j) = 0$, for $m$ large enough. Using Corollary 1.4.5 from [Den], we see that $(p^{s+1} - 1)Z_{\Phi_{y,p}}(p, \chi_{\mathrm{triv}}, s, f_j)$ and $Z_{\Phi_{y,p}}(p, \chi, s, f_j)$, for $\chi \neq \chi_{\mathrm{triv}}$, cannot have any poles. Because the resolution $(Y_j, h_j)$ of $f_j$ has good reduction modulo $p$, for $p > N'$, and $C_{\overline{f}_j} \cap \mathrm{Supp}(\overline{\Phi}_{y,p}) \subset \overline{f}_j^{-1}(0)$, for $y \in B_{j,p}$, the Theorem 3.2.5 applies. By combining it with Proposition 3.2.3, we get that for all $p > N'$ and $y \in B_{j,p}$, the sum $E_{m,p}^y(f_j)$ equals

$$Z_{\Phi_{y,p}}(p, \chi_{\mathrm{triv}}, 0, f_j) + \mathrm{Coeff}_{t^{m-1}}\left(\frac{(t-p)Z_{\Phi_{y,p}}(p, \chi_{\mathrm{triv}}, s, f_j)}{(p-1)(1-t)}\right) \quad (3.6.0.2)$$
$$+ \sum_{\substack{\chi \neq \chi_{\mathrm{triv}}, \\ c(\chi)=1}} g_{\chi^{-1}}\chi(u)\mathrm{Coeff}_{t^{m-1}}(Z_{\Phi_{y,p}}(p, \chi, s, f_j)).$$

Since $Z_{\Phi_{y,p}}(p, \chi, s, f_j)$ does not have any poles for $\chi \neq \chi_{\mathrm{triv}}$, we can see that, for $m$ big enough, $\mathrm{Coeff}_{t^{m-1}}(Z_{\Phi_{y,p}}(p, \chi, s, f_j))$ will not depend on $m$. Also the total expression 3.6.0.2 is independent of $m$, for $m$ big enough (because it is equal to 0). Therefore the part $\mathrm{Coeff}_{t^{m-1}}\left(\frac{(t-p)Z_{\Phi_{y,p}}(p, \chi_{\mathrm{triv}}, s, f_j)}{(p-1)(1-t)}\right)$ must be independent of $m$ as well, for $m$ big enough. This can only be the case if $\frac{(t-p)Z_{\Phi_{y,p}}(p,\chi_{\text{triv}},s,f_j)}{(p-1)(1-t)}$, as a function in $t$, has at most two poles, one pole at $t = 1$ of order 1 and one pole at $t = 0$. However, the explicit formula of $Z_{\Phi_{y,p}}(p, \chi_{\text{triv}}, s, f_j)$ implies that it can not have poles at $t = 0$. So $\frac{(t-p)Z_{\Phi_{y,p}}(p,\chi_{\text{triv}},s,f_j)}{(p-1)(1-t)}$ has at most one pole, and this pole (if it exists) must be of order 1 at $t = 1$.

According to 4.1.1 from [Den], the degree of $Z_{\Phi_p}(p, \chi, s, f_j) \leq 0$ (as a rational function in $t$), for all $p > N'$ and all characters $\chi$ with conductor $c(\chi) = 1$. This implies that $\frac{(t-p)Z_{\Phi_{y,p}}(p,\chi_{\text{triv}},s,f_j)}{(p-1)(1-t)}$ is of the form $c + \frac{d}{1-t}$, for certain $c, d \in \mathbb{C}$, and that $Z_{\Phi_{y,p}}(p, \chi, s, f_j)$ is equal to a constant function, for $\chi \neq \chi_{\text{triv}}$. Now we can easily see that for all $m \geq 2$, $\text{Coeff}_{t^{m-1}}\left(\frac{(t-p)Z_{\Phi_{y,p}}(p,\chi_{\text{triv}},s,f_j)}{(p-1)(1-t)}\right)$ and $\text{Coeff}_{t^{m-1}}(Z_{\Phi_{y,p}}(p, \chi, s, f_j))$, for $\chi \neq \chi_{\text{triv}}$, are indepent of $m$. We conclude that $E^y_{m,p}(f_j) = 0$, for all $m \geq 2$, $p > N'$ and $y \in B_{j,p}$.

The last case is the one where $y \in A_{j,p}$, for $1 \leq j \leq d$. We will show that in this case there exists a constant $C_j$ and a natural number $N_j$ (only depending on $j$, not on $y$), such that for all $p > N_j$, $m \geq 2$ and $y \in A_{j,p}$, we have

$$|E^y_{m,p}(f_j)| \leq C_j m^{n-1} p^{-m\sigma_{y,p}}. \tag{3.6.0.3}$$

By taking $N := \max\{N', N_1, \ldots, N_d\}$ and $C := \max\{C_1, \ldots, C_d\}$ (both independent of $y$), the formula 3.6.0.1 will hold for all $y \in \mathbb{Z}^n$, $p > N$ and $m \geq 1$. In what follows, we will show how to adapt the proofs of both Sections 3.3 and 3.4, to obtain the formula 3.6.0.3.

### 3.6.1 Adapting Section 3.4

If we want to be able to use the method of proof that was outlined in Section 3.4, then we need to show the following result, for all $j > 0$, $y \in A_{j,p}$ and $p > N'$:

$$a_{y,p}(f) = \min_{E: \overline{y} \in \overline{h}_j(\overline{E}(\mathbb{F}_p))} \left\{\frac{\nu}{N}\right\}, \tag{3.6.1.1}$$

where $E$ is an irreducible component of $h_j^{-1}(Z(f_j))$ with numerical data $(N, \nu)$. When we compare this to the Formula 3.4.0.1, we see that, by replacing $c_0(f)$ by $a_{y,p}(f)$, we can adapt the results of Section 3.4 to $f_j$ with $\Phi_p = \mathbf{1}_{y+(p\mathbb{Z}_p)^n}$. Indeed the condition $C_{\overline{f}_j} \cap \text{Supp}(\overline{\Phi}_p) \subset \overline{f}_j^{-1}(0)$ is satisfied. This proves the Formula 3.6.0.3 and by Remark 3.4.3 we know that the constant $C_j$ and the natural number $N_j$ only depend on $f_j$ and the chosen resolution $(Y_j, h_j)$.

All that is left, is to prove Equation 3.6.1.1 for $y \in A_{j,p}$ and $j > 0$. We remark that if $y' \in y + (p\mathbb{Z}_p)^n$ is not a critical point of $f$, then $c_{y'}(f(x) - f(y')) = 1$. If $y' \in y+(p\mathbb{Z}_p)^n$ is a critical point of $f$, then we know by Lemma 3.5.1 that $f(y') = z_j$, hence $f_j(y') = f(y') - f(y') = 0$. Since $(Y_j, h_j)$ has good reduction modulo $p$, for $p > N'$, we know that, after possibly enlarging $N'$ as we did for 3.4.0.1, we have

$$c_{y'}(f(x) - f(y')) = c_{y'}(f_j) = \min_{E: \overline{y'} \in \overline{h}_j(\overline{E}(\mathbb{F}_p))} \left\{\frac{\nu}{N}\right\} = \min_{E: \overline{y} \in \overline{h}_j(\overline{E}(\mathbb{F}_p))} \left\{\frac{\nu}{N}\right\} \leq 1.$$

If $y \in A_{j,p}$, then $y + (p\mathbb{Z}_p)^n$ contains at least one critical point of $f$, in which case Equation 3.6.1.1 holds.

### 3.6.2 Adapting Section 3.3

For $j > 0$ and $y \in A_{j,p}$, we will split the exponential sum $E^y_{m,p}(f_j)$ into three subsums in exactly the same way as in Section 3.3. In each of the Lemmas 3.3.9, 3.3.12, 3.3.13 and in the proof of the Main Theoreom 3.1.2 from Section 3.3 we need to make some changes.

**Lemma 3.6.1.** *Let $f \in \mathbb{Z}[x_1, \ldots, x_n]$ be a nonconstant polynomial and let $z_j \in V_f$ be a critical value of $f$. There exists a natural number $N_0 > N'$, such that for all $m \geq 1$, for all $p > N_0$ and for all $y \in A_{j,p}$, we have*

$$\sum_{\substack{\overline{x} \in \overline{y} + (p\mathbb{Z}/p^m\mathbb{Z})^n, \\ \operatorname{ord}_p(f_j(x)) \leq m-2}} \exp\left(\frac{2\pi i f_j(x)}{p^m}\right) = 0.$$

*Remark that if $A_{j,p} \neq \emptyset$, then $z_j \in V_{f,p} \subset \mathbb{Z}_p$, so the term $\exp\left(\dfrac{2\pi i f_j(x)}{p^m}\right)$ is well-defined.*

To prove this lemma, we adapt the proof of Lemma 3.3.9 as follows. We replace the formula $\phi$ by

$$\phi_j(x_1, \ldots, x_n, z, \xi_1, \ldots, \xi_n, m) = \bigwedge_{i=1}^n (\overline{x}_i = \xi_i) \wedge (\operatorname{ord}(z - z_j) \leq m - 2) \wedge (\operatorname{ord}(z - f(x_1, \ldots, x_n)) \geq m),$$

where $x_i, z$ are in the valued field-sort, $\xi_i$ are in the residu field-sort and $m$ is in the value group-sort. This is an $\mathcal{L}_\mathbb{Z} \cup \{z_j\}$-formula, with $z_j$ a constant symbol in the valued field-sort. We remark that the function $\mathcal{O}_K \to k_K : x \mapsto \overline{x} = (x \mod \mathcal{M}_K)$ is definable in $\mathcal{L}_{\mathrm{DP}}$.

Now $\phi_j$ induces a definable subassignment $X_j \subset h[n+1, n, 1]$ and constructible functions $F_j := \mathbf{1}_{X_j}$ and $G_j := \pi_!(F_j)$, where $\pi : h[n+1, n, 1] \to h[1, n, 1]$ is the projection onto the last $n+2$ coordinates. For each prime $p$, for each uniformiser $\varpi_p$ of $\mathbb{Q}_p$ and for each $y \in A_{i,p}$, we have the following interpretation of $G_j$ in $\mathbb{Q}_p$:

$$G_{j,\varpi_p}(z, \overline{y}, m) = \#\{\overline{x}^{(m)} \in \overline{y}^{(m)} + (p\mathbb{Z}/p^m\mathbb{Z})^n \mid f(x) \equiv z \mod p^m\},$$

if $\operatorname{ord}_p(z - z_j) \leq m - 2$, and

$$G_{j,\varpi_p}(z, \overline{y}, m) = 0,$$

if $\operatorname{ord}_p(z - z_j) \geq m - 1$. Here the notation $\overline{x}^{(m)}$ means the class of $(x \mod p^m)$. Note however that $G_{j,\varpi_p}$ actually only depends on $(y \mod p)$, i.e., on $\overline{y}$. We remark that if $A_{j,p} \neq \emptyset$, then $z_j \in \mathbb{Q}_p$, which makes it possible to interprete $\operatorname{ord}(z - z_j)$ (and other formulas that contain the symbol $z_j$) in $\mathbb{Q}_p$. We apply Corollary 3.3.8 to $G_j$ to obtain a cell decomposition where the centers $c_i$ are given by $\mathcal{L}_\mathbb{Z} \cup \{z_j\}$-formulas $\theta_j(z, \xi, \eta, \gamma, m)$. By elimination of quantifiers, $\theta_i$ is equivalent to the formula

$$\bigvee_k \left(\zeta_{ik}\big(\overline{\mathrm{ac}}g_1(z), \ldots, \overline{\mathrm{ac}}g_s(z), \xi, \eta\big) \wedge \nu_{ik}\big(\operatorname{ord} g_1(z), \ldots, g_s(z), \gamma, m\big)\right),$$

where $\zeta_{ik}$ is an $\mathcal{L}_{\text{ring}}$-formula and $\nu_{ik}$ an $\mathcal{L}_{\text{oag}}$-formula, and $g_1, \ldots, g_s \in (\mathbb{Z}[z_j][[t]])[z]$. The rest of the proof of Lemma 3.3.9 still applies if we replace $\text{ord}(z)$ by $\text{ord}(z - z_j)$ everywhere. By going over the proof, we can see that the natural number $N_0$ that is obtained in the proof, only depends on $j$.

**Lemma 3.6.2.** *Let $f \in \mathbb{Z}[x_1, \ldots, x_n]$ be a nonconstant polynomial and let $z_j \in V_f$ be a critical value of $f$. There exists, for each integer $m > 1$, a natural number $N_m > N'$ and a positive constant $D_m$, such that for all $p > N_m$ and for all $y \in A_{j,p}$, we have*

$$\left| \sum_{\substack{\overline{x} \in \overline{y} + (p\mathbb{Z}/p^m\mathbb{Z})^n, \\ \text{ord}_p(f_j(x)) = m-1}} p^{-mn} \exp\left(\frac{2\pi i f_j(x)}{p^m}\right) \right| \leq D_m p^{-m\sigma_{y,p}}.$$

To prove this lemma, we adapt the proof of Lemma 3.3.12 as follows. We replace the formulas $\phi$ and $\overline{\phi}$ by

$$\phi_j(x_1, \ldots, x_n, z, \xi_1, \ldots, \xi_n, m) = \\ \bigwedge_{i=1}^n (\overline{x}_i = \xi_i) \wedge (\text{ord}(z - z_j) = m - 1) \wedge (\text{ord}(z - f(x_1, \ldots, x_n)) \geq m),$$

$$\overline{\phi}_j(x_1, \ldots, x_n, \xi, \xi_1, \ldots, \xi_n, m) = \\ \bigwedge_{i=1}^n (\overline{x}_i = \xi_i) \wedge (\text{ord}(f(x_1, \ldots, x_n) - z_j) = m - 1) \wedge \overline{\text{ac}}(f(x_1, \ldots, x_n) - z_j) = \xi),$$

where $x_i, z$ are in the valued field-sort, $\xi_i, \xi$ are in the residue field-sort and $m$ is in the value group-sort. These are also $\mathcal{L}_\mathbb{Z} \cup \{z_j\}$-formulas. Most of the other modifications in the proof of Lemma 3.3.12 are the same as we discussed above for Lemma 3.6.1.

The only moment that we have to be more careful, is when estimating $\#\{\overline{x}^{(m)} \in \overline{y}^{(m)} + (p\mathbb{Z}_p/p^m\mathbb{Z}_p)^n \mid \text{ord}_p(f_j(x)) = m - 1\}$. From Section 3.6.1 we know that if $p > N'$ and if $y \in A_{j,p}$, then there exists $y' \in y + (p\mathbb{Z}_p)^n$, such that $a_{y,p}(f) = c_{y'}(f_j)$. By Corollary 3.6 from [Mus02], we have

$$a_{y,p}(f) = c_{y'}(f_j) \leq \frac{(m-1)n - \dim_{\mathbb{F}_p}(\tilde{A}_{p,m,y})}{m - 1},$$

where $A_{p,m,y} := \{\overline{x}^{(m)} \in \overline{y}^{(m)} + (p\mathbb{Z}_p/p^m\mathbb{Z}_p)^n \mid \text{ord}_p(f_j(x)) = m - 1\}$, viewed as a subvariety of $\mathbb{F}_p^{mn}$, and where $\tilde{A}_{p,m,y}$ is the image of $A_{p,m,y}$ under the projection $\pi_m : (\mathbb{Z}_p/p^m\mathbb{Z}_p)^n \to (\mathbb{Z}_p/p^{m-1}\mathbb{Z}_p)^n$, viewed as a subvariety of $\mathbb{F}_p^{mn-n}$. Then $\#A_{p,m,y} \leq \#\tilde{A}_{p,m,y} \cdot p^n$. By the Lang-Weil estimate, there exists a constant $D'_{m,y}$, not depending on $p$, such that

$$\#\tilde{A}_{p,m,y} = D'_{m,y} p^{\dim_{\mathbb{F}_p}(\tilde{A}_{p,m,y})} + O(p^{\dim_{\mathbb{F}_p}(\tilde{A}_{p,m,y}) - \frac{1}{2}}).$$

By looking at the arcspace of $Z(f_j)$, we can see that, for each $m$, there are finitely many schemes $Z_1^{(m)}, \ldots, Z_{k_m}^{(m)}$, such that for all $p$ and $y$, $\tilde{A}_{p,m,y} \cong Z_i^{(m)}(\mathbb{F}_p)$ for some $i \in \{1, \ldots, k_m\}$. This means that the constant $D'_{m,y}$, which we know already to be independent of $p$, only depends on the set of schemes $\{Z_1^{(m)}, \ldots, Z_{k_m}^{(m)}\}$. Hence there exists a constant $D'_{m,j}$, such that $D'_{m,j} \geq D'_{m,y}$ for all $y \in A_{j,p}$. By going over the rest of the proof of Lemma 3.3.12, we can see that the natural number $N_m$ and the constant $D_m$, that are obtained in the proof, only depend on $m$ and $j$.

We need to make similar adjustments in the proof of Lemma 3.3.13, to obtain the following lemma.

**Lemma 3.6.3.** *Let $f \in \mathbb{Z}[x_1, \ldots, x_n]$ be a nonconstant polynomial and let $z_j \in V_f$ be a critical value of $f$. There exists, for each integer $m > 1$, a natural number $N_m > N'$ and a positive constant $D_m$, such that for all $p > N_m$ and for all $y \in A_{j,p}$, we have*

$$\Big| \sum_{\substack{\overline{x} \in \overline{y} + (p\mathbb{Z}/p^m\mathbb{Z})^n, \\ \mathrm{ord}_p(f_j(x)) \geq m}} p^{-mn} \exp\left(\frac{2\pi i f_j(x)}{p^m}\right) \Big| \leq D_m p^{-m\sigma_{y,p}}.$$

The final step after these three lemmas, is to modify the proof of the Main Theorem 3.1.2 at the end of Section 3.3. According to Corollary 3.2.6 and its proof, there exist natural numbers $s_j, M_j, N''_j$, such that for all $p > N''_j$, $m > M_j$ and $y \in A_{j,p}$, we have

$$E^y_{m,p}(f_j) = \sum_{i=1}^{s_j} a_{i,p,y} m^{\beta_{ij}} p^{-\lambda_{ij} m} \mathbf{1}_{A_{ij}}(m). \tag{3.6.2.1}$$

We can easily see that $\beta_{ij}$, $\lambda_{ij}$ and $A_{ij}$ only depend on $f_j$ and not on $y$. By going through the proof of Claim 3.3.14 we obtain a constant $C_0$ and natural numbers $\tilde{M}, \tilde{N}$ (that depend on $\beta_{ij}$, $\lambda_{ij}$ and $A_{ij}$, but not on $a_{i,p,y}$), such that for all $m > \tilde{M}$, $p > \tilde{N}$, $y \in A_{j,p}$ and $1 \leq i \leq s_j$, we have

$$|a_{i,p,y} p^{-\lambda_{ij} m}| \leq C_0 p^{-\sigma_{y,p} m}.$$

Now 3.6.0.3 follows easily.

# Chapter 4

# Conjectures on exponential sums and conjectures on numerical data

*This chapter is joint work with Raf Cluckers, see [CNc].*

## Contents

## Abstract

We relate questions in the spirit of Igusa's question on finite exponential sums to relations among numerical data associated to an embedded resolution from Hironaka's construction.

## 4.1 Introduction

On the one hand we study finite exponential sums like

$$S_{f,a} := \sum_{x \in (\mathbb{Z}/a\mathbb{Z})^n} \exp(\frac{2\pi i f(x)}{a})$$

where $f$ is a non-constant polynomial over $\mathbb{Z}$ in $n$ variables, and $a > 0$ any integer. On the other hand, we study relations among numerical data associated to an embedded resolution with normal crossings of $\mathbb{A}_{\mathbb{C}}^n$ of $f = b$ for complex $b$, using [Hir64]. We formulate conjectures about these sums and numerical data, and prove that they are equivalent. The conjectures on exponential sums originate in work by Igusa [Igu78] and were varied upon by Denef, Sperber [DS01] and by Veys and the first author [CV16]. The conjectures that we put forward on the numerical data are new. We also give local variants, and show that they are all equivalent. By the Chinese remainder theorem, in order to estimate $|S_{f,a}|$ in terms of $a$, it is enough to consider integers $a$ which are positive powers of prime numbers $p$, and then one can rewrite $S_{f,a}$ as a $p$-adic integral; these $p$-adic integrals are the subject of a conjecture we now recall.

### 4.1.1 A conjecture on sums

We recall a variant by Veys and the first author from [CV16] of Igusa's conjecture on exponential sums. Let $\mathcal{O}$ be a ring of integers. Consider the (non-archimedean) integrals

$$E_{f,K,\psi} := \int_{x \in \mathcal{O}_K^n} \psi(f(x))|dx| \qquad (4.1.1.1)$$

where $f$ is a non-constant polynomial in $n$ variables with coefficients in $\mathcal{O}$, $K$ is a local field [1] over $\mathcal{O}$, $\mathcal{O}_K$ the ring of integers of $K$, $\psi : K \to \mathbb{C}^\times$ a nontrivial additive character [2] on $K$, and $|dx|$ the Haar measure on $K^n$ normalized so that $\mathcal{O}_K^n$ has measure 1. Write $q_K$ for the number of elements in the residue field $k_K$ of $K$.

**Definition 4.1.1.** For a nontrivial additive character $\psi$ on $K$, write $m_\psi$ for the unique integer $m$ such that $\psi$ is trivial on $\varpi_K^{-m}\mathcal{O}_K$ and nontrivial on $\varpi_K^{-m-1}\mathcal{O}_K$, where $\varpi_K$ is a uniformizer of $\mathcal{O}_K$.

Note that when $K = \mathbb{Q}_p$, then for each integer $m > 0$ there is a nontrivial additive character $\psi$ on $K$ with $m = m_\psi$ and with $S_{f,p^m} = E_{f,K,\psi}$.

**Definition 4.1.2.** Let $\sigma(f)$ be the minimum over all $b \in \mathbb{C}$ of the log canonical threshold of $f - b$ (see Section 4.1.2).

**Conjecture S.1** (Sum conjecture [CV16])**.** *Given a non-constant polynomial $f$ with coefficients in $\mathcal{O}$, there exist $M > 0$ and $M' > 0$ such that*

$$|E_{f,K,\psi}| < M m_\psi^{n-1} q_K^{-\sigma(f)m_\psi}$$

*for all local fields $K$ over $\mathcal{O}$ whose residue field characteristic is at least $M'$ and for all nontrivial additive characters $\psi$ on $K$ satisfying $m_\psi \geq 2$.*

The main point of the conjecture goes back to Igusa [Igu78, first page of the introduction] and is the independence of $M$ on $K$. This independence is useful for reasons related to adèlic integrability, see [Igu78] .

---

1. By a local field $K$ over a commutative ring $R$ with unit we mean a finite field extension of $\mathbb{Q}_p$ or of $\mathbb{F}_p((t))$ for some prime number $p$ such that moreover there is a unit-preserving homomorphism of rings $R \to K$.

2. By an additive character is meant a continuous group homomorphism from the additive group on $K$ to $\mathbb{C}^\times$.

In this paper we formulate a conjecture purely on numerical data associated to Hironaka's embedded resolutions (see Conjecture N.1), and we show the equivalence with Conjecture S.1.

Denef and Sperber formulated a local variant of Igusa's conjecture; we generalize these variants to Conjecture S.2.

### 4.1.2 A numerical conjecture

Let $f(x)$ be a nonconstant polynomial in $\mathbb{C}[x]$ for a tuple of variables $x = (x_1, \ldots, x_n)$ for some $n$. Write $D(f)$ for $\text{Spec}(\mathbb{C}[x]/(f(x)))$, and $Z(f)$ for the reduced scheme associated to $D(f)$. Let $h : Y \to \mathbb{A}_K^n$ be an embedded resolution with normal crossings for $D(f)$, obtained by successive blowing-ups at irreducible nonsingular schemes, as follows from [Hir64, page 142, Main Theorem II]. In particular, we have:
— $Y$ is a closed nonsingular subscheme of $\mathbb{P}_{\mathbb{A}_{\mathbb{C}}^n}^k$ for some $k \geq 0$.
— $h$ is a proper birational morphism which is an isomorphism outside $h^{-1}(\text{Sing}(Z(f)))$, where $\text{Sing}(Z(f))$ denotes the singular locus of $Z(f)$.
— the divisor of $h \circ f$ on $Y$ equals $\sum_{j \in J} N_j E_j$ for a finite set $J$ and some positive integers $N_j$, where each $E_j$ is an irreducible component of the reduced scheme of $h^{-1}(D(f))$.
— $\bigcup_{j \in J} E_j$ has only normal crossings as subscheme of $Y$, namely, for any closed point $a \in Y$, there is an affine neighborhood $V$ of $a$ so that

$$f \circ h = u \prod_{i \in I} y_i^{N_i} \tag{4.1.2.1}$$

holds in $\mathcal{O}(V)$, $I \subset J$ is such that $i \in I$ if and only if $a \in E_i$, $(y_i)_{i \in I}$ forms a regular sequence of parameters for the stalk $\mathcal{O}_{Y,a}$ at $a$, and where $y_i$ is a representation of the divisor $E_i$ in $V$ for each $i \in I$ and $u$ a unit in $\mathcal{O}(V)$.
— the divisor of $h^*(dx_1 \wedge \ldots \wedge dx_n)$ equals $\sum_{j \in J} (\nu_j - 1) E_j$ for some positive integers $\nu_j$.

For any subset $I \subset J$, define

$$E_I := \bigcap_{i \in I} E_i$$

if $I$ is nonempty and put $E_I = Y$ if $I$ is the empty set.

For each $j \in J$, put

$$\sigma_j := \nu_j/N_j.$$

For $V$ a closed subvariety of $\mathbb{A}_{\mathbb{C}}^n$ write

$$\text{lct}_V(f)$$

for the log canonical threshold of $f$ along $V$, namely, the minimum of the values $\sigma_j$ over $j \in J$ with $h(E_j) \cap V \neq \emptyset$, and where the minimum over the empty set is taken to be 1. If $V$ is a complex point $P$, we write $c_P(f) := \text{lct}_P(f)$ for the log canonical threshold of $f$ along $P$. Write $\text{lct}(f)$ for $\text{lct}_{\mathbb{A}_{\mathbb{C}}^n}(f)$, which is called the log canonical threshold of $f$. Finally write

$$\sigma(f)$$

for the minimum over all $b \in \mathbb{C}$ of the values $\mathrm{lct}(f - b)$, namely, the log-canonical threshold of the polynomial $f - b$. Obviously, if $P$ is a complex point and $V$ is a closed subvariety of $\mathbb{A}^n_\mathbb{C}$ containing $P$, then one has

$$1 \geq c_P(f) \geq \mathrm{lct}_V(f) \geq \mathrm{lct}(f) \geq \sigma(f).$$

**Definition 4.1.3** (Power condition for $(f, h)$). Let $f$ be in $\mathbb{C}[x]$ and fix an embedded resolution $h: Y \to \mathbb{A}^n_\mathbb{C}$ for $Z(f)$ with notation as above. Choose $I \subset J$ such that $E_I$ is nonempty, let $W$ be an open subvariety of $E_I$, $g$ be in $\mathcal{O}_W(W)$, and let $d > 1$ be an integer. Say that $(f, h)$ satisfies the power condition, witnessed by $(I, W, g, d)$, if $d | N_i$ for all $i \in I$ and there is an open subvariety $V \subset Y$ such that we can write $f \circ h = u \prod_{i \in I} y_i^{N_i}$ as in (4.1.2.1) on $V$, such that $W = E_I \cap V$, and such that

$$u|_W = g^d. \tag{4.1.2.2}$$

We write shortly that the power condition holds for $(f, h)$ if there exists $(I, W, g, d)$ witnessing the power condition for $(f, h)$. (As usual, $\mathcal{O}_W$ stands for the structure sheaf of $W$.)

We introduce the following conjecture, purely about numerical data.

**Conjecture N.1** (Numerical conjecture). *Let $f$ be in $\mathbb{C}[x]$ and $h: Y \to \mathbb{A}^n_\mathbb{C}$ be an embedded resolution for $Z(f)$. Suppose that $(f, h)$ satisfies the power condition, say, witnessed by $(I, W, g, d)$. Then the following inequality holds*

$$\sigma(f) - \frac{1}{2} \leq \sum_{i \in I} N_i(\sigma_i - \sigma(f)). \tag{4.1.2.3}$$

### 4.1.3  The general conjecture

Let $f$ be a non-constant polynomial in $\mathcal{O}[x]$ and use notation from section 4.1.1. Let $Z$ be a closed subvariety of $\mathbb{A}^n_\mathcal{O}$ such that $f$ vanishes on $Z(\mathbb{C})$. If $K$ is a local field over $\mathcal{O}$, we consider the Schwartz-Bruhat function

$$\Phi_{Z,K} = \mathbf{1}_{\{x \in \mathcal{O}^n_K | \sup_{z \in Z(K)} \mathrm{ord}(x-z) > 0\}} \tag{4.1.3.1}$$

and put

$$E^Z_{f,K,\psi} := \int_{x \in \mathcal{O}^n_K} \Phi_{Z,K}(x) \psi(f(x)) |dx|.$$

We generalize conjecture S.1 on exponential sums as follows.

**Conjecture S.2** (General sum conjecture). *Let $f$ be a nonconstant polynomial in $n$ variables and with coefficients in $\mathcal{O}$ and let $Z$ be a closed subvariety of $\mathbb{A}^n_\mathcal{O}$. Suppose that $f$ vanishes on $Z(\mathbb{C})$. Then there exist $M > 0$ and $M' > 0$ such that*

$$|E^Z_{f,K,\psi}| < M m_\psi^{n-1} q_K^{-\mathrm{lct}_Z(f) m_\psi}$$

*for all local fields $K$ over $\mathcal{O}$ whose residue field characteristic is at least $M'$ and for all nontrivial additive characters $\psi$ on $K$ satisfying $m_\psi \geq 2$.*

Similarly, we formulate a natural variant of conjecture N.1.

**Conjecture N.2** (General numerical conjecture). *Let $f$ be a non-constant polynomial in $\mathbb{C}[x]$ and let $Z$ be a closed subvariety of $\mathbb{A}_{\mathbb{C}}^n$ such that $f$ vanishes on $Z$. Suppose that the power condition holds for $(f,h)$, witnessed by some $(I,W,g,d)$, such that moreover*
$$h(W) \subset Z.$$
*Then the following inequality holds*
$$\operatorname{lct}_Z(f) - \frac{1}{2} \leq \sum_{i \in I} N_i(\sigma_i - \operatorname{lct}_Z(f)). \tag{4.1.3.2}$$

The variant of Conjecture 2 with $Z = \{0\}$ is known is the Denef-Sperber conjecture, from [DS01]. The exponent $\sigma(f)$ or $\operatorname{lct}_Z(f)$ in the sum conjectures S.1, S.2 is not always optimal (namely, when they equal 1 they may not be optimal). As optimal exponent, one can use the motivic oscillation index, or, the maximum of the real parts of the non-trivial poles of Igusa local zeta functions, similar as in [Clu08b].

### 4.1.4 Main result

We are ready to state the main result of this chapter.

**Theorem 4.1.4.** *Let $f$ be a non-constant polynomial in $\mathcal{O}[x]$ with $x = (x_1, \ldots, x_n)$.*

1. *Conjecture N.1 holds for $f$ and all choices of $h$, if and only if Conjecture N.1 holds for $f$ and any specific choice of $h$, if and only if Conjecture S.1 holds for $f$.*
2. *Conjecture S.2 holds for $f$ and a closed subvariety $Z$ of $\mathbb{A}_{\mathcal{O}}^n$ if and only if Conjecture N.2 holds for $f$ and $Z$ (and all, resp. any, choice of $h$).*

Known special cases are the cases of Conjectures S.1 and S.2 with $Z = 0$ in the case that moreover
$$\sigma(f) \leq 1/2$$
(resp. $\operatorname{lct}_0(f) \leq 1/2$), see [CNb]. Also the non-degenerate cases [DS01, Clu08a, Clu10, CNa], and the homogeneous cases in at most 3 variables [Lic13, Lic16, Wri] are known, as well as the nonsingular homogeneous case [Igu78], but will not be used in this paper.

**Remark 4.1.5.** On top of (4.1.3.2), it may be interesting to study when the inequality
$$\operatorname{lct}_Z(f) - \frac{1}{d} \leq \sum_{i \in I} N_i(\sigma_i - \operatorname{lct}_Z(f))$$
holds. The numerical conjectures may relate to Veys' conjecture on poles of the Igusa local Zeta function with maximal order, and the proof techniques from Nicaise and Xu in [NX16].

## 4.2 Embedded resolutions of singularities

We recall some properties and terminology for embedded resolutions in the case of affine hypersurfaces over a subfield $K$ of $\mathbb{C}$. Let $f(x)$ be a nonconstant polynomial in $K[x]$ for a tuple of variables $x = (x_1, \ldots, x_n)$ for some $n$, where $K$ is a subfield of $\mathbb{C}$. Write $D(f)_K$ for $\mathrm{Spec}(K[x]/(f(x)))$, $D(f)$ for $\mathrm{Spec}(\mathbb{C}[x]/(f(x)))$, $Z(f)_K$ for the reduced scheme associated to $D(f)_K$, and $Z(f)$ for the reduced scheme associated to $D(f)$. Let $h: Y \to \mathbb{A}_K^n$ be an embedded resolution with normal crossings for $D(f)_K$, obtained by successive blowing-ups at irreducible non-singular schemes, as follows from [Hir64, page 142, Main Theorem II]. In particular, we have:

— $Y$ is a closed nonsingular subscheme of $\mathbb{P}_{\mathbb{A}_K^n}^k$ for some $k \geq 0$.
— $h$ is a proper birational morphism which is an isomorphism outside $h^{-1}(S)$, where $S$ is the singular locus of $Z(f)_K$.
— the divisor of $h \circ f$ on $Y$ equals $\sum_{j \in T_K} N_j E_j$ for a finite set $T_K$ and some positive integers $N_j$, where each $E_j$ is an irreducible component of the reduced scheme of $h^{-1}(D(f)_K)$.
— $\bigcup_{j \in T_K} E_j$ has only normal crossings as subscheme of $Y$, namely, for any closed point $a \in Y$, there is an affine neighborhood $V$ of $a$ so that

$$f \circ h = u \prod_{i \in I} y_i^{N_i}$$

holds in $\mathcal{O}(V)$, $I \subset T_K$ is such that $i \in I$ if and only if $a \in E_i$, $(y_i)_{i \in I}$ forms a regular sequence of parameters in the stalk $\mathcal{O}_{Y,a}$ at $a$, and where $y_i$ is a representation of the divisor $E_i$ in $V$ for each $i \in I$ and $u$ a unit in $\mathcal{O}(V)$.
— the divisor of $h^*(dx_1 \wedge \ldots \wedge dx_n)$ equals $\sum_{j \in T_K} (\nu_j - 1) E_j$ for some positive integers $\nu_j$.

For any subset $I \subset T_K$, define $E_I := \bigcap_{i \in I} E_i$ if $I$ is nonempty and put $E_I = Y$ if $I$ is the empty set. Further, write

$$E_I^\circ := E_I \setminus (\bigcup_{i \notin I} E_i).$$

For each $j \in T_K$, put $\sigma_j = \nu_j/N_j$.

By the functoriality of embedded resolutions for extensions of the base field, $h$ will induce an embedded resolution $h: Y \otimes K' \to \mathbb{A}_{K'}^n$ of $\mathrm{Spec}(K'[x]/(f(x)))$ for any field $K'$ containing $K$. We remark that each blowing-up center $C$ of $h$ can be split into finitely many irreducible components $C_i$ over $K'$ and we can replace the blowing-up with center $C$ by a composition of blowing-ups at the $C_i$. Similarly, an irreducible component $E_i$ can be split into a union of finitely many irreducible components $E_{ij}$ over $K'$, with $(i, j) \in T_{K'}$ for a corresponding finite set $T_{K'}$, and where one always has $N_i = N_{ij}$, $\nu_i = \nu_{ij}$, $\sigma_i = \sigma_{ij}$. When $K = \mathbb{C}$, we write $J$ or $J(h)$ for $T_K$.

## 4.3 Igusa's local zeta function and exponential sums

Let $K$ be a number field, $\mathcal{O}$ the ring of algebraic integers of $K$ and $\mathfrak{p}$ any maximal ideal of $\mathcal{O}$. We denote the completions of $K$ and $\mathcal{O}$ with respect to $\mathfrak{p}$ by $K_\mathfrak{p}$ and $\mathcal{O}_\mathfrak{p}$. Let $q = p^m$ be the cardinality of the residue field $k_\mathfrak{p}$ of the local ring $\mathcal{O}_\mathfrak{p}$, thus one has $k_\mathfrak{p} = \mathbb{F}_q$. Fix a uniformizer $\varpi$ of $\mathcal{O}_\mathfrak{p}$. For $x \in K_\mathfrak{p}$, denote by $\text{ord}(x) \in \mathbb{Z} \cup \{+\infty\}$ the $\mathfrak{p}$-valuation of $x$, write $|x| = q^{-\text{ord}(x)}$ and write $\text{ac}(x) = x\varpi^{-\text{ord}(x)}$. Let $\chi : \mathcal{O}_\mathfrak{p}^\times \to \mathbb{C}^\times$ be a multiplicative character, namely, a continuous group homomorphism on the group of units $\mathcal{O}_\mathfrak{p}^\times$ of $\mathcal{O}_\mathfrak{p}$. Automatically $\chi$ has finite image. By the *order* of such a character we mean the number of elements in its image. The *conductor* $c(\chi)$ of the character is the smallest $c \geq 1$ for which $\chi$ is trivial on $1 + \mathfrak{p}^c$, where we also write $\mathfrak{p}$ for the maximal ideal of $\mathcal{O}_\mathfrak{p}$. We formally put $\chi(0) = 0$. Let $f(x) \in \mathcal{O}[x]$ be a non-constant polynomial in $n$ variables $x = (x_1, \ldots, x_n)$, and let $\Phi : K_\mathfrak{p}^n \to \mathbb{C}$ be a *Schwartz-Bruhat function*, i.e., a locally constant function with compact support. We say that $\Phi$ is *residual* if $\text{Supp}(\Phi) \subset \mathcal{O}_\mathfrak{p}^n$ and $\Phi(x)$ only depends on $x$ mod $\mathfrak{p}$. If $\Phi$ is residual, it induces a function $\overline{\Phi} : k_\mathfrak{p}^n \to \mathbb{C}$. Now we associate to these data the (twisted) *Igusa's local zeta function*

$$Z_\Phi(K_\mathfrak{p}, \chi, s, f) := \int_{K_\mathfrak{p}^n} \Phi(x) \chi\big(\text{ac}(f(x))\big) |f(x)|^s |dx|. \tag{4.3.0.1}$$

In [Igu78], Igusa showed that $Z_\Phi(K_\mathfrak{p}, \chi, s, f)$ is a rational function in $t = q^{-s}$, starting the study of a now vast subject.

Let $p$ be a prime. For each $x \in \mathbb{Q}_p \setminus \mathbb{Z}_p$, there exist unique co-prime natural numbers $m$ and $n$ such that $1 \leq m \leq p^n - 1$ and

$$x - \frac{m}{p^n} \in \mathbb{Z}_p.$$

Write

$$\{x\}_p = \frac{m}{p^n}$$

and call it the *p-adic fractional part* of $x$. For $x \in \mathbb{Z}_p$ put $\{x\}_p = 0$. It is easy to see that $\{x+y\}_p - \{x\}_p - \{y\}_p \in \mathbb{Z}$ so the map

$$\exp : \mathbb{Q}_p \to S^1 : x \mapsto \exp^{2\pi i \{x\}_p}$$

is a group homomorphism which is moreover continuous. We call this map the standard additive character of $\mathbb{Q}_p$.

By the standard additive character $\Psi$ on $K_\mathfrak{p}$ we mean the map sending $z \in K_\mathfrak{p}$ to

$$\Psi(z) := \exp(2\pi i \, \text{Tr}_{K_\mathfrak{p}/\mathbb{Q}_p}(z)),$$

where $\text{Tr}_{K_\mathfrak{p}/\mathbb{Q}_p}$ denotes the trace map. We set

$$E_\Phi(f, K_\mathfrak{p}, \psi_z) := \int_{K_\mathfrak{p}^n} \Phi(x) \Psi(zf(x)) |dx|, \tag{4.3.0.2}$$

where $\psi_z$ stands for the additive character on $K_\mathfrak{p}$ sending $t$ to $\Psi(zt)$. (Note that $E_\Phi(f, K_\mathfrak{p}, \psi_z)$ is a finite exponential sum.)

Whenever $\Phi = \mathbf{1}_{\mathcal{O}_\mathfrak{p}^n}$ or $\Phi = \mathbf{1}_{(\mathfrak{p}\mathcal{O}_\mathfrak{p})^n}$, where $\mathbf{1}$ stands for the characteristic function, we will simply denote this function by

$$E_{f,K_\mathfrak{p},\psi_z} \text{ resp. } E^{(0)}_{f,K_\mathfrak{p},\psi_z}.$$

Note that for any $\psi$ as in conjecture S.1 there exists unique $z \in K_\mathfrak{p}^\times$ such that $\psi = \psi_z$. Moreover, one then has $m_\psi = -\operatorname{ord} z$. We will write $m_z$ for $-\operatorname{ord} z$.

We recall the following proposition from [Den], relating exponential sums to Igusa's local zeta functions.

**Proposition 4.3.1** ([Den], Proposition 1.4.4)**.** *Consider $u \in \mathcal{O}_\mathfrak{p}^\times$ and $m \in \mathbb{Z}$. Then $E_\Phi(f, K_\mathfrak{p}, \psi_{u\varpi^{-m}})$ is equal to*

$$Z_\Phi(K_\mathfrak{p}, \chi_{\mathrm{triv}}, 0, f) + \operatorname{Coeff}_{t^{m-1}}\left(\frac{(t-q)Z_\Phi(K_\mathfrak{p}, \chi_{\mathrm{triv}}, s, f)}{(q-1)(1-t)}\right)$$
$$+ \sum_{\chi \neq \chi_{\mathrm{triv}}} g_{\chi^{-1}} \chi(u) \operatorname{Coeff}_{t^{m-c(\chi)}}\left(Z_\Phi(K_\mathfrak{p}, \chi, s, f)\right),$$

*where $g_\chi$ is the Gaussian sum*

$$g_\chi = \frac{q^{1-c(\chi)}}{q-1} \sum_{\overline{v} \in (\mathcal{O}_\mathfrak{p}/\mathfrak{p}^{c(\chi)})^\times} \chi(v)\Psi(v/\varpi^{c(\chi)})$$

*and where $\chi_{\mathrm{triv}}$ is the trivial character.*

Now we will describe Denef's formula for Igusa's local zeta function from [Den], based on resolution of singularities.

Let $\mathcal{O}$, $K$, and $f$ be as in the beginning of this section. We use the notation of section 4.2. In particular, $X = \operatorname{Spec} K[x]$ and $D(f)_K = \operatorname{Spec}(K[x]/(f))$ and there is an embedded resolution $(Y, h)$ with respect to $D(f)_K$ over $K$. Denote the critical locus of $f$ by $C_f$, namely the locus given by $\operatorname{grad} f = 0$. If $Z$ is a closed subscheme of $Y$, for any any local field $L$ over $\mathcal{O}$ of sufficiently large residue field characteristic we denote the reduction modulo $\mathcal{M}_L$ of $Z$ by $\overline{Z}$ (see [Den87, Definition 2.2] and [Shi55]). We say that the resolution $(Y, h)$ of $f$ has *good reduction modulo* $\mathcal{M}_L$ if $\overline{Y}$ and all $\overline{E}_i$ are smooth, $\cup_{i \in T} \overline{E}_i$ has only normal crossings, and the schemes $\overline{E}_i$ and $\overline{E}_j$ have no common components whenever $i \neq j$. There exists a finite subset $S$ of $\operatorname{Spec} \mathcal{O}$, such that for all $\mathfrak{p} \notin S$, $f \not\equiv 0 \mod \mathfrak{p}\mathcal{O}_\mathfrak{p}[x]$ and the resolution $(Y, h)$ for $f$ has good reduction mod $\mathfrak{p}$ (see [Den87], Theorem 2.4). Hence, there exists $M'$ such that for all primes $p > M'$ and all finite field $\mathbb{F}_q$, where $q = p^r$, which is an algebra over $\mathcal{O}$, $h$ induces an embedded resolution $\overline{h}$ of $\operatorname{Spec}(\mathbb{F}_q[x]/(\overline{f}))$, where $\overline{f}$ is the image of $f$ via the structure map $\mathcal{O}[x] \to \mathbb{F}_q[x]$. For all $p > M'$, and assuming that $M'$ is large, $h$ also induces an embedded resolution $\overline{h}$ of $\operatorname{Spec}(\mathbb{F}_p^a[x]/(\overline{f}))$, where $\mathbb{F}_p^a$ is the algebraic closure of $\mathbb{F}_p$ and there is a bijection between $J$ and the set of irreducible components of $\overline{h}^{-1}(Z(\overline{f}))$, where $Z(\overline{f})$ is the reduced scheme associated to $\operatorname{Spec}(\mathbb{F}_p^a[x]/(\overline{f}))$.

For $\mathfrak{p} \notin S$ and $I \subset T_K$ it is easy to prove that $\overline{E}_I = \cap_{i \in I} \overline{E}_i$. We put

$$\overset{\circ}{\overline{E}}_I := \overline{E}_I \setminus \cup_{j \notin I} \overline{E}_j.$$

Let $a$ be a closed point of $\overline{Y}$ and $T_{a,K} = \{i \in T_K \mid a \in \overline{E}_i\}$. In the local ring of $\overline{Y}$ at $a$ we can write

$$\overline{f} \circ \overline{h} = \overline{u} \prod_{i \in T_{a,K}} \overline{g}_i^{N_i},$$

where $\overline{u}$ is a unit, $(\overline{g}_i)_{i \in T_{a,K}}$ is a part of a regular system of parameters and $N_i$ is as above.

In two cases, depending on the conductor $c(\chi)$ of the character $\chi$, we will give a more explicit description of Igusa's zeta function $Z_\Phi(\mathfrak{p}, \chi, s, f)$. In the first case we consider a character $\chi$ on $\mathcal{O}_\mathfrak{p}^\times$ of order $d$ which is trivial on $1 + \mathfrak{p}\mathcal{O}_\mathfrak{p}$, i.e., $c(\chi) = 1$. Then $\chi$ induces a character (denoted also by $\chi$) on $k_\mathfrak{p}^\times$. We define a map

$$\Omega_\chi : \overline{Y}(k_\mathfrak{p}) \to \mathbb{C}$$

as follows. Let $a$ be in $\overline{Y}(k_\mathfrak{p})$. If $d | N_i$ for all $i \in T_a$, then we put $\Omega_\chi(a) = \chi(\overline{u}(a))$, otherwise we put $\Omega_\chi(a) = 0$. This definition is independent of the choice of the $\overline{g}_i$. In the following theorem we recall the formula of Igusa's local zeta function.

**Theorem 4.3.2** ([Den91], Theorem 2.2 or [Den], Theorem 3.4). *Let $\chi$ be a character on $\mathcal{O}_\mathfrak{p}^\times$ of order $d \geq 1$ which is trivial on $1 + \mathfrak{p}\mathcal{O}_\mathfrak{p}$. Suppose that $\mathfrak{p} \notin S$ and that $\Phi$ is residual. Then we have*

$$Z_\Phi(K_\mathfrak{p}, \chi, s, f) = q^{-n} \sum_{\substack{I \subset T_K, \\ \forall i \in I: d | N_i}} c_{I,\Phi,\chi} \prod_{i \in I} \frac{(q-1)q^{-N_i s - \nu_i}}{1 - q^{-N_i s - \nu_i}}, \quad (4.3.0.3)$$

*where, for any $I$, we write*

$$c_{I,\Phi,\chi} = \sum_{a \in \overset{\circ}{\overline{E}}_I(k_\mathfrak{p})} \overline{\Phi}(\overline{h}(a)) \Omega_\chi(a)$$

If $\overline{W}$ is an irreducible open subset of $\overline{E}_I$, we will set $\overset{\circ}{\overline{W}} = \overline{W} \cap \overset{\circ}{\overline{E}}_I$ and

$$c_{I,\overline{W},\Phi,\chi} = \sum_{a \in \overset{\circ}{\overline{W}}(k_\mathfrak{p})} \overline{\Phi}(\overline{h}(a)) \Omega_\chi(a). \quad (4.3.0.4)$$

Because $E_I$ is smooth, any two irreducible components of $\overline{E}_I$ are disjoint. Hence it follows that

$$c_{I,\Phi,\chi} = \sum_{\overline{W}} c_{I,\overline{W},\Phi,\chi} \quad (4.3.0.5)$$

where the sum taken over the set of all irreducible components of $E_I$.

**Remark 4.3.3.** Note that the sum in (4.3.0.3) only involves $I$ such that $d \mid N_i$ for each $i \in I$. Therefore the number of characters $\chi$ for which $c(\chi) = 1$ and $c_{I,\Phi,\chi} \neq 0$ for some $I \subset T_K$ will have an upper bound $\mathcal{N}$ which will only depend on $(Y, h)$.

The next result treats characters of conductor larger than one.

**Theorem 4.3.4** ([Den91], Theorem 2.1 or [Den], Theorem 3.3). *Given $f$ and $h$ as above, we can take $M'$ large enough so that the following hold for any local field $L$ over $\mathcal{O}$ with residue field characteristic at least $M'$ and any character $\chi$ on $\mathcal{O}_L^\times$ which is non-trivial on $1 + \mathcal{M}_L$. Suppose that $\Phi$ is residual and that $C_{\overline{f}}(\mathbb{F}_q) \cap \mathrm{Supp}(\overline{\Phi}) \subset \overline{f}^{-1}(0)$, where $C_{\overline{f}}$ is the locus given by $\mathrm{grad}\,\overline{f} = 0$. Then one has*

$$Z_\Phi(L, \chi, s, f) = 0.$$

Note that there is also another expression of $c_{I,\Phi,\chi}$, based on Grothendieck's trace formula given in detail in [Den91]. We will not use this in our paper.

## 4.4 Proof of theorem 4.1.4

In this section, we will prove theorem 4.1.4. The proof of the second part is similar to the proof of the first part and will be treated more shortly at the end of the section.

We recall some notation from section 4.3. Let $f$ be as in conjecture S.1. Let $K$ be a number field containing all coefficients of $f$, $\mathcal{O}$ be the ring of integers of $K$; by supposition we have $f \in \mathcal{O}[x]$. We fix an embedded resolution $(Y, h)$ of $D(f)_K$ and a number $M'$ such that $h$ has good reduction and induces an embedded resolution of $D(f)_{\mathbb{F}_q}$ whenever $\mathrm{char}(\mathbb{F}_q) > M'$ and $\mathbb{F}_q$ is an algebra over $\mathcal{O}$. Moreover we suppose that for the numerical data $(\nu_i, N_i)_{i \in T_K}$ of $(Y, h)$ we have $N_i < M'$ and that $f$ modulo $\mathcal{M}_L$ is nonzero. From now on we consider local fields $L$ over $\mathcal{O}$ and with residue field characteristic larger than $M'$. We recall Lemma 5.1 from [CNb], see also Lemma 3.5.1 from chapter 3.

**Lemma 4.4.1** ([CNb]). *Let $L$ be a local field as above and $V_{f,L}$ be the set of critical values $z$ of $f$ in $L$. Then $\#(V_{f,L})$ has an upper bound $e$ which does not depend on $L$. Furthermore, up to enlarging $M'$, the following holds:*

1. *for all $z \in V_{f,L}$, we have $\mathrm{ord}_L(z) = 0$;*
2. *for any two distinct points $z_1, z_2$ in $V_{f,L}$, we have $\mathrm{ord}_L(z_1 - z_2) = 0$;*
3. *if $x \in \mathcal{O}_L^n$ is such that $\mathrm{ord}_L(f(x) - z) = 0$ for all $z \in V_{f,L}$, then $x$, resp. $\overline{x}$, is a regular point of $f$, resp. $\overline{f} := (f \mod \mathcal{M}_L)$,*

Write $V_f = \{z_1, \ldots, z_e\}$ for the set of critical values of $f$ over $\mathbb{C}$ (in particular, this value of $e$ can serve for Lemma 4.4.1). Let $L$ be a local field as above. Up to reindexing the elements of $V_f$, we can assume that $V_{f,L} = \{z_1, \ldots, z_r\}$ with $r \leq e$. For each $1 \leq i \leq r$, put

$$\Phi_{i,L} := \mathbf{1}_{\{x \in \mathcal{O}_L^n \mid \mathrm{ord}_L(f(x) - z_i) > 0\}} : L^n \to \mathbb{C}$$

for the characteristic function of $\{x \in \mathcal{O}_L^n \mid \mathrm{ord}_L(f(x) - z_i) > 0\}$. Because $f \in \mathcal{O}_L[x_1, \ldots, x_n]$ and by Lemma 4.4.1 we see that $\Phi_{i,L}$ is residual for all $1 \leq i \leq r$, and that $\text{Supp}(\Phi_{i,L}) \cap \text{Supp}(\Phi_{j,L}) = \emptyset$, if $i \neq j$. Write $\Phi_{0,L} := \mathbf{1}_{\mathcal{O}_L^n} - \sum_{i=1}^{r} \Phi_{i,L}$. Clearly $\Phi_{0,L}$ is also residual. Now we have the following formula

$$E_{f,L,\psi_z} = \int_{\mathcal{O}_L^n} \Psi(zf(x))|dx| \tag{4.4.0.1}$$

$$= \sum_{i=0}^{r} \int_{\mathcal{O}_L^n} \Phi_{i,L}(x)\Psi(zf(x))|dx| \tag{4.4.0.2}$$

$$= \sum_{i=1}^{r} \int_{\mathcal{O}_L^n} \Phi_{i,L}(x)\Psi(z(f(x)-z_i))\Psi(zz_i)|dx| + E_{\Phi_{0,L}}(f,L,\psi_z) \tag{4.4.0.3}$$

$$= \sum_{i=1}^{r} \Psi(zz_i) E_{\Phi_{i,L}}(f_i, L, \psi_z) + E_{\Phi_{0,L}}(f, L, \psi_z), \tag{4.4.0.4}$$

where $f_i(x) = f(x) - z_i$ for $1 \leq i \leq r$. For all local field $L$ as in above, one has

$$E_{\Phi_{0,L}}(f, L, \psi_z) = 0 \text{ if } \text{ord}(z) < -1 \,. \tag{4.4.0.5}$$

Indeed, this follows from Lemma 4.4.1 and from $C_{\overline{f}} \cap \text{Supp}(\overline{\Phi}_{0,L}) = \emptyset$ (see remark 4.5.3, [Den]). The implication from conjecture N.1 for a polynomial $f$ to conjecture S.1 for $f$ relies on the following result.

**Proposition 4.4.2.** *Suppose that $z_1 = 0$ is a critical value of $f$. If conjecture N.1 holds for an embedded resolution $(Y,h)$ of $D(f)_K$ then there exists a positive number $M$ such that*
$$|E_{\Phi_{1,L}}(f,L,\psi_z)| \leq M m_z^{n-1} q_L^{-\sigma(f)m_z}$$
*for all local fields $L$ as above and all $z \in L$ with $m_z = -\text{ord}(z) \geq 2$.*

We recall a variant of the Lang-Weil estimates from [Mus, Prop. 6.6], and show some corollaries.

**Proposition 4.4.3.** *Let $k = \mathbb{F}_q$ be a finite field and $X \hookrightarrow \mathbb{P}_k^n$ be an irreducible closed subvariety of degree $d > 0$ and dimension $r$. Suppose that $X$ splits into $m$ irreducible components over $k^a$, the algebraic closure of $k$. Then there are positive constants $c_X$ and $c'_X$ such that for every $l \geq 1$ we have*

$$|\#X(\mathbb{F}_{q^l}) - mq^{lr}| \leq \frac{(d-m)(d-2m)}{m} q^{l(r-\frac{1}{2})} + c_X q^{l(r-1)}, \text{ if } m|l \text{ and}$$

$$\#X(\mathbb{F}_{q^l}) \leq c'_X q^{l(r-1)} \text{ if } m \nmid l$$

*In particular, if $X$ is smooth then we can take $c'_X = 0$ and $c_X$ only depends on $n, d, r$. Note that bounds from the case $m \nmid l$ even hold when $X$ is a general (not necessarily closed) subvariety of $\mathbb{P}_k^n$.*

**Corollary 4.4.4.** *Let $\overline{W}$ be an irreducible component of $\overline{E}_I$ over $\mathbb{F}_q$. If $\overline{W}$ is not geometrically irreducible then $\overline{W}(\mathbb{F}_q) = \emptyset$. In particular, for all characters $\chi$ of $\mathcal{O}_L^\times$ and all residual Schwartz-Bruhat functions $\Phi$ we have $c_{I,\overline{W},\Phi,\chi} = 0$.*

*Proof.* This claim follows from the smoothness of $\overline{W}$ and by the final part of proposition 4.4.3. □

We also need the following result from [RL05].

**Lemma 4.4.5.** *Let $\mathcal{O}$ be a ring of integers. Let $d > 1$ be an integer, $Z \subset \mathbb{A}_{\mathcal{O}}^n$ be a subscheme such that $Z \otimes \mathbb{C}$ is an irreducible subvariety of $\mathbb{A}_{\mathbb{C}}^n$ of dimension $r$, and $F : Z \to \mathbb{G}_{m,\mathcal{O}}$ be a regular morphism. Suppose that there does not exist a regular morphism $g : Z \otimes \mathbb{C} \to \mathbb{G}_{m,\mathbb{C}}$ such that $g^d$ equals the restriction of $F \otimes \mathbb{C}$. Then there exists constants $c$ and $M'$ such that for all finite fields $\mathbb{F}_q$ of characteristic at least $M'$ such that $\mathbb{F}_q$ is an algebra over $\mathcal{O}$ and all nontrivial characters of $\mathbb{F}_q^\times$ of order $d$, one has*

$$|\sum_{x \in Z(\mathbb{F}_q)} \chi(F(x))| \leq c q^{r-1/2}. \tag{4.4.0.6}$$

*Proof.* See [RL05, Theorem 1.1] □

We fix a uniformizer $\varpi_L$ of $\mathcal{O}_L$ for each local field $L$ over $\mathcal{O}$. For each $z \in L^\times$ with $m_z = -\operatorname{ord} z \geq 2$ we can write $z = \varpi_L^{m_z} u_z$, where $u_z \in \mathcal{O}_L^\times$. We recall another lemma from [CNb].

**Lemma 4.4.6** ([CNb], lemma 4.1). *There exist a positive constant $C$ such that for all local fields $L$ as above and for all $z \in L^\times$ with $m_z \geq 2$ we have*

$$|Z_{\Phi_{1,L}}(L, \chi_{\mathrm{triv}}, 0, f) + \operatorname{Coeff}_{t^{m_z-1}} \frac{(t - q_L) Z_{\Phi_{1,L}}(L, \chi_{\mathrm{triv}}, s, f)}{(q_L - 1)(1 - t)}| \leq C m_z^{n-1} q_L^{-m_z \sigma(f)}.$$

We will now estimate the other terms of Proposition 4.3.1. By Theorems 4.3.2 and 4.3.4 we may focus on non-trivial characters $\chi$ with conductor $c(\chi)$ one.

Let $d > 1$ be an integer and $I \subset T_L$ such that $\overline{E}_I \neq \emptyset$ and $d | N_i$ for all $i \in I$. Let $\chi$ be a character of $\mathcal{O}_L^\times$ of order $d$ and $\mathcal{F}_\chi$ as in section 4.3. Let $\overline{W}$ be an irreducible component of $\overline{E}_I$ such that $\overline{W}(k_L) \neq \emptyset$. By proposition 4.4.3 and since $\overline{E}_I$ is smooth one has that $\overline{W}$ is geometrically irreducible and of dimension $m = n - \#(I)$. As in section 4.3 we write $\overset{\circ}{\overline{W}} = \overline{W} \cap \overset{\circ}{\overline{E}}_I$. We have the following proposition.

**Proposition 4.4.7.** *Choose $I \subset T_L$ such that $\overline{E}_I \neq \emptyset$. Then there exists positive constants $M'$ and $C_I$ depending only on $I$, $h$, and $f$, such that for all local fields $L$ over $\mathcal{O}$ and with residue field characteristic at least $M'$, and with $\overline{W}$ any irreducible component of $\overline{E}_I$, one has*

$$|g_{\chi^{-1}} \operatorname{Coeff}_{t^{m-1}}(q_L^{-n} c_{I,\overline{W},\Phi_{1,L},\chi} \prod_{i \in I} \frac{(q_L - 1) t^{N_i} q_L^{-\nu_i}}{1 - t^{N_i} q_L^{-\nu_i}})| \leq C_I q_L^{-m\sigma(f) + \sigma_I} m^{\#(I)-1} \tag{4.4.0.7}$$

*for all $m \geq 2$, with*

$$\sigma_I = \sigma(f) - \frac{1}{2} - \sum_{i \in I} N_i(\sigma_i - \sigma(f)).$$

*If furthermore the power condition is not satisfied for $\overline{W}$, i.e., there exists no open $W'$ of $W \otimes \mathbb{C}$ such that the induced resolution $(f, h_\mathbb{C})$ satisfies the power condition witnessed by $(I, W', g, d)$ for some $g$ and where $W$ reduces to $\overline{W}$, then we can choose $c_I$ such that*

$$|g_{\chi^{-1}} \operatorname{Coeff}_{t^{m-1}}(q_L^{-n} c_{I,\overline{W},\Phi_{1,L},\chi} \prod_{i \in I} \frac{(q_L - 1) t^{N_i} q_L^{-\nu_i}}{1 - t^{N_i} q_L^{-\nu_i}})| \leq C_I q_L^{-m\sigma(f) + \sigma(f) - 1} m^{\#(I)-1} \tag{4.4.0.8}$$

*for all $m \geq 2$.*

*Proof.* We have

$$\text{Coeff}_{t^{m-1}} \prod_{i \in I} \frac{t^{N_i} q_L^{-\nu_i}}{1 - t^{N_i} q_L^{-\nu_i}} = \sum_{(a_i)_{i \in I} \in A_{I,m}} q_L^{-\sum_{i \in I} \nu_i(a_i+1)},$$

where $A_{I,m} = \{(a_i)_{i \in I} \in \mathbb{N}^{\#(I)} \mid \sum_{i \in I} N_i(a_i + 1) = m - 1\}$. For each $(a_i)_{i \in I} \in A_{I,m}$ we have

$$-\sum_{i \in I} \nu_i(a_i + 1) = -\sum_{i \in I} N_i \sigma_i(a_i + 1)$$
$$= -\sum_{i \in I} N_i(a_i + 1)(\sigma_i - \sigma(f)) - (m-1)\sigma(f)$$
$$\leq -(m-1)\sigma(f),$$

since $\sigma_i \geq \sigma(f)$ for all $i$. Since $\#(A_{I,m}) \leq m^{\#(I)-1}$ we deduce that

$$|\text{Coeff}_{t^{m-1}} \prod_{i \in I} \frac{t^{N_i} q_L^{-\nu_i}}{1 - t^{N_i} q_L^{-\nu_i}}| \leq m^{\#(I)-1} q_L^{-(m-1)\sigma(f)}$$

We additionally observe that for each $(a_i)_{i \in I} \in A_{I,m}$ we have

$$-\sum_{i \in I} \nu_i(a_i + 1) = -\sum_{i \in I} N_i \sigma_i(a_i + 1)$$
$$= -\sum_{i \in I} N_i(a_i + 1)(\sigma_i - \sigma(f)) - (m-1)\sigma(f)$$
$$\leq -m\sigma(f) + \sigma_I,$$

where $\sigma_I = \sigma(f) - \frac{1}{2} - \sum_{i \in I} N_i(\sigma_i - \sigma(f))$. By the Lang-Weil estimates of Proposition 4.4.3, there exists a positive constant $C_I$ depending only on $\overline{W}$ such that

$$|c_{I,\overline{W},\Phi_{1,L},\chi}| \leq \#(\overset{\circ}{\overline{W}}(k_L)) \leq C_I q_L^{\dim(\overline{W})} = C_I q_L^{n-\#(I)}.$$

The absolute value of the Gauss sum $g_{\chi^{-1}}$ is $q_L^{-1/2}$. Combining, we find

$$|g_{\chi^{-1}} \text{Coeff}_{t^{m-1}}(q_L^{-n} c_{I,W,\Phi_{1,L},\chi} \prod_{i \in I} \frac{(q_L - 1)t^{N_i} q_L^{-\nu_i}}{1 - t^{N_i} q_L^{-\nu_i}})| \leq C_I q_L^{-m\sigma(f) + \sigma_I} m^{\#(I)-1}$$

for all $m \geq 2$. This proves the first inequality (4.4.0.7). Using Corollary 4.4.5 instead of the Lang-Weil estimates (which is allowed by the absence of the power condition), we see that there exists a positive constant $C_I$ depending only on $I$ such that

$$|g_{\chi^{-1}} \text{Coeff}_{t^{m-1}}(q_L^{-n} c_{I,\overline{W},\Phi_{1,L},\chi} \prod_{i \in I} \frac{(q_L - 1)t^{N_i} q_L^{-\nu_i}}{1 - t^{N_i} q_L^{-\nu_i}})| \leq C_I q_L^{-m\sigma(f) + \sigma(f) - 1} m^{\#(I)-1},$$

as desired for the inequality from (4.4.0.8). □

Let $(Y, h)$ be an embedded resolution of $D(f)_K$ as in sections 4.2, 4.3. Then $h$ induces an embedded resolution $h$ of $D(f)_{K'}$ for all fields $K' \supset K$ and an embedded resolution $\overline{h}$ of $D(\overline{f})_{\mathbb{F}_q}$ and of $D(\overline{f})_{\mathbb{F}_q^a}$ for all finite fields $\mathbb{F}_q$ which is an algebra over $\mathcal{O}$ and $\mathrm{char}(\mathbb{F}_q) > M'$.

As in section 4.1.2, when $F$ is a field and if $h$ induces an embedded resolution $h_F$ of $D(f)_F$ ( or $D(\overline{f})_F$) over $F$, we will denote by $\{E_i\}_{i \in T_F}$ the set of irreducible components of the reduced scheme of $h_F^{-1}(D(f)_F)$ (resp. $h_F^{-1}(D(\overline{f})_F)$). For all finite fields $\mathbb{F}_q$ with $\mathrm{char}(\mathbb{F}_q) > M'$, one has $T_{\mathbb{C}} = T_{\mathbb{F}_q^a} = J$ and there exists a natural correspondence between irreducible component $E_i$ of the reduced scheme of $h_{\mathbb{C}}^{-1}(D(f)_{\mathbb{C}})$ with its reduction $\overline{E}_i$ over $\mathbb{F}_q^a$ which is an irreducible component of the reduced scheme of $h_{\mathbb{F}_q^a}^{-1}(D(\overline{f})_{\mathbb{F}_q^a})$. Moreover for each $I \subset J$ we also have a natural bijection between

$$\{\text{irreducible components of } E_I\} \longleftrightarrow \{\text{irreducible components of } \overline{E}_I\} \quad (4.4.0.9)$$

given by $W$ corresponding to its reduction $\overline{W}$. (Since $I$ is a subset of $J$, we consider the components over $\mathbb{C}$, or over $\mathbb{F}_q^a$.)

**Lemma 4.4.8.** *We can find a finite extension $K'$ of $K$ with the following properties*

1. *For all $I \subset J$ and $W$ is an irreducible component of $E_I$ then $W$ is already defined over $K'$.*

2. *For any $I \subset J$ and $W$ any irreducible component of $E_I$, we can find finitely many affine open subsets $(V_j)_{j \in J_W}$ of $Y$ such that $V_j$ is defined over $K'$ for all $j$ and $W \subset \cup_{j \in J_W} V_j$. Moreover, we can write $f \circ h|_{V_j} = u_j \prod_{i \in I} y_{ij}^{N_i}$ as in (4.1.2.1) such that $u_j, y_{ij}$ are defined over $K'$ for all $i \in I$ and $j \in J_W$. If moreover $(f, h)$ satisfies the power condition, then there are witnesses defined over $K'$. Namely, the power condition is witnessed by some $(I, W', g, d)$ where $W' = V_j \cap W$ for some $j$ and $g^d = u_j$, and $g$ is defined over $K'$.*

*Proof.* This lemma follows from the algebraically closed case, by taking a field of definition over $K$. □

**Corollary 4.4.9.** *Let $\mathcal{O}'$ be the ring of integers of $K'$ for a well-chosen finite field extension $K'$ of $K$ as in Lemma 4.4.8, and use notation of that lemma. Up to enlarging $M'$, we can suppose that for each finite field $\mathbb{F}_q$ which is an algebra over $\mathcal{O}'$ and with $\mathrm{char}(\mathbb{F}_q) > M'$ we have that the reductions of the data from Lemma 4.4.8 are defined over $\mathbb{F}_q$. Namely,*

1. *For all $I \subset J = T_{\mathbb{F}_q^a}$ and for any irreducible component $W$ of $E_I$, the reduction $\overline{W}$ of $W$, via (4.4.0.9), is defined over $\mathbb{F}_q$*

2. *For each $I \subset J$ and for any irreducible component $W$ of $E_I$, let finitely many affine open subsets $(V_j)_{j \in J_W}$ of $Y$ be given such that $V_j$ is defined over $K'$ for all $j$, such that $W \subset \cup_{j \in J_W} V_j$ and such that $f \circ h|_{V_j} = u_j \prod_{i \in I} y_{ij}^{N_i}$ as in (4.1.2.1) with moreover the $u_j, y_{ij}$ defined over $K'$ for all $i \in I$ and $j \in J_W$. Then the reductions $(\overline{V}_j)_{j \in J_W}, \overline{u}_j,$ and $\overline{y}_{ij}$ are defined over $\mathbb{F}_q$.*

   *Suppose moreover for some $j \in J_W$ with $W'_j = W \cap V_j \neq \emptyset$ there exists an integer $d > 1$ such that $d | N_i$ for all $i \in I$. Let $\overline{W}'_j$ be the reduction of $W'_j$ via (4.4.0.9), suppose that there is a regular function $\overline{g}$ in $\mathcal{O}(\overline{W}'_j)$ such that*

   $$\overline{u}|_{\overline{W}'_j} = \overline{g}_j^d$$

with $\overline{g}_j$ defined over $\mathbb{F}_q^a$. Then there exists a regular function $g_j$ in $\mathcal{O}(W'_j)$ such that
$$u_j|_{W'_j} = g_j^d$$
with $g_j$ defined over $K'$ and $\overline{g}_j$ the reduction of $g_j$.

*Proof of proposition 4.4.2.* Let $L$ be a local field as in proposition 4.4.2. Then there exist finitely many character $\chi$ of $\mathcal{O}_L^\times$ with $c(\chi) = 1$ and with $1 < \text{order}(\chi)|N_i$ for some $i \in T_L$. Denote by $\Pi$ the set of all such characters, thus $\#\Pi \leq \mathcal{N}$ with $\mathcal{N}$ as in remark 4.3.3. Fix $I$ and let $\overline{W}$ be an irreducible component of $\overline{E}_I$.

If $\overline{W}$ is not geometrically irreducible then by Corollary 4.4.6 we have $\overline{W}(k_L) = \emptyset$ and $c_{I,\overline{W},\Phi_{1,L},\chi} = 0$ for all characters $\chi$ of $\mathcal{O}_L^\times$. Suppose now that $\overline{W}$ is geometrically irreducible. Let $\chi \in \Pi$ be such that $1 < d = \text{order}(\chi)$ and such that $d|N_i$ for all $i \in I$. Take $C_I$ as in Proposition 4.4.7, as in (4.4.0.7), resp. in (4.4.0.8), according to whether the power condition is satisfied or not.

**Case 1:** The power condition is not satisfied for $\overline{W}$ (see Proposition 4.4.7).

**Case 2:** The power condition is satisfied and witnessed by some $W'$ which reduces to an open of $\overline{W}$.

If conjecture N.1 holds for $f$ then we have in case 2 that
$$\sigma_I \leq 0.$$

Hence, in both cases 1 and 2 there is a positive constant $C_I$ such that (4.4.0.8) holds, that is, such that
$$|g_{\chi^{-1}}\text{Coeff}_{t^{m-1}}(q_L^{-n} c_{I,\overline{W},\Phi_{1,L},\chi} \prod_{i \in I} \frac{(q_L - 1)t^{N_i} q_L^{-\nu_i}}{1 - t^{N_i} q_L^{-\nu_i}})| \leq C_I q_L^{-m\sigma(f) + \sigma(f) - 1} m^{\#(I) - 1}$$
(4.4.0.10)
for all $m \geq 2$ and uniformly in $L$ as needed. On the other hand, we have $\sigma(f) \leq 1$. So the Proposition 4.4.2 follows from (4.4.0.7) of Proposition 4.4.7 by using the formula (4.3.0.5), Lemma 4.4.1 and Proposition 4.3.1. This proves proposition 4.4.2. □

We have in fact obtained one direction of the first item of Theorem 4.1.4, recapitalised as follows.

*Proof that conjecture N.1 for $f$ and some $h$ implies conjecture S.1 for $f$.* By Proposition 4.4.2, (4.4.0.1) and (4.4.0.5), it follows that conjecture N.1 for a polynomial $f$ implies conjecture S.1 for $f$. □

We now focus on the inverse implication, namely, that conjecture S.1 for $f$ implies conjecture N.1 for $f$ and all $h$.

Let $\{z_1, ..., z_e\}$ be the set of critical values of $f$ and put $f_i = f - z_i$ for $1 \leq i \leq e$. For each $1 \leq i \leq e$ we take an embedded resolution $(Y_i, h_i)$ of $D(f_i)_K$. Let $K'$ be a finite extension of $K$ such that all critical values of $f$ belong to $K'$ and such that lemma 4.4.8 is verified for each of the $(Y_i, h_i)$. Let $M'$ be such that lemma 4.4.1 and corollary 4.4.9 are verified each of the $(Y_i, h_i)$. Let
$$J_i = J(h_i)$$

be the index set of irreducible components of $h_i^{-1}(Z(f_i))$. For $i$ with $1 \leq i \leq e$, consider the set
$$\mathcal{A}_i$$
consisting of all the irreducible components $W$ of the $E_I$ for $I \subset J_i$ which witness the power condition for $(f_i, h_i)$ (namely, there are some $d$ and some $g$ such that $(I, W, g, d)$ witnesses the power condition for $(f_i, h_i)$). For each $W \in \mathcal{A}_i$, let $d_W > 1$ be the maximal integer which witnesses the power condition on $W$ for $(f_i, h_i)$ for some $g$, and set

$$\rho_W = \sigma(f) - \frac{1}{2} - \sum_{i \in I} N_i(\sigma_i - \sigma(f))$$

and

$$N_W = \sum_{i \in I} N_i,$$

where $W$ is an irreducible components of $E_I$ and $I \subset J_i$. Put

$$\rho = \max\{\rho_W \mid W \in \mathcal{A}_i,\ i = 1, \ldots, e\}$$

and

$$N = \min\{N_W \mid W \in \mathcal{A}_i, \rho_W = \rho,\ i = 1, \ldots, e\}.$$

Finally put

$$\mathcal{B}_i = \{W \in \mathcal{A}_i \mid \rho_W = \rho \text{ and } N_W = N\}.$$

*Proof that conjecture S.1 for $f$ implies conjecture N.1 for $f$ and any $h$.* Suppose that conjecture S.1 is true for $f$ but conjecture N.1 is false for $f$ and $h$. This implies that

$$\rho > 0. \tag{4.4.0.11}$$

We need to find a contradiction. For $W \in \mathcal{A}_i$ and for all local fields $L$ containing $K'$ and with $\mathrm{char}(k_L) > M'$, and for all characters $\chi$ of $\mathcal{O}_L^\times$ with conductor $c_\chi = 1$, we have by lemma 4.4.8, corollary 4.4.9 and formula 4.3.0.4 that

$$c_{I,\overline{W},\Phi_{i,L},\chi} = \#\overset{\circ}{\overline{W}}(k_L)$$

whenever $\chi^{d_W} = 1$, and, that

$$c_{I,\overline{W},\Phi_{i,L},\chi} = 0$$

when $\chi^{d_W} \neq 1$. By Dirichlet's theorem on arithmetic progressions, we can find $L$ with arbitrarily large residue field characteristic and such that $d_W | p_L - 1$ for all $W \in \mathcal{A}_i$ and all $i$. For each $i$ and each $W \in \mathcal{A}_i$ we call $G_{W,L}$ the subgroup of order $d_W$ of the group of characters of $\mathcal{O}_L^\times$ which are trivial on $1 + \mathcal{M}_L$.

Put $m = N + 1$, fix $z$ of order $-m$. Write $z = u\varpi_L^{-m}$ for some uniformizer $\varpi_L$ of $\mathcal{O}_L$. By 4.3.1, 4.4.6, 4.4.7, (4.4.0.1), and the Lang-Weil estimates 4.4.3, we have

$$|E_{f,L,\psi_z}| = |\sum_{i=1}^{e} \Psi_L(zz_i) \sum_{W \in \mathcal{B}_i} \sum_{\chi \in G_{W,L} \setminus \{1\}} \chi(u) q_L^{\frac{1}{2}} g_{\chi^{-1}} |q_L^{-m\sigma(f)} q_L^\rho + o(m^{n-1} q_L^{-m\sigma(f)} q_L^\rho),$$
$$\tag{4.4.0.12}$$

where $\Psi_L = \exp(2\pi i \operatorname{Tr}_{L/\mathbb{Q}_{p_L}})$. Let us put

$$A = |\sum_{i=1}^{e} \Psi_L(zz_i) \sum_{W \in \mathcal{B}_i} \sum_{\chi \in G_{W,L} \setminus \{1\}} \chi(u) q_L^{\frac{1}{2}} g_{\chi^{-1}}|.$$

By Corollary 4.3 of [CGH16], there is a constant $c > 0$ such that for infinitely many $L$ of aribitrarily large residue field characteristic one has for some $z$ of order $-m$ that $c \leq A$. Combined with (4.4.0.12) and (4.4.0.11) this gives a contradiction to S.1 for $f$. Hence, conjecture N.1 must hold for $f$ and $h$. But then conjecture N.1 must hold for $f$ and any $h$, since we can always assume that $h$ is defined over the algebraic closure of $\mathbb{Q}$ (for example by the first order completeness of the theory of algebraically closed fields of characteristic zero). □

*Proof of theorem 4.1.4.* We have already proved the first item of theorem 4.1.4. For the second item, if $Z$ is a closed subset of $Z(f)$ then 4.4.6 is true if we replace $\Phi_{1,L}$ by $\Phi_{Z,L}$ and $\sigma(f)$ by $\operatorname{lct}_Z(f)$. If $W, I, d$ are as in the formulation of conjecture N.2, in the case $\overline{Z} \cap \overline{h}(\overline{W}) = \emptyset$ one has $c_{I,\overline{W},\Phi_{Z,L},\chi} = 0$ and there is nothing to do. If $\overline{W}$ is an irreducible component of $\overline{E}_I$ with $\overline{W} \neq \emptyset$ and $\overline{h}(\overline{W}) \subset \overline{Z}$, then 4.4.7 is still true if we replace $\Phi_{1,L}$ by $\Phi_{Z,L}$ and $\sigma(f)$ by $\operatorname{lct}_Z(f)$. If $\overline{Z} \cap \overline{h}(\overline{W}) \neq \emptyset$ but $\overline{h}(\overline{W}) \not\subseteq \overline{Z}$, we note that $\dim(\overline{h}^{-1}(\overline{Z}) \cap \overline{W}) \leq n - |I| - 1$, so by the Lang-Weil estimates 4.4.3, there exists a constant $C_I$ only depend on $(Y, h)$ such that

$$|c_{I,\overline{W},\Phi_{Z,L},\chi}| \leq \#(\overline{h}^{-1}(\overline{Z}) \cap \overline{W}) \leq C_I q_L^{n-|I|-1},$$

so for such $\overline{W}$ we get

$$|g_{\chi^{-1}} \operatorname{Coeff}_{t^{m-1}}(q_L^{-n} c_{I,\overline{W},\Phi_{Z,L},\chi} \prod_{i \in I} \frac{(q_L-1)t^{N_i} q_L^{-\nu_i}}{1-t^{N_i} q_L^{-\nu_i}})| \leq C_I q_L^{-m \operatorname{lct}_Z(f) + \operatorname{lct}_Z(f) - \frac{3}{2}} m^{\#(I)-1}$$

for all $m \geq 2$. Now the discussion for the proof of first item of Theorem 4.1.4 can be repeated, with $\operatorname{lct}_Z(f) \leq 1$ instead of $\sigma(f)$. □

The following corollary relates to the uniform aspect in (1.2.2) of Conjecture 1.2 of [CV16]. In fact, such uniform behavior of bounds on exponential sums follows much more generally from the conjectures on numerical data and the results of this paper.

**Corollary 4.4.10.** *If conjecture N.2 is true for each polynomial $f_y(x) = f(x) - f(y)$ where $y \in \mathbb{C}^n$ and with $Z = \{y\}$, then conjecture S.2 for $f$ and with $Z = \{y\}$ holds uniformly in $y \in (\mathcal{O}^{\operatorname{int}})^n$, where $\mathcal{O}^{\operatorname{int}}$ is the integral closure of $\mathcal{O}$ in $\mathbb{C}$. Namely, one can choose $M$ and $M'$ as in conjecture S.2 for $f$, $Z = \{y\}$ but with $\operatorname{lct}_y(f_y)$ replaced by $a_{y,L} = \min_{y' \in y + (\mathcal{M}_L \cap \mathcal{O}^{\operatorname{int}})^n} \operatorname{lct}_{y'}(f_{y'})$, where $M$ and $M'$ do not depend on $y$ and $L$ is a local field over $\mathcal{O}$ of residue field characteristic at least $M'$.*

*Proof.* Write $\{z_1, ..., z_e\}$ for the set of critical values of $f$. If one has $\operatorname{ord}(f(y) - z_i) = 0$ for each $i$, then there is nothing to do, by a reasoning as for (4.4.0.5). Since the formula 3.6.1.1 and section 3.6 in chapter 3, we may in fact assume that $f(y) = z_i$ for some $i$. In the proof of the implication from N.2 to S.2, the choices of $M'$ and of $M$ can be made depending only on a fixed embedded resolution $(Y_i, f_i)$ where $f_i(x) = f(x) - z_i$. This gives the desired uniformity. □

# Bibliography

[A'C73]   N. A'Campo. Sur la monodromie des singularités isolées d'hypersurfaces complexes. *Inventiones mathematicae*, 20:147–170, 1973.

[AGZV12]  V. Arnold, S. Gusein-Zade, and A. Varchenko. *Singularities of differentiable maps. Volume 2. Monodromy and asymptotics of integrals.* Birkhäuser Classics, Springer, New York, 2012. translated from Russian by H. Porteous and revised by the authors and J. Montaldi.

[BS66]    S. E. Borevich and I. R. Shafarevich. *Number Theory.* Academic Press. 1966.

[CCL10]   R. Cluckers, G. Comte, and F. Loeser. Lipschitz continuity properties for p-adic semi-algebraic and subanalytic functions. *GAFA*, 20(1):68–87, 2010.

[CCL15]   R. Cluckers, G. Comte, and F. Loeser. Non-archimedean Yomdin-Gromov parametrizations and points of bounded height. *Forum of Mathematics, Pi*, 3(e5):60 pages, 2015.

[CCM$^+$] R. Cluckers, G. Comte, D. J. Miller, J.-P. Rolin, and T. Servi. Integration of oscillatory and subanalytic functions. *ArXiv e-prints*.

[CGH]     R. Cluckers, J. Gordon, and I. Halupczok. Uniform analysis on local fields and applications to orbital integrals. *ArXiv e-prints*.

[CGH14a]  R. Cluckers, J. Gordon, and I. Halupczok. Integrability of oscillatory functions on local fields: transfer principles. *Duke Math. J.*, 163(8):1549–1600, 2014.

[CGH14b]  R. Cluckers, J. Gordon, and I. Halupczok. Local integrability results in harmonic analysis on reductive groups in large positive characteristic. *Ann. Sci. Ecole Norm. Sup*, 47(6):1163–1195, 2014.

[CGH14c]  R. Cluckers, J. Gordon, and I. Halupczok. Motivic functions, integrability, and applications to harmonic analysis on p-adic groups. *Electronic Research Announcements in Math.*, 21:137–152, 2014.

[CGH16]   R. Cluckers, J. Gordon, and I. Halupczok. *Transfer principles for Bounds of motivic exponential functions.* Springer, 2016. chapter in Families of automorphic forms and the trace formula, Müller, Werner, Shin, Sug Woo, Templier, Nicolas (Eds.).

[CH]      R. Cluckers and I. Halupczok. Integration of functions of motivic exponential class, uniform in all non-archimedean local fields of characteristic zero. *ArXiv e-prints*.

[CHL11]   R. Cluckers, T. Hales, and F. Loeser. *Transfer Principle for the Fundamental Lemma*. On the Stabilization of the Trace Formula, edited by L. Clozel, M. Harris, J.-P. Labesse, B.-C. Ngô. International Press of Boston, 2011.

[CHLR]    R. Cluckers, I. Halupczok, F. Loeser, and M. Raibaut. Distributions and wave front sets in the uniform non-archimedean setting. *ArXiv e-prints*.

[CKN16]   P. Cubides Kovacsics and K. H. Nguyen. An example of a p-minimal structure without skolem functions. *The Journal of Symbolic Logic*, 82(2):147–170, 2016.

[CL]      R. Cluckers and L. Lipshitz. Strictly convergent analytic structures. *J. Eur. Math. Soc.* to appear.

[CL05]    R. Cluckers and F. Loeser. *Ax-Kochen-Ershov Theorems for p-adic integrals and motivic integration*. Geometric methods in algebra and number theory, Y. Tschinkel and F. Bogomolov (Eds.). 2005. Proceedings of the Miami Conference 2003.

[CL07]    R. Cluckers and F. Loeser. b-minimality. *Journal of Mathematical Logic*, 7(2):195–227, 2007.

[CL08]    R. Cluckers and F. Loeser. Constructible motivic functions and motivic integration. *Inventiones Mathematicae*, 173(1):23–121, 2008.

[CL10]    R. Cluckers and F. Loeser. Constructible exponential functions, motivic fourier transform and transfer principle. *Annals of Mathematics*, 171(2):1011–1065, 2010.

[CL11]    R. Cluckers and L. Lipshitz. Fields with analytic structure. *J. Eur. Math. Soc*, 13:1147–1223, 2011.

[CL15]    R. Cluckers and F. Loeser. Motivic integration in all residue field characteristics for henselian discretely valued fields of characteristic zero. *J. Reine und Angew. Math.*, 701:1–31, 2015.

[CLR06]   R. Cluckers, L. Lipshitz, and Z. Robinson. Analytic cell decomposition and analytic motivic integration. *Ann. Sci. École Norm. Sup*, 39(4):535–568, 2006.

[Clu04]   R. Cluckers. Analytic *p*-adic cell decomposition and integrals. *Trans. Amer. Math. Soc.*, 356(4):1489–1499, 2004.

[Clu08a]  R. Cluckers. Igusa and Denef-Sperber conjectures on nondegenerate *p*-adic exponential sums. *Duke Math. J.*, 141(1):205–216, 2008.

[Clu08b]  R. Cluckers. The modulo $p$ and $p^2$ cases of the Igusa conjecture on exponential sums and the motivic oscillation index. *Int. Math. Res. Not*, 2008, 2008.

[Clu10]   R. Cluckers. Exponential sums: questions by Denef, Sperber, and Igusa. *Trans. Amer. Math. Soc.*, 362(7):3745–3756, 2010.

[CNa]     W. Castryck and K. H. Nguyen. New bounds for exponential sums with a non-degenerate phase polynomial. *ArXiv e-prints*.

[CNb] S. Chambille and K. H. Nguyen. Proof of Cluckers-Veys's conjecture on exponential sums for polynomials with log-canonical threshold at most a half. *ArXiv e-prints*.

[CNc] R. Cluckers and K. H. Nguyen. Conjectures on exponential sums and on numerical data. preprint.

[Coh69] P. J. Cohen. Decision procedures for real and *p*-adic fields. *Comm. Pure Appl. Math.*, 22(1):131–151, 1969.

[CPW] R. Cluckers, J. Pila, and A. Wilkie. Uniform parameterization of subanalytic sets and diophantine applications. *ArXiv e-prints*.

[CV16] R. Cluckers and W. Veys. Bounds for *p*-adic exponential sums and log-canonical thresholds. *Amer. J. Math.*, 138(1):61–80, 2016.

[DD88] J. Denef and L. van den Dries. *p*-adic and real subanalytic sets. *Annals of Math.*, 128(1):79–138, 1988.

[Del74] P. Deligne. La conjecture de Weil i. *Publications Mathématiques de l'IHÉS*, 43:273–307, 1974.

[Del77] P. Deligne. *Cohomologie étale*, volume 569 of *Lecture Notes in Math.* Springer, Berlin, 1977.

[Del80] P. Deligne. La conjecture de Weil II. *Publications Mathématiques de l'IHÉS*, 52:137–252, 1980.

[Den] J. Denef. Report on Igusa's local zeta function. *Séminaire Bourbaki,*, 77.

[Den84] J. Denef. The rationality of the Poincaré series associated to the *p*-adic points on a variety. *Inventiones Mathematicae*, 77:1–23, 1984.

[Den85] J. Denef. On the evaluation of certain *p*-adic integrals. *Séminaire de théorie des nombres (Paris 1983-84)*, 59:25–47, 1985. edited by C. Goldstein, Prog. in math.

[Den87] J. Denef. On the degree of Igusa's local zeta function. *Amer. J. Math.*, 109(6):991–1008, 1987.

[Den91] J. Denef. Local zeta functions and euler characteristics. *Duke Math. J.*, 63(3):713–721, 1991.

[Den00] J. Denef. *Arithmetic and geometric applications of quantifier elimination for valued fields*, volume 39. Cambridge Univ. Press, Cambridge, 2000. Model theory, algebra, and geometry, Math. Sci. Res. Inst. Publ.

[DH01] J. Denef and K. Hoornaert. Newton polyhedra and Igusa's local zeta function. *Journal of Number Theory*, 89:31–64, 2001.

[DK] P. Deligne and N. Katz. *Séminaire de Géométrie Algébrique du Bois Marie - 1967-69 - Groupes de monodromie en géométrie algébrique*, volume 2 of *Lecture Notes in Mathematics (in French). 340.* Berlin; New York: Springer-Verlag.

[DK01] J.-P. Demailly and J. Kollár. Semi-continuity of complex singularity exponents and Kähler-Einstein metrics on Fano orbifolds. *Annales Scientifiques de l'École Normale Supérieure*, 34(4):525–556, 2001.

[DL98]    J. Denef and F. Loeser. Motivic Igusa functions. *Journal of Algebraic Geometry*, 7:505–537, 1998.

[DL99a]   J. Denef and F. Loeser. Germs of arcs on singular algebraic varieties and motivic integration. *Invent. Math.*, 135(1):201–232, 1999.

[DL99b]   J. Denef and F. Loeser. Motivic exponential integrals and a motivic Thom-Sebastiani theorem. *Duke Mathematical Journal*, 99:285–309, 1999.

[DL01]    J. Denef and F. Loeser. Definable sets, motives and p-adic integrals. *Journal of the Amer. Math. Soc.*, 14:429–469, 2001.

[DL02]    J. Denef and F. Loeser. Motivic integration and the grothendieck group of pseudo-finite fields. *Proceedings of the International Congress of Mathematicians,*, II:13–23, 2002.

[Dri92]   L. van den Dries. Analytic Ax-Kochen-Ersov theorems. *Proceedings of the International Conference on Algebra, Part 3 (Novosibirsk, 1989)*, pages 379–398, 1992.

[Dri98]   L. van den Dries. *Tame topology and o-minimal structures*. Lond. Math .Soc. Lect. Note Series. 1998.

[dS93]    M. P. F. du Sautoy. Finitely generated groups, p-adic analytic groups and Poincaré series. *Annals of Mathematics Second Series*, 137(3):639–670, 1993.

[DS01]    J. Denef and S. Sperber. Exponential sums mod $p^n$ and Newton polyhedra. *Bulletin Belg. Math. Soc.–Simon Stevin*, suppl.:55–63, 2001.

[DV95]    J. Denef and W. Veys. On the holomorphy conjecture for Igusa's local zeta function. *Proc. Amer. Math. Soc.*, 123(1):2981–2988, 1995.

[ELT]     A. Esterov, A. Lemahieu, and K. Takeuchi. On the monodromy conjecture for non-degenerate hypersurfaces. *ArXiv e-prints*.

[FEM10]   T. D. Fernex, L. Ein, and M. Mustaţă. Shokurov's ACC conjecture for log canonical thresholds on smooth varieties. *Duke Math. J.*, 152(1):93–114, 2010.

[Fre]     E. Frenkel. *Langlands Correspondence for Loop Groups*. Cambridge Studies in Advanced Mathematics. 9780521854436.

[Ful98]   W. Fulton. *Intersection theory*, volume 2 of *Ergebnisse der Mathematik und ihrer Grenzgebiete. 3. Folge. A Series of Modern Surveys in Mathematics [Results in Mathematics and Related Areas. 3rd Series. A Series of Modern Surveys in Mathematics]*. Springer-Verlag, Berlin, second edition, 1998.

[Gro]     A. Grothendieck. *Séminaire de Géométrie Algébrique du Bois Marie - 1967-69 - Groupes de monodromie en géométrie algébrique*, volume 1 of *Lecture Notes in Mathematics (in French). 288*. Berlin; New York: Springer-Verlag.

[Har77]   R. Hartshorne. *Algebraic geometry*, volume 52 of *Graduate Texts in Mathematics*. Springer, 1977.

[HHM06]  D. Haskel, H. Hrushovski, and D. Macpherson. Definable sets in algebraically closed valued fields: elimination of imaginaries. *J. Reine und Angew. Math.*, 597:175–236, 2006.

[HHM13]  D. Haskel, H. Hrushovski, and D. Macpherson. Unexpected imaginaries in valued fields with analytic structure. *J. Symbolic Logic*, 78(2):523–542, 2013.

[Hir64]  H. Hironaka. Resolution of singularities of an algebraic variety over a field of characteristic zero: I. *Annals of Mathematics*, 79(1):109–203, 1964.

[HK06]  E. Hrushovski and D. Kazhdan. Integration in valued fields. *Algebraic and Number Theory, Progress in Mathematics*, 253(1):261–405, 2006.

[HK08]  E. Hrushovski and D. Kazhdan. The value ring of geometric motivic integration, and the Iwahori hecke algebra of sl2. *Geom. Funct. Anal*, 17:1924–1967, 2008. With an appendix by Nir Avni.

[HL15]  E. Hrushovski and F. Loeser. Monodromy and the lefschetz fixed point formula. *Ann. Sci. École Norm. Sup.*, 48:313–349, 2015.

[HLP14]  E. Hrushovski, F. Loeser, and B. Poonen. Berkovich spaces embed in Euclidean spaces. *L'Enseignement Mathématique*, 60:273–292, 2014.

[HMR]  E. Hrushovski, B. Martin, and S. Rideau. Definable equivalence relations and zeta functions of groups. *ArXiv e-prints*. with an appendix by R. Cluckers.

[Hoo02]  K. Hoornaert. Newton polyhedra and the poles of Igusa's local zeta function. *Bulletin of the Belgian Mathematical Society – Simon Stevin*, 9(4):589–606, 2002.

[Igu74]  J. I. Igusa. Complex powers and asymptotic expansions. *I. J. Reine angew. Math.*, 268-269:110–130, 1974.

[Igu78]  J. I. Igusa. *Lectures on forms of higher degree*, volume 59 of *Lectures on mathematics and physics*. Tata Institute of Fundamental Research. Narosa Pub. House, 1978. (notes by S. Raghavan).

[Kat99]  N. M. Katz. Estimates for "singular" exponential sums. *Internat. Math. Res. Notices*, (16):875–899, 1999.

[Kol97]  J. Kollár. Singularities of pairs. *Algebraic Geometry, Santa Cruz 1995, Proceedings of Symposia in Pure Mathematics*, (62):221–287, 1997.

[Kou76]  A. G. Kouchnirenko. Polyèdres de Newton et nombres de Milnor. *Invent. Math.*, 32(1):875–899, 1976.

[Kow]  E. Kowalski. The large sieve, monodromy and zeta functions of curves. *Crelle Jounal*.

[Laz04]  R. Lazarsfeld. *Positivity in algebraic geometry II: positivity for vector bundles, and multiplier ideals*, volume 49 of *Ergebnisse der Mathematik und ihrer Grenzgebiete . 3. Folge. A Series of Modern Surveys in Mathematics*. Springer, 2004.

[Le15]  Quy Thuong Le. Proofs of the integral identity conjecture over algebraically closed fields. *Duke Math. J.*, 164(1):157–194, 2015.

[Lic13]   B. Lichtin. On a conjecture of Igusa. *Mathematica*, 59(2):399–425, 2013.

[Lic16]   B. Lichtin. On a conjecture of Igusa II. *American Journal of Mathematics*, 138(1):201–249, 2016.

[LS03]    F. Loeser and J. Sebag. Motivic integration on smooth rigid varieties and invariants of degeneration. *Duke Mathematical Journal*, 119:315–344, 2003.

[LW54]    S. Lang and A. Weil. Number of points of varieties in finite fields. *Amer. J. Math.*, 76(4):819–827, 1954.

[Mac76]   A. Macintyre. On definable subsets of $p$-adic fields. *J. Symb. Logic*, 41:605–610, 1976.

[Mac90]   A. Macintyre. Rationality of $p$-adic Poincaré series: uniformity in $p$. *Ann. Pure Appl. Logic*, 49(1):31–74, 1990.

[Mal83]   B. Malgrange. Le polynôme de Bernstein-Sato et cohomologie évanescente. *Astérisque*, 101-102:233–267, 1983.

[Mal90]   B. Malgrange. Intégrales asymptotiques et monodromie. *Annales scientifiques de l'École Normale Supérieure*, Série 4 : Tome 7(3):405–430, 1990.

[Mar02]   D. Marker. *Model Theory: An introduction*, volume 217 of *Graduate texts in mathematics*. Springer-Verlag, 2002.

[Meu81]   D. Meuser. On the rationality of certain generating functions. *Math. Ann.*, 256:303–310, 1981.

[Mil68]   J. Milnor. *Singular Points of Complex Hypersurfaces*, volume 61 of *Annals of Mathematics Studies*. Princeton University Press, 1968.

[MN15]    M. Mustaţă and J. Nicaise. Weight functions on non-archimedean analytic spaces and the Kontsevich-Soibelman skeleton. *Algebraic Geometry*, 2(3):365–404, 2015.

[Mus]     M. Mustaţă. *Zeta functions in algebraic geometry*. Lecture notes.

[Mus02]   M. Mustaţă. Singularities of pairs via jet schemes. *J. Amer. Math. Soc.*, 15(3):599–615, 2002.

[Mus12]   M. Mustaţă. Impanga lecture notes on log canonical thresholds. *In: Contributions to algebraic geometry, EMS Ser. Congr. Rep.*, pages 407–442, 2012. Notes by Tomasz Szemberg.

[Ngua]    K. H. Nguyen. Motivic oscillation index conjecture on exponential sums of polynomials in three or four variables. preprint.

[Ngub]    K. H. Nguyen. On some criterions of Poisson summation formula of Siegel-Weil type. preprint.

[Nguc]    K. H. Nguyen. The rational property of the Poincaré series of geometric definable equivalent relations on some minimal structures. preprint.

[Ngud]    K. H. Nguyen. Uniform rationality of the Poincaré series of definable, analytic equivalence relations on local fields. *to appear in Trans. Amer. Math. Soc, series B*, page 26 pages.

[NP08]  J. Nicaise and S. Payne. A tropical motivic Fubini theorem with applications to Donaldson-Thomas theory. *J. Algebra*, 319:1585–1610, 2008.

[NS07]  J. Nicaise and J. Sebag. Motivic Serre invariants of curves. *Manuscripta Math.*, 123(2):105–132, 2007.

[NS08a]  J. Nicaise and J. Sebag. Motivic Serre invariants and Weil restriction. *J. Algebra*, 319:1585–1610, 2008.

[NS08b]  J. Nicaise and J. Sebag. Motivic Serre invariants, ramification, and the analytic milnor fiber. *Inventiones Mathematicae*, 168(1):133–173, 2008.

[NX16]  J. Nicaise and C. Xu. Poles of maximal order of motivic zeta functions. *Duke Mathematical Journal*, 165(2):217–243, 2016.

[Oes82]  J. Oesterlé. Réduction modulo $p^n$ des sous-ensembles analytiques fermés de $\mathbb{Z}_p^N$. *Invent. math.*, 66:325–341, 1982.

[Pas89]  J. Pas. Uniform *p*-adic cell decomposition and local zeta functions. *Journal für die reine und angewandte Mathematik*, 399:137–172, 1989.

[RL05]  Antonio Rojas-León. Estimates for singular multiplicative character sums. *Int. Math. Res. Not.*, (20):1221–1234, 2005.

[Seg06]  D. Segers. Lower bounds for the poles of Igusa's *p*-adic zeta functions. *Mathematische Annalen*, 336(3):659–669, 2006.

[Ser81]  J-P. Serre. Quelques applications du théorème de densité de Chebotarev. *Inst. Hautes études Sci. Publ. Math.*, 54:323–401, 1981.

[Shi55]  G. Shimura. Reduction of algebraic varieties with respect to a discrete valuation of the basis field. *Amer. J. Math.*, 77(1):134–176, 1955.

[ST71]  M. Sebastiani and R. Thom. Un résultat sur la monodromie. *Inventiones mathematicae*, 13:90–96, 1971.

[Vey91]  W. Veys. Relations between numerical data of an embedded resolution. *Amer. J. Math.*, 113(4):573–592, 1991.

[Vey98]  W. Veys. More congruences for numerical data of an embedded resolution. *Composition Math.*, 112(3):313–331, 1998.

[Vol10]  C. Voll. Functional equations for zeta functions of groups and rings. *Ann. Math.*, 172:1181–1218, 2010.

[Wei48]  A. Weil. On some exponential sums. *Proc. N. A. S.*, 34(5):204–207, 1948.

[Wri]  J. Wright. Exponential sums and polynomial congruences in two variables: the quasi-homogeneous case. *ArXiv e-prints*.